# THE POSTAGE STAMP PROBLEM

## Formulae and proof for the case of 3 denominations

M.F. Challis

Hempnall House
Lundy Green
Hempnall
NORWICH NR15 2NU
England

Abstract

$A_k = \{1, a_2, ..., a_k\}$ is an h-basis for n if every positive integer not exceeding n can be expressed as the sum of no more than h values $a_i$; an extremal h-basis $A_k$ is one for which n is as large as possible. Computing such extremal bases has become known as the Postage Stamp Problem.

In 1968 Gerd Hofmeister published formulae for solutions to the Postage Stamp Problem with three denominations (h = 3). This paper presents an alternative approach developed independently by the author 25 years later.

## Revision History

### Issue 1.0 December 31st 1990

This version was sent to Professor Selmer in May 1992.

### Issue 2.0 November 1st 1997

This version contains corrections and additions made after a complete re-read during 1995 and 1996. Here is a summary of the changes:

Some clarifications have been added, and a number of typographical or other trivial errors have been corrected.

An error in Theorem 300 (and hence in Tables 100 and 300) has been corrected: there is an additional formula for OSG(n,0) for even n. This first came to light in February 1991 when I considered an alternative approach to the proof of this Theorem; this alternative is now included as an addendum to the Theorem.

An error in the proof of Theorem 500 (where I failed to show that PP(s)<4 for all s>=40) has been corrected. This came to light as I was rereading my proof: the original argument shows only that PP(40)<4 and PP(s) -> 3.574 ... as s tends to infinity, but does not show that PP(s) decreases monotonically as s increases.

In March 1995, Selmer (to whom I had sent Issue 1 of the proof in 1992) observed that Conjecture 1 of section 2.14 - "all underlying stride generators are canonical" - is equivalent to "$h_1$<=$h_0$ for h-bases $A_3$". This had already been proved elsewhere, but Selmer wondered if a simpler proof might be found using my techniques. By early 1996 I had succeeded in proving both conjectures from section 2.14 using stride generators and thread diagrams, and wrote up the results as a separate document entitled "A proof that $h_1$,$h_2$<=$h_0$ for any h-basis $A_3$". My proof remains complicated, but throws some further light on topics mentioned in this document, and I have added forward references to link the two proofs together.

### Issue 3.0 November 5th 2013

Added an Abstract and a list of relevant references.

CONTENTS

# SECTION 1

# INTRODUCTION

Section 1.1 introduces the postage stamp problem and the concept of the "stride generator", and suggests how the two may be related.

Section 1.2 is the top level proof of the theorem: maximal sets are optimal stride generators of order 1 at least for all s>=81.

Section 1.3 contains tables of formulae for optimal stride generators and M(3,s) and the corresponding maximal sets.

## 1.1 Background and basic terminology

The general postage stamp problem can be stated informally as follows:

> The Post Office decides to rationalise postal services by providing standard envelopes which have space for at most s stamps, and charging one penny per ounce for postage. They decide to issue just d different stamp values.
>
> What are the best choices for these values, and what is the heaviest letter that can be posted?

The following definition describes the problem more formally, and introduces some of the terminology used throughout this document.

Definition 100

We define the stamp values or "denominations" as the set:

$$A = \{a_i : i=1:d\}$$

Without loss of generality we assume $a_i<a_{i+1}$, and clearly $a_1=1$.

We say that the set A "generates" a value x with at most s stamps if there exist non-negative coefficients $\{c_i : i=1:d\}$ such that:

$$c_1a_1 + \ldots + c_ia_i + \ldots + c_da_d = x$$

and $c_1 + \ldots + c_i + \ldots c_d \le s$

We define the "cover" X = C(A,d,s) of the set $A = \{a_i : i=1:d\}$ using at most s stamps as:

i) A generates x for all x<=X

and ii) A does not generate the value X+1

The maximum cover achieved by any set A is called M(d,s):

M(d,s) = Max[over all sets A] C(A,d,s)

Any set whose cover is equal to M(d,s) is called a "maximal" set.

The "general postage stamp problem" is to construct formulae for M(d,s) and the corresponding maximal sets for arbitrary s and d.

[ As an example, consider the set A = {1,3,6} with s=3.

The following values can be generated:

```
1 = 1
2 = 1+1
3 = 1+1+1  or  3
4 = 3+1
5 = 3+1+1
6 = 3+3  or  6
7 = 3+3+1  or  6+1
8 = 6+1+1
```

```
     9 = 3+3+3  or  6+3
    10 = 6+3+1
    12 = 6+3+3  or  6+6
    13 = 6+6+1
    15 = 6+6+3
    18 = 6+6+6
```

We see that the cover C({1,3,6},3,3) = 10.

In fact, {1,4,5} is the maximal set for s=3, and M(3,3) = 15. ]

It is easy to determine formulae for s=1, d=1, and d=2; this document defines and proves formulae for d=3. No general formulae are known for other cases.

The cover of any set A = {1, $a_2$, $a_3$} can be subdivided as follows:

```
|            |            |    ....    |            |            |
 0            a3           2a3          ka3          (k+1)a3   X
                                                     <--- Y --->
```

where:

C(A,3,s) = X  is the cover, and:

$X = (k+1)a_3 + Y, \quad 0 <= Y < a_3-1$*

[*see Theorem 201 for proof that Y = $a_3$-1 is not possible]

We call each set of values $ja_3 <= x < (j+1)a_3$ the jth "stride"; the cover is made up of k+1 complete strides and one incomplete stride.

We call the kth stride - that is, the last complete stride - the "final stride".

Any value x may have more than one possible generation.

[Example:  A = {1, 4, 6},  s = 4,  x = 16:

16 = 2.6 + 1.4 + 0.1
16 = 0.6 + 4.4 + 0.1   ]

We define the "canonical generation" of x to be that which uses the maximum number of $a_3$ stamps - and if there are several with the same number of $a_3$ stamps, then choose the one with most $a_2$ stamps.

[Thus the first generation above is the canonical one]

Now consider the generation

$x = c_3a_3 + c_2a_2 + c_1$   where  $x = b_3a_3 + b_2$,  $0 <= b_2 < a_3$

We say that this generation of x is "of order ($b_3$-$c_3$)"; clearly, of all the generations of x, the canonical generation has the smallest order, and this is what we mean by "the order of the value x".

[In the example above, the canonical generation of 16 has order 0.

However, the only generation of 20 is:

20 = 2.6 + 2.4 + 0.1    -  which is of order 1]

We define the "order" of a sequence of values $x_1 <= x < x_2$ to be the maximum order of any value x belonging to the sequence.

We now show an important result:

Theorem 100

Suppose the set A = {1, $a_2$, $a_3$} generates all values $ka_3 <= x < (k+1)a_3$, with maximum order p.

Then A will also generate all values $pa_3 <= x < ka_3$.

[In other words, if A generates the final stride with order p, then it will also generate strides p, p+1, ... k-1 ]

Proof:

Consider a value x belonging to the kth stride and its canonical generation of order i:

$$x = c_3a_3 + c_2a_2 + c_1 \quad \text{where} \quad c_1+c_2+c_3 <= s \quad \text{and} \quad c_3 = k-i,\ c_1,c_2,c_3 >= 0$$

By definition, p is the maximum value of i over all values x in the kth stride - and because of our choice of canonical generation, we ensure that p is as small as possible. Clearly, p <= k (since $c_3$>=0).

We now show we can generate all strides p, p+1, ... k.

For consider p <= q <= k, and y belonging to the qth stride.

Then $y+(k-q)a_3$ belongs to the kth stride, and so can be generated as:

$$y+(k-q)a_3 = c_3a_3 + c_2a_2 + c_1 \quad \text{where} \quad c_1+c_2+c_3 <= s \quad \text{and} \quad c_3 >= k-p,\ c_1,c_2,c_3 >= 0$$

$$=> \quad y = c_3'a_3 + c_2a_2 + c_1$$

$$\begin{aligned} \text{where} \quad c_3' &= c_3-(k-q) \\ &>= (k-p)-(k-q) \\ &= q-p \\ &>= 0 \end{aligned}$$

$$\begin{aligned} \text{and} \quad c_3'+c_2+c_1 &= c_3-(k-q)+c_2+c_1 \\ &>= c_3+c_2+c_1 \\ &>= 0 \end{aligned}$$

This is a simple, but powerful result, because in practice p is quite small: indeed, we shall show later that for maximal sets p=1.

In any case, the result strongly suggests that we should look in more detail at ways in which values in the final stride can be generated.

Possibilities are:

0) $ka_3 + c_2a_2 + c_1$ $k+c_2+c_1 <= s$ (order 0 generation)

1) $(k-1)a_3 + c_2a_2 + c_1$ $(k-1)+c_2+c_1 <= s$ (order 1 generation)

......

p) $(k-p)a_3 + c_2a_2 + c_1$ $(k-p)+c_2+c_1 <= s$ (order p generation)

Another way of looking at this is to write n = (s-k) and consider $0 <= x < a_3$ and see if:

0) x can be generated using at most n of $\{1, a_2\}$

1) $x+a_3$ can be generated using at most (n+1) of $\{1, a_2\}$

......

p) $x+pa_3$ can be generated using at most (n+p) of $\{1, a_2\}$

We define SG(n,p) = $\{1, a_2, a_3\}$ to be an "n-stride generator of order p" if:

i) Each $0 < x < a_3$ satisfies one of the rules (0) ... (p) above

ii) At least one value $0 < x < a_3$ satisfies rule (p) above but does not satisfy any rule (0) ... (p-1)

iii) There is at least one value $a_3 < x < 2a_3$ that does not satisfy any of the rules (0) ... (p), and also does not satisfy rule (-1):

-1) $x-a_3$ can be generated using at most (n-1) of $\{1, a_2\}$

$a_3$ is known as the "length" of the n-stride generator SG(n,p).

We show in Theorem 202 that every cover C(A,3,s) has an underlying (s-k)-stride generator A = SG(s-k,p) for some p <= k: conditions (i) and (ii) follow from the fact that all values in the final stride have generations, and condition (iii) from the fact that at least one value in the (k+1)th stride cannot be generated.

The converse is not true: that is, a stride generator A = SG(n,p) does not necessarily define a cover C(A,3,s) for any s >= n+p. Conditions (i) and (ii) coupled with Theorem 100 show that strides p ... s-n can all be generated, and condition (iii) shows that there is at least one value in stride s-n+1 that cannot be generated; however, it may still be the case that some values in strides 0 ... p-1 cannot be generated.

Finally, we define an "optimal n-stride generator of order p" OSG(n,p) = $\{1, a_2, a_3\}$ to be an n-stride generator of order p with maximal length: that is, there is no SG(n,p) = $\{1, a_2', a_3'\}$ such that $a_3' > a_3$.

Optimal n-stride generators are likely to be of interest because if they do indeed provide a cover for some value of s this cover is likely to exceed that provided by any sub-optimal n-stride generator because of the influence of the factor $a_3$ in the formula for X.

It seems plausible that maximal sets will be optimal stride generators, and so computer programs were written to list both OSG(n,1) for various n and to determine M(3,s) and the corresponding maximal sets for various values of s. It transpires that the computed optimal stride generators of order 1 are precisely the computed maximal sets (for sufficiently large n):

| e.g. | n | OSG(n,1) | Value(s) of s for which OSG(n,1) is the maximal set |
|---|---|---|---|
| | 30 | 1, 34, 352 | 44 |
| | 31 | 1, 33, 374 | 45, 46 |
| | 32 | 1, 35, 397 | 47, 48 |
| | 33 | 1, 37, 420 | 49 |
| | 34 | 1, 36, 444 | 50 |
| | 35 | 1, 38, 469 | 51, 52 |

Furthermore, examination of the computed results suggests that the structure of the OSG(n,1) sequence as n varies is simpler than that of the maximal set sequence as s varies: instead of a cycle of 9 with respect to s, it exhibits a cycle of 3 with respect to n.

Clearly, we are going in the right direction, and one might hope that the proof would be reasonably simple; unfortunately it turns out to be much more complex than expected!

The rest of this document is structured as follows:

The formulae for M(3,s) and the corresponding maximal sets are given in Section 1.3; other formulae relating to stride generators are repeated here for easy reference.

The proof itself is presented at the top level as Theorem 101 in Section 1.2; this contains references to further theorems which deal with each major step in turn. Rather than present each of these theorems in sequence, the remainder of the document is structured by topic, with each section leading up to one or more of the major results.

Section 2 is an investigation into the properties of stride generators. Many of the results of this section are used elsewhere, and some further development of stride generator properties occurs in Section 5.2.

Section 3 groups together results about optimal and near optimal stride generators of low order.

Section 4 deals with optimising the choice of certain classes of stride generator for a given value of s.

Section 5 contains results which limit the order of the stride generators which need to be taken into consideration; the final result in this section shows that the order must be equal to 1 if s>=81.

During the course of the proof, a number of results - both numerical and algebraic - are generated by computer; Section 6 includes further details of the relevant programs.

Theorem 101 gives formulae valid only for s>=81; Section 7 looks at smaller values of s.

Finally, Section 8 includes summary lists of the theorems and definitions together with a glossary for reference.

Throughout the proof, notes, comments and examples appear within square brackets: [ ... ]; such items do not form part of the formal proof, but are intended to aid understanding or to provide pointers to other interesting aspects arising from this work.

## 1.2 Overview of proof

This section details the "top-level" proof, giving references to other sections where detailed proofs of the necessary intermediate results can be found.

Theorem 101

For any s>=81:

The maximal set A is an optimal (s-k)-stride generator OSG(s-k,1) of order 1.

The formulae for k and M(3,s) = C(A,3,s) depend on the value of r = s mod 9, and are given as $k_{opt}$ and $X_{opt}$ in Table 103.

Formulae for A = $\{1, a_2, a_3\}$ are given in Table 105.

Proof

Let A = $\{1, a_2, a_3\}$ have a non-trivial cover X = C(A,3,s); then by Theorem 201 we can write:

$$X = (k+1)a_3 + Y \quad \text{for some } 0<=k<s \text{ and } 0<=Y<a_3-1$$

[ "non-trivial" (cf "admissible" in the literature) means that $X>=a_3$ which ensures k>=0; the proof of Theorem 201 shows why $Y=a_3-1$ is not possible. ]

Our job is to choose A in such a way as to maximise X; for this we want both large $a_3$ and k, and, to a lesser extent, large Y.

For a given value of k, maximal $a_3$ is obtained by choosing A to be an optimal stride generator OSG(s-k,p) of some order p; but what value of p will give the highest value of $a_3$? Somewhat surprisingly, the answer is p=1, provided that both s and s-k are sufficiently large.

1) Definition 208 defines the "potential cover" of a stride generator SG(n,p) with respect to a value s>=n as:

$$(s-n+1)a_3 + Y$$

where Y is a function of the stride generator alone (that is, Y does not depend on s), and satisfies $0<=Y<a_3-1$.

[ Y=y-1 where $0<y<a_3$ is the first "break" in the stride generator; see Definition 200 for details. ]

Then Theorems 202 and 231 prove that every non-trivial cover C(A,3,s) defines an underlying (s-k)-stride generator A=SG(s-k,p) such that:

$$C(A,3,s) = (k+1)a_3 + Y \quad 0<=p<=k<s, \quad 0<=Y<a_3-1$$

where Y is the same function of the stride generator as above.

[ The potential cover P of a stride generator A = SG(n,p) with respect to s is the largest value >= $pa_3$ that has an s-generation of order i for some i<=p.

The cover X = C(A,3,s) of a set A is the largest value >= 0 that has an s-generation of *any* order.

Now suppose that SG(n,p) is the stride generator underlying the cover C(A,3,s). Clearly, X>=P, but there is no a-priori reason why X=P: there might be an s-generation for P+1 of order i>p. ]

In other words, potential cover - which is a property of any stride generator - is equal to actual cover whenever that stride generator happens to underlie a real cover.

The proof proceeds by:

a) identifying a good lower bound $X_{opt}$ for M(3,s);

b) identifying those stride generators which could potentially underlie covers C(A,3,s) >= $X_{opt}$;

c) proving that no such stride generator has a potential cover >= $X_{opt}$.

The significance of the "potential cover" concept lies in the fact that we do not need to determine whether a given stride generator in (b) actually underlies a cover or not; it is sufficient to prove in (c) that its potential cover is less than $X_{opt}$.

2) Theorem 301 determines formulae for OSG(n,1), the optimal n-stride generators of order 1; these, together with corresponding formulae for y=Y+1 are given in Table 101.

3) In Theorem 400, the potential cover X of OSG(s-k,1) is maximised for a given value of s; the corresponding optimal values $k_{opt}$ and $X_{opt}$ - valid for s>=18 - are given in Table 103.

Then Theorem 401 shows that A = OSG(s-$k_{opt}$,1) is indeed the underlying stride generator for C(A,3,s): in other words, this potential cover is a real one. Formulae for A = {1, $a_2$, $a_3$} are given in Table 105.

[The potential cover is defined as:

$$X = (k+1)a_3 + Y$$

As k increases, so s-k decreases and so does the corresponding value $a_3$ of OSG(s-k,1); it turns out that the value X reaches a maximum when $k_{opt}$ is approximately equal to s/3.

e.g. consider s=54:

| k | s-k | OSG(s-k,1) | Y | (k+1)$a_3$ | X | |
|---|---|---|---|---|---|---|
| 3 | 51 | {1, 55, 954} | 914 | 3816 | 4730 | |
| ... | | ... | ... | ... | | |
| 15 | 39 | {1, 43, 574} | 542 | 9184 | 9726 | |
| 16 | 38 | {1, 41, 547} | 517 | 9299 | 9816 | |
| 17 | 37 | {1, 39, 520} | 492 | 9360 | 9852 | <- optimum |
| 18 | 36 | {1, 40, 494} | 464 | 9386 | 9850 | |
| 19 | 35 | {1, 38, 469} | 441 | 9380 | 9821 | |

For s=54 we see that $X_{opt}$=9852, and $k_{opt}$=17; note that $k_{opt}$ is one less than the value which gives maximal (k+1)$a_3$!

It is also easy to show that the potential cover of 9852 supplied by OSG(37,1) is an actual cover - that is, C({1,39,520},3,54) = 9852.

On the other hand, the potential cover 4730 of OSG(51,1) is not realised: C({1,55,954},3,54) = 108. ]

4) We now have established $X_{opt}$ as a (good!) lower bound for M(3,s); we make use of this fact in Theorem 500 to show that for all s>=40 the order of the stride generator underlying the maximal cover must be <=3.

Of course, $X_{opt}$ might not be the maximal value because:

a) There might be a stride generator SG(n,p) for p = 0, 2 or 3 whose length is greater than or equal to that of OSG(n,1), or:

b) There might be a "sub-optimal" stride generator SG(n,p) whose length is less than that of OSG(n,1) but which results in a better value of X because of an improved value for Y.

[e.g. Suppose {1, 42, 519} were an SG(37,1) with a Y value of 516.

Then $(k+1)a_3$ + Y = 18*519 + 516 = 9858 would be an improvement on the optimum value 9852 delivered by the optimal generator OSG(37,1) = {1, 39, 520} ]

5) Sub-optimal n-stride generators can be classified according to their length as follows (see Definition 400):

OSG(n,1) has maximal length - say $A_3$

SG1(n,p) has length $A_3$-1

...

SGi(n,p) has length $A_3$-i

[e.g. OSG(37,1) = {1, 39, 520}

SG1(37,1) = {1, 42, 519}
SG1(37,2) = {1, 54, 519}

SG2(37,3) = {1, 77, 518} ]

Consider A = SGi(s-k,p) for i>=2 which therefore has length $a_3$<=$A_3$-2.

Then if SGi(s-k,p) is the underlying stride generator for a cover C(A,3,s) we know that:

$$\begin{aligned} C(A,3,s) &= (k+1)a_3 + Y \\ &< (k+2)a_3 \\ &<= (k+2)(A_3-2) \end{aligned}$$

Theorem 402 shows that for any s>=36 this value is < $X_{opt}$ for all 0<=k<s.

[e.g. continuing the s=54 example:

| k | s-k | SG2(s-k,p) | p | X' |
|---|---|---|---|---|
| 15 | 39 | {1, 85, 572} | 3 | 9724 |
| 16 | 38 | {1, 59, 545} | 2 | 9810 |
| 17 | 37 | {1, 77, 518} | 3 | 9842 |
| 18 | 36 | {1, 75, 492} | 4 | 9840 |
| 19 | 35 | {1, 84, 467} | 4 | 9807 |

All values of X' = $(k+2)a_3$ are less than $X_{opt}$ = 9852.

(As it happens, each of the stride generators cited above generates only a trivial cover for s=54 because $a_2$>55, but this may not always be the case for SG2 generators) ]

This means that we now need only consider stride generators whose length is >=$A_3$-1 of order p<=3.

6) Theorem 302 establishes formulae for SG1(n,1), the immediately sub-optimal n-stride generators of order 1; these, together with corresponding formulae for y=Y+1, are given in Table 102.

7) Theorem 403 shows that for any s>=9 the potential cover of SG1(s-k,1) is less than $X_{opt}$.

[e.g. continuing the s=54 example:

| k | s-k | SG1(s-k,1) | Y | X' | X |
|---|---|---|---|---|---|
| 15 | 39 | {1, 40, 573} | 546 | 9741 | 9714 |
| 16 | 38 | none exists | | | |
| 17 | 37 | {1, 42, 519} | 488 | 9861 | 9830 |
| 18 | 36 | {1, 37, 493} | 468 | 9860 | 9835 |
| 19 | 35 | none exists | | | |

Note that the upper bound test, successful for SG2 and beyond, fails with SG1: X' = (k+2)$a_3$ > $X_{opt}$ for both k=17 and k=18. ]

Now the only stride generators left to consider are those of order 0, 2 or 3 whose length >= $A_3$-1.

8) Theorem 510 shows that for any n>=52, no SG(n,p) for p = 0, 2 or 3 exists whose length is >= $A_3$-1.

This leaves us with the possibility that for some value(s) of s such a stride generator SG(n,p) for some n<52 might improve on $X_{opt}$.

9) Let A = SG(s-k,p) be a stride generator of order 0, 2 or 3 whose length is $a_3$>=$A_3$-1; by (8), s-k<52.

Suppose A is the stride generator underlying some cover C(A,3,s); then:

$$C(A,3,s) = (k+1)a_3 + Y < (k+2)a_3$$

We show in Theorem 512 that for any s>=81:

$$(k+2)a_3 < X_{opt} \quad \text{for any such A.}$$

This concludes the proof:

(3) above shows that A = OSG(s-$k_{opt}$,1) defines a cover C(A,3,s) = $X_{opt}$, with formulae for $k_{opt}$ and $X_{opt}$ given in Table 103 and formulae for A = {1, $a_2$, $a_3$} given in Table 105.

(4) to (9) above show that no other underlying stride generator for any set A can equal or better this cover for sufficiently large s.

The most onerous constraint - s>=81 - is imposed by (9).

## 1.3 Formulae and results

In this section, we summarise the major results and formulae.

Theorem 300 gives formulae for $\{1, a_2, a_3\}$ = OSG(n,0) and for the first break y.

The form depends on r = n mod 2 as follows:

| r | $a_2$ | $a_3$ | y |
|---|---|---|---|
| 0 | (n+2)/2 | $(n^2+6n+4)/4$ | $(n^2+4n)/4$ |
| | (n+4)/2 | $(n^2+6n+4)/4$ | $(n^2+4n-4)/4$ |
| 1 | (n+3)/2 | $(n^2+6n+5)/4$ | $(n^2+4n-1)/4$ |

Table 100 - OSG(n,0)

Theorem 301 gives formulae for $\{1, a_2, a_3\}$ = OSG(n,1) and for the first break y.

The form depends on r = n mod 3 as follows:

| r | $a_2$ | $a_3$ | y |
|---|---|---|---|
| 0 | n+4 | $(n^2+5n+6)/3$ | $(n^2+3n-9)/3$ |
| 1 | n+2 | $(n^2+5n+6)/3$ | $(n^2+3n-1)/3$ |
| 2 | n+3 | $(n^2+5n+7)/3$ | $(n^2+3n-4)/3$ |

Table 101 - OSG(n,1)

Theorem 302 gives formulae for $\{1, a_2, a_3\}$ = SG1(n,1) and for the first break y.

The form (and existence) depends on r = n mod 3 as follows:

| r | $a_2$ | $a_3$ | y |
|---|---|---|---|
| 0 | n+1 | $(n^2+5n+3)/3$ | $(n^2+3n)/3$ |
| 1 | n+5 | $(n^2+5n+3)/3$ | $(n^2+3n-16)/3$ |
| 2 | | No SG1(n,1) exists | |

Table 102 - SG1(n,1)

Theorem 400 maximises X = C(OSG(s-k,1),3,s) for given s to obtain the formulae for $X_{opt}$ and $k_{opt}$ in the following table. Theorem 101 proves that $X_{opt}$ = M(3,s) for s>=81.

The form depends on r = s mod 9 as follows, where s = 9t+r:

| r | $k_{opt}$ | $n_{opt}$=s-$k_{opt}$ | $X_{opt}$ |
|---|---|---|---|
| 0 | 3t-1 | 6t+1 | $36t^3 + 54t^2 + 22t$ |
| 1 | 3t | 6t+1 | $36t^3 + 66t^2 + 36t + 4$ |
| 2 | 3t | 6t+2 | $36t^3 + 78t^2 + 53t + 8$ |
| 3 | 3t+1 | 6t+2 | $36t^3 + 90t^2 + 71t + 15$ |
| 4 | 3t+1 | 6t+3 | $36t^3 + 102t^2 + 92t + 22$ |
| 5 | 3t+1 | 6t+4 | $36t^3 + 114t^2 + 116t + 36$ |
| 6 | 3t+1 | 6t+5 | $36t^3 + 126t^2 + 143t + 49$ |
| 7 | 3t+2 | 6t+5 | $36t^3 + 138t^2 + 173t + 68$ |
| 8 | 3t+2 | 6t+6 | $36t^3 + 150t^2 + 204t + 86$ |

Table 103 - M(3,s)

There are just six different optimal stride generators defined in table 103; their values can be worked out from table 101, and are as follows:

Stride generators corresponding to optimal covers

| r | $n_{opt}$ | $a_2$ | $a_3$ |
|---|---|---|---|
| 0,1 | 6t+1 | 6t+3 | $12t^2 + 14t + 4$ |
| 2,3 | 6t+2 | 6t+5 | $12t^2 + 18t + 7$ |
| 4 | 6t+3 | 6t+7 | $12t^2 + 22t + 10$ |
| 5 | 6t+4 | 6t+6 | $12t^2 + 26t + 14$ |
| 6,7 | 6t+5 | 6t+8 | $12t^2 + 30t + 19$ |
| 8 | 6t+6 | 6t+10 | $12t^2 + 34t + 24$ |

Table 104

The following formulae can be derived directly from tables 103 and 104 above.

Alternative formulae for M(3,s) and maximal sets

| r | $a_2$ | $a_3$ | $X_{opt}$ = M(3,s) |
|---|---|---|---|
| 0 | 6t+3 | $(2t+1)a_2 + (2t+1)$ | $(3t+0)a_3 + (2t+0)a_2 + (4t+0)$ |
| 1 | 6t+3 | $(2t+1)a_2 + (2t+1)$ | $(3t+1)a_3 + (2t+0)a_2 + (4t+0)$ |
| 2 | 6t+5 | $(2t+1)a_2 + (2t+2)$ | $(3t+1)a_3 + (2t+0)a_2 + (4t+1)$ |
| 3 | 6t+5 | $(2t+1)a_2 + (2t+2)$ | $(3t+2)a_3 + (2t+0)a_2 + (4t+1)$ |
| 4 | 6t+7 | $(2t+1)a_2 + (2t+3)$ | $(3t+2)a_3 + (2t+0)a_2 + (4t+2)$ |
| 5 | 6t+6 | $(2t+2)a_2 + (2t+2)$ | $(3t+2)a_3 + (2t+1)a_2 + (4t+2)$ |
| 6 | 6t+8 | $(2t+2)a_2 + (2t+3)$ | $(3t+2)a_3 + (2t+1)a_2 + (4t+3)$ |
| 7 | 6t+8 | $(2t+2)a_2 + (2t+3)$ | $(3t+3)a_3 + (2t+1)a_2 + (4t+3)$ |
| 8 | 6t+10 | $(2t+2)a_2 + (2t+4)$ | $(3t+3)a_3 + (2t+1)a_2 + (4t+4)$ |

Table 105

# SECTION 2

## STRIDE GENERATORS

The primary aim of this section is to prove that every non-trivial cover has an underlying stride generator whose potential cover is the same as the actual cover. To do this we have to investigate properties of stride generators in some depth, and the proof of the theorem mentioned above is not complete until Section 2.10.

We first take the informal definition of a stride generator (from Section 1.1) and derive equivalent formal definitions; Section 2.1 also includes some examples. Section 2.3 shows that every cover has an underlying stride generator with certain properties and proves the result required for "canonical" stride generators which were previously introduced in Section 2.2.

Section 2.4 introduces the idea of a stride generator as a collection of "threads" which can be represented in a "thread diagram". Examination of these diagrams suggests further properties that are developed in Sections 2.5 to 2.7.

The background work is now complete, and Section 2.8 is able to provide a full characterisation of zero order stride generators. Section 2.9 continues with non-zero order stride generators, and we show that for non-canonical stride generators the first break always has the highest order. This makes sense of the concept of "potential cover" in the non-canonical case, thus enabling Section 2.10 to complete the theorem partially proved (for canonical stride generators only) in Section 2.3.

It is perhaps not obvious that one set can be two or more stride generators SG(n,p) for different values of n,p. Section 2.11 investigates some of the properties of such series of stride generators. This section is not part of the proof, and is included for interest only.

Section 2.12 establishes some limits on the values $a_2$ and $a_3$ for a stride generator SG(n,p). These are useful when developing programs to list all stride generators with certain properties.

Sections 2.13 and 2.14 - again, not part of the proof - include some further results of interest and some as yet unproven speculations. The most significant conjecture is that all underlying stride generators are canonical; if a simple proof could be found, most of the content of Sections 2.4 to 2.9 inclusive would become redundant!

## 2.1 Equivalent definitions

In this section we first restate as Definition 200 the informal definition of a stride generator given in the introductory section 1.1, and then derive an alternative working definition as three conditions (SG1), (SG2) and (SG3). It is these conditions that are referenced in proofs throughout the sections that follow.

Definition 200

A set $A = \{1, a_2, a_3\}$ is an n-stride generator SG(n,p) of order p if:

(i) Every value $0<x<a_3$ is such that:

$x + ia_3$ can be generated using at most $(n+i)$ of $\{1, a_2\}$ for some $i<=p$*

(ii) At least one value $0<x<a_3$ requires $i=p$ in (i) above; in other words:

$x + pa_3$ can be generated using at most $(n+p)$ of $\{1, a_2\}$,

but $x + ia_3$ cannot be generated using at most $(n+i)$ of $\{1, a_2\}$ for any $i<p$

(iii) At least one value $a_3<x<2a_3$ is such that:

$x + ia_3$ cannot be generated using at most $(n+i)$ of $\{1, a_2\}$ for any $-1<=i<=p$

$a_3$ is known as the "length" of the stride generator.

p is known as the "order" of the stride generator.

$x-a_3$ (where $a_3<x<2a_3$ satisfies (iii)) is known as a "break" in the stride generator.

[*"i<=p" means "0<=i<=p": the lower bound is stated explicitly only if it is different from zero]

[Note that in any generation

$$x + ia_3 = c_2a_2 + c_1 \quad c_1,c_2>=0$$

we can assume without loss of generality that $c_1<=a_2-1$, since any other generation with $c_1,c_2>=0$ can only have a higher coefficient sum $c_1+c_2$.]

We now derive an equivalent working definition of a stride generator in terms of three conditions which must all hold.

Theorem 200

A is an n-stride generator SG(n,p) of order p iff:

For all $0<x<a_3$:

$x + ia_3 = c_2a_2 + c_1 \quad c_2+c_1<=n+i$ is soluble for some $i<=p$ (SG1)

There exists $0<x<a_3$ such that (SG1) is soluble for $i=p$,
but not for any $i<p$ (SG2)

There exists $0<y<a_3$ such that:

$y + ja_3 = c_2a_2 + c_1 \quad c_2+c_1<=n+j-1$ is not soluble for any $j<=p+1$ (SG3)

[As always, $c_2$ and $c_1$ must be >=0]

Proof

(SG1) and (SG2) are simply rewrites of (i) and (ii):

"$x + ia_3$ can be generated using at most (n+i) of {1, $a_2$}"

<=> "there exist $c_2$, $c_1$ such that $x + ia_3 = c_2a_2 + c_1$ with $c_2+c_1<=n+i$"

Similarly, (iii) <=>

There exists $a_3<x<2a_3$ such that:

$x + ia_3 = c_2a_2 + c_1 \quad c_2+c_1<=n+i$ is not soluble for any $-1<=i<=p$

We derive (SG3) directly by writing $y = x-a_3$ and $j = i+1$:

There exists $0<y<a_3$ such that:

$y + a_3 + ia_3 = y + ja_3 = c_2a_2 + c_1 \quad c_2+c_1<=n+j-1$
is not soluble for any $0<=j<=p+1$

[Example

We show that SG(4,2) = {1,6,13} is an order 2 stride generator with two break values at y=9 and y=10.

| must be: | | | <=4 | | | | <=5 | | | | <=6 |
|---|---|---|---|---|---|---|---|---|---|---|---|
| x | $c_2$ | $c_1$ | $c_2+c_1$ | $x+a_3$ | $c_2$ | $c_1$ | $c_2+c_1$ | $x+2a_3$ | $c_2$ | $c_1$ | $c_2+c_1$ |
| 1 | 0 | 1 | 1 | | | | | | | | |
| 2 | 0 | 2 | 2 | | | | | | | | |
| 3 | 0 | 3 | 3 | | | | | | | | |
| 4 | 0 | 4 | 4 | | | | | | | | |
| 5 | 0 | 5 | 5 | 18 | 3 | 0 | 3 | | | | |
| 6 | 1 | 0 | 1 | | | | | | | | |
| 7 | 1 | 1 | 2 | | | | | | | | |
| 8 | 1 | 2 | 3 | | | | | | | | |
| 9 | 1 | 3 | 4 | | | | | | | | |
| 10 | 1 | 4 | 5 | 23 | 3 | 5 | 7 | 36 | 6 | 0 | 6 |
| 11 | 1 | 5 | 6 | 24 | 4 | 0 | 4 | | | | |
| 12 | 2 | 0 | 2 | | | | | | | | |

The left-most group of columns above shows the order 0 generations:

$x + 0.a_3 = c_2a_2 + c_1 \qquad c_2+c_1 <= n+0 = 4$

This shows that 1,2,3,4,6,7,8,9,12 are all valid order 0 generations.

The middle group looks at the remaining values to see if they are order 1 generations:

$$x + 1.a_3 = c_2a_2 + c_1 \qquad c_2+c_1 <= n+1 = 5$$

We see that 5 and 11 are valid order 1 generations.

The right-most group shows that the remaining value 10 is an order 2 generation.

This shows that conditions (SG1) and (SG2) for a stride generator SG(4,2) are satisfied: all values 0<x<13 can be generated, and at least one value (10) can only be generated with an order p=2 generation.

Condition (SG3) requires that there exist $0<y<a_3$ such that:

$$y + ja_3 = c_2a_2 + c_1 \quad c_2+c_1<=n+j-1 \quad \text{is not soluble for any } j<=p+1$$

The following four groups of columns show that there are two such values: 9 and 10.

| | must be: | | <=3 | | | | <=4 | | | | <=5 | | | | <=6 |
|---|---|---|---|---|---|---|---|---|---|---|---|---|---|---|---|
| x | $c_2$ | $c_1$ | $c_2+c_1$ | $x+a_3$ | $c_2$ | $c_1$ | $c_2+c_1$ | $x+2a_3$ | $c_2$ | $c_1$ | $c_2+c_1$ | $x+3a_3$ | $c_2$ | $c_1$ | $c_2+c_1$ |
| 1 | 0 | 1 | 1 | | | | | | | | | | | | |
| 2 | 0 | 2 | 2 | | | | | | | | | | | | |
| 3 | 0 | 3 | 3 | | | | | | | | | | | | |
| 4 | 0 | 4 | 4 | 17 | 2 | 5 | 7 | 30 | 5 | 0 | 5 | | | | |
| 5 | 0 | 5 | 5 | 18 | 3 | 0 | 3 | | | | | | | | |
| 6 | 1 | 0 | 1 | | | | | | | | | | | | |
| 7 | 1 | 1 | 2 | | | | | | | | | | | | |
| 8 | 1 | 2 | 3 | | | | | | | | | | | | |
| 9 | 1 | 3 | 4 | 22 | 3 | 4 | 7 | 35 | 5 | 5 | 10 | 48 | 8 | 0 | 8 |
| 10 | 1 | 4 | 5 | 23 | 3 | 5 | 8 | 36 | 6 | 0 | 6 | 49 | 8 | 1 | 9 |
| 11 | 1 | 5 | 6 | 24 | 4 | 0 | 4 | | | | | | | | |
| 12 | 2 | 0 | 2 | | | | | | | | | | | | |

Thus SG(4,2) = {1,6,13} has two breaks at y=9 and y=10.

In fact, it turns out that any stride generator $\{1, a_2, a_3\}$ where $a_3>=2a_2$ has at most two breaks; this follows from Theorems 220 and 227 proved in section 2.9.]

## 2.2 Canonical stride generators

This section defines the concepts of "canonical break" and "canonical stride generator", and gives some examples.

In later sections we prove that every set A defines a unique canonical stride generator, and hypothesise that the stride generator which underlies a cover C(A,3,s) is always canonical.

Definition 201

For any break y in a stride generator SG(n,p) condition (SG3) shows that:

$y + ja_3 = c_2a_2 + c_1 \quad c_2+c_1<=n+j-1$ is not soluble for any $j<=p+1$

If this equation is in fact not soluble for any value of j, we say that the break is "canonical", and consider the "order of the break" to be infinite.

Otherwise, the "order of the break" is defined to be the smallest value $k>p+1$ for which there is a solution of:

$y + ka_3 = c_2a_2 + c_1 \quad c_2+c_1<=n+k-1$

Definition 202

A stride generator SG(n,p) is said to be "canonical" if every one of its breaks is canonical.

In other words, every break y in a canonical stride generator is such that:

$y + ja_3 = c_2a_2 + c_1 \quad c_2+c_1<=n+j-1$ is not soluble for any j

We show in section 2.9 (Theorem 230) that either all or none of the breaks in a stride generator are canonical; in other words, none of the breaks in a non-canonical stride generator is canonical.

[Examples

The first example is of a canonical stride generator: we demonstrate that the two breaks in the stride generator SG(4,2) = {1,6,13} are both canonical, and hence the stride generator itself is canonical.

Consider y = 9:

We know that $y + 3a_3 = 8a_2 + 0$, and so $y + (p+1)a_3$ already uses more than n+p (=7) $a_2$'s.

$y + (p+2)a_3$ will use at least 2 more $a_2$'s, since 13 = 2.6 + 1, whereas the limitation on the sum $c_2+c_1$ can only increase by 1 to n+p+1.

So matters can only get worse, and there is no value j for which a solution exists.

This is illustrated by the following table of $y+ja_3$ for j = 0,1 ..., where we see $c_2$ is increasing faster than n+j-1, and overtaking it at j=3:

| | j | $y+ja_3$ | $c_2$ | $c_1$ | $c_2+c_1$ | $n+j-1$ | |
|---|---|---|---|---|---|---|---|
| | 0 | 9 | 1 | 3 | 4 | 3 | |
| | 1 | 22 | 3 | 4 | 7 | 4 | |
| | 2 | 35 | 5 | 5 | 10 | 5 | |
| p+1 = | 3 | 48 | 8 | 0 | 8 | 6 | <- $c_2$ (on its own) > n+j-1 |
| | 4 | 61 | 10 | 1 | 11 | 7 | from here on in |
| | 5 | 74 | 12 | 2 | 14 | 8 | ... |

The second break at y = 10 can be shown to be canonical in the same way.

The second example is of a non-canonical stride generator:

SG(8,2) = {1,14,33}

with a single break y=22 with break order 4.

We demonstrate this is so using tables in the same way as in the example in section 2.1:

| | must be: | | <=8 | | | | <=9 | | | | <=10 |
|---|---|---|---|---|---|---|---|---|---|---|---|
| x | $c_2$ | $c_1$ | $c_2+c_1$ | $x+a_3$ | $c_2$ | $c_1$ | $c_2+c_1$ | $x+2a_3$ | $c_2$ | $c_1$ | $c_2+c_1$ |
| 1 | 0 | 1 | 1 | | | | | | | | |
| 2 | 0 | 2 | 2 | | | | | | | | |
| 3 | 0 | 3 | 3 | | | | | | | | |
| 4 | 0 | 4 | 4 | | | | | | | | |
| 5 | 0 | 5 | 5 | | | | | | | | |
| 6 | 0 | 6 | 6 | | | | | | | | |
| 7 | 0 | 7 | 7 | | | | | | | | |
| 8 | 0 | 8 | 8 | | | | | | | | |
| 9 | 0 | 9 | 9 | 42 | 3 | 0 | 3 | | | | |
| 10 | 0 | 10 | 10 | 43 | 3 | 1 | 4 | | | | |
| 11 | 0 | 11 | 11 | 44 | 3 | 2 | 5 | | | | |
| 12 | 0 | 12 | 12 | 45 | 3 | 3 | 6 | | | | |
| 13 | 0 | 13 | 13 | 46 | 3 | 4 | 7 | | | | |
| 14 | 1 | 0 | 1 | | | | | | | | |
| 15 | 1 | 1 | 2 | | | | | | | | |
| 16 | 1 | 2 | 3 | | | | | | | | |
| 17 | 1 | 3 | 4 | | | | | | | | |
| 18 | 1 | 4 | 5 | | | | | | | | |
| 19 | 1 | 5 | 6 | | | | | | | | |
| 20 | 1 | 6 | 7 | | | | | | | | |
| 21 | 1 | 7 | 8 | | | | | | | | |
| 22 | 1 | 8 | 9 | 55 | 3 | 13 | 16 | 88 | 6 | 4 | 10 |
| 23 | 1 | 9 | 10 | 56 | 4 | 0 | 4 | | | | |
| 24 | 1 | 10 | 11 | 57 | 4 | 1 | 5 | | | | |
| 25 | 1 | 11 | 12 | 58 | 4 | 2 | 6 | | | | |
| 26 | 1 | 12 | 13 | 59 | 4 | 3 | 7 | | | | |
| 27 | 1 | 13 | 14 | 60 | 4 | 4 | 8 | | | | |
| 28 | 2 | 0 | 2 | | | | | | | | |
| 29 | 2 | 1 | 3 | | | | | | | | |
| 30 | 2 | 2 | 4 | | | | | | | | |
| 31 | 2 | 3 | 5 | | | | | | | | |
| 32 | 2 | 4 | 6 | | | | | | | | |

This shows that all values 0<x<33 have valid generations of order <=2, with one value (x=22) requiring a generation of order 2: so conditions (SG1) and (SG2) for a stride generator SG(8,2) are satisfied.

Condition (SG3) requires that there exist $0<y<a_3$ such that:

$$y + ja_3 = c_2a_2 + c_1 \quad c_2+c_1<=n+j-1 \quad \text{is not soluble for any } j<=p+1$$

The following four groups of columns show that there is just one such value: y=22.

| must be: | | | <=7 | | | | <=8 | | | | <=9 | | | | <=10 |
|---|---|---|---|---|---|---|---|---|---|---|---|---|---|---|---|
| x | $c_2$ | $c_1$ | $c_2+c_1$ | $x+a_3$ | $c_2$ | $c_1$ | $c_2+c_1$ | $x+2a_3$ | $c_2$ | $c_1$ | $c_2+c_1$ | $x+3a_3$ | $c_2$ | $c_1$ | $c_2+c_1$ |
| 1 | 0 | 1 | 1 | | | | | | | | | | | | |
| 2 | 0 | 2 | 2 | | | | | | | | | | | | |
| 3 | 0 | 3 | 3 | | | | | | | | | | | | |
| 4 | 0 | 4 | 4 | | | | | | | | | | | | |
| 5 | 0 | 5 | 5 | | | | | | | | | | | | |
| 6 | 0 | 6 | 6 | | | | | | | | | | | | |
| 7 | 0 | 7 | 7 | | | | | | | | | | | | |
| 8 | 0 | 8 | 8 | 41 | 2 | 13 | 15 | 74 | 5 | 4 | 9 | | | | |
| 9 | 0 | 9 | 9 | 42 | 3 | 0 | 3 | | | | | | | | |
| 10 | 0 | 10 | 10 | 43 | 3 | 1 | 4 | | | | | | | | |
| 11 | 0 | 11 | 11 | 44 | 3 | 2 | 5 | | | | | | | | |
| 12 | 0 | 12 | 12 | 45 | 3 | 3 | 6 | | | | | | | | |
| 13 | 0 | 13 | 13 | 46 | 3 | 4 | 7 | | | | | | | | |
| 14 | 1 | 0 | 1 | | | | | | | | | | | | |
| 15 | 1 | 1 | 2 | | | | | | | | | | | | |
| 16 | 1 | 2 | 3 | | | | | | | | | | | | |
| 17 | 1 | 3 | 4 | | | | | | | | | | | | |
| 18 | 1 | 4 | 5 | | | | | | | | | | | | |
| 19 | 1 | 5 | 6 | | | | | | | | | | | | |
| 20 | 1 | 6 | 7 | | | | | | | | | | | | |
| 21 | 1 | 7 | 8 | 54 | 3 | 12 | 15 | 87 | 6 | 3 | 9 | | | | |
| 22 | 1 | 8 | 9 | 55 | 3 | 13 | 16 | 88 | 6 | 4 | 10 | 121 | 8 | 9 | 17 |
| 23 | 1 | 9 | 10 | 56 | 4 | 0 | 4 | | | | | | | | |
| 24 | 1 | 10 | 11 | 57 | 4 | 1 | 5 | | | | | | | | |
| 25 | 1 | 11 | 12 | 58 | 4 | 2 | 6 | | | | | | | | |
| 26 | 1 | 12 | 13 | 59 | 4 | 3 | 7 | | | | | | | | |
| 27 | 1 | 13 | 14 | 60 | 4 | 4 | 8 | | | | | | | | |
| 28 | 2 | 0 | 2 | | | | | | | | | | | | |
| 29 | 2 | 1 | 3 | | | | | | | | | | | | |
| 30 | 2 | 2 | 4 | | | | | | | | | | | | |
| 31 | 2 | 3 | 5 | | | | | | | | | | | | |
| 32 | 2 | 4 | 6 | | | | | | | | | | | | |

Thus SG(8,2) = {1,14,33} has a break at y=22; this is not a canonical break, because:

$$y + 4a_3 = 22 + 132 = 154 = 11a_2 + 0 \quad \text{and} \quad 11+0<=8+4-1$$

So the break order of y=22 is 4. ]

## 2.3 Every cover has an underlying stride generator

Why are stride generators of such interest? This section shows the close relationship between stride generators and covers which provides the key to the overall proof: we show that every cover C(A,3,s) has an underlying stride generator SG(s-k,p) such that:

$$C(A,3,s) = (k+1)a_3 + Y$$

where y = Y-1 is one of the breaks in SG(s-k,p).

We also show that in the case that SG(s-k,p) is canonical, y is the first break in the stride generator. This is significant, because it shows that Y in the formula above is a property of the stride generator alone, and is independent of s.

It turns out to be much more difficult to prove this result for non-canonical stride generators as well: sections 2.4 to 2.9 prepare the ground work, and the result is finally proved as Theorem 231 in section 2.10.

Definition 203

A set A = {1, $a_2$, $a_3$} has a "non-trivial" cover iff C(A,3,s) >= $a_3$.

We first prove that the formula for any non-trivial cover can be expressed in a certain way.

Theorem 201

Let A = {1, $a_2$, $a_3$} be a set with a non-trivial cover C(A,3,s); then we can write:

$$C(A,3,s) = (k+1)a_3 + Y \quad \text{where} \quad 0 <= Y < a_3-1 \quad \text{and} \quad 0 <= k < s$$

[Note that the case Y=$a_3$-1 is specifically excluded: this is important for the proof of the following theorem.]

Proof

a) k>=0 because the cover is non-trivial.

k<s because it is not possible to generate the value (s+1)$a_3$ using no more than s stamps.

b) We define (k+1) to be the integer result of dividing C(A,3,s) by $a_3$, and so Y is the remainder. By definition, 0 <= Y <= $a_3$-1; so we have only to show that Y=$a_3$-1 is not possible.

Suppose Y=$a_3$-1.

By definition, the value C(A,3,s)+1 has no generation:

$$C(A,3,s)+1 = (k+1)a_3 + a_3-1 + 1 = (k+2)a_3$$

This means that s<k+2, for otherwise there would be a trivial generation.

But the value (k+1)$a_3$ can be generated, and since s<=k+1 this can only be generated as k+1 $a_3$ stamps; therefore s=k+1.

But $(k+1)a_3+1$ can also be generated, say as:

$$(k+1)a_3+1 = c_3a_3 + c_2a_2 + c_1 \quad c_3+c_2+c_1<=s=k+1$$

$c_3=k+1$ is not a solution, so $c_3<k+1$

So $(k+1-c_3)a_3+1 = c_2a_2 + c_1 \quad c_2+c_1 <= k+1-c_3$ for some $c_3<k+1$

But $c_2a_2 + c_1 < (c_2+c_1)a_2$ - since $a_2>1$
$<= (k+1-c_3)a_2$
$< (k+1-c_3)a_3$ - since $a_3>a_2$

So $(k+1)a_3+1$ cannot be generated if s=k+1.

This contradicts the assumption that C(A,3,s) is a cover, and so Y cannot be equal to $a_3$-1, and the theorem is proved.

[ Note that Y=0 is possible:

C({1,3,4},3,3) = 12 = $3a_3+0$ ]

Next we prove that every non-trivial cover defines an underlying stride generator, and relate the value of the cover to a certain break in the stride generator.

Theorem 202

Let A = {1, $a_2$, $a_3$} be a set with a non-trivial cover C(A,3,s).

Then A is also an (s-k)-stride generator SG(s-k,p) of order p<=k for some 0<=k<s, and:

$$C(A,3,s) = (k+1)a_3 + Y$$

where y=Y+1 is the first break in SG(s-k,p) with break order > k+1.

[Note that any canonical break is considered as having an "infinite" break order, and so will always meet this condition.]

Definition 204

We say that SG(s-k,p) is the stride generator "underlying" the cover C(A,3,s).

Proof

We can write:

$$C(A,3,s) = (k+1)a_3 + Y \qquad 0<=Y<a_3-1, \quad 0<=k<s \qquad \text{by Theorem 201}$$

and, writing y=Y+1, we have:

$$C(A,3,s) = (k+1)a_3 + y-1 \qquad 0<y<a_3, \quad 0<=k<s$$

so we now have only to show that A = SG(s-k,p) for some p<=k, and that y is a break in A with break order > k+1.

Consider any $0<x<a_3$; then $x+ka_3$ is a value in the final stride and so has a generation:

$$x + ka_3 = c_3a_3 + c_2a_2 + c_1a_1 \quad c_3+c_2+c_1<=s$$

$$=> x + (k-c_3)a_3 = c_2a_2 + c_1 \quad c_2+c_1<=s-c_3$$

Writing $i = k-c_3$ we have:

$$x + ia_3 = c_2a_2 + c_1 \quad c_2+c_1<=(s-k)+i \quad 0<=i<=k, \text{ since } 0<=c_3<=k$$

Let $p<=k$ be the maximum value of i needed for any x; then:

For all $0<x<a_3$:

$$x + ia_3 = c_2a_2 + c_1 \quad c_2+c_1<=(s-k)+i \text{ is soluble for some } i<=p$$

and there exists $0<x<a_3$ such that this equation is soluble for i=p but not for any i<p.

Thus conditions (SG1) and (SG2) for A = SG(s-k,p) are met, with p<=k.

Now consider y = Y+1; then $y + (k+1)a_3$ is the first value that cannot be generated by A, and so:

There exists $0<y<a_3$ (since $0<=Y<a_3-1$) such that no $c_3$, $c_2$, $c_1$ exist such that:

$$y + (k+1)a_3 = c_3a_3 + c_2a_2 + c_1 \quad c_3+c_2+c_1<=s$$

$$=> \quad y + (k+1-c_3)a_3 = c_2a_2 + c_1 \quad c_2+c_1<=s-c_3$$

Writing $j = k+1-c_3$, we see that there is no j<=k+1 such that:

$$y + ja_3 = c_2a_2 + c_1 \quad c_2+c_1<=(s-k)+j-1$$

Since p<=k, this shows that y satisfies (SG3) which requires j<=p+1, and so is a break with break order > k+1 in A = SG(s-k,p).

We have now proved that A is, indeed, a stride generator SG(s-k,p), and that y=Y+1 is one of its breaks with break order > k+1; we now show that y must be the first such break.

Suppose y is any break in SG(s-k,p) with break order > k+1.

By (SG3) and definition 201, we have:

$$y + ja_3 = c_2a_2 + c_1 \quad c_2+c_1<=(s-k)+j-1 \text{ is not soluble for } j<=k+1$$

Writing $j = k+1-c_3$ we have:

$$y + (k+1-c_3)a_3 = c_2a_2 + c_1 \quad c_2+c_1<=s-c_3 \text{ is not soluble for } c_3<=k+1$$

$$=> y + (k+1)a_3 = c_3a_3 + c_2a_2 + c_1 \quad c_3+c_2+c_1<=s \text{ is not soluble for } c_3<=k+1$$

and so y cannot be generated.

The cover is determined by the smallest value that cannot be generated, and so y must be the first break in SG(s-k,p) with break order > k+1.

Finally, we show that in the case of a canonical underlying stride generator, it is always the first break that counts.

Theorem 203

If the stride generator SG(s-k,p) underlying the cover C(A,3,s) is canonical, then:

$$C(A,3,s) = (k+1)a_3 + Y$$

where y = Y+1 is the first break in SG(s-k,p).

[This means that the value Y is independent of s - it is a property of the canonical stride generator alone.]

Proof

All breaks in a canonical stride generator are canonical by Definition 202, and so the result follows directly from Theorem 202.

This is as far as we can go without a lot more ground work: it is only in section 2.10 with Theorem 231 that we are able to extend Theorem 203 above to cover the case of a non-canonical underlying stride generator.

[Conjecture:

All underlying stride generators are canonical.

This appears to be the case, but I have so far failed to prove it.

Note that a simple proof of this result would remove the need for much of what follows!]

[Example:

It is easy to show that:

C({1,6,13},3,6) = 47 = 3.13 + 8

ie C({1,6,13},3,6) = (k+1).13 + y-1  where k=2 and y=9

Thus the stride generator underlying this cover must be SG(4,p) for some p<=2 with a break at y=9 with break order > 3.

In fact, a previous example has shown that {1,6,13} = SG(4,2) is a canonical stride generator with breaks at y=9 and y=10. ]

## 2.4 Threads and thread diagrams

In this section, we introduce the idea of a stride generator being made up of a number of "threads" which can be displayed in the form of a two dimensional diagram.

Later sections make use of the thread concept and of the diagrams both to illustrate results and within some of the proofs themselves.

[To help my early investigations I developed an interactive computer program to display diagrams of this kind which allows me to see how the threads move as the parameters $a_2$, $a_3$, n and p are altered. The program also indicates whether a given set of parameters define a valid stride generator or not and, if so, where the breaks occur.

Many of the insights behind this proof were the result of experiments with this program.]

An example of such a "thread diagram" is given in Figure 200 for the stride generator SG(8,2) = {1,14,33}; this should be compared with the tables given in the examples at the end of section 2.2.

The x-axis represents values to be generated, and normally includes the range 0 to $a_3$; the y-axis represents the orders of generation, and usually runs from 0 to p.

Suppose that:

$$x_1 + ia_3 = c_2a_2 \qquad c_2 <= n+i$$

and $$x_2 + ia_3 = c_2a_2 + c_1 \qquad c_2+c_1 = n+i$$

This means that all values from $x_1$ to $x_2$ inclusive can be generated with an order i generation, and so we draw a line at y=i from x=$x_1$ to x=$x_2$ inclusive, and label it $c_2$.

For example:

18 + 2.33 = 6.14      6 < 8+2
22 + 2.33 = 6.14 + 4      10 = 8+2

and so a line labelled 6 runs from x=18 to x=22 at y=2.

Such a line is called a "thread", and the completed diagram is called a "thread diagram".

It is easy to see from such a diagram whether a particular value can be generated and, if so, in what ways, by simply looking to see which threads "cover" the value:

For example we see that 15 has two possible generations of order 0 and 1:

15 + 0.33 = 1.14 + 1      1+1 <= 8+0
15 + 1.33 = 3.14 + 6      3+6 <= 8+1

whereas 12 has a single generation of order 1:

12 + 1.33 = 3.14 + 3      3+3 <= 8+1

We can also easily check that all values $0<x<a_3$ (33) have at least one generation by simply making sure that the threads of order $0<=i<=p$ (2) together cover the whole range.

If the thread diagram's scope is extended to include values $a_3<x<2a_3$ and threads of order -1 (which start at $x=a_3$), then we can also identify any breaks in the stride generator: they will be indicated by values which are not covered by any thread of order $-1<=j<=p$.

Figure 201 is an example, where the value 55 has been shaded to show that it is the only value in the range $0<x<2a_3$ that is not covered by any thread; therefore $y=55-a_3 = 22$ is the only break in SG(8,2) = {1,14,33}.

Another substantial example of a thread diagram is Figure 202, which shows the breaks at y=53 and y=59 in the stride generator SG(12,3) = {1,30,82}.

Definition 205

The notation for threads is defined formally as follows:

$T_i(k)$ is a thread k of order i and runs from:

$$ST_i(k) = ka_2 - ia_3$$

$$\text{to } ET_i(k) = ka_2 - ia_3 + (n+i-k)$$

$$\text{and has length } L_i(k) = (n+i-k)+1$$

[ For example, $T_2(6)$ runs from 18 to 22 inclusive:

$ST_2(6) = 18$
$ET_2(6) = 22$
and $L_2(6) = 5$ ]

We say that thread $T_i(k)$ "covers" a value x if $ST_i(k)<=x<=ET_i(k)$, and that $T_i(k)$ "crosses" a value x if $ST_i(k)<=x<ET_i(k)$.

In following sections we may describe a thread of order i as an "i-thread", and will sometimes refer loosely to "a thread $T_i$ starting at $ST_i$" when there is no ambiguity as to which i-thread is being referenced.

{1,14,33} = SG(8,2)  Break at 22

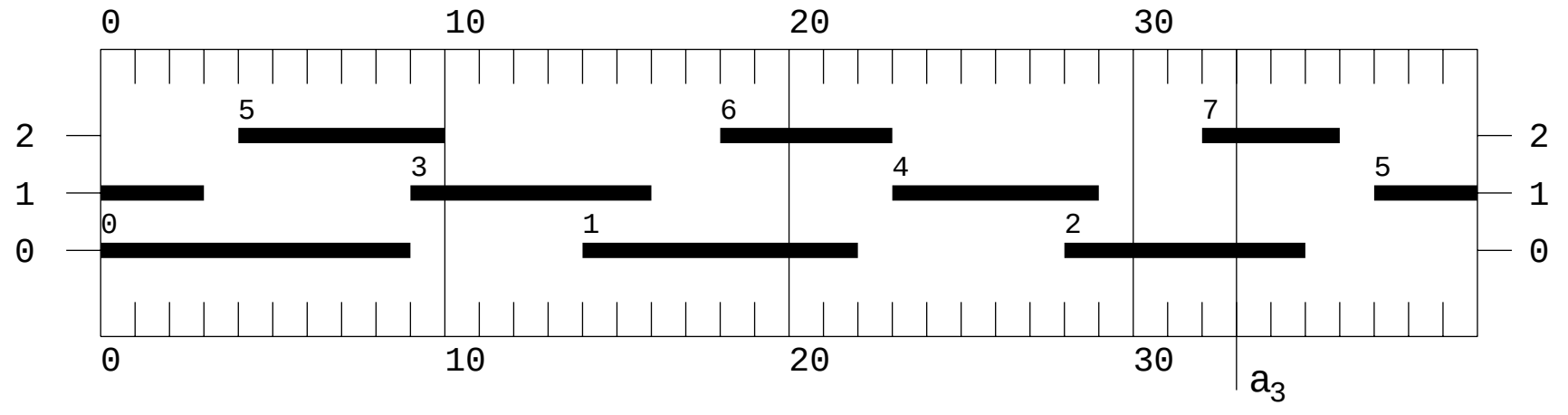

Figure 200

{1,14,33} = SG(8,2)  Break at 22

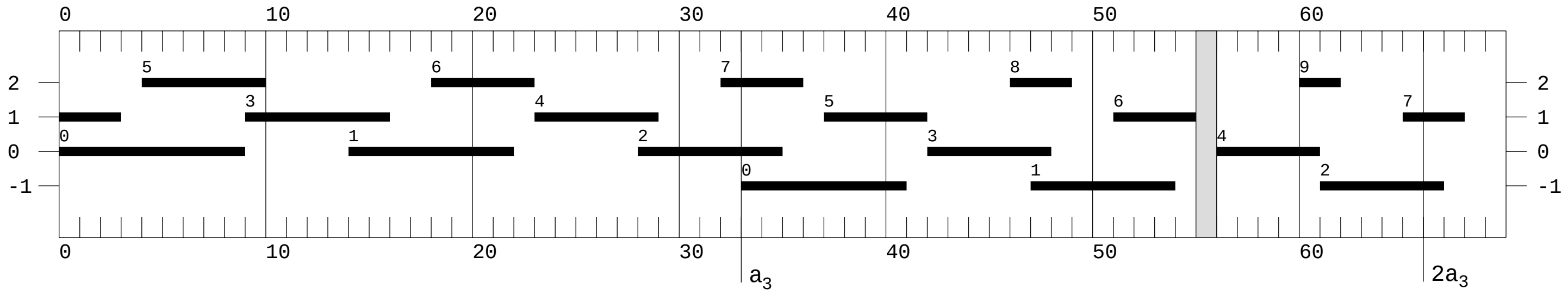

Figure 201

{1,30,82} = SG(12,3)  Breaks at 53, 59

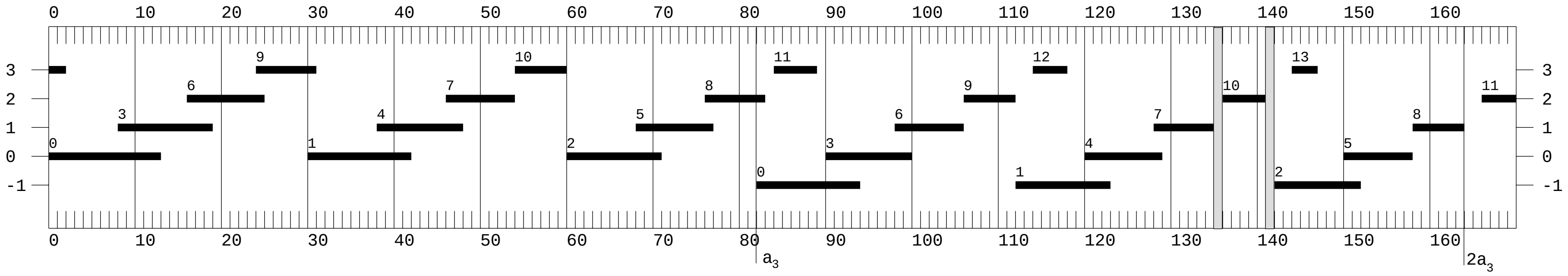

Figure 202

## 2.5 Stride generators as collections of threads

In this section we further characterise stride generators in terms of their thread diagrams.

We show that any break y must lie at the end of a thread, and that (unless $y=a_3-1$) y+1 is at the start of a thread - and that no thread of order less than or equal to p+1 can cross this boundary. This is illustrated in Figure 203.

We first show that any break y lies at the end of a thread; this can be viewed as an alternative but equivalent definition of a break in terms of threads.

Theorem 204

If $A = \{1, a_2, a_3\}$ is a stride generator SG(n,p), then y is a break in A iff:

a) y is at the end of an i-thread for some $i<=p$

i.e. $y + ia_3 = c_2a_2 + c_1$  $c_2+c_1=n+i$  is soluble for some $i<=p$

and b) y is not crossed by any thread of order $<=p+1$

i.e. $y + ia_3 = c_2a_2 + c_1$  $c_2+c_1<n+i$  is not soluble for any $i<=p+1$

Proof

Condition (b) =>

$y + ja_3 = c_2a_2 + c_1$  $c_2+c_1<=n+j-1$  is not soluble for any $j<=p+1$

which is (SG3); so if A is a stride generator, then (b) => y is a break.

On the other hand, (SG1) =>

$y + ia_3 = c_2a_2 + c_1$  $c_2+c_1<=n+i$  is soluble for some $i<=p$

and (SG3) for y to be a break =>

$y + ia_3 = c_2a_2 + c_1$  $c_2+c_1<=n+i-1$  is not soluble for any $i<=p+1$

These two conditions can only both be true if:

$y + ia_3 = c_2a_2 + c_1$  $c_2+c_1=n+i$  is soluble for some $i<=p$
and $y + ia_3 = c_2a_2 + c_1$  $c_2+c_1<n+i$  is not soluble for any $i<=p+1$

which shows that y is a break => conditions (a) and (b).

Note that Theorem 204 together with (SG2) shows that in any stride generator SG(n,p) n and p are "mutually minimal":

a) Given n, (SG2) requires that p be as small as possible - so p is minimal with respect to n.

b) Given p, Theorem 204 shows that (SG3) requires there to be a value y that is covered by the end of a thread: so that if n were to be made any smaller, y would no longer be covered by any thread of order <=p at all. So n is minimal with respect to p.

Next we show that for any break $y<a_3-1$, y+1 lies at the start of a thread.

Theorem 205

If A = {1, $a_2$, $a_3$} is a stride generator SG(n,p) and $y<a_3-1$ is a break in A, then (y+1) is at the start of an i-thread for some i<=p.

i.e. $(y+1) + ia_3 = c_2a_2$ $c_2<=n+i$ is soluble for some i<=p

[Note: The case of $y=a_3-1$ - i.e. a break right at the end of the stride generator - is treated separately as a special case when characterising stride generators in sections 2.8 and 2.9.

An alternative approach might be to weaken this theorem by allowing -1<=i<=p which would then cover the $y=a_3-1$ case as well, since $(y+1) - a_3 = 0 = 0.a_2$ ]

Proof

(SG1) for SG(n,p) =>

$(y+1) + ia_3 = c_2a_2 + c_1$ $c_2+c_1<=n+i$ is soluble for some i<=p

Suppose $c_1>0$; we can then write $c_1' = c_1-1 >= 0$:

=> $y + ia_3 = c_2a_2 + c_1'$ $c_2+c_1'<n+i$ is soluble for some i<=p

which is contrary to (SG3).

Therefore $c_1 = 0$ and the theorem is proved.

The above two theorems together show that any break $y<a_3-1$ will appear as a "step" in the thread diagram as shown in Figure 203. See, for example, Figure 200, where the break y=22 lies at the end of $T_2(6)$ which is contiguous with $T_1(4)$.

We now provide an equivalent definition of break order in terms of threads.

Theorem 206

Let y be a break in a stride generator SG(n,p).

Then the break order of y is given by the order of the first thread, if any, which crosses y (i.e. covers both y and y+1).

If no such thread exists, the break is canonical.

[In other words, the first thread, if any, that "crosses the boundary" determines the break order; this is illustrated in Figure 204.]

Proof

We first show that if y has break order q then y and y+1 are both covered by the same q-thread:

y has break order q

=> $y + qa_3 = c_2a_2 + c_1 \quad c_2+c_1<=n+q-1$ is soluble

=> $y + qa_3 = c_2a_2 + c_1 \quad c_2+c_1<=n+q$
and $(y+1) + qa_3 = c_2a_2 + c_1+1 \quad c_2+c_1+1<=n+q$ are both soluble

=> y and y+1 are both covered by the thread $T_q(c_2)$

Next we show that if y and y+1 are both covered by the same r-thread, then the break order q of y must be <=r:

y and y+1 are both covered by a thread $T_r(c_2)$

=> $y + ra_3 = c_2a_2 + c_1 \quad c_2+c_1<=n+r$
and $(y+1) + ra_3 = c_2a_2 + c_1+1 \quad c_2+c_1+1<=n+r$ are both soluble

=> $y + ra_3 = c_2a_2 + c_1 \quad c_2+c_1<=n+r-1$ is soluble

=> y has break order q<=r

Therefore the break order of any non-canonical break y must be equal to the order of the first thread that covers both y and y+1.

"No thread exists <=> the break is canonical" follows as a direct corollory.

Finally, we give an equivalent definition of a stride generator in terms of the corresponding thread diagram.

Theorem 207

SG(n,p) is a stride generator iff its thread diagram shows that:

a) Every value $0<x<a_3$ is covered by at least one thread of order <=p.

b) At least one value $0<x<a_3$ is covered only by a thread of order p.

c) At least one value $0<y<a_3$ lies at the end of a thread of order <=p and is not crossed by any thread of order <=p+1.

Proof

(a) and (b) are conditions equivalent to (SG1) and (SG2).

(c) follows from Theorem 204.

[ Examples

We know that breaks occur at the end of the first of two contiguous threads, provided the boundary is not crossed by another thread of order <=p+1.

A straightforward example is shown in Figure 200:

$T_2(6)$ and $T_1(4)$ are contiguous threads corresponding to the break at y=22.

However, there is no break that corresponds to the contiguous threads $T_0(0)$ and $T_1(3)$ because their boundary is covered by the thread $T_2(5)$.

The following somewhat pathological example again shows how breaks appear as steps in a thread diagram, and how the corresponding break orders are determined by higher order "boundary-crossing" threads.

Figure 205 is a thread diagram for SG(8,3) = {1,30,38} that has been extended to include threads of order > p=3; these additional threads are represented as shaded instead of solid lines to indicate that they do not belong to the stride generator itself.

We now analyse each pair of contiguous threads (ie each "step") that appears in the stride generator:

$T_3(4)$, $T_2(3)$ - This step is first crossed by $T_6(8)$. Since 6>p+1, the step defines a break y=13 with break order 6.

$T_2(3)$, $T_1(2)$ - This step similarly defines a break y=21 with break order 5.

$T_1(2)$, $T_0(1)$ - This step is first crossed by $T_4(6)$. Since 4=p+1, this step does not define a break at all. ]

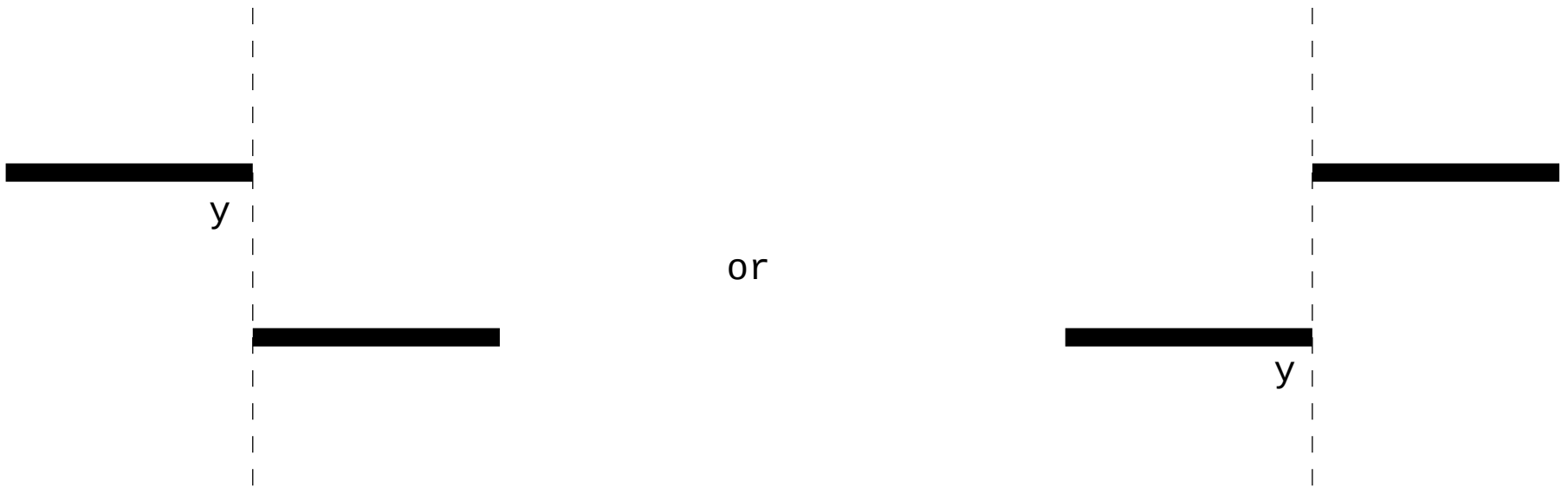

Figure 203

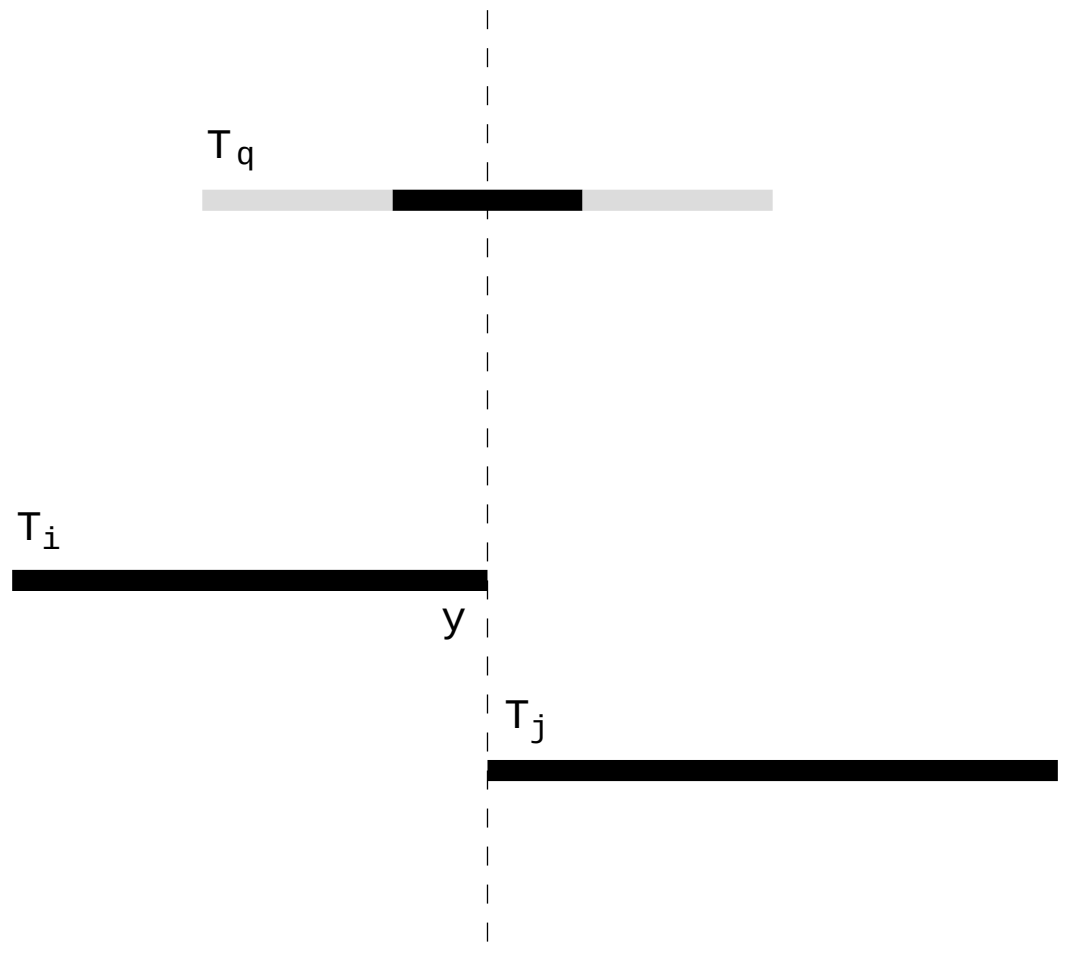

Figure 204

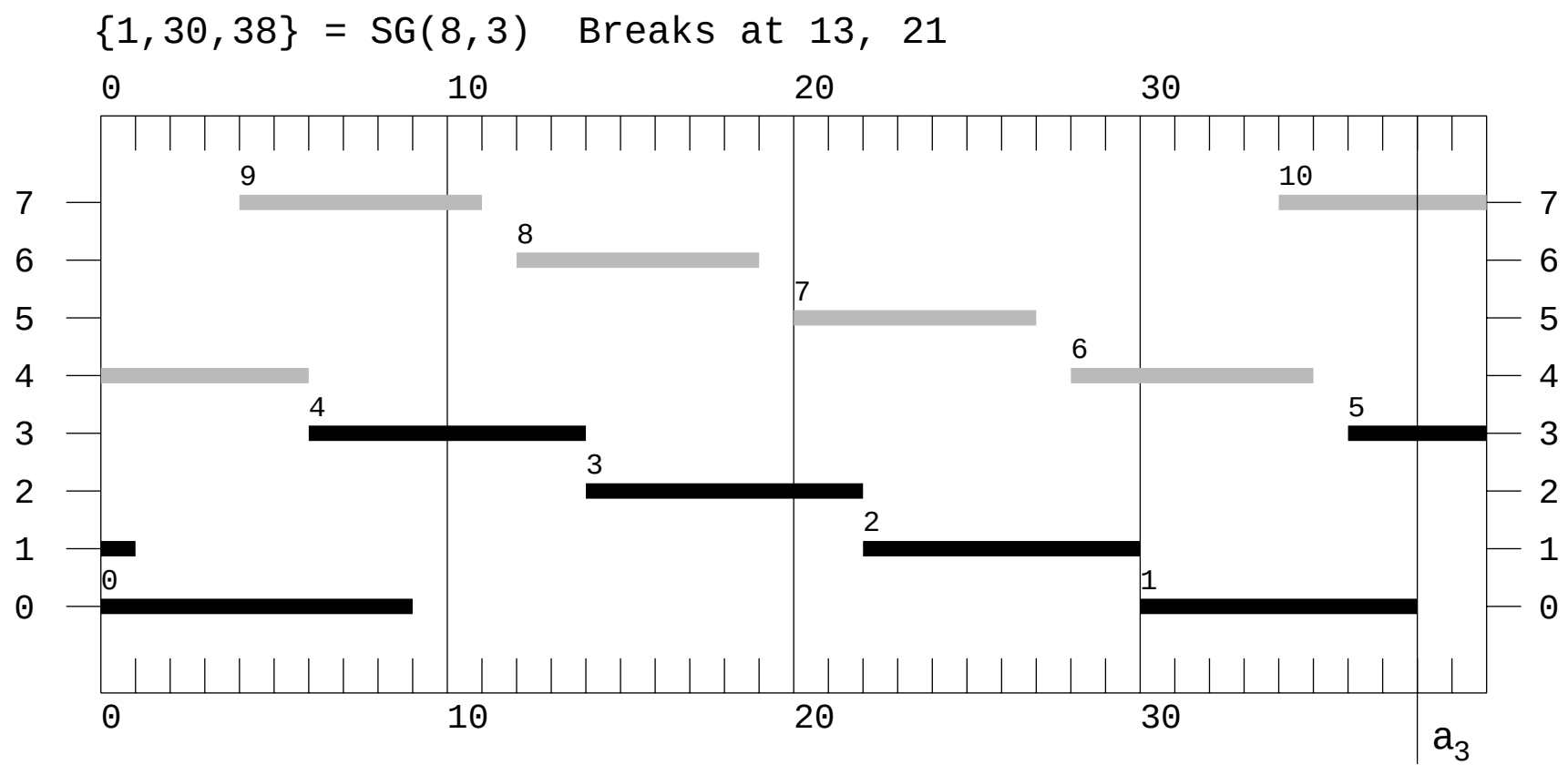

Figure 205

## 2.6 Thread series

When we look at a thread diagram (such as that shown in Figure 202, for example) it is immediately apparent that the threads group together in patterns or "series". In this section we identify three such series and provide formal proofs which show how the threads in each series are related to each other.

Given a thread $T_i$ at one position in a stride generator, we can use one of the series theorems to predict the existence of a related thread $T_j$ at some other position. This is a powerful technique which we use in later sections to help determine the properties of stride generators.

However, the series are not all infinite, and so we must take care when deriving $T_j$ from $T_i$ to show that $T_i$ is not the last member of the series.

Fortunately, we are able to show in a later section (Theorem 224) that all possible positions for threads in a stride generator are filled. This means that if $T_i$ and $T_j$ belong to the same series and are both part of the same stride generator then we know that they both must exist: there is no need to prove explicitly that any constraints required by the series theorems (such as L>1 or $c_2$>=$C_2$) are met.

The first series consists of all threads of the same order:

$T_i(0)$, $T_i(1)$, ... $T_i(K)$

Figure 206 shows all* threads associated with SG(19,4) = {1,34,148} within the range 0 to 199, and all the shaded threads $T_2(9)$ ... $T_2(14)$ in the top section are an example of such a series.

[*Note that threads of order >4 are included in the diagram although these are not part of the stride generator itself.]

The start points of successive threads in the series differ by $a_2$, and each thread is one shorter than its predecessor; this derivation rule is shown in the top section of Figure 207.

The final thread $T_i(K)$ in the series is of length 1 and satisfies:

$K = n+i$
$ST_i(n+i) = ET_i(n+i) = (n+i)a_2 - ia_3$

Theorem 208

Let $T_i(c_2)$ be a thread of length L starting at x.

Then $T_i(c_2+1)$ is a thread of length L-1 starting at $x+a_2$, provided L>1.

Proof

$T_i(c_2)$ is a thread of length L starting at x

=> $x + ia_3 = c_2a_2$    $c_2<=n+i$    $L=n+i-c_2+1$

=> $(x+a_2) + ia_3 = (c_2+1)a_2$

This is a thread of order i of length $n+i-(c_2+1)+1 = L-1$ provided L>1, since:

$L>1$

$\Rightarrow n+i-c_2+1>1$

$\Rightarrow c_2<n+i$

$\Rightarrow (c_2+1)<=n+i$

Corollory

Let $T_i(c_2)$ be a thread of length L starting at x.

Then $T_i(c_2-1)$ is a thread of length L+1 starting at $x-a_2$ provided $c_2>0$.

Proof

Similar to above.

The second series consists of threads of decreasing order whose start points are separated by $a_3$ mod $a_2$ - that is, by the remainder obtained when $a_3$ is divided by $a_2$.

To be precise, let $a_3 = C_2a_2 + C_1$; then the series is:

$T_i(n+i), T_{i-1}(n+i-C_2), \ldots T_{i-k}(n+i-kC_2)$

Each thread is $(C_2-1)$ longer than the preceding thread, and the series terminates at $T_{i-k}$ when:

a) $i-k = 0$

or b) $n+i-(k+1)C_2 < 0$

The shaded threads in the middle section of Figure 206 are an example of such a series, and the corresponding derivation rule is shown in the middle of Figure 207.

Theorem 209

Let $T_i(c_2)$ be a thread of length L starting at x, and let $a_3 = C_2a_2 + C_1$.

Then $T_{i-1}(c_2-C_2)$ is a thread of length $L+(C_2-1)$ starting at $x+C_1$, provided that $i>0$ and $c_2>=C_2$.

Proof

$T_i(c_2)$ is a thread of length L starting at x

$\Rightarrow x + ia_3 = c_2a_2 \quad c_2<=n+i \quad L=n+i-c_2+1$

$\Rightarrow x + (i-1)a_3 = c_2a_2 - C_2a_2 - C_1$

$\Rightarrow (x+C_1) + (i-1)a_3 = (c_2-C_2)a_2$

This is a thread of order i-1 with length $n+i-1-c_2+C_2+1 = L+(C_2-1)$ provided $i>0$ and $c_2>=C_2$ because:

$c_2<=n+i$

$\Rightarrow c_2-C_2<=n+(i-1)$ since $a_3>a_2$ and so $C_2>=1$

Corollory

Let $T_i(c_2)$ be a thread of length L starting at x, and let $a_3 = C_2a_2 + C_1$.

Then $T_{i+1}(c_2+C_2)$ is a thread of length $L-(C_2-1)$ starting at $x-C_1$, provided that $L>(C_2-1)$.

Proof

Similar to above.

[There is another series analogous to the second series in which the threads are of increasing order and separated by $a_2-C_1$:

$T_i(c_2)$ at x of length L <=> $T_{i+1}(c_2+C_2+1)$ at $x+(a_2-C_1)$ of length $L-C_2$

This is easily derived by first taking a step back along the second series followed by a step forwards along the first:

$T_i(c_2)$ at x of length L

-> $T_{i+1}(c_2+C_2)$ at $x-C_1$ of length $L-(C_2-1)$ - by series 2

-> $T_{i+1}(c_2+C_2+1)$ at $x+a_2-C_1$ of length $L-C_2$ - by series 1]

The third series consists of threads of decreasing order whose start points are separated by $a_3-a_2$; all threads in the series have the same length.

... $T_i(m)$, $T_{i-1}(m-1)$, ... $T_{i-k}(m-k)$

The series has no start, and finishes when either the order reaches zero (i=k) or when k>m.

The bottom sections of Figures 206 and 207 illustrate an example of part of the series and the corresponding derivation rule respectively.

This series can, in fact, be derived from the first two - in other words, it is not an independent series - but it is useful in its own right and easy to prove directly.

Theorem 210

Let $T_i(c_2)$ be a thread of length L starting at x.

Then $T_{i-1}(c_2-1)$ is a thread of length L starting at $x+(a_3-a_2)$ provided that $i>0$ and $c_2>0$.

Proof

$T_i(c_2)$ is a thread of length L starting at x

=> $x + ia_3 = c_2a_2$ $\quad c_2<=n+i \quad L=n+i-c_2+1$

=> $x+(a_3-a_2) + ia_3 = c_2a_2 + (a_3-a_2)$

=> $x+(a_3-a_2) + (i-1)a_3 = (c_2-1)a_2$

This is a thread of length $n+(i-1)-(c_2-1)+1 = n+i-c_2+1 = L$, provided $i>0$ and $c_2>0$, since:

$c_2 <= n+i$

$=> (c_2-1) <= n+(i-1)$

Corollory

Let $T_i(c_2)$ be a thread of length L starting at x.

Then $T_{i+1}(c_2+1)$ is a thread of length L starting at $x-(a_3-a_2)$.

Proof

Similar to above.

SG(19,4) = {1,34,148}

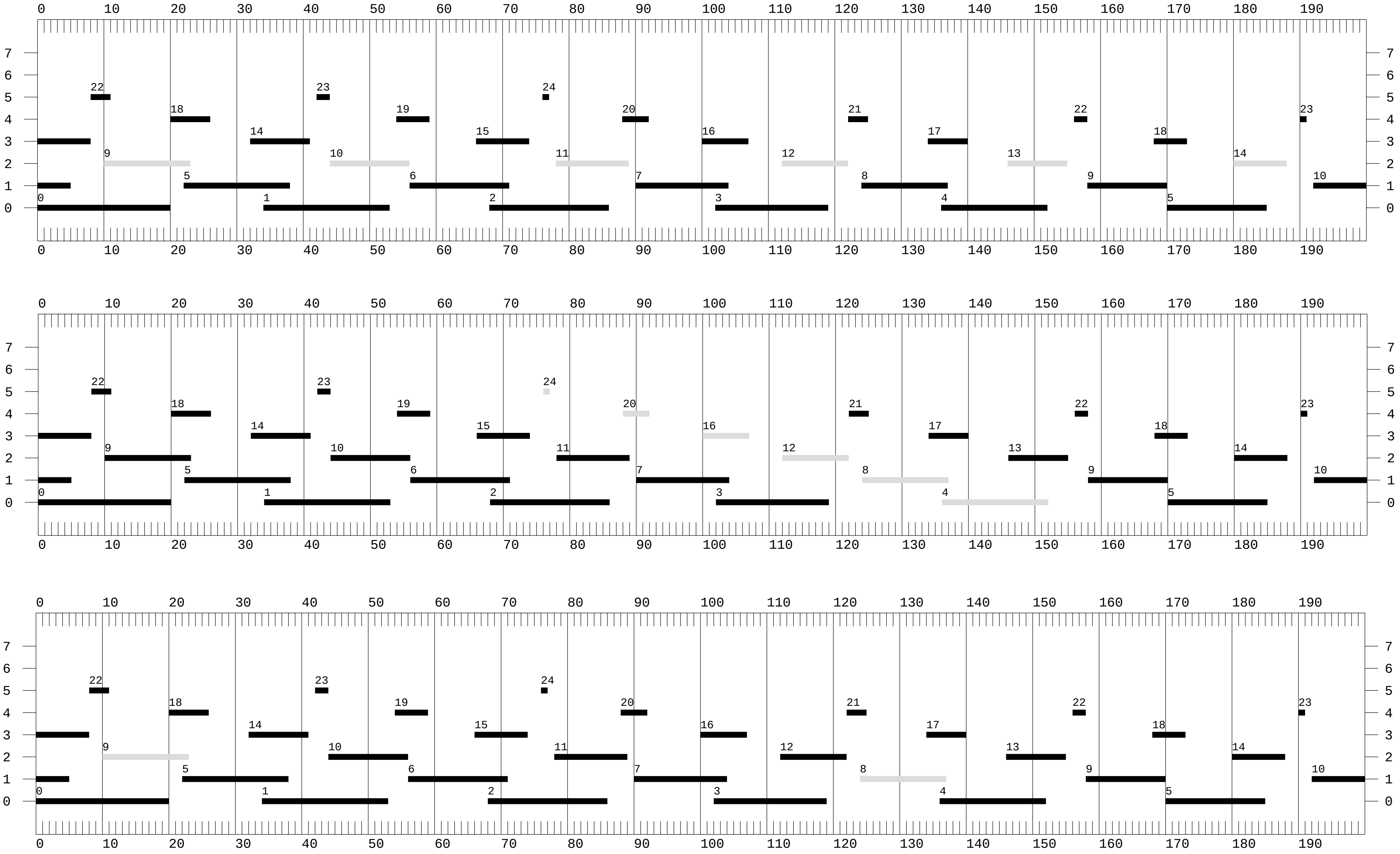

Figure 206

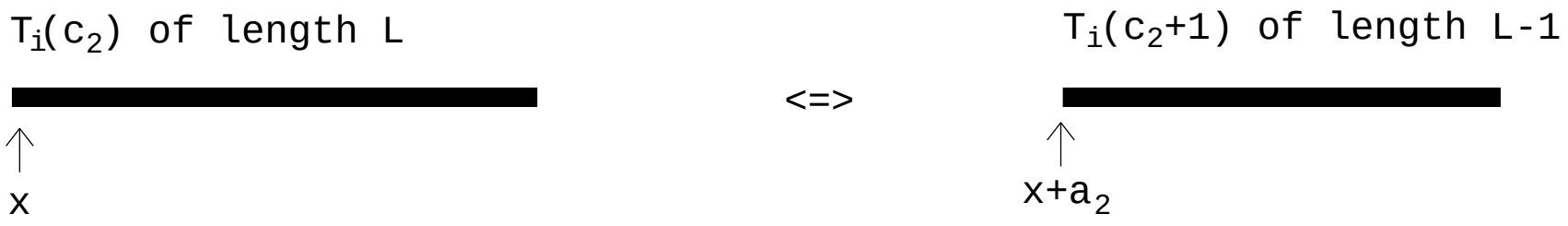

Theorem 208

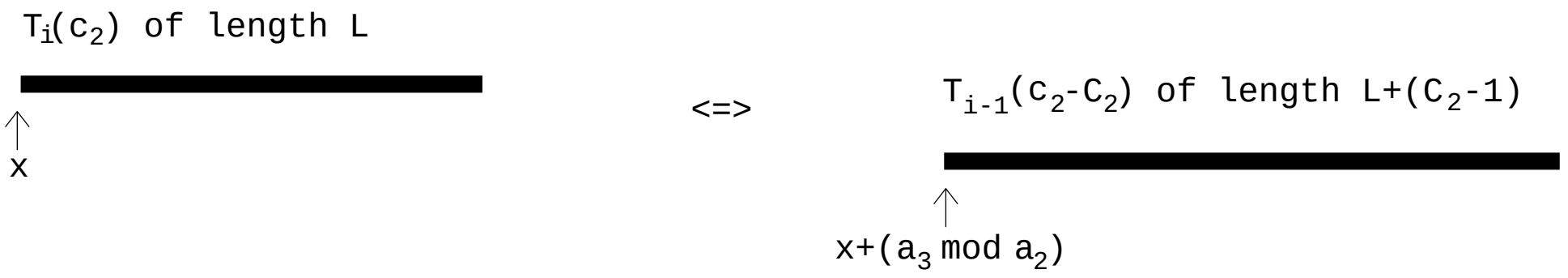

Theorem 209

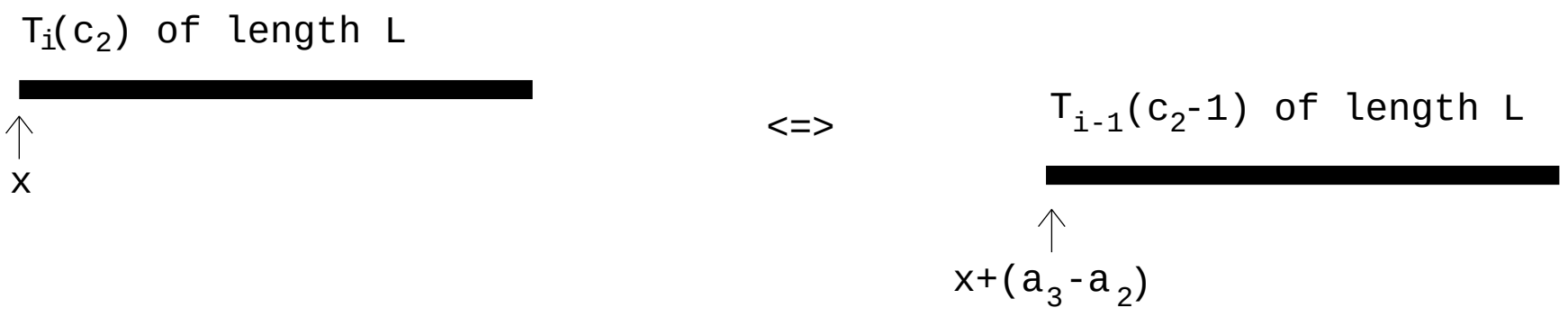

Theorem 210

Figure 207

## 2.7 Covered threads

One thread can be completely covered by another thread - in which case it is in some sense redundant: that is, all values generated by the first thread can also be generated by the covering thread.

We show in this section that such coverings cannot exist amongst the threads of a stride generator. This important result helps prove that certain stride generators are canonical, as well as assisting in the characterisation of stride generators.

### Definition 206

A thread $T_i(b_2)$ is said to be "covered" by another thread of different order $T_j(c_2)$ iff:

$$ST_i(b_2) >= ST_j(c_2) \quad \text{and} \quad ET_i(b_2) <= ET_j(c_2)$$

[This is illustrated in Figure 208.]

It is intuitively obvious that any thread $T_i$ that is above and to the right of a thread $T_j$ must be shorter than $T_j$, and so an i-thread can only be covered by a j-thread if $j<i$. The proof follows for completeness.

### Theorem 211

$T_i(b_2)$ and $T_j(c_2)$ are two threads of lengths $L_i$ and $L_j$ respectively such that $i>j$ and $ST_i(b_2)>=ST_j(c_2)$.

Then $L_i<L_j$, and $b_2>c_2+1$.

Proof

Write $x = ST_j(c_2)$

Then $x + ja_3 = c_2a_2$ with $L_j=n+j-c_2+1$

and the i-thread of the same length $L_j$ starts at y where:

$$y + ia_3 = c_2'a_2$$

where $n+i-c_2'+1 = n+j-c_2+1$

$\Rightarrow c_2' = (i-j)+c_2$

So $(x-y) + (j-i)a_3 = (c_2-c_2')a_2 = (j-i)a_2$

$\Rightarrow (x-y) = (i-j)(a_3-a_2)$

Now $a_3>a_2$ and $i>j$, so $x>y$.

So the i-thread of equal length to the j-thread starts before the j-thread, and so any i-thread that starts at or beyond the j-thread must be shorter (by Theorem 208).

Now $L_i = n+i-b_2+1$ and $L_j = n+j-c_2+1$, so:

$$L_i<L_j \Rightarrow n+i-b_2+1 < n+j-c_2+1$$
$$\Rightarrow i-b_2 < j-c_2$$
$$\Rightarrow b_2 > c_2+(i-j)$$
$$\Rightarrow b_2 > c_2+1$$

We now show that once one i-thread is covered by a j-thread, then all k-threads for k>=i are covered by other threads of order <i.

Theorem 212

Suppose thread $T_i$ is covered by another thread $T_j$.

Then any thread $T_k$ for k>=i is covered by some thread $T_m$ for m<i.

Proof

Consider $T_i$ and $T_j$ as members of two series of the same kind, and let $T_{i'}$ and $T_{j'}$ be the next members of each series.

Then: $ST_{i'} = ST_i + K$
$ST_{j'} = ST_j + K$

and: $L_{i'} = L_i + L$
$L_{j'} = L_j + L$

where K, L are order-independent constants according to the kind of series [e.g. $K=a_2$, L=-1 for the first kind of series].

Thus if $T_j$ covers $T_i$, then $T_{j'}$ covers $T_{i'}$, and hence by induction each member of the i-series is covered by a corresponding member of the j-series.

Using the first series (Theorem 208) we have:

If $T_k(c_2)$ is covered by $T_m(b_2)$ then $T_k(c_2+1)$* is covered by $T_m(b_2+1)$ and $T_k(c_2-1)$* is covered by $T_m(b_2-1)$; so if one thread of order k is covered by some thread of order m, then every thread of order k is covered by some corresponding thread of order m.

[* if it exists]

Using the third series (Theorem 210) we have:

If $T_k(c_2)$ is covered by $T_m(b_2)$ then $T_{k+1}(c_2+1)$ will be covered by $T_{m+1}(b_2+1)$; so if threads of order k are covered by threads of order m, then threads of order k+1 are covered by threads of order m+1.

So by induction from "$T_i$ is covered by $T_j$" we see first that "i-threads are covered by j-threads" and thence "k-threads are covered by (k-(i-j))-threads" for all k>=i.

But if k-(i-j)>=i, then (k-(i-j))-threads will themselves be covered by (k-2(i-j))-threads and so on, until threads of order <i are finally reached.

So "$T_i$ is covered by $T_j$" => "every k-thread is covered by some m-thread for m<i for all k>=i".

We can now show that no stride generator can include a thread that is covered by another thread.

Theorem 213

Let $T_i(b_2)$ and $T_j(c_2)$ be two threads of different orders i>j belonging to some stride generator SG(n,p) = {1, $a_2$, $a_3$}.

Then $T_i$ is not covered by $T_j$.

Proof

Suppose that $T_i$ is covered by $T_j$; then by Theorem 212 every thread $T_k$ for k>=i is covered by some thread $T_m$ for m<i.

This means that any value x is either covered by no thread, or by a thread of order < i, and so either has no generation or has a generation whose order is less than i.

This means that the order of the stride generator must be less than i - which is contrary to hypothesis; hence no thread is ever completely covered by any other thread in a stride generator.

[ We can define a thread as "redundant" in a stride generator if each of the values covered by the thread is also covered by another thread.

Such redundant threads do exist; an example is $T_2(5)$ in SG(8,2) = {1,14,33}, illustrated in Figure 200. ]

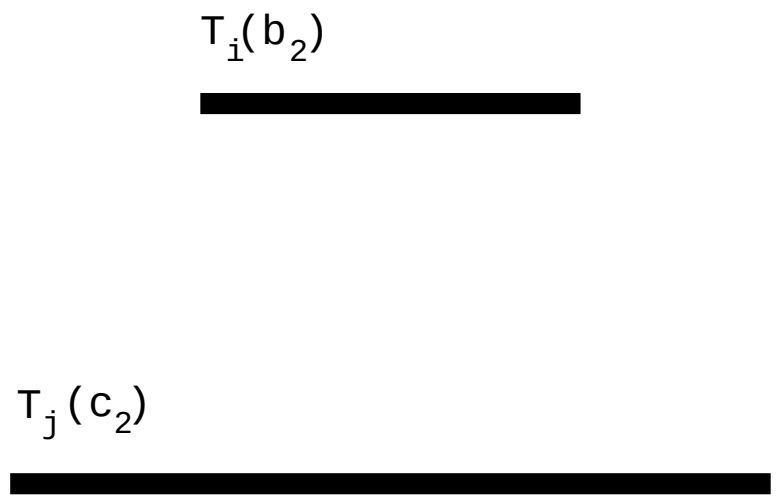

Figure 208

## 2.8 Order zero stride generators

This section identifies all order 0 stride generators SG(n,0). In particular, we show that sets of the form {1, $a_2$, $C_2a_3$} - i.e. those where $a_3$ is an exact multiple of $a_2$ - are always order 0 stride generators, and so can be ignored in future sections when only higher order stride generators are considered.

We also show that a set that is an order 0 stride generator cannot also be a stride generator of higher order.

Theorem 214

Let $a_3 = C_2a_2 + C_1$

Then {1, $a_2$, $a_3$} is an order 0 stride generator SG(n,0) provided:

$C_1 = 0$ or $a_2-C_2 <= C_1 < a_2$

in which case:

$n = a_2+C_2-2$

and the (first) break is at:

$y = C_2a_2-1$

[ Note (added in 1995): "$a_2-C_2 <= C_1$ or $C_1 = 0$" is exactly the requirement for A to be a "pleasant" basis - see Selmer's notes in the literature. ]

Proof

Theorems 207 and 204 show that the threads of a stride generator satisfy:

- the threads together cover the range $0<x<a_3$, and
- each break value lies at the end of a thread

Since we have only 0-threads, the arrangement must be as shown in Figure 209 where $T_0(k-1)$ and $T_0(k)$ are contiguous and so there is a single gap between $T_0(k)$ and $T_0(k+1)$.

Clearly all values from 0 to $(k+1)a_2-2$ inclusive are covered by 0-threads, and so $a_3<=(k+1)a_2-1$.

The (potential) break points are at $y=ka_2-1$ and $y=(k+1)a_2-2$, and so $a_3>=ka_2$.

Since $a_3=C_2a_2+C_1$, this shows that $k=C_2$ for an order 0 stride generator.

We can determine the value of n from the length of $T_0(C_2-1)$, which must be $a_2$ since it is contiguous with $T_0(C_2)$:

$length = n+0-(C_2-1)+1 = a_2$

$=> n = a_2+C_2-2$

We now have only to determine the conditions under which the potential breaks are actual breaks.

We start with $y = C_2a_2-1$:

y is a break <=>
$y + ja_3 = c_2a_2 + c_1$  $c_2+c_1<=n+j-1$  is not soluble for any $j<=p+1$

i.e. $y = c_2a_2 + c_1$  $c_2+c_1<=n-1$  is not soluble
and $y + a_3 = c_2a_2 + c_1$  $c_2+c_1<=n$  is not soluble

The first equation is not soluble because y is at the end of $T_0(C_2-1)$.

We now consider the second equation:

$$y + a_3 = (C_2a_2-1) + (C_2a_2+C_1)$$

$$= 2C_2a_2+C_1-1$$

There are two cases to consider:

a) $C_1=0$ - in which case the way of expressing $y+a_3$ as $c_2a_2+c_1$ with $0<=c_1<a_2$ is as:

$$y + a_3 = (2C_2-1)a_2+a_2-1$$

This is not a solution provided that:

$$(2C_2-1)+(a_2-1) > n$$

$$=> 2C_2+a_2-2 > a_2+C_2-2$$

$$=> C_2 > 0$$

This is always true since $a_3>a_2$, and so $C_1=0$ makes y a break and gives an order 0 stride generator.

b) $C_1>0$ - in which case the way of expressing $y+a_3$ as $c_2a_2+c_1$ with $0<=c_1<a_2$ is as:

$$y + a_3 = 2C_2a_2+(C_1-1)$$

This is not a solution provided that:

$$2C_2+C_1-1 > n$$

$$=> 2C_2+C_1-1 > a_2+C_2-2$$

$$=> C_1 > a_2-C_2-1$$

$$=> C_1 >= a_2-C_2$$

So $a_2-C_2 <= C_1 < a_2$ also gives a break and an order 0 stride generator.

Next we consider the second case, where $y=(C_2+1)a_2-2$, and so the only possible value for $a_3$ is $(C_2+1)a_2-1$; this means that $C_1=a_2-1 >= a_2-C_2$ since $C_2>=1$.

We have just shown that under these circumstances we have a zero order stride generator with a break at $y=C_2a_2-1$, and so the theorem is proved.

[The next theorem but one shows that $y=(C_2+1)a_2-2$ is, indeed, a second break under these circumstances.]

Theorem 215

All sets of the form {1, $a_2$, $a_3$} where $a_3$ >= $(a_2-1)a_2$ are order 0 stride generators.

Proof

This is a direct corollory of the theorem above.

In such cases, $a_2-C_2<=1$, and so either $C_1=0$ or $C_1>=(a_2-C_2)$ whatever the value of $C_1$.

Theorem 216

Order 0 stride generators SG(n,0) of the form {1, $a_2$, $(C_2+1)a_2-1$} have a second (and final) break at $a_3-1$.

No other order 0 stride generator has more than one break.

Proof

The proof of Theorem 214 above identified $y=a_3-1 = (C_2+1)a_2-2$ as a potential second break for such a stride generator.

To show that y is, indeed, a break we must show that:

$y = c_2a_2 + c_1 \quad c_2+c_1<=n-1$ is not soluble
and $y + a_3 = c_2a_2 + c_1 \quad c_2+c_1<=n$ is not soluble

The first equation is not soluble because y is at the end of $T_0(C_2)$.

We now consider the second equation:

$$y + a_3 = C_2a_2+a_2-2 + C_2a_2+a_2-1$$

$$= (2C_2+2)a_2-3$$

There are two cases to consider:

a) $a_2>=3$ - in which case the way of expressing $y+a_3$ as $c_2a_2+c_1$ with $0<=c_1<a_2$ is as:

$$y + a_3 = (2C_2+1)a_2+(a_2-3)$$

This is not a solution provided that:

$$2C_2+1+a_2-3 > n$$

$$=> 2C_2+1+a_2-3 > a_2+C_2-2 \quad \text{- by Theorem 214}$$

$$=> C_2 > 0$$

This is always true, and so y is a break in this case.

b) $a_2=2$ - in which case the way of expressing $y+a_3$ as $c_2a_2+c_1$ with $0<=c_1<a_2$ is as:

$$y + a_3 = 2C_2a_2+1$$

This is not a solution provided that:

$2C_2+1 > n$

$\Rightarrow 2C_2+1 > C_2$

$\Rightarrow C_2 > -1$

This is also always true.

So $y = (C_2+1)a_2-2 = a_3-1$ is the second break in this case.

Theorem 217

Every order 0 stride generator is canonical.

Proof

We prove this result by showing that there is a 1-thread $T_1$ that is covered by some 0-thread $T_0$.

Then by Theorem 212 we deduce that every i-thread for i>0 is covered by some 0-thread.

From this we deduce that every break y is canonical and so the stride generator itself is canonical:

Consider any break y.

By definition, y lies at the end of a 0-thread, and so if y is covered by any k-thread for k>0 it must also lie at the end of that k-thread (since any k-thread is covered by a 0-thread).

In other words, any solution to:

$y + ka_3 = c_2a_2 + c_1 \quad c_2+c_1<=n+k \qquad k>0$

can only be:

$y + ka_3 = c_2a_2 + c_1 \quad c_2+c_1=n+k \qquad k>0$

and therefore no solution is possible for:

$y + ka_3 = c_2a_2 + c_1 \quad c_2+c_1<=n+k-1 \quad k>0$

and so y is a canonical break.

It now remains only to show that there is a 1-thread covered by a 0-thread; there are two cases to consider:

a) $a_2-C_2<=C_1<a_2$ - in which case we show that $T_1(2C_2)$ is covered by $T_0(C_2-1)$:

$ST_1 = 2C_2a_2-a_3$

$= 2C_2a_2-C_2a_2-C_1$

$= C_2a_2-C_1$

$ST_0 = (C_2-1)a_2$

$ET_1 = C_2a_2-C_1+n+1-2C_2$

$= C_2a_2-C_1+a_2+C_2-2+1-2C_2$ - since $n=a_2+C_2-2$ by Theorem 214

$= (C_2+1)a_2-C_2-C_1-1$

$ET_0 = (C_2-1)a_2+n-C_2+1$

$= (C_2-1)a_2+a_2+C_2-2-C_2+1$

$= C_2a_2-1$

We have to show:

$ST_1>=ST_0 \Rightarrow ST_1-ST_0>=0$

and $ET_1<=ET_0 \Rightarrow ET_0-ET_1>=0$

$ST_1-ST_0 = C_2a_2-C_1-(C_2-1)a_2$

$= a_2-C_1$

$> 0$ because $0<=C_1<a_2$

$ET_0-ET_1 = C_2a_2-1-(C_2+1)a_2+C_2+C_1+1$

$= -a_2+C_2+C_1$

$= C_1-(a_2-C_2)$

$>= 0$ because $a_2-C_2<=C_1$

b) $C_1=0$ - in which case we show that $T_1(2C_2-1)$ is covered by $T_0(C_2-1)$:

$ST_1 = (2C_2-1)a_2-C_2a_2$

$= (C_2-1)a_2$

So $ST_1=ST_0$.

$ET_1 = (C_2-1)a_2+n+1-2C_2+1$

$= (C_2-1)a_2+a_2+C_2-2+1-2C_2+1$

$= C_2a_2-C_2$

So $ET_0-ET_1 = C_2a_2-1-C_2a_2+C_2$

$= C_2-1$

$>= 0$ because $C_2>=1$

Theorem 218

No set A that is an order 0 stride generator can also be a stride generator of higher order.

Proof

The proof of Theorem 217 above shows that in every order 0 stride generator every i-thread for i>0 is covered by some corresponding 0-thread.

It is easy to show that this relationship remains true regardless of the value of n, since a change in n alters the length of the threads equally, and does not alter their relative positions:

Suppose $T_i$ is covered by $T_0$ for SG(n,0) and let $T_i'$, $T_0'$ be the corresponding threads for the value n'; then:

$$ST_i' = ST_i \qquad ST_0' = ST_0$$
$$ET_i' = ET_i+(n'-n) \qquad ET_0' = ET_0+(n'-n)$$

So:

$$ST_i >= ST_0 \;\; => \;\; ST_i' >= ST_0'$$
$$\text{and } ET_i <= ET_0 \;\; => \;\; ET_i' <= ET_0'$$

Therefore there is no stride generator A = SG(n',p) of order p>0, since any attempt to create such an object would result in a situation where i-threads are covered by 0-threads - which is not permissible by Theorem 213.

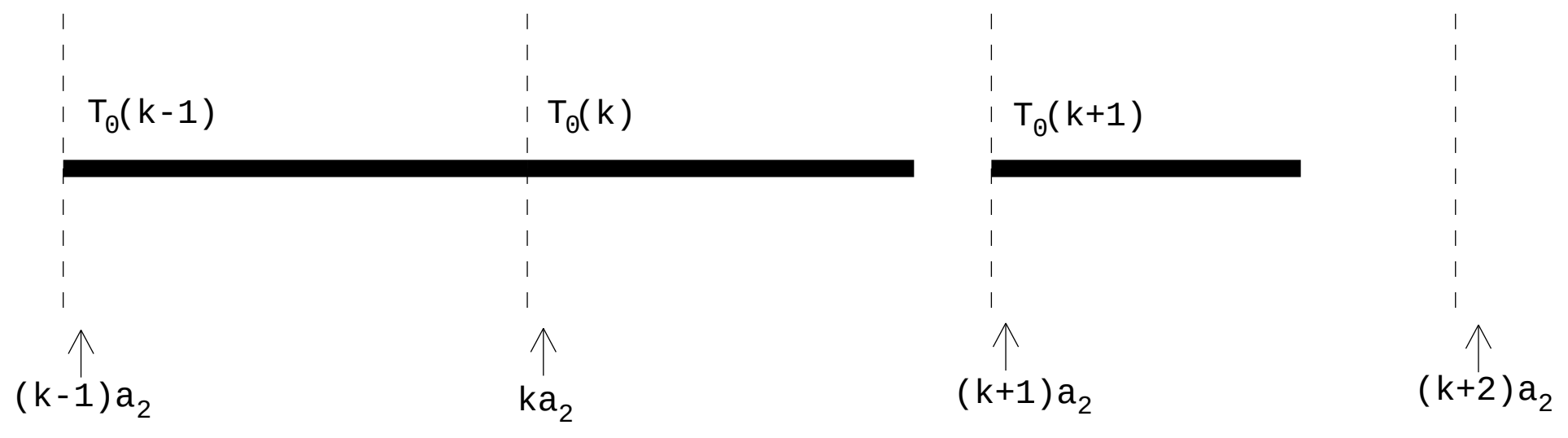

Figure 209

## 2.9 Stride generators of non-zero order

In this section we characterise the breaks in stride generators of non-zero order. The results are summarised as:

For any stride generator of order greater than zero:

Either: All breaks are canonical

or: The first break has highest break order

It is this result that makes the concept of "potential cover" (see section 2.10) worthwhile.

Note that throughout this section we consider only stride generators of non-zero order, and so we know by Theorems 214 and 218 that $a_3$ cannot be an exact multiple of $a_2$.

We first show that any break value must be $\geq a_3-a_2$:

Theorem 219

If y is any break in a stride generator $A = SG(n,p) = \{1, a_2, a_3\}$, then:

$$y \geq a_3-a_2$$

Proof

Consider $y' = y+a_2$

If $y' < a_3$, then:

$y' + ia_3 = c_2a_2 + c_1 \quad c_2+c_1\leq n+i$ is soluble for some $i\leq p$ by (SG1)

$\Rightarrow y + ia_3 = (c_2-1)a_2 + c_1 \quad c_2+c_1\leq n+i$ is soluble for some $i\leq p$

Write $c_2' = c_2-1$:

$\Rightarrow y + ia_3 = c_2'a_2 + c_1 \quad c_2'+c_1\leq n+i-1$ is soluble for some $i\leq p$

This contradicts (SG3): so the assumption that $y'<a_3$ is incorrect.

Therefore $y' \geq a_3 \Rightarrow y \geq a_3-a_2$.

We now refine this limit:

Theorem 220

Let $A = \{1, a_2, a_3\}$ be a stride generator SG(n,p).

a) If $a_3>2a_2$ then any break y satisfies:

$$(a_3-a_2) \leq y < (a_3-a_2)+n$$

b) If $a_3<2a_2$ then any break y is part of a series of breaks $y_0, y_1, \ldots$ where $y_i = y_{i-1}+C_1$ (where $C_1=a_3-a_2$) and $y_0$ satisfies:

$$(a_3-a_2) \leq y_0 < (a_3-a_2)+n$$

[ There may, of course, be more than one such series; we show later in Theorem 227 that there are at most two. ]

Proof

y is a break =>

$y + ja_3 = c_2a_2 + c_1 \quad c_2+c_1<=n+j-1$ is not soluble for $0<=j<=p+1$

Write $y' = y-(a_3-a_2)$ (note that $y'>=0$ by Theorem 219 above):

$y'+(a_3-a_2) + ja_3 = c_2a_2 + c_1 \quad c_2+c_1<=n+j-1$ is not soluble for $0<=j<=p+1$

=> $y' + (j+1)a_3 = (c_2+1)a_2 + c_1 \quad c_2+c_1<=n+j-1$ is not soluble for $0<=j<=p+1$

Write $j'=j+1$, $c_2'=c_2+1$:

=> $y' + j'a_3 = c_2'a_2 + c_1 \quad c_2'+c_1<=n+j'-1$ is not soluble for $1<=j'<=p+2$

y' would be a break if this were true for j'=0 as well, so y' is a break iff:

$y' = c_2'a_2 + c_1 \quad c_2'+c_1<=n-1$ is not soluble - (1)

Now $a_3-a_2 <= y < a_3$ (by Theorem 219 above) and so $0 <= y' < a_2$.

Furthermore, if y' is also a break, then:

$a_3-a_2 <= y' < a_3$ - by Theorem 219

and so $a_3-a_2 <= y' < a_2$

=> $a_3 <= 2a_2$ - (2)

We now split the proof into two cases:

a) $a_3>2a_2$

In this case y' cannot be a break by (2).

Therefore by (1):

$y' = c_2'a_2 + c_1 \quad c_2'+c_1<=n-1$ is soluble

Since $y'<a_2$, the solution must be:

$y' = c_1 \quad c_1<=n-1$

=> $y' <= n-1$

Therefore $y = y'+(a_3-a_2) < (a_3-a_2)+n$

b) $a_3<2a_2$

In this case y' is a break iff (1) holds, which is true iff $y'>=n$ (see (a) above).

In other words:

y is a break => $y-(a_3-a_2)$ is a break provided $y>=(a_3-a_2)+n$

Thus given any break $y_i$ there is a series of breaks $y_i$, $y_{i-1}$, ... $y_0$ such that:

$$y_{i-1} = y_i-(a_3-a_2)$$

and this series terminates with $y_0<(a_3-a_2)+n$.

By Theorem 219 each break $y_i>=(a_3-a_2)$, and so there exists a break $y_0$ such that:

$$(a_3-a_2) <= y_0 < (a_3-a_2)+n$$

We now show in two analogous steps that there is a value in the range $(a_3-a_2) <= x <= (a_3-a_2)+n$ that requires an order p generation.

Theorem 221

If A = {1, $a_2$, $a_3$} is a stride generator SG(n,p), then at least one value:

$$a_3-a_2 <= x < a_3$$

requires an order p generation.

[In other words, the range $a_3-a_2<=x<a_3$ always includes at least part of a p-thread]

Proof

Suppose this is not the case; then for all $a_3-a_2<=x<a_3$:

$$x + ia_3 = c_2a_2 + c_1 \quad c_2+c_1<=n+i \quad \text{is soluble for } i<p$$

Now consider any value $0<y<a_3-a_2$; there exists k>0 such that:

$$a_3-a_2 <= y+ka_2 < a_3$$

So writing $x = y+ka_2$ we have:

$$y + ka_2 + ia_3 = c_2a_2 + c_1 \quad c_2+c_1<=n+i \quad \text{is soluble for } i<p$$

$$=> y + ia_3 = (c_2-k)a_2 + c_1 \quad c_2+c_1<=n+i \quad \text{is soluble for } i<p$$

Write $c_2' = c_2-k$:

$$=> y + ia_3 = c_2'a_2 + c_1 \quad c_2'+c_1<=n+i \quad \text{is soluble for } i<p$$

Clearly, $c_2'>=0$ (since y>0, i>=0, $c_1<a_2$).

So all values $0<x<a_3$ satisfy (SG1) for i<p - which is contrary to(SG2) for SG(n,p).

Therefore the supposition is wrong, and the theorem is proved.

Theorem 222

Let SG(n,p) be a stride generator {1, $a_2$, $a_3$}.

Then there exists a value x which requires an order p generation such that:

$$(a_3-a_2) <= x <= (a_3-a_2)+n$$

Proof

By Theorem 221, we know that at least one value:

$$(a_3-a_2) <= x < a_3$$

requires an order p generation.

We complete the proof by showing that no value:

$$(a_3-a_2)+n < x < a_3$$

requires an order p generation.

Suppose the contrary; then there exists $(a_3-a_2)+n < x < a_3$ such that:

$x + ia_3 = c_2a_2 + c_1 \quad c_2+c_1<=n+i$ is soluble for i=p
and not soluble for any i<p

By the corollory to Theorem 210, we have $x' = x-(a_3-a_2)$ such that:

$x' + (i+1)a_3 = (c_2+1)a_2 + c_1 \quad c_2+1+c_1<=n+i+1$ is soluble for i=p
and not soluble for any i<p

Write $i'=i+1$, $c_2'=c_2+1$:

$x' + i'a_3 = c_2'a_2 + c_1 \quad c_2'+c_1<=n+i'$ is soluble for i'=p+1
and not soluble for any 1<=i'<=p

But by (SG1), there must exist some i'<=p for which the above equation is soluble - and the only possibility is i'=0; so:

$$x' = c_2'a_2 + c_1 \quad c_2'+c_1<=n$$

But by hypothesis:

$$n < x' < a_2 \;\; => \;\; x' = 0.a_2 + c_1 \;\; \text{where } c_1>n$$

So such a value x' has no valid generation of order <=p, and so our hypothesis that $x > (a_3-a_2)+n$ is false.

The following two theorems show that all possible positions for threads in a stride generator are actually filled. This means that we do not need to check any constraints when deriving one thread $T_j$ from another one $T_i$ by one of the "series" theorems 208, 209 and 210 provided only that both $T_j$ and $T_i$ lie within the stride generator.

Theorem 223

In every stride generator SG(n,p) there is exactly one thread $T_k$ such that:

$$a_3-a_2 <= ST_k < a_3$$

for all k<=p.

[In other words, there is always a k-thread for every k<=p that starts in the last $a_2$ positions of a stride generator SG(n,p).]

Proof

The range $a_3$-$a_2$<=x<$a_3$ contains exactly $a_2$ values, and threads of order k start at intervals of $a_2$; so there is exactly one potential starting position for a k-thread in the range.

We now show that such a k-thread exists: that is, that its length is greater than zero.

By Theorem 221, we know that there is a p-thread $T_p(c_2)$ part of which lies within the range; we now show that this means that there is a p-thread that starts within the range:

a) $ST_p(c_2)$>=$a_3$-$a_2$ - nothing to prove

b) $ST_p(c_2)$<$a_3$-$a_2$

In this case:

$L_p(c_2)$ >= ($a_3$-$a_2$)-$ST_p(c_2)$ + 1 > 1, since part of the thread lies within the range.

Therefore $T_p(c_2+1)$ starting at $ST_p(c_2)$+$a_2$ of length $L_p(c_2)$-1 > 0 exists by Theorem 208.

Either $T_p(c_2+1)$ starts within the range, or we repeat step (b).

So there is a p-thread that starts within the range: let this be $T_p(c_2)$.

Now consider the k-thread $T_k(b_2)$ closest but before $T_p(c_2)$:

0 < $ST_p(c_2)$-$ST_k(b_2)$ < $a_2$

By Theorem 211, $L_k(b_2)$>1, since $L_p(c_2)$>=1 by definition.

If $ST_k(b_2)$>=$a_3$-$a_2$, the theorem is proved.

If not, consider $T_k(b_2+1)$ of length $L_k(b_2)$-1 > 0 derived by Theorem 208. This thread starts at $ST_k(b_2)$+$a_2$ > $ST_p(c_2)$ >= $a_3$-$a_2$.

Theorem 224

If x + k$a_3$ = $c_2a_2$ and x<$a_3$, k<=p, then there is a thread $T_k(c_2)$ starting at x.

[In other words, all possible positions for threads in a stride generator are filled. Note that it is convenient to allow x<0 to accommodate the case of a k-thread with $ST_k(c_2)$<0, $ET_k(c_2)$>=0.]

Proof

Theorem 223 shows that the highest possible starting position for a k-thread is filled, and repeated use of Theorem 208 shows that all lower positions must also be filled.

Definition 207

A break y in a stride generator SG(n,p) is called a "fundamental" break if it lies in the range $(a_3-a_2)<=y<(a_3-a_2)+n$.

[Theorem 220 shows that every stride generator has at least one fundamental break, and that any stride generator where $a_3>2a_2$ has only fundamental breaks]

Theorems 204 and 205 show that any break $y<a_3-1$ must lie at the end of the first of two contiguous threads (see Figure 203). The next theorem shows that for a fundamental break one of those two threads must be of order p.

Theorem 225

Let $y<a_3-1$ be a break in a stride generator SG(n,p) = {1, $a_2$, $a_3$}.

Then Theorems 204 and 205 show that:

$y + ia_3 = c_2a_2 + c_1$  $c_2+c_1=n+i$  is soluble for some i<=p
and $(y+1) + ja_3 = c_2'a_2$  $c_2'<=n+j$  is soluble for some j<=p;

in other words, $y = ET_i$ and $y+1 = ST_j$ for some i,j<=p

This theorem states that either i=p or j=p for any fundamental break y.

[Note that the case of $y=a_3-1$ is dealt with separately as Theorem 226 below.]

Proof

Let $a_3 = C_2a_2 + C_1$  =>  $a_3-a_2 = (C_2-1)a_2 + C_1$

We first consider the threads covering the region around $(C_2-1)a_2$.

First, we have $T_0(C_2-1)$ of length $n-(C_2-1)+1$:

$ST_0(C_2-1) = (C_2-1)a_2$
$ET_0(C_2-1) = (C_2-1)a_2 + n - (C_2-1)$

We now show formally the intuitive fact that this thread must be met or overlapped at both ends by other threads of non-zero order as illustrated in Figure 210.

If $T_0(C_2-1)$ has length $a_2$ or greater, then it is easy to see that at least all values from 0 to $(C_2+1)a_2-2$ inclusive are covered by zero order threads. This means that the stride generator itself is of order zero, contrary to our assumption that p>0.

So we may assume that $ET_0<C_2a_2-1$; we write $x = ET_0(C_2-1)+1$ and $y = ST_0(C_2)-1$. We know that neither of these values (which may be the same) is covered by a zero order thread.

So x must be covered by a j-thread $T_j(d_2)$ which meets or overlaps the end of $T_0(C_2-1)$.

Similarly, y must be covered by an i-thread $T_i(b_2+1)$ which meets or overlaps the start of $T_0(C_2)$; hence, by Theorem 208, the i-thread $T_i(b_2)$ - whose length is one greater than that of $T_i(b_2+1)$ - must overlap the start of $T_0(C_2-1)$.

So we know that threads surround $T_0(C_2-1)$; we now choose as $T_i$ and $T_j$ the two that are "closest" to $(C_2-1)a_2$:

$T_i(b_2)$ of length $L_i$ such that $ST_0(C_2-1) - ST_i(b_2)$ is minimal
$T_j(d_2)$ of length $L_j$ such that $ST_j(d_2) - ST_0(C_2-1)$ is minimal

Referring to Figure 210 again, we note that $M_i$ and $M_j$ are non-zero:

e.g. Suppose $M_i=0$

Then $T_i$ is completely covered by $T_0$, which by Theorem 213 is not permissible in a stride generator.

We now consider the impact of these threads on the region around $a_3-a_2 = (C_2-1)a_2 + C_1$.

By Theorems 209 and 224*, we can derive threads $T_{i-1}$ and $T_{j-1}$ from $T_i$ and $T_j$, but no thread can be derived from $T_0$:

$T_{i-1}(b_2-C_2)$ of length $L_{i-1} = L_i + (C_2-1)$ with $ST_{i-1} = ST_i + C_1$
$T_{j-1}(d_2-C_2)$ of length $L_{j-1} = L_j + (C_2-1)$ with $ST_{j-1} = ST_j + C_1$

This is illustrated in Figure 211; note that it is possible for $T_{i-1}$ and $T_0$ to be one and the same thread.

[*Note:

Theorem 224 is applicable only if threads $T_{i-1}$ and $T_{j-1}$ both start within the stride generator - that is, if $ST_{i-1}<a_3$ and $ST_{j-1}<a_3$; fortunately, this is the case:

We have $ST_{i-1} = ST_i + C_1 < ST_0 + C_1 = (C_2-1)a_2 + C_1 = a_3-a_2 < a_3$

and: $ST_{j-1} = ST_j + C_1 <= (ET_0+1) + C_1 < ET_0 + 2 + C_1$

$= (C_2-1)a_2 + n - (C_2-1) + 2 + C_1$

$= C_2a_2 - a_2 + n - C_2 + 3 + C_1$

So: $ST_{j-1} < a_3 \iff C_2a_2 - a_2 + n - C_2 + 3 + C_1 < C_2a_2 + C_1$

$\iff -a_2 + n - C_2 + 3 < 0$

$\iff n < a_2+C_2-3$

This is true because $n >= a_2+C_2-2 \Rightarrow T_0(C_2-1)$ has length $>=a_2$, and so all values $0<x<a_3$ are covered by 0-threads and the stride generator must be of order 0 (cf Theorem 214).]

We now show that with the possible exception of a gap between the end of $T_{i-1}$ and the start of $T_{j-1}$ these two threads cover the range $(a_3-a_2)-1$ to $(a_3-a_2)+n+1$ inclusive:

$ST_{i-1} = ST_i + C_1$

$< (C_2-1)a_2 + C_1 = a_3-a_2$

$ET_{j-1} = ET_j + (C_2-1) + C_1$

$> ET_0 + (C_2-1) + C_1$

$= (C_2-1)a_2 + n - (C_2-1) + (C_2-1) + C_1$

$= (C_2-1)a_2 + n + C_1$

$= (a_3-a_2) + n$

So: $ST_{i-1} <= (a_3-a_2)-1 < (a_3-a_2)+n+1 <= ET_{j-1}$

Suppose the two threads $T_{i-1}$ and $T_{j-1}$ overlap, and consider any break value y.

By Theorem 220, $a_3-a_2 <= y < (a_3-a_2)+n$, and so y is covered by either $T_{i-1}$ or $T_{j-1}$.

By Theorem 204, y lies at the end of a thread, and since $ET_{j-1}>(a_3-a_2)+n$ y must lie at the end of $T_{i-1}$.

Further, if y is at the end of $T_{i-1}$ it is also covered by $T_{j-1}$ (since the threads overlap) and so:

either: $y + (i-1)a_3 = c_2a_2 + c_1 \quad c_2+c_1<n+(i-1)$

or: $y + (j-1)a_3 = c_2a_2 + c_1 \quad c_2+c_1<n+(j-1)$

In other words, y cannot be a break, and so $T_{i-1}$ and $T_{j-1}$ do not overlap.

Suppose the two threads are contiguous:

$ET_{i-1} + 1 = ST_{j-1}$

In this case, $T_{i-1}$ and $T_{j-1}$ together cover the range $a_3-a_2 <= y <= (a_3-a_2)+n$ and so by Theorem 222 either i-1=p or j-1=p. This would mean that either i=p+1 or j=p+1, which is contrary to the assumption that $T_i$ and $T_j$ are part of the stride generator.

So $T_{i-1}$ and $T_{j-1}$ are not contiguous.

So there must be a gap between $ET_{i-1}$ and $ST_{j-1}$ and hence a (greater or equal) gap between $ET_i$ and $ST_j$. This latter gap is, of course, covered by $T_0$; we now show that no other thread of order <=p can cover any part of this gap.

Any such thread $T_k$ would have to start before $ST_i$, because $T_i$ and $T_j$ were chosen to be the threads closest to the start of $T_0$; but since no thread can cover another thread (Theorem 213), $T_k$ must finish before $ET_i$ and so cannot cover any part of the gap.

But the gap between $ET_{i-1}$ and $ST_{j-1}$ must be covered by one (or more) threads. Suppose that one of these threads is $T_{k-1}$ of order k-1<p.

Using the corollory to Theorem 209 together with Theorem 224, we can derive $T_k$, k<=p, such that:

$ST_k = ST_{k-1} - C_1$
$ET_k = ET_{k-1} - C_1 - (C_2-1)$

We now show that $T_k$ must cover some part of the gap between $ET_i$ and $ST_j$; there are three cases to consider:

a) $ET_{k-1}<ST_{j-1}$ (see Figure 212)

In this case, the value $x = ET_{k-1}$ must lie in the gap between $ET_{i-1}$ and $ST_{j-1}$, and we show that the corresponding value $x' = x-C_1-(C_2-1)$ is part of $T_k$ and lies in the gap between $ET_i$ and $ST_j$:

$x' = ET_{k-1}-C_1-(C_2-1) = ET_k$ - and so lies at the end of $T_k$

$x < ST_{j-1} \Rightarrow x' < ST_{j-1}-C_1-(C_2-1) = ST_j-(C_2-1) < ST_j$

$x > ET_{i-1} \Rightarrow x' > ET_{i-1}-C_1-(C_2-1) = ET_i$

b) $ET_{k-1}>=ST_{j-1}$, $ST_{k-1}>ET_{i-1}$ (see Figure 213)

In this case, the value $x = ST_{k-1}$ must be in the gap between $ET_{i-1}$ and $ST_{j-1}$, and we show that the corresponding value $x' = x-C_1 = ST_k$ lies in the gap between $ET_i$ and $ST_j$:

$x < ST_{j-1} \Rightarrow x' < ST_{j-1}-C_1 = ST_j$

$x > ET_{i-1} \Rightarrow x' > ET_{i-1}-C_1 = ET_i+(C_2-1) >= ET_i$

c) $ET_{k-1}>=ST_{j-1}$, $ST_{k-1}<=ET_{i-1}$ (see Figure 214)

This case cannot occur, because it would mean that the entire range $(a_3-a_2)<=y<(a_3-a_2)+n$ was covered by overlapping threads and so would contain no break, contrary to Theorem 220.

So we have shown that if any part of the gap between $ET_{i-1}$ and $ST_{j-1}$ is covered by a thread of order <p, then some part of the gap between $ET_i$ and $ST_j$ is covered by a thread of order <=p.

But we showed above that this is not possible, and so the gap between $ET_{i-1}$ and $ST_{j-1}$ must be covered by a single thread of order p.

In other words, the region $(a_3-a_2)<=y<(a_3-a_2)+n$ in which any fundamental break y must occur is covered by the threads $T_{i-1}$, $T_p$ and $T_{j-1}$ where:

$ST_{i-1} <= (a_3-a_2)-1$
$ET_{j-1} >= (a_3-a_2)+n+1$

$ST_{i-1} < ST_p < ST_{j-1} < a_3$

and there is no other thread $T_k$ such that $ST_{i-1} < ST_k < ST_{j-1}$

But breaks $y < a_3-1$ only occur at the junctions of contiguous threads, so:

either a) $T_{i-1}$, $T_p$ are contiguous and $y = ET_{i-1}$, $y+1 = ST_p$

or b) $T_p$, $T_{j-1}$ are contiguous and $y = ET_p$, $y+1 = ST_{j-1}$

as illustrated in Figure 215; this concludes the proof.

[ Note:

The theorems above which give upper bounds of the form $(a_3-a_2)+n$ may not be as good as the earlier ones where an upper bound $(a_3-a_2)$ was given. This is because it is possible that $n>a_2$; as an example, consider A = {1, 13, 84} which is an order 1 stride generator SG(16,1).

Great care has to be taken in the proofs above to make sure that assumptions such as "$(a_3-a_2)+n<=a_3$" are not made; an example of this is the need to show that $ST_{j-1}<a_3$. ]

Next we need to cover the special case of $y=a_3-1$; the following theorem shows this case not to be relevant to our argument.

Theorem 226

No stride generator of non-zero order has a fundamental break at $y=a_3-1$.

Proof

Since the stride generator is of non-zero order, there must exist a gap between 0-threads that is covered by one or more threads of higher order.

Suppose the first gap appears between $T_0(c_2)$ and $T_0(c_2+1)$ and suppose the thread that covers the value $ST_0(c_2+1)-1$ is $T_i(b_2)$, as illustrated in Figure 216.

We now show that $c_2>=2$:

y is a fundamental break

=> $y < (a_3-a_2)+n$    by Definition 207

=> $a_3-1 < a_3-a_2+n$

=> $n > a_2-1$

So we have:

$ST_0(0) = 0$    $ET_0(0) = n >= a_2$
$ST_0(1) = a_2$    $ET_0(1) = a_2+n-1 >= 2a_2-1$
$ST_0(2) = 2a_2$    $ET_0(2) = 2a_2+n-2 >= 3a_2-2$
$ST_0(3) = 3a_2$    ...

which shows that no gap can occur between $T_0(c_2)$ and $T_0(c_2+1)$ while $c_2<2$.

Now consider the threads $T_0(0)$, $T_0(1)$ and $T_i(b_2-c_2)$ as illustrated in Figure 217.

These are derived by repeated applications of Theorem 208, and so each derived thread's start position is $c_2a_2$ before its original, and each derived thread is $c_2$ longer than its original.

This means that $T_i(b_2-c_2)$ must overlap $T_0(1)$ by at least $c_2>=2$ places.

So $T_i(b_2-c_2)$ covers at least $ST_0(1)-1 = a_2-1$ and $ST_0(1) = a_2$.

We can now derive the thread $T_{i-1}(b_2-c_2-1)$ by Theorems 210 and 224, since i>0.

This thread is the same length as $T_i(b_2-c_2)$ but is displaced from it by $a_3-a_2$; so $T_{i-1}(b_2-c_2-1)$ will cover both $a_3-1$ and $a_3$.

Thus the values y and y+1 are both covered by the same thread of order <= p, which means that y cannot be a break by Theorem 204.

So no stride generator of non-zero order has a fundamental break at $y=a_3-1$.

[Notes

First, we know that stride generators of order 0 do have breaks at $y=a_3-1$, and we can also characterise them precisely:

Theorems 214 and 216 show that order zero stride generators can have breaks at $y=a_3-1$, and that this happens iff the stride generator has one of the following forms:

a) $SG(a_2+C_2-2,0) = \{1, a_2, C_2a_2\}$
which has one break at $y=C_2a_2-1 = a_3-1$

b) $SG(a_2+C_2-2,0) = \{1, a_2, (C_2+1)a_2-1\}$
which has two breaks $y=C_2a_2-1$ and $y=(C_2+1)a_2-2 = a_3-1$

In the case of $a_2=34$, this gives:

$C_2=1$, (b): {1,34,67} = SG(33,0) with breaks at 33, 66
$C_2=2$, (a): {1,34,68} = SG(34,0) with a break at 67
$C_2=2$, (b): {1,34,101} = SG(34,0) with breaks at 67, 100
$C_2=3$, (a): {1,34,102} = SG(35,0) with a break at 101

Secondly, it is not true to say that no stride generator of order > 0 has any break at $y=a_3-1$. Since such a break cannot be fundamental, the stride generator must be such that $a_3<2a_2$, by Theorem 220.

There are many such examples, as the following exhaustive list for $a_2=34$ shows:

{1,34,35} = SG(1,32) has breaks at 1*,2,3,4, ... ,33,34
{1,34,36} = SG(2,16) has breaks at 3*,5,7,9, ... ,33,35
{1,34,37} = SG(3,10) has breaks at 3*,6,9,12, ... ,33,36
{1,34,38} = SG(3,16) has breaks at 5*,9,13,17, ... 33,37
{1,34,42} = SG(5,16) has breaks at 9*,17,25,33,41
{1,34,45} = SG(11,2) has breaks at 11*,22,33,44
{1,34,49} = SG(6,8) has breaks at 15*,18*,30,33,45,48
{1,34,50} = SG(9,16) has breaks at 17*,33,49
{1,34,51} = SG(17,1) has breaks at 33*,50
{1,34,61} = SG(9,4) has breaks at 27*,33*,54,60
{1,34,63} = SG(9,6) has breaks at 29*,33*,58,62
{1,34,66} = SG(17,16) has breaks at 33*,65

where * indicates a fundamental break. ]

We are now on the homeward straight!

We first show that every stride generator has at most two fundamental breaks, and that if two are present they are both canonical. This shows that for any stride generator with $a_3>2a_2$, the first break has highest order.

Then we look at the series of breaks which may occur when $a_3<2a_2$, and show that if the first break is canonical then so are the rest, and that if the first has break order q, then all subsequent ones have a lower break order: once again, the first break has highest order.

[In summary:

$a_3>2a_2$: either one or two breaks
if two, both are canonical

$a_3<2a_2$: either one or two series of breaks
if two, all breaks are canonical
if one, either: all breaks are canonical
or: break order decreases as y increases

Thus for any stride generator:

either: all breaks are canonical
or: the first break has maximum break order.]

Theorem 227

There are at most two fundamental breaks in any stride generator.

Proof

This is a corollory of theorem 225, where the proof shows that the region in which fundamental breaks can occur is covered by just three threads $T_{i-1}$, $T_p$ and $T_{j-1}$. (Note that the case of $y=a_3-1$ is not possible by Theorem 226.)

Figure 215 illustrates the possible cases: any fundamental break y must be either at the end of $T_{i-1}$ or $T_p$, and so there can be at most two of them.

Theorem 228

If a stride generator has two fundamental breaks, then both of these are canonical breaks.

Proof

By Theorem 227, SG(n,p) has either one or two fundamental breaks, and if it has exactly two breaks then the threads $T_{i-1}$, $T_p$ and $T_{j-1}$ are contiguous. Let L be the length of the thread $T_p$.

By the corollory to Theorem 209 we may be able to derive a thread $T_{p+1}$ of length $L-(C_2-1)$ which will lie between $T_i$ and $T_j$ and hence over $T_0$ as illustrated in Figure 218.

This will not be possible if $L<=(C_2-1)$; but if this is the case, then no threads of order > p can start after $(C_2-1)a_2$, and so there is no possibility that either break can be covered by a thread of order > p - and hence both must be canonical.

So suppose $L>(C_2-1)$, and hence the thread $T_{p+1}$ exists. We now show formally that $T_{p+1}$ is covered by $T_0$:

$ST_p = ET_{i-1} + 1$

$= ET_i + (C_2-1) + C_1 + 1$ where $a_3=C_2a_2+C_1$ (see proof of Theorem 225)

$ST_{p+1} = ST_p - C_1$

$= ET_i + (C_2-1) + 1$

$= ET_i + C_2$

$> ET_i >= ST_0-1$

$ET_p = ST_{j-1} - 1$

$= ST_j + C_1 - 1$

$ET_{p+1} = ET_p - (C_2-1) - C_1$

$= ST_j + C_1 - 1 - (C_2-1) - C_1$

$= ST_j - C_2$

$< ST_j <= ET_0+1$

So: $ST_0-1 < ST_{p+1} <= ET_{p+1} < ET_0+1$

=> $ST_0 <= ST_{p+1} <= ET_{p+1} <= ET_0$

and so $T_{p+1}$ is completely covered by $T_0$.

By Theorem 212, this means that any thread of order >= p+1 is covered by some thread of order <= p.

Now suppose that one of the fundamental breaks y in SG(n,p) is not canonical, but has break order q; then:

$y + qa_3 = c_2a_2 + c_1$  $c_2+c_1<=n+q-1$  is soluble for q>=p+1

This means that y is at least one from the end of the thread $T_q(c_2)$:

$ST_q(c_2) <= y <= ET_q(c_2)-1$

Now we have shown above that $T_q(c_2)$ must be covered by some thread $T_j(b_2)$ where j<=p, and so:

$ST_j(b_2) <= y <= ET_j(b_2)-1$

=> $y + ja_3 = c_2a_2 + c_1$  $c_2+c_1<=n+j-1$  is soluble for some j<=p

=> y is not a break

This is a contradiction, and so y must be a canonical break.

Theorem 220 showed that when $a_3<2a_2$ there is a possibility of a series of breaks, where each is separated from the next by $a_3-a_2$. The next theorem shows that in this case either all of the breaks in the series are canonical, or the break order decreases as the break value increases.

Theorem 229

Let $y_0$, $y_1$, ... be a series of breaks in a stride generator SG(n,p) where $a_3<2a_2$, with $y_0$ a fundamental break.

a) If $y_i$ is canonical, then so is $y_{i+1}$

b) If $y_i$ has break order q, then $y_{i+1}$ has break order q-1

Proof

Let $y = y_i$ and $y' = y_{i+1}$ be two consecutive members of the series:

$y' = y + a_3-a_2$

We first show that:

$y' + ja_3 = c_2a_2 + c_1$

$\iff y'-(a_3-a_2) + ja_3 = c_2a_2 + c_1 - (a_3-a_2)$

$\iff y + (j+1)a_3 = (c_2+1)a_2 + c_1$

Now suppose y is canonical and y' is not.

Then:

$y' + ja_3 = c_2a_2 + c_1$  $c_2+c_1<=n+j-1$  is soluble for some $j>p+1$

=> $y + (j+1)a_3 = (c_2+1)a_2 + c_1$  $c_2+c_1<=n+j-1$  is soluble for some $j>p+1$

Write $j'=j+1$, $c_2'=c_2+1$:

=> $y + j'a_3 = c_2'a_2 + c_1$  $c_2'+c_1<=n+j'-1$  is soluble for some $j'>p+2$

=> y is not canonical

This is a contradiction, so y is canonical => y' is canonical.

Now suppose y has break order q.

Then:

$y + qa_3 = c_2a_2 + c_1$  $c_2+c_1<=n+q-1$  is soluble

=> $y' + (q-1)a_3 = (c_2-1)a_2 + c_1$  $c_2+c_1<=n+q-1$  is soluble

Write $c_2'=c_2-1$:

=> $y' + (q-1)a_3 = c_2'a_2 + c_1$  $c_2'+c_1<=n+(q-1)-1$  is soluble

=> y' has break order $q'<=q-1$.

Now suppose that $q'<q-1$.

Then:

$y' + q'a_3 = c_2a_2 + c_1$  $c_2+c_1<=n+q'-1$  is soluble

=> $y + (q'+1)a_3 = (c_2+1)a_2 + c_1$  $c_2+c_1<=n+q'-1$  is soluble

Write $c_2'=c_2+1$:

=> $y + (q'+1)a_3 = c_2'a_2 + c_1$  $c_2'+c_1<=n+(q'+1)-1$  is soluble

=> y has a break order $<= q'+1 < q$

This is a contradiction, and so:

y has a break order q => y' has a break order q-1.

This final theorem brings together all the results proved above, and provides the rationale for the "potential cover" concept introduced in section 2.10.

Theorem 230

The breaks in a stride generator SG(n,p) are either all canonical or all non-canonical.

If they are non-canonical, then the first break has a higher break order than any succeeding break.

Proof

This follows straightforwardly from theorems already proved; there are three cases to consider:

a) p=0: (Theorem 217)

Every order 0 stride generator is canonical.

b) p>0, $a_3>2a_2$: (Theorems 220, 227, 228)

Such a stride generator has either one or two breaks; if it has two breaks, both are canonical.

c) p>0, $a_3<2a_2$: (Theorems 220, 227, 228, 229)

In this case, the stride generator has either one or two series of breaks, where each series is of the kind defined by Theorem 220 with each break separated from the next by $C_1 = a_3-a_2$.

If there are two series, all the breaks are canonical.

If there is just one series, then either all the breaks are canonical, or none are and each break has a lower break order than the preceding one.

[ Further clarification for SG(n,p), p>0:

Every stride generator has one or two fundamental breaks $y_0$ (i.e. which satisfy $a_3-a_2 \le y_0 < (a_3-a_2)+n$). [Theorem 227]

If $a_3>2a_2$, these are the only breaks. [Theorem 220]

If $a_3<2a_2$, each fundamental break $y_0$ gives rise to a series of breaks $y_0$, $y_1$, ... where $y_{i+1} = y_i + (a_3-a_2)$. [Theorem 220]

If a stride generator has two fundamental breaks, both are canonical. [Theorem 228]

If $a_3<2a_2$ and $y_i$, $y_{i+1}$ are two consecutive breaks in the same series, then:

$y_i$ is canonical => $y_{i+1}$ is canonical
$y_i$ has break order q => $y_{i+1}$ has break order q-1 ]

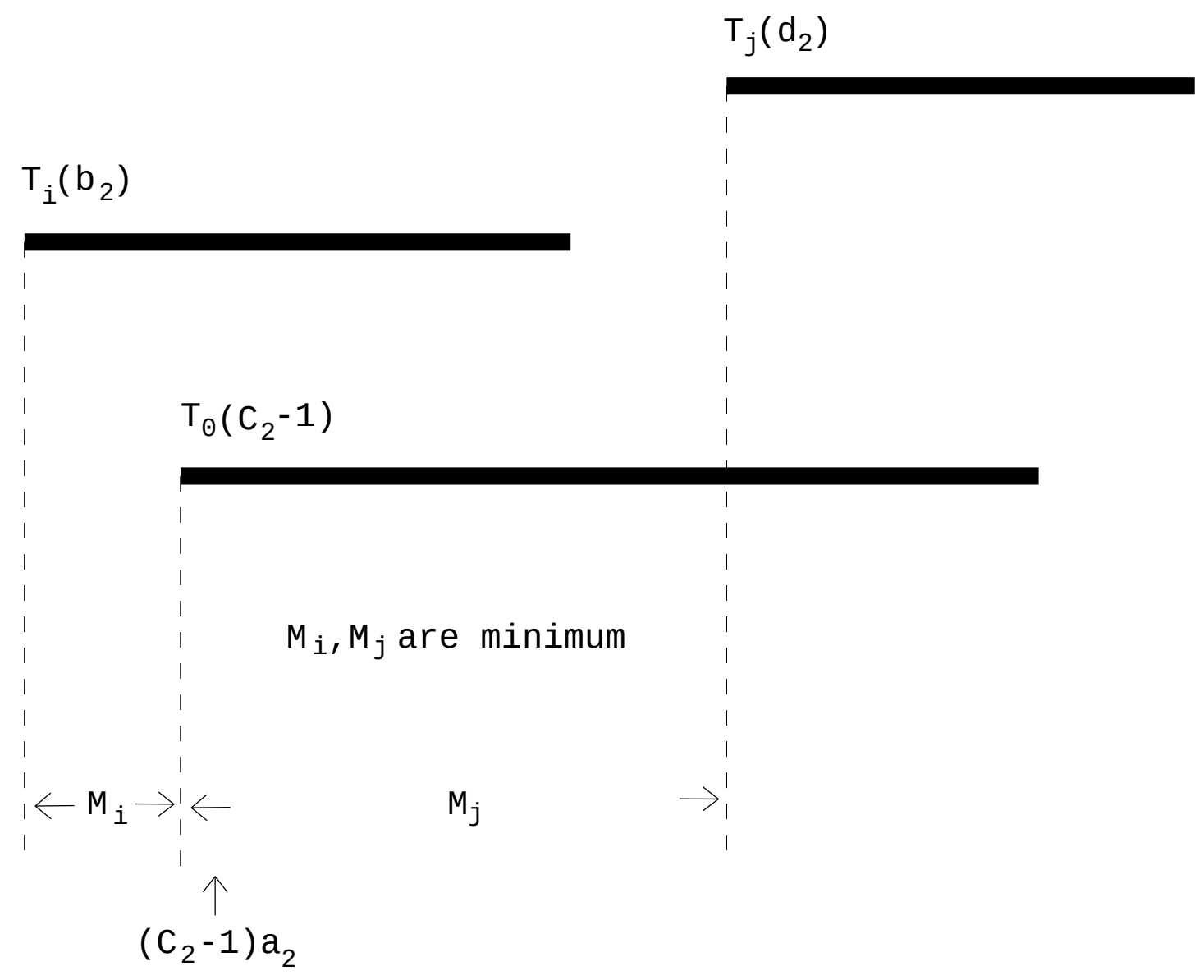

Figure 210

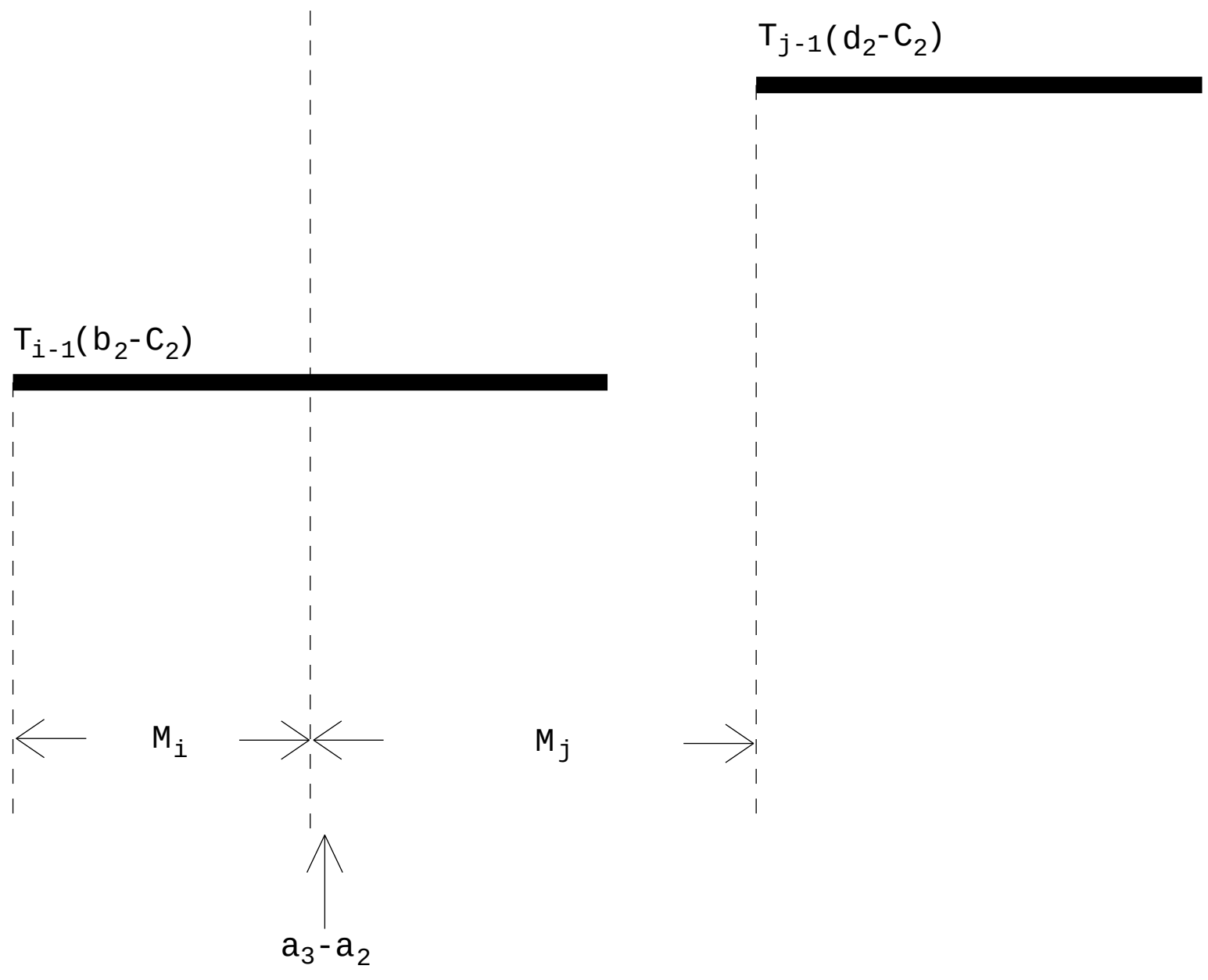

Figure 211

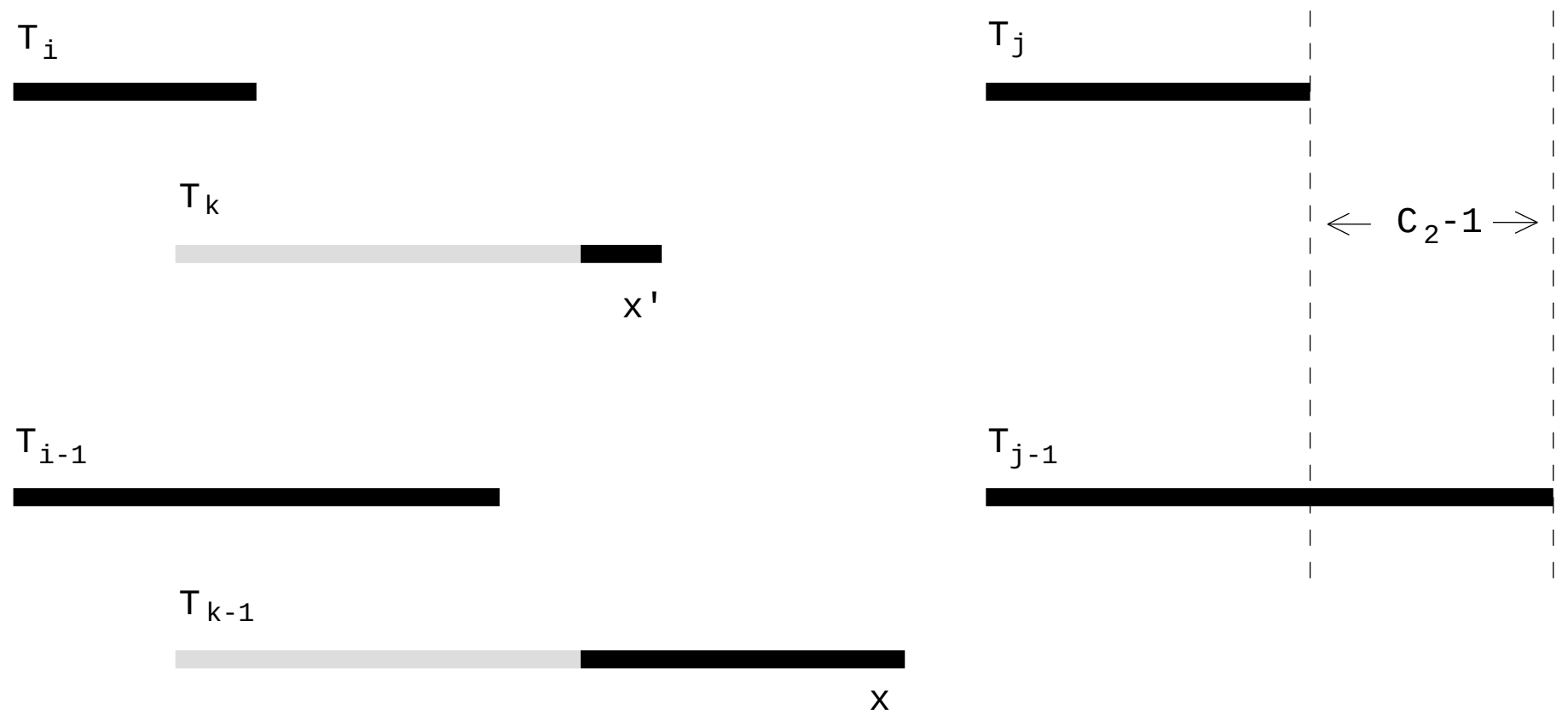

Figure 212

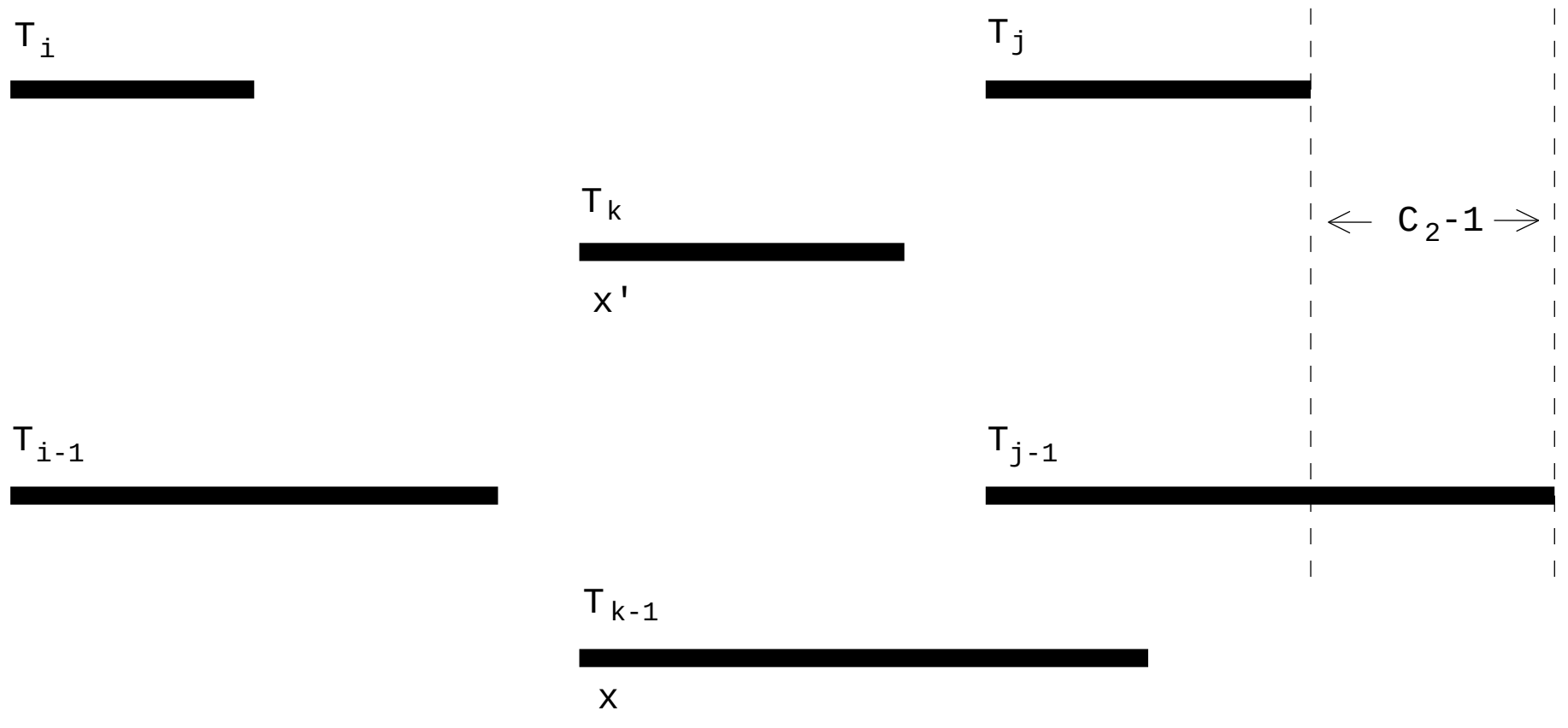

Figure 213

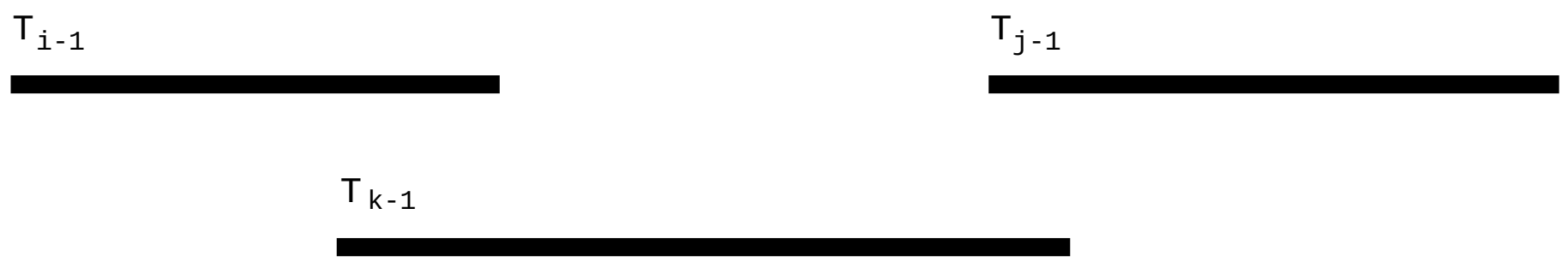

Figure 214

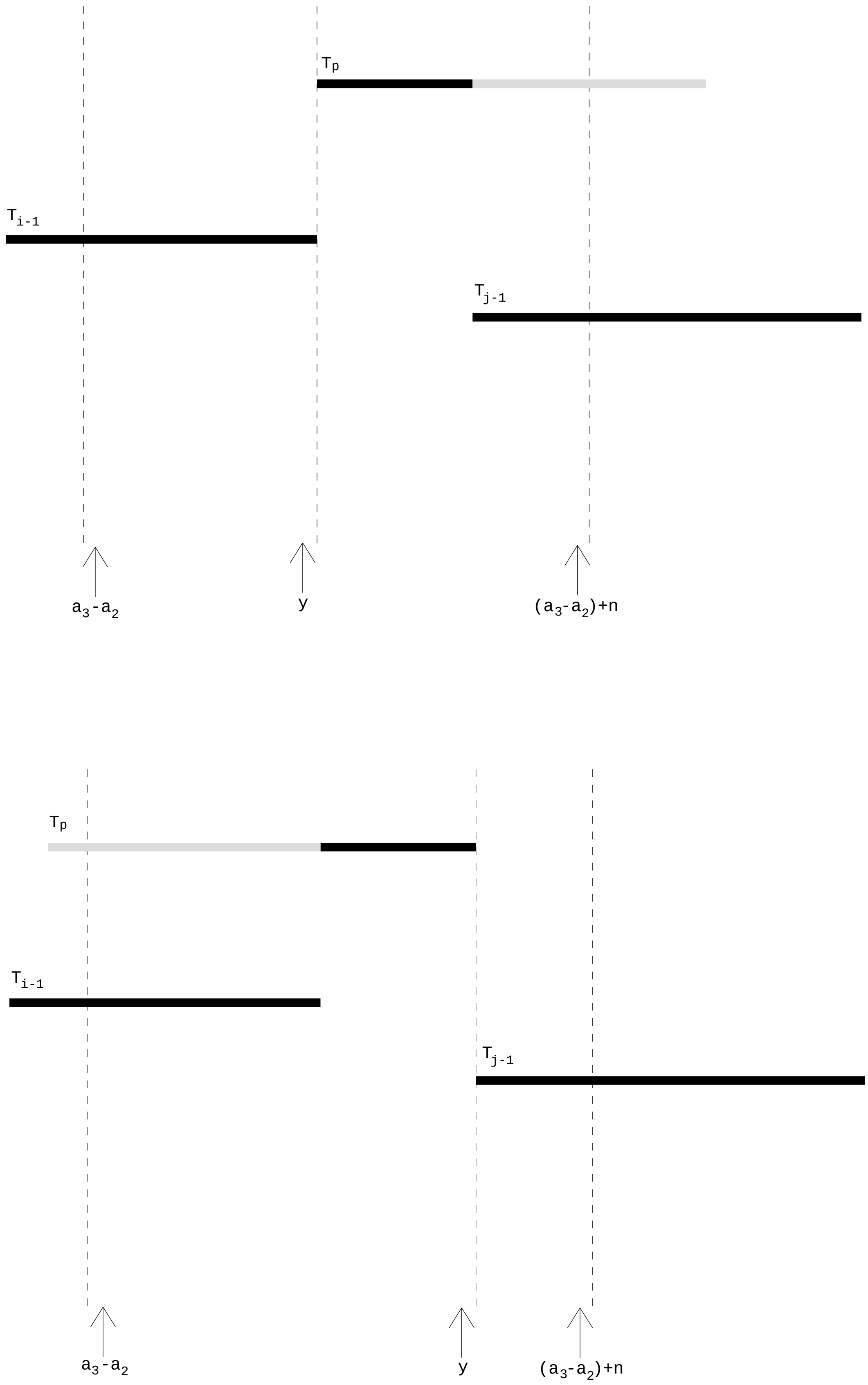

Figure 215

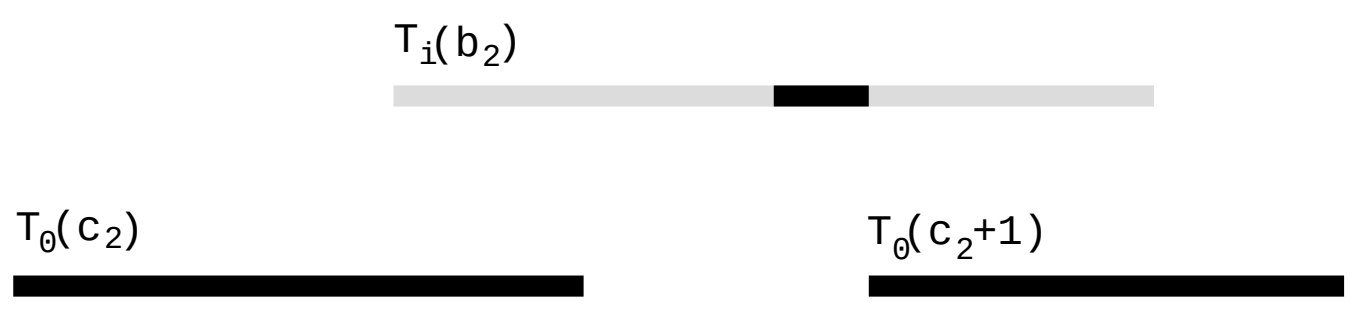

Figure 216

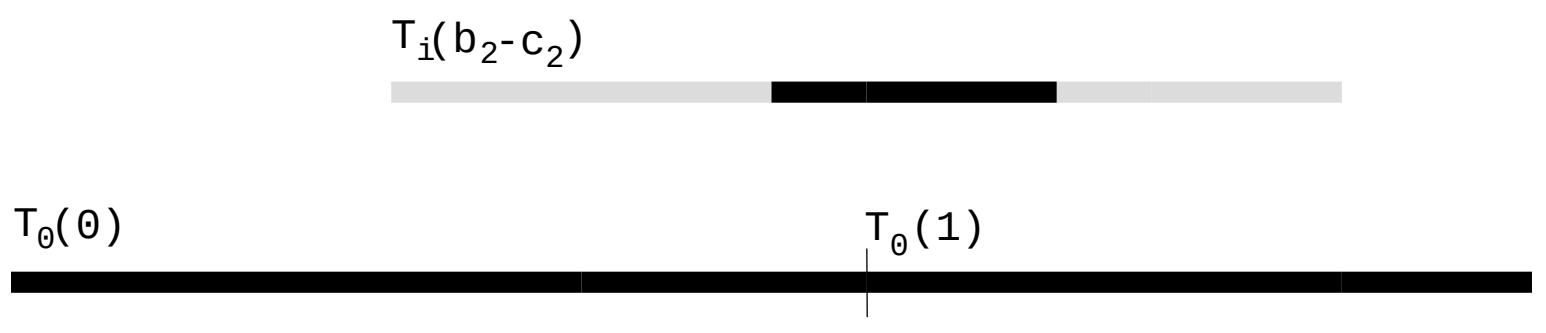

Figure 217

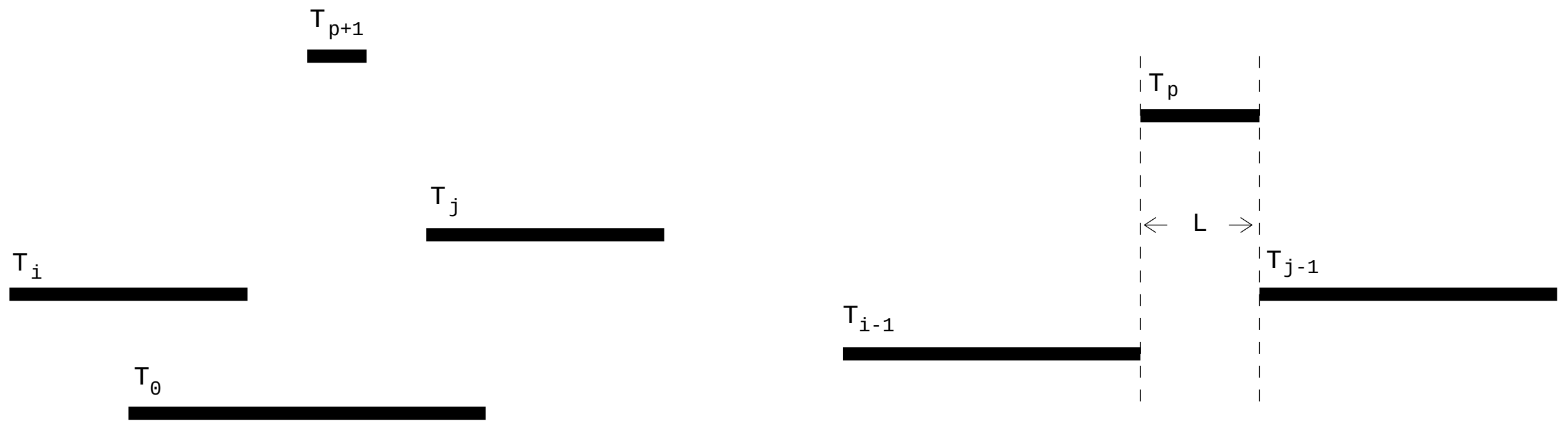

Figure 218

## 2.10 Potential covers

This section tidies up the loose end left in Section 2.3 by extending the scope of Theorem 203 to cover non-canonical stride generators as well.

Definition 208

The "potential cover" of a stride generator SG(n,p) with respect to a value s>=n is defined to be:

$$(s-n+1)a_3 + y-1$$

where y is the first break in SG(n,p).

We now show that if SG(n,p) is the stride generator underlying a cover C(A,3,s), then its potential cover is equal to the actual cover.

Theorem 231

Let SG(s-k,p) be the stride generator underlying a non-trivial cover C(A,3,s), with potential cover P defined as:

$$P = (k+1)a_3 + y-1$$

where y is the first break in SG(s-k,p).

Then:

$$C(A,3,s) = P$$

Proof

Theorem 203 proves this result for canonical stride generators, so we need only consider the case where SG(s-k,p) is non-canonical.

By Theorem 202 we know that:

$$C(A,3,s) = (k+1)a_3 + y-1$$

where y is the first break in SG(s-k,p) with break order > k+1.

But Theorem 230 proves that the first break in any non-canonical stride generator always has higher break order than any subsequent breaks, and so y must be the first break (of any order) in SG(s-k,p).

Therefore the actual cover C(A,3,s) is equal to the potential cover in all cases.

## 2.11 Stride generator series

It is perhaps not obvious that some sets define more than one stride generator.

[For example, {1,38,97} defines three stride generators as follows:

SG(19,2) with a break at 71 of order 4
SG(15,4) with a break at 67 of order 6
SG(14,6) with a canonical break at 67

A more pathological example - illustrated in Figure 219 - is {1,30,38} which defines three stride generators each with multiple breaks:

SG(8,3) with breaks at 13(6)*, 21(5)
SG(6,6) with breaks at 11(10), 19(9), 27(8)
SG(4,10) with canonical breaks at 9,11,17,19,25,27,33,35

* the break order is given in brackets after the break value]

In this section we develop some of the properties of such series of stride generators. The section is not relevant to the overall proof, and can be omitted if desired.

We present the summary theorem first.

Theorem 232

Every set A has a sequence of one or more stride generators:

$SG(n_1,p_1)$
$SG(n_2,p_2)$
...

$SG(n_t,p_t)$

where $n_{i+1}<n_i$, $p_{i+1}>p_i+1$ for $1<i<=t$ and only $SG(n_t,p_t)$ is canonical.

Proof

Theorem 237 shows that every set A defines a unique canonical stride generator $SG(n_t,p_t)$, so if A defines no other stride generators then the theorem is proven.

Suppose A defines other non-canonical stride generators $SG(n_i,p_i)$ for $i<t$.

Then Theorem 235 shows that these can be ordered so that $n_{i+1}<n_i$ and $p_{i+1}>p_i+1$, and Theorem 238 shows that $n_t<n_i$ and $p_t>p_i+1$ for all $i<t$.

Finally, we know that any such series of stride generators must be finite in length: if the order of the canonical stride generator is $p_t$, there can be no more than $p_t/2+1$ members since the order of each member is at least two less than the order of the following member.

There now follow the theorems referred to above.

Theorem 233

Every set A defines at least one stride generator SG(n,p).

Proof

Let A = {1, $a_2$, $a_3$} and choose s = $a_3$; then C(A,3,s) is a non-trivial cover since all values from 1 to $a_3$ can be generated using unit stamps only.

By Theorem 202, A is the stride generator SG(s-k,p) underlying the cover C(A,3,$a_3$) where $0<=k<a_3$ and $p<=k<a_3$.

Theorem 234

Given $a_2$, $a_3$, n and any $0<y<a_3$, there is a value J such that:

$$y + ja_3 = c_2a_2 + c_1 \quad c_2+c_1<=n+j-1 \quad \text{is not soluble for any } j>=J$$

Proof

Since $a_3>a_2$, $c_2$ increases faster than j, so we can find a value J such that:

$$y + Ja_3 > c_2a_2 \quad \text{and } c_2>n+J-1$$

Clearly $y + ja_3 > c_2a_2$ and $c_2>n+j-1$ for all j>=J as well.

We can identify a lower bound for J as follows:

We need to find integer j and integer $c_2$ such that:

$$y + ja_3 = c_2a_2 + c_1 \quad \text{with } c_2>n+j-1$$

This will certainly be so if:

$$y + ja_3 >= (n+j)a_2$$

and this is true for all y if:

$$ja_3 >= (n+j)a_2$$

$$=> j(a_3-a_2)/a_2 >= n$$

$$=> j >= na_2/(a_3-a_2)$$

Suppose $na_2/(a_3-a_2) = j'+x$ where j' is an integer and $0<=x<1$

Then J'=j'+1 is an integer that can serve as a lower bound for J:

$$J' = \text{intpt}[\ na_2/(a_3-a_2) + 1\ ]$$

$$=> J' = \text{intpt}[\ ((n-1)a_2+a_3)/(a_3-a_2)\ ]$$

[ We can check this result backwards as follows:

Let $j = ((n-1)a_2+a_3)/(a_3-a_2)$ => $j-1 < J' <= j$

Then:

$$y + J'a_3 > J'a_3 > (j-1)a_3 = na_2a_3/(a_3-a_2)$$

and:

$$(n+J'-1)a_2 <= (n+j-1)a_2 = ((n-1)(a_3-a_2) + (n-1)a_2 + a_3)a_2/(a_3-a_2)$$

$$= ((n-1)a_3+a_3)a_2/(a_3-a_2)$$

$$= na_2a_3/(a_3-a_2)$$

So $y + J'a_3 > na_2a_3/(a_3-a_2) >= (n+J'-1)a_2$ ]

Theorem 235

If A defines two stride generators SG(n,p) and SG(n',p') with n'<n, then p'>p+1.

[Note that n'=n => p'=p because by (SG2) "p is minimal for n"]

Proof

(SG3) for SG(n,p) to be a stride generator => there exists $0<y<a_3$ such that:

$y + ja_3 = c_2a_2 + c_1$ $c_2+c_1<=n+j-1$ is not soluble for j<=p+1

=> $y + ja_3 = c_2a_2 + c_1$ $c_2+c_1<=n'+j$ is not soluble for j<=p+1 since n'<n

(SG1) for SG(n',p') to be a stride generator =>

$y + ia_3 = c_2a_2 + c_1$ $c_2+c_1<=n'+i$ is soluble for some i<=p'

These two statements can only both be true if p'>p+1.

Theorem 236

Every set A has a canonical stride generator.

Proof

Let SG(n,p) be a stride generator underlying some cover C(A,3,s); note that by Theorem 233 we know that at least one such stride generator exists for every set A.

If SG(n,p) is canonical, the theorem is proven; so suppose that it is not.

By Theorem 202, there exists a break $0<y<a_3$ such that:

$y + ja_3 = c_2a_2 + c_1$ $c_2+c_1<=n+j-1$ is not soluble for any j<=s-n+1

This means that s-n+1 < break order of y.

Let q be the maximum break order of all breaks in SG(n,p) - in other words, the order of the first break by Theorem 229; then:

$s-n+1 < q$

$\Rightarrow s < n+q-1$

In other words, if the cover C(A,3,s) defines a non-canonical stride generator SG(n,p), then s<n+q-1.

We now choose s'>=n+q-1.

By Theorems 201 and 202:

$C(A,3,s') = (k'+1)a_3 + Y$ defines a stride generator SG(n',p') for n'=s'-k'

We now show that k'>k, where k=s'-n:

a) We know we can generate strides 0 ... s-n using at most s<s' stamps, because the cover C(A,3,s) defines SG(n,p).

   So we can generate all values in stride (s-n)+i using at most s+i stamps - by simply adding i $a_3$ stamps to the s-n generations.

   Therefore we can generate strides 0 ... s'-n  =>  k'>=k.

b) Now consider stride k+1.

   Because SG(n,p) is not canonical, we know that for all $0<y<a_3$:

$$y + ja_3 = c_2a_2 + c_1 \quad c_2+c_1<=n+j-1 \quad \text{for some } j<=q$$

$$\Rightarrow y + (k+1)a_3 = (k+1-j)a_3 + c_2a_2 + c_1$$

   This is a generation if:

   i) $k+1-j >= 0$

   and ii) $k+1-j+c_2+c_1 <= s'$

   (i): $k+1-j = s'-n+1-j$

   $>= s'-n+1-q = s'-(n+q-1) >= 0$

   (ii): $k+1-j+c_2+c_1 = s'-n+1-j+c_2+c_1$

   $<= s'-n+1-j+n+j-1 = s'$

   So all values in stride k+1 can be generated  =>  k'>k.

Now k'>k  =>  s'-n' > s'-n

=>  n'<n

=>  p'>p+1     by Theorem 235

This shows that given a non-canonical stride generator SG(n,p) underlying some cover C(A,3,s) we can derive another stride generator SG(n',p') underlying some cover C(A,3,s') for s'>s where:

n'<n  and  p'>p+1

If SG(n',p') is canonical, then the theorem is proved.

If not, we repeat the process until we derive some stride generator that is canonical.

How do we know that the process does not simply generate a (finite or infinite) sequence of non-canonical stride generators?

Firstly, we know that the sequence is finite:

a) $C(A,3,s) <= sa_3$, and so $k<s \Rightarrow n>=1$.

b) Each iteration causes n to decrease.

Secondly, the final member of the sequence must be canonical - for if it were not, we could choose $s>=n+q-1$ and complete a further iteration to derive another member.

Corollary

If SG(n,p) is a canonical stride generator underlying the cover C(A,3,s), then SG(n,p) will also be the stride generator underlying all covers C(A,3,s') for s'>s.

Proof

By Theorems 201 and 202:

$C(A,3,s') = (k'+1)a_3 + Y$  defines a stride generator SG(n',p') for n'=s'-k'

Let k=s'-n; we now show that k'=k and hence n'=n:

a) We know we can generate strides 0 ... s-n using at most s<s' stamps, and so we can generate strides 0 ... s'-n using at most s' stamps (see proof of Theorem above).

Therefore k'>=k.

b) Consider stride k+1.

SG(n,p) is canonical, and so there exists $0<y<a_3$ such that:

$y + ja_3 = c_2a_2 + c_1 \quad c_2+c_1<=n+j-1$  is not soluble for any j

$\Rightarrow y + (k+1)a_3 = (k+1-j)a_3 + c_2a_2 + c_1 \quad c_2+c_1<=n+j-1$ is not soluble for any j

But:  $k+1-j+c_2+c_1 <= s'-n+1-j+n+j-1 = s'$

so $y + (k+1)a_3 = c_3a_3 + c_2a_2 + c_1 \quad c_3+c_2+c_1<=s'$  is not soluble for any $c_3$

So stride k+1 cannot be generated and k'<k+1.

So k'=k => n'=n => p'=p  because of (SG2) ["p is minimal for n"].

Theorem 237

Each set A defines a unique canonical stride generator SG(n,p).

Proof

Suppose this is not so, and that A = SG(n,p) and SG(n',p') where (n,p) is different from (n',p') and both stride generators are canonical.

Suppose that n' and n are not equal, and choose n'>n:

Then by (SG3) for SG(n',p') there exists $0<y<a_3$ such that:

$y + ja_3 = c_2a_2 + c_1$  $c_2+c_1<=n'+j-1$  is not soluble for any j
(because SG(n',p') is canonical)

n<n' => n+j<=n'+j-1, so:

$y + ja_3 = c_2a_2 + c_1$  $c_2+c_1<=n+j$  is not soluble for any j

But by (SG1) for SG(n,p):

$y + ia_3 = c_2a_2 + c_1$  $c_2+c_1<=n+i$  is soluble for some i<=p

These two statements are in contradiction, so n'=n => p'=p because of (SG2).

Theorem 238

Suppose A defines two stride generators SG(n,p) and SG(n',p'), the first of which is canonical, and the other of which is not.

Then n<n' and p>p'+1.

Proof

Consider any break y of the canonical stride generator SG(n,p); by (SG3):

$y + ja_3 = c_2a_2 + c_1$  $c_2+c_1<=n+j-1$  is not soluble for any j    - (1)

Now suppose n'<n, and consider any value $0<y<a_3$; by (SG1):

$y + ia_3 = c_2a_2 + c_1$  $c_2+c_1<=n'+i$  is soluble for some i<=p'

=> $y + ia_3 = c_2a_2 + c_1$  $c_2+c_1<=n+i-1$  is soluble for some i<=p' since n>n'

This contradicts (1) and so the assumption that n'<n is wrong.

We know from Theorem 237 that n is not equal to n', and so n<n' and hence p>p'+1 by Theorem 235.

[Although every canonical stride generator underlies some cover (this follows from Theorems 236 and 237), this is not true of non-canonical ones.

For example, consider A = {1,8,11} = SG(3,2) and SG(2,4), where only the latter is canonical.

The smallest non-trivial cover C(A,3,7) has SG(2,4) as its underlying stride generator, and by the corollory to Theorem 236 all covers C(A,3,s) for s>7 have the same underlying stride generator.

So no cover has {1,8,11} = SG(3,2) as its underlying stride generator.]

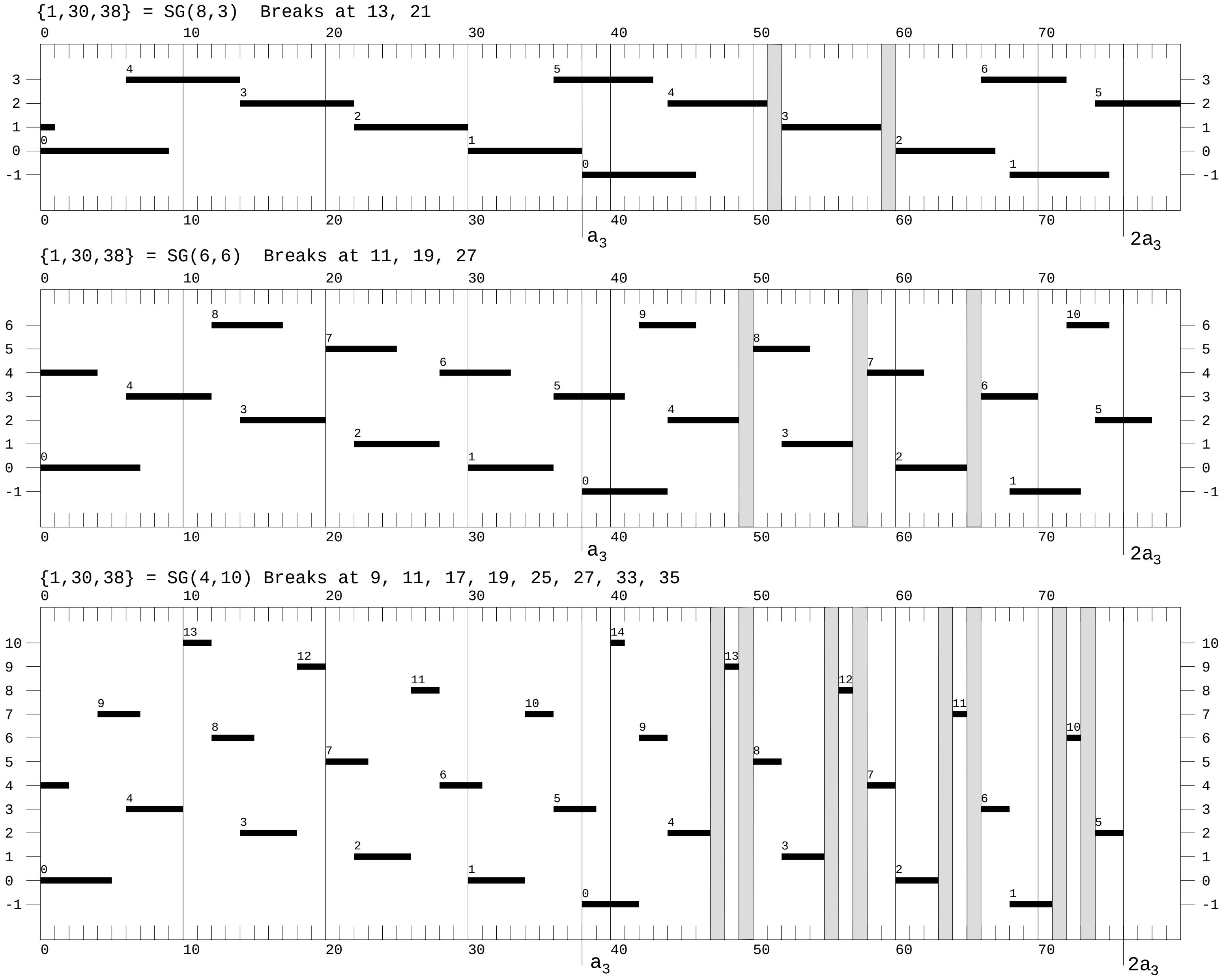
{1,30,38} = SG(8,3)  Breaks at 13, 21
$a_3$
$2a_3$
{1,30,38} = SG(6,6)  Breaks at 11, 19, 27
$a_3$
$2a_3$
{1,30,38} = SG(4,10) Breaks at 9, 11, 17, 19, 25, 27, 33, 35
$a_3$
$2a_3$

Figure 219

## 2.12 Limits

This section presents theorems that give limits on the values $a_2$ and $a_3$ for a stride generator SG(n,p). These are referenced from other sections of the proof (for example, where we show that the order of the stride generator underlying a maximal cover is less than 3). They are also useful when developing programs to produce exhaustive lists of stride generators of a given order or length, although they are quite "weak" limits.

Theorem 239

If $A = \{1, a_2, a_3\}$ is a stride generator SG(n,p), then:

$$a_2 >= a_3(p+1)/(n+p+1)$$

Proof

By Theorem 221, at least one value $a_3-a_2<=x<a_3$ has an order p generation:

$$x + pa_3 = c_2a_2 + c_1 \quad c_2+c_1<=n+p$$

$$=> c_2a_2 = x + pa_3 - c_1$$

$$=> c_2a_2 >= a_3-a_2+pa_3 \quad \text{since } c_1>=0$$

$$=> (c_2+1)a_2 >= (p+1)a_3$$

$$=> a_2 >= (p+1)a_3/(c_2+1)$$

The rhs is smallest when $c_2+1$ is largest, so:

$$=> a_2 >= a_3(p+1)/(n+p+1) \quad \text{since } c_2<=n+p$$

Theorem 240

If $A = \{1, a_2, a_3\}$ is a stride generator SG(n,p), then:

$$a_2 <= n(p+1)+1$$

Proof

Consider how values $0<=x<a_2$ can be generated by SG(n,p): how many could have order i generations?

The distance between the start points of two consecutive threads of order i is $a_2$:

$$ST_i(k+1) - ST_i(k) = (k+1)a_2 - ia_3 - ka_2 + ia_3 = a_2$$

So the best that we can hope for is that for some value of k $T_i(k)$ lies wholly within the range $0<=x<a_2$, in which case:

$$L_i(k) = n+i-k+1$$

values would have order i generations.

What can we say about k?

$ST_0(0) = 0$ and has length n+1; for i>0, $ST_i(k)$ must therefore be greater than 0 if it is to contribute effectively to the cover:

$ST_i(k) = ka_2 - ia_3 > 0$

$\Rightarrow ka_2 > ia_3$

$\Rightarrow k > ia_3/a_2$

$\Rightarrow k > i$ since $a_3 > a_2$

$\Rightarrow k >= i+1$ since k is an integer

Therefore we know that at best:

$n+i-k+1 <= n$

values in the range $0<=x<a_2$ have order i generations for i>0.

Now the best conceivable cover for this range would occur if the threads of order 0, 1, ... p all fitted neatly together one after the other (in some order): in other words, every value in the range had exactly one generation.

In this case, the number of values covered would be <= n(p+1)+1, since there is one thread of length n+1 followed by p threads each of which covers at most n values.

This proves that $a_2 <= n(p+1)+1$.

[An example of a stride generator which exhibits such an efficient cover for $a_2$ is SG(2,3) = {1,9,11}; this is shown in Figure 220.]

Theorem 241

If A = {1, $a_2$, $a_3$} is a stride generator SG(n,p), then:

$a_3 <= n(n+p+1) + (n+p+1)/(p+1)$

Proof

By Theorem 239:

$a_2 >= a_3(p+1)/(n+p+1)$

$\Rightarrow a_3 <= (n+p+1)a_2/(p+1)$

Using the upper bound for $a_2$ from Theorem 240 we have:

$a_3 <= (n+p+1)(n(p+1)+1)/(p+1)$

$\Rightarrow a_3 <= n(n+p+1) + (n+p+1)/(p+1)$

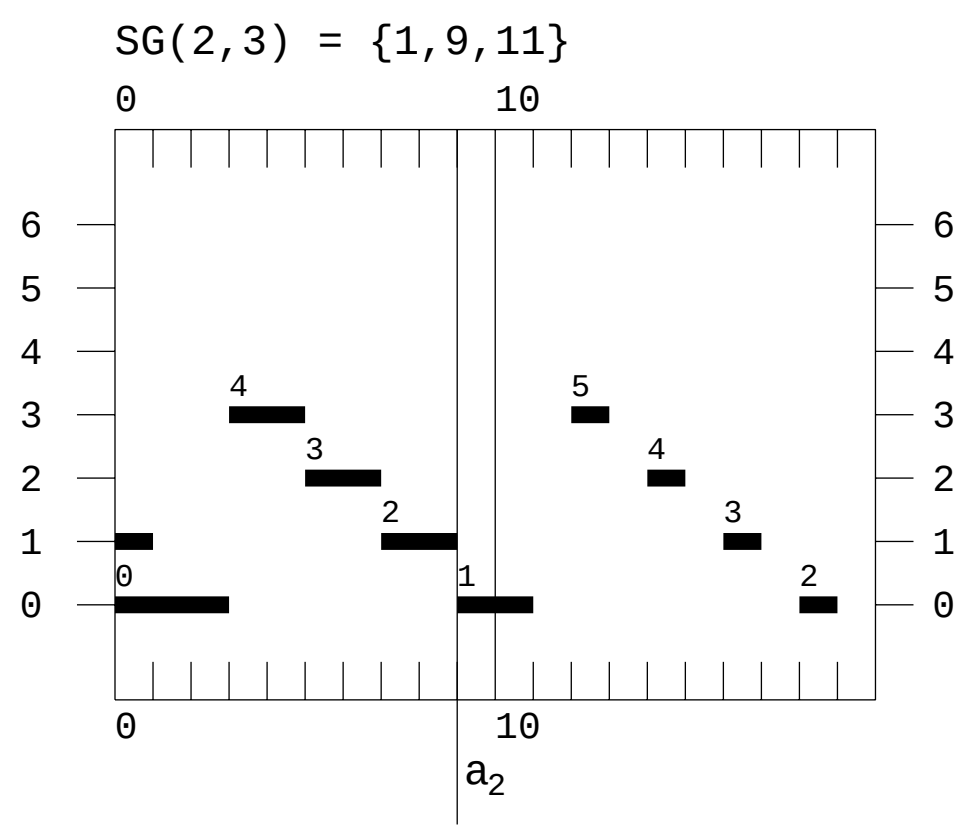

Figure 220

## 2.13 Miscellaneous results

The following result is unrelated to the proof as a whole, but is quite interesting: we show that we can construct an n-stride generator that is as long as we wish it to be.

### Theorem 242

Given any n and any X, we can find p such that the length of SG(n,p) is greater than or equal to X.

### Proof

We will show that A = {1, (k-1)n, kn} is an n-stride generator SG(n,p) for some p<=k-1, with k>1 and n>1.

Given n and X, we then choose k such that kn>=X, and the theorem is proven.

[Figure 221 shows the thread diagram for an example of such a stride generator with n=5 and k=20]

Consider the thread $T_i(i+1)$; this runs from $ST_i(i+1)$ to $ET_i(i+1)$ inclusive where:

$$ST_i(i+1) = (i+1)a_2 - ia_3$$

$$= (i+1)(k-1)n - ikn$$

$$= ((k-1)-i)n$$

$$ET_i(i+1) = ST_i(i+1) + (n+i-(i+1))$$

$$= ((k-1)-i)n + (n-1)$$

$$= (k-i)n - 1$$

Now let i range from k-1 down to 0:

| | | | | |
|---|---|---|---|---|
| $T_{k-1}(k)$ | covers | 0 | ... | n - 1 |
| $T_{k-2}(k-1)$ | covers | n | ... | 2n - 1 |
| | | .... | | |
| $T_1(2)$ | covers | (k-2)n | ... | (k-1)n - 1 |
| $T_0(1)$ | covers | (k-1)n | ... | kn - 1 |

So these threads together cover all values $0<x<a_3$ using generations of order less than or equal to k-1, thus showing that conditions (SG1) and (SG2) for an n-stride generator of some order p<=k-1 are met.

Now consider the value y = (k+2)n-1; clearly kn < y < 2kn since k>1. We will show that $y-a_3$ is a break for SG(n,p).

Consider any thread $T_i(m)$ where:

$$ST_i(m) = ma_2 - ia_3$$

$$= mn(k-1) - ikn$$

$$= (mk-m-ik)n \qquad - \quad (1)$$

$ET_i(m) = ST_i(m) + (n+i-m)$ - (2)

This thread can only cover y if $ET_i(m) >= y$:

=> (mk-m-ik)n + (n+i-m) >= (k+2)n-1

=> m(kn-n-1) >= i(kn-1) + (k+1)n-1

=> m >= i(kn-1)/(kn-n-1) + (kn+n-1)/(kn-n-1)

=> m > i+1 - since n>1

So from (1) and (2), any thread that covers y must start on an integral multiple of n, say qn, and runs as far as:

qn + (n+i-m) < qn + n+i-(i+1) = qn+n-1

But y is of the form qn+n-1 - so there is no thread of any order that covers y, and so $y-a_3$ is a break in SG(n,p).

{1,95,100} = SG(5,18)  First break at 9

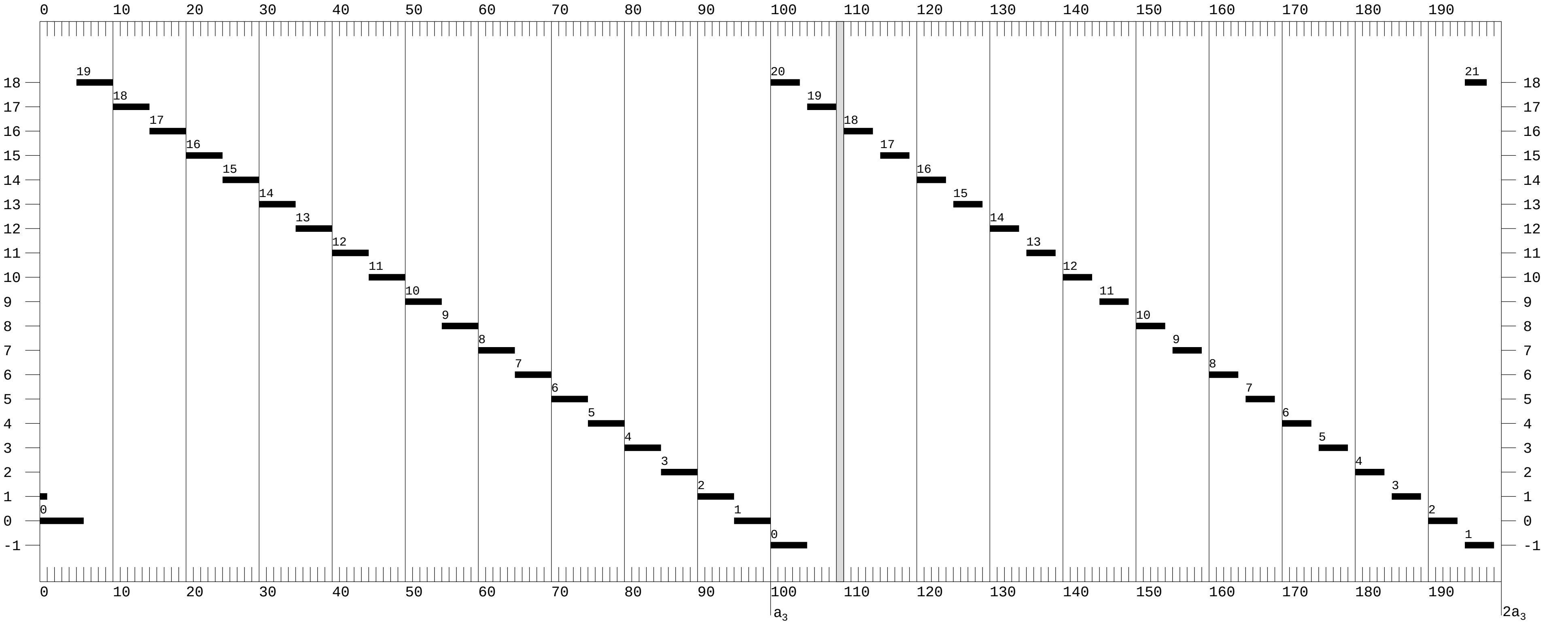

Figure 221

## 2.14 Areas for further research

This section includes notes on outstanding issues of interest.

### Conjecture 1

The stride generator underlying a cover is always canonical.

### Notes

We know that if it is not canonical, s<n+q-1 where q is the break order of the first break (see proof of Theorem 236).

With unit stamps we can make 1...s, but not s+1 - so for a cover, $s+1>=a_2$.

Now s<n+q-1 => s+1<n+q; so if we can show $n+q<=a_2$, we're there - and it does appear to be true.

Note that it also appears that q<=2p, but it is certainly not the case that $n+2p<=a_2$ - indeed, $n-a_2$ can become arbitrarily large for zero order stride generators where $n=a_2+C_2-2$.

It is obvious that for a series of stride generators for the same A, all threads start in the same positions regardless of n and p: they just get longer as n increases and shorter as n decreases. This may help.

[ I finally proved this in late 1995/early 1996 - see "A proof that $h_1,h_2<=h_0$ for any h-basis $A_3$" - and the proof does, indeed, proceed by showing that $n+q<=a_2$. Unfortunately, it is a complex proof with many special cases to be considered. ]

### Conjecture 2

If A defines SG(n,p) and SG(n',p') with p'>p, then the first break y for SG(n,p) is always >= the first break y' for SG(n',p').

[ This follows easily when we examine what happens in a thread diagram as n is reduced from $n_1$ through to $n_t$. The document referenced above defines $SG(n_1,p_1)$ as the "fundamental stride generator" of the set A, and contains a series of thread diagrams for A = {1, 38, 97} which illustrate how this Conjecture can be proved. ]

# SECTION 3

## OPTIMAL AND NEAR-OPTIMAL STRIDE GENERATORS

This section develops formulae for OSG(n,0), OSG(n,1) and SG1(n,1).

[See Definition 400 for an explanation of SG1(n,p)]

## 3.1 Optimal stride generators of order 0

Theorem 300

Formulae for OSG(n,0), the optimal n-stride generators of order 0, are as given in the following table:

The form depends on r = n mod 2 as follows:

| r | $a_2$ | $a_3$ | y |
|---|---|---|---|
| 0 | (n+2)/2 | $(n^2+6n+4)/4$ | $(n^2+4n)/4$ |
| | (n+4)/2 | $(n^2+6n+4)/4$ | $(n^2+4n-4)/4$ |
| 1 | (n+3)/2 | $(n^2+6n+5)/4$ | $(n^2+4n-1)/4$ |

Formulae for $\{1, a_2, a_3\}$ = OSG(n,0) and for the first break y

Table 300

Proof

Theorem 214 shows that $(1, a_2, a_3\}$ = SG(n,0) iff:

$n = a_2+C_2-2$

$C_1=0$ or $a_2-C_2<=C_1<a_2$ where $a_3 = C_2a_2+C_1$

with the first break at $y = C_2a_2-1$

To determine OSG(n,0) we must maximise $a_3=C_2a_2+C_1$ given n; we start by choosing $C_1$ as large as possible ($C_1=a_2-1$) giving:

$a_3 = (C_2+1)a_2 - 1$

Now $n = a_2+C_2-2$

$\Rightarrow C_2 = (n+2)-a_2$

$\Rightarrow a_3 = (n+3-a_2)a_2 - 1 = (n+3)a_2 - a_2^2 - 1$

This will be a maximum when $da_3/da_2 = 0$

$\Rightarrow a_2 = (n+3)/2$

There are two cases to consider.

a) n = 2m+1 (n odd)

Maximum occurs when $a_2=m+2$ and $C_2=m+1$, and so:

$a_3 = (C_2+1)a_2 - 1$

$= (m+2)(m+2) - 1$

$= m^2 + 4m + 3$

So in terms of n we have (for odd n):

$$OSG(n,0) = \{1, (n+3)/2, (n^2+6n+5)/4\}$$

and the first break occurs at:

$$y = (m+1)(m+2)-1 = m^2+3m+1 = (n^2+4n-1)/4$$

b) n = 2m (n even)

In this case the maximum will occur either when

$a_2 = m+1$ and $C_2 = m+1$
or $a_2 = m+2$ and $C_2 = m$

The first case gives:

$$a_3 = (m+2)(m+1) - 1$$
$$= m^2 + 3m + 1$$

with $y = (m+1)(m+1) - 1 = m^2+2m$

The second gives:

$$a_3 = (m+1)(m+2) - 1$$
$$= m^2 + 3m + 1$$

with $y = m(m+2) - 1 = m^2+2m-1$

Both possibilities are optimal stride generators, with the first having a higher initial break value.

In terms of n, we have (for even n):

$$OSG(n,0) = \{1, (n+2)/2, (n^2+6n+4)/4\} \text{ with } y = (n^2+4n)/4$$

or:

$$OSG(n,0) = \{1, (n+4)/2, (n^2+6n+4)/4\} \text{ with } y = (n^2+4n-4)/4$$

[ Note:

My original proof of this Theorem contained an error, and did not reveal the second possible form of OSG(n,0) for even n. I discovered this problem in February 1991 when I was considering an alternative (and more intuitive) proof of the theorem, but did not resolve the issue until November 1997.

The alternative proof runs roughly as follows:

If $A = \{1, a_2, a_3\}$ is an order 0 stride generator, then for every $0<x<a_3$ there is a generation:

$$x = c_2a_2 + c_1 \quad \text{with } c_2+c_1<=n$$

and this means that $C(\{1,a_2\},2,n) >= a_3-1$.

For A to be an OSG(n,0), we need $a_3$ to be as large as possible, and this is achieved by choosing $a_2$ such that $C(\{1,a_2\},2,n)$ is as large as possible. In other words, choose $A_2$ to be a maximal set for M(2,n) and let $a_3$=M(2,n)+1.

The formulae for maximal sets for d=2 are as follows:

a) s even, s=2t:

$M(2,s) = t(t+3)$ with $a_2 = t+1$ or $a_2 = t+2$

b) s odd, s=2t+1:

$M(2, s) = t(t+4)+2$ with $a_2 = t+2$

So OSG(n,0) should be:

{1, t+1, t(t+3)+1} or {1, t+2, t(t+3)+1} for even n

{1, t+2, t(t+4)+3} for odd n

These yield the formulae of Theorem 300 when rewritten in terms of n. ]

## 3.2 Optimal stride generators of order 1

### Theorem 301

Formulae for OSG(n,1), the optimal n-stride generators of order 1, are as given in the following table:

The form depends on r = n mod 3 as follows:

| r | $a_2$ | $a_3$ | y |
|---|---|---|---|
| 0 | n+4 | $(n^2+5n+6)/3$ | $(n^2+3n-9)/3$ |
| 1 | n+2 | $(n^2+5n+6)/3$ | $(n^2+3n-1)/3$ |
| 2 | n+3 | $(n^2+5n+7)/3$ | $(n^2+3n-4)/3$ |

Formulae for $\{1, a_2, a_3\}$ = OSG(n,1) and for the first break y

Table 301

### Proof

Any stride generator of order 1 must consist of a sequence of 0-threads where at least two are separated by a gap which is covered by a 1-thread; this is illustrated in Figure 300, where $a_3 = C_2a_2+C_1$:

$T_0(C_2-1)$ and $T_0(C_2)$ are the last two 0-threads in the stride generator (since $(C_2+1)a_2>a_3$), and by the corollory to Theorem 209 we know that:

$T_0(C_2)$ starting at $C_2a_2$

=> $T_1(2C_2)$ starting at $C_2a_2-C_1$

which shows that the first 1-thread must be $T_1(2C_2)$.

Theorem 220 shows that there must exist a break y satisfying:

$(a_3-a_2) <= y < (a_3-a_2)+n$

=> $(C_2-1)a_2 + C_1 <= y < (C_2-1)a_2 + C_1 + n$

Looking at Figure 300, we see that this break can only occur at the end of $T_0(C_2-1)$, $T_1(2C_2)$ or $T_0(C_2)$; we first prove that the last possibility cannot occur:

Suppose $y = ET_0(C_2)$

=> $y = C_2a_2 + n - C_2$

But $y < (a_3-a_2)+n$

=> $C_2a_2 + n - C_2 < (C_2-1)a_2 + C_1 + n$

=> $-C_2 < -a_2+C_1$

=> $C_1 > a_2-C_2$

=> $\{1, a_2, a_3\}$ is an order 0 stride generator by Theorem 214

This leaves us with two possibilities (A) and (B) as shown in Figure 301.

To find the formulae for OSG(n,1) we have to maximise $a_3$ for fixed n. We do this in two stages: first we assume fixed n and $a_2$, and determine formulae for maximum $a_3$ as functions of $a_2$, n; and then we allow $a_2$ to vary to determine maximum $a_3$ across all $a_2$.

So first we suppose n and $a_2$ are given, and we wish to make $a_3 = C_2a_2+C_1$ as large as possible. We do this by first making $C_2$ as large as possible; then if there are any alternatives for $C_1$ we will choose the largest!

Another way of saying "make $C_2$ as large as possible" is "make the length of $T_1(2C_2)$ as small as possible", since $LT_1(2C_2) = n+2-2C_2$ and so decreases as $C_2$ increases for fixed n.

Looking at both cases (A) and (B), $T_1(2C_2)$ is shortest when it fits exactly between $T_0(C_2-1)$ and $T_0(C_2)$ - ie when $LT_1(2C_2)=C_1$, and there are two fundamental breaks; for this to happen we must have:

$$ST_1(2C_2) = ET_0(C_2-1) + 1 \quad \text{and} \quad ET_1(2C_2) = ST_0(C_2) - 1$$

Now $ST_1(2C_2) = ET_0(C_2-1) + 1$

$\Rightarrow C_2a_2 - C_1 = (C_2-1)a_2 + n - (C_2-1) + 1$

$\Rightarrow n = C_2-C_1+a_2-2$ - (1)

And $ET_1(2C_2) = ST_0(C_2) - 1$

$\Rightarrow C_2a_2 - C_1 + n + 1 - 2C_2 = C_2a_2 - 1$

$\Rightarrow n = 2C_2+C_1-2$ - (2)

(1)+(2)

$\Rightarrow 2n = 3C_2+a_2-4$

$\Rightarrow 3C_2 = (2n+4)-a_2$ - (3)

2(1)-(2)

$\Rightarrow n = -3C_1+2a_2-2$

$\Rightarrow 3C_1 = (2a_2-n)-2$ - (4)

(3) gives an integer value for $C_2$ iff 3 divides $(2n-a_2)+1$, and it is easy to see that under these conditions $C_1$ will also have an integer value by (1).

But suppose $(2n-a_2)+1$ is not divisible by 3, and so there is no SG(n,1) in which the 1-thread fits exactly between the two 0-threads: in this case we have to look at the next shortest 1-thread, taking the two cases (A) and (B) separately.

(A): $ST_1(2C_2) = ET_0(C_2-1) \Rightarrow n = C_2-C_1+a_2-1$

$ET_1(2C_2) = ST_0(C_2)-1 \Rightarrow n = 2C_2+C_1-2$

$\Rightarrow 3C_2 = (2n+3)-a_2, \quad 3C_1 = 2a_2-n$

(B): $ST_1(2C_2) = ET_0(C_2-1)+1 \Rightarrow n = C_2-C_1+a_2-2$

$ET_1(2C_2) = ST_0(C_2) \Rightarrow n = 2C_2+C_1-1$

$\Rightarrow 3C_2 = (2n+3)-a_2, \quad 3C_1 = (2a_2-n)-3$

We see that both possibilities give integer values for $C_2$ and $C_1$ iff 3 divides $(2n-a_2)$, and that both also give the same value for $C_2$; but the value for $C_1$ is better in case (A) and so this is the one to choose to obtain maximum $a_3$.

But suppose neither $(2n-a_2)+1$ nor $(2n-a_2)$ is divisible by 3: we now consider the case where the 1-thread is two greater than the minimal length.

(A): $ST_1(2C_2) = ET_0(C_2-1)-1 \Rightarrow n = C_2-C_1+a_2$

$ET_1(2C_2) = ST_0(C_2)-1 \Rightarrow n = 2C_2+C_1-2$

$\Rightarrow 3C_2 = (2n+2)-a_2, \quad 3C_1 = (2a_2-n)+2$

(B): $ST_1(2C_2) = ET_0(C_2-1)+1 \Rightarrow n = C_2-C_1+a_2-2$

$ET_1(2C_2) = ST_0(C_2)+1 \Rightarrow n = 2C_2+C_1$

$\Rightarrow 3C_2 = (2n+2)-a_2, \quad 3C_1 = (2a_2-n)-4$

Both possibilities give integer values for $C_1$, $C_2$ provided 3 divides $(2n-a_2)+2$, and the best is again (A) because of its superior $C_1$ value.

So in summary, given fixed n and $a_2$, the best $a_3$ is given by one of the following formulae according to the value of $r = (2n-a_2) \bmod 3$:

| r | $3C_2$ | $3C_1$ | $3a_3 = 3(C_2a_2+C_1)$ | Type |
|---|---|---|---|---|
| 2 | $(2n+4)-a_2$ | $(2a_2-n)-2$ | $(2n+6)a_2 - a_2^2 - (n+2)$ | A,B |
| 0 | $(2n+3)-a_2$ | $(2a_2-n)$ | $(2n+5)a_2 - a_2^2 - n$ | A |
| 1 | $(2n+2)-a_2$ | $(2a_2-n)+2$ | $(2n+4)a_2 - a_2^2 - (n-2)$ | A |

Table 302

The next step is to optimise these formulae over $a_2$ to determine the $a_2$ value that yields the highest $a_3$ value for any given n.

For a given value of n there are three series of $a_2$ values according to the value of r.

[ For example, if n=9 we have the series:

| | $a_2$ | $C_2$ | $C_1$ | $a_3=C_2a_2+C_1$ | |
|---|---|---|---|---|---|
| r=0: | [ 3 | 6 | -1 | 17 ] | $C_1<0$ |
| | 6 | 5 | 1 | 31 | |
| | 9 | 4 | 3 | 39 | |
| | 12 | 3 | 5 | 41 * | Best |
| | 15 | 2 | 7 | 37 | |
| | [ 18 | 1 | 9 | 27 ] | $C_1>=a_2$ |
| r=1: | [ 2 | 6 | -1 | 11 ] | $C_1<0$ |
| | 5 | 5 | 1 | 26 | |
| | 8 | 4 | 3 | 35 | |
| | 11 | 3 | 5 | 38 * | Best |
| | 14 | 2 | 7 | 35 | |
| | [ 17 | 1 | 9 | 26 ] | $C_1>=a_2$ |
| r=2: | [ 4 | 6 | -1 | 23 ] | $C_1<0$ |
| | 7 | 5 | 1 | 36 | |
| | 10 | 4 | 3 | 43 | |
| | 13 | 3 | 5 | 44 * | Best |
| | 16 | 2 | 7 | 39 | |
| | [ 19 | 1 | 9 | 28 ] | $C_1>=a_2$ ] |

So we must find the best value of $a_2$ in each series, and then choose the best of the three optimal values.

[ In the above example, we have:

| r | $a_2$ | $a_3$ |
|---|---|---|
| 0 | 12 | 41 |
| 1 | 11 | 38 |
| 2 | 13 | 44 |

and so $a_2$=13, $a_3$=44 is the optimal solution. ]

The $a_3$ values in each series are determined by a quadratic formula, and simple differentiation shows that the maximum occurs as follows:

| r | $a_2$ |
|---|---|
| 0 | n+2.5 |
| 1 | n+2 |
| 2 | n+3 |

So in each case we need to examine permitted values around n+2, n+3 until we are certain we have identified the best one.

[ For the r=1 case above:

| | $a_2$ | $a_3$ |
|---|---|---|
| n-1 | 8 | 35 |
| n+2 | 11 | 38 |
| n+5 | 14 | 35 |

Clearly $a_2$=11 must be the best one, since the two cases on either side in the series are both inferior. ]

We have three separate cases to consider, according to the value of n mod 3.

a) n mod 3 = 0

| | $a_2$ | $3C_2$ | $3C_1$ | $3a_3=3C_2a_2+3C_1$ |
|---|---|---|---|---|
| r=0: | n | n+3 | n | $n^2+4n$ |
| | n+3 | n | n+6 | $n^2+4n+6$ |
| | n+6 | n-3 | n+12 | $n^2+4n-6$ |
| r=1: | n-1 | n+3 | n | $n^2+3n-3$ |
| | n+2 | n | n+6 | $n^2+3n+6$ |
| | n+5 | n-3 | n+12 | $n^2+3n-3$ |
| r=2: | n+1 | n+3 | n | $n^2+5n+3$ |
| | n+4 | n | n+6 | $n^2+5n+6$ |
| | n+7 | n-3 | n+12 | $n^2+5n-9$ |

Table 303

So $OSG(n,1) = \{1, n+4, (n^2+5n+6)/3\}$ if n mod 3 = 0

b) n mod 3 = 1

| | $a_2$ | $3C_2$ | $3C_1$ | $3a_3=3C_2a_2+3C_1$ |
|---|---|---|---|---|
| r=0: | n-2 | n+5 | n-4 | $n^2+4n-14$ |
| | n+1 | n+2 | n+2 | $n^2+4n+4$ |
| | n+4 | n-1 | n+8 | $n^2+4n+4$ |
| | n+7 | n-4 | n+14 | $n^2+4n-14$ |
| r=1: | n | n+2 | n+2 | $n^2+3n+2$ |
| | n+3 | n-1 | n+8 | $n^2+3n+5$ |
| | n+6 | n-4 | n+14 | $n^2+3n-10$ |
| r=2: | n-1 | n+5 | n-4 | $n^2+5n-1$ |
| | n+2 | n+2 | n+2 | $n^2+5n+6$ |
| | n+5 | n-1 | n+8 | $n^2+5n+3$ |

Table 304

So $OSG(n,1) = \{1, n+2, (n^2+5n+6)/3\}$ if n mod 3 = 1

c) n mod 3 = 2

| | $a_2$ | $3C_2$ | $3C_1$ | $3a_3=3C_2a_2+3C_1$ |
|---|---|---|---|---|
| r=0: | n-1 | n+4 | n-2 | $n^2+4n-6$ |
| | n+2 | n+1 | n+4 | $n^2+4n+6$ |
| | n+5 | n-2 | n+10 | $n^2+4n$ |
| r=1: | n-2 | n+4 | n-2 | $n^2+3n-10$ |
| | n+1 | n+1 | n+4 | $n^2+3n+5$ |
| | n+4 | n-2 | n+10 | $n^2+3n+2$ |
| r=2: | n | n+4 | n-2 | $n^2+5n-2$ |
| | n+3 | n+1 | n+4 | $n^2+5n+7$ |
| | n+6 | n-2 | n+10 | $n^2+5n-2$ |

Table 305

So $OSG(n,1) = \{1, n+3, (n^2+5n+7)/3\}$ if n mod 3 = 2

The next step is to determine the value of the first break: note that in all three cases the optimal solution is of the "r=2" type which corresponds to the situation where the 1-thread fits exactly between the two 0-threads, and so each OSG(n,1) has two fundamental breaks, and the first of these is at:

$$y = ST_1(2C_2)-1 = C_2a_2-C_1-1$$

The corresponding formulae are:

| n mod 3 | y |
|---|---|
| 0 | $(n^2+3n-9)/3$ |
| 1 | $(n^2+3n-1)/3$ |
| 2 | $(n^2+3n-4)/3$ |

To complete the proof, we must show that each of the three formulae does, indeed, define a stride generator of order 1.

To do this we must show:

a) All values $0<x<a_3$ are covered by some thread

b) The value y is, indeed, a break

a) It is clear by definition that all values up to and including $ET_0(C_2)$ are covered; we show that $ET_0(C_2)+1>=a_3$:

$$ET_0(C_2)+1-a_3 = C_2a_2 + n - C_2 + 1 - C_2a_2 - C_1$$

$$= (n+1) - (C_1+C_2)$$

| n mod 3 | $3(C_1+C_2)$ | $3((n+1)-(C_1+C_2))$ | Hypothesis true iff |
|---|---|---|---|
| 0 | 2n+6 | n-3 | n>=3 |
| 1 | 2n+4 | n-1 | n>=1 |
| 2 | 2n+5 | n-2 | n>=2 |

So all values are covered for n>=1.

b) For y to be a break we must show that:

$y + ja_3 = c_2a_2 + c_1 \quad c_2+c_1<=n+j-1$ is not soluble for j<=2

We know that this equation has no solution for j=0 and j=1, because y is not covered by a 1-thread and lies at the end of a 0-thread.

So we have only to show that:

$y + 2a_3 = c_2a_2 + c_1 \quad c_2+c_1<=n+1$ is not soluble

Now:

$$y + 2a_3 = C_2a_2 - C_1 - 1 + 2C_2a_2 + 2C_1$$

$$= 3C_2a_2 + (C_1-1)$$

| n mod 3 | $3C_2$ | $C_1$-1 | $a_2$ |
|---|---|---|---|
| 0 | n | (n+3)/3 | n+4 |
| 1 | n+2 | (n-1)/3 | n+2 |
| 2 | n+1 | (n+1)/3 | n+3 |

We consider each case separately:

n mod 3 = 0:

$y + 2a_3 = c_2a_2 + c_1 \Rightarrow c_2 = n, \quad c_1 = (n+3)/3 < a_2$

So $c_2+c_1 = n + (n+3)/3$

$> n+1$ when $n>=1$

n mod 3 = 1:

We have $c_2 = n+2, \quad c_1 = (n-1)/3 >= 0$ provided $n>=1$

Clearly, $c_2+c_1 > n+1$

n mod 3 = 2:

We have $c_2 = n+1, \quad c_1 = (n+1)/3 >= 1$ provided $n>=2$

So $c_2+c_1 > n+1$

So in each case the equation is not soluble for any permitted n>=1, and so y is indeed a break and the OSG(n,1) are indeed order 1 n-stride generators.

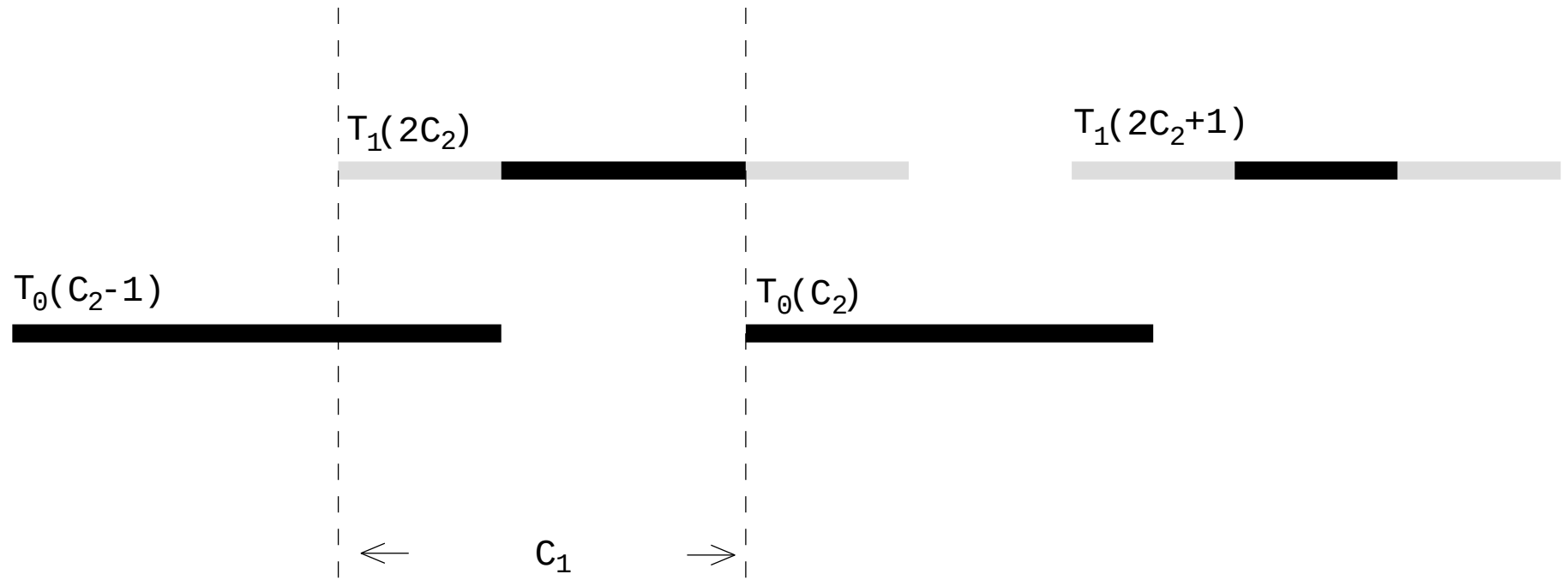

Figure 300

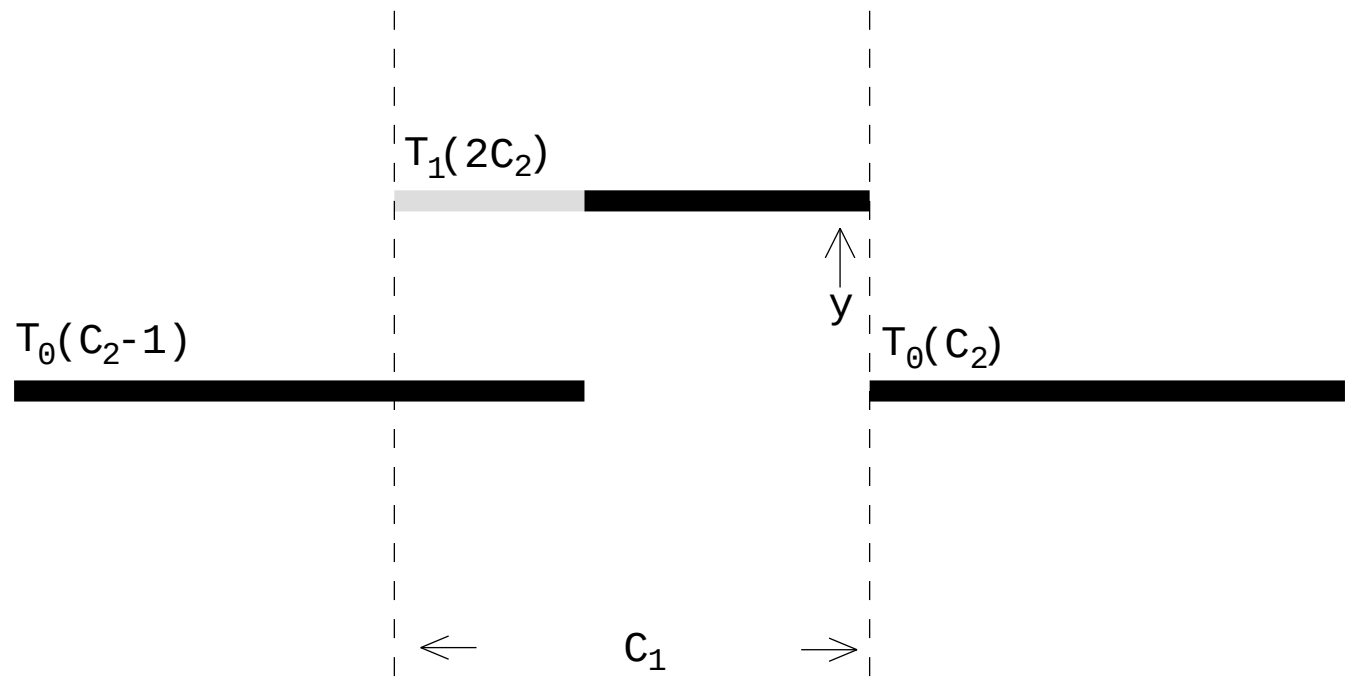

(A)

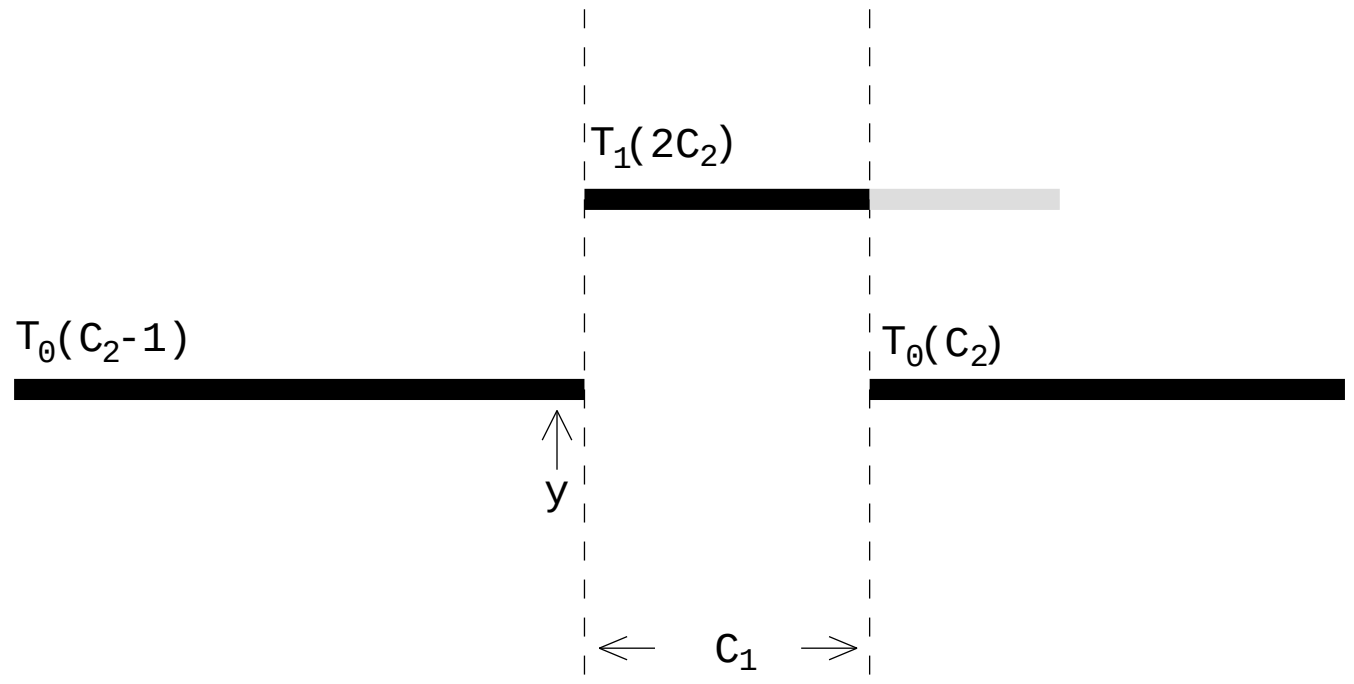

(B)

Figure 301

## 3.3 Immediately sub-optimal stride generators of order 1

Theorem 302

Formulae for SG1(n,1), the immediately sub-optimal n-stride generators of order 1, are as given in the following table for n>=3:

The form (and existence) depends on r = n mod 3 as follows:

| r | $a_2$ | $a_3$ | y |
|---|---|---|---|
| 0 | n+1 | $(n^2+5n+3)/3$ | $(n^2+3n)/3$ |
| 1 | n+5 | $(n^2+5n+3)/3$ | $(n^2+3n-16)/3$ |
| 2 | No SG1(n,1) exists | | |

Formulae for $\{1, a_2, a_3\}$ = SG1(n,1) and for the first break y

Table 306

There are two additional special cases:

SG1(2,1) = {1,4,6} and SG1(3,1) = {1,6,9}

Proof

Tables 303, 304 and 305 in the proof of Theorem 301 include details of both the best and the next best SG(n,1) for each of the three series. Examination of these tables shows that for general n the only possibilities for SG1(n,1) are as follows:

SG1(n,1) = $\{1, n+1, (n^2+5n+3)/3\}$ if n mod 3 = 0

SG1(n,1) = $\{1, n+5, (n^2+5n+3)/3\}$ if n mod 3 = 1

There is no SG1(n,1) if n mod 3 = 2

Both SG1(n,1) are of the "r=2" type, and so have two fundamental breaks with the first given by:

$$y = C_2a_2 - C_1 - 1$$

So for n mod 3 = 0 we have:

$$y = ((n+1)(n+3)-n-3)/3 = (n^2+3n)/3$$

and for n mod 3 = 1 we have:

$$y = ((n+5)(n-1)-(n+8)-3)/3 = (n^2+3n-16)/3$$

Next we must show that the two SG1(n,1) above are, indeed, order 1 stride generators (for more details of the method used, see the proof of Theorem 301):

a) We show $ET_0(C_2)+1>=a_3$:

| n mod 3 | $3(C_1+C_2)$ | $3((n+1)-(C_1+C_2))$ | Hypothesis true iff |
|---|---|---|---|
| 0 | 2n+3 | n | n>=0 |
| 1 | 2n+7 | n-4 | n>=4 |

So all permitted values of n are covered with the exception of n=1.

b) We show that y is, indeed, a break.

| n mod 3 | $3C_2$ | $C_1-1$ | $a_2$ |
|---|---|---|---|
| 0 | n+3 | (n-3)/3 | n+1 |
| 1 | n-1 | (n+5)/3 | n+5 |

So for n mod 3 = 0 we have:

$y + 2a_3 = c_2 + c_1 \Rightarrow c_2=n+3, \; c_1=(n-3)/3 >= 0$ for n>=3

So $c_2+c_1 > n+1$ for all permitted n

And for n mod 3 = 1 we have:

$y + 2a_3 = c_2 + c_1 \Rightarrow c_2=n-1, \; c_1=(n+5)/3$

Now $0<c_1<a_2$, since (n+5)/3 < (n+5), so:

$c_2+c_1 = n-1 + (n+5)/3 > n+1$ for n>=4

So y is, indeed, a break for all permitted values of n other than n=1.

We now look at the n=1 case.

Substituting in the formula, we get SG1(1,1) = {1,6,3} - which is certainly not a stride generator since $a_2>a_3$.

There remain only the special cases to deal with.

Examination of Tables 303, 304 and 305 show that the following are potential special cases for small n:

n mod 3 = 0: {1, n+3, $(n^2+4n+6)/3$} for n=3 => {1,6,9}

n mod 3 = 1: {1, n+1, $(n^2+4n+4)/3$} for n=1 => {1,2,3}
{1, n+4, $(n^2+4n+4)/3$} for n=1 => {1,5,3}*
{1, n+3, $(n^2+3n+5)/3$} for n=1 => {1,4,3}*

n mod 3 = 2: {1, n+2, $(n^2+4n+6)/3$} for n=2 => {1,4,6}

Those marked with an asterisk are immediately rejected because $a_2>a_3$; analysis of the remainder shows that {1,2,3} is an order 0 stride generator, but that the other two are valid.

In summary, SG1(n,1) are as given in Table 306 for n>=3 with two additional special cases:

SG1(2,1) = {1,4,6}
SG1(3,1) = {1,6,9}

# SECTION 4

## MAXIMISING (s-k)-STRIDE GENERATORS OVER k

This Section is concerned with optimising the choice of underlying stride generator for a given value of s.

Section 4.1 identifies the best optimal stride generator OSG($n_{opt}$,1) for a given value of s>=18.

Section 4.2 shows that no sub-optimal stride generator SGi(n,p) for any i>=2 can improve on OSG($n_{opt}$,1) for s>=36.

Section 4.3 shows that no immediately sub-optimal stride generator SG1(n,1) can improve on OSG($n_{opt}$,1) for s>=9.

## 4.1 Optimal OSG(s-k,1)

Theorem 400

Let A = {1, $a_2$, $a_3$} be an optimal stride generator OSG(s-k,1) of order 1 with potential cover:

$$X = (k+1)a_3 + y-1$$

Then given s, the maximum value of X and the corresponding value of 0<=k<s are given by $X_{opt}$ and $k_{opt}$ in Table 103 for all s>=18.

Proof

We define:

$$C_0(k) = X = (k+1)a_3 + y_0-1$$

where $a_3$ and $y_0$ are functions of s-k.

Table 101 gives the values $a_3$ and $y=y_0$ for OSG(n,1) as functions of n, and we substitute as follows:

(s-k) mod 3 = 0:

$$C_0(k) = (k+1)(s-k+3)(s-k+2)/3 + ((s-k)(s-k+3)-12)/3$$

(s-k) mod 3 = 1:

$$C_0(k) = (k+1)(s-k+3)(s-k+2)/3 + ((s-k)(s-k+3)-4)/3$$

(s-k) mod 3 = 2:

$$C_0(k) = (k+1)((s-k+3)(s-k+2)+1)/3 + ((s-k)(s-k+3)-7)/3$$

We now define $C_0'(k)$ as follows:

$$C_0'(k) = (k+1)((s-k+3)(s-k+2)+1)/3 + ((s-k)(s-k+3)-4)/3$$

$C_0'(k)$ is very similar to $C_0(k)$, but has been chosen so that it is always greater than $C_0(k)$ and is independent of (s-k) mod 3.

Expansion gives:

$$C_0'(k) = (k^3 - (2s+3)k^2 + (s^2+s-1)k + 2s^2+8s+3)/3$$

This function has the form shown in Figure 400, where the values $k_1$, $k_2$ are the solutions of $dC_0'(k)/dk = 0$:

i.e. $3k^2 - (4s+6)k + s^2+s-1 = 0$

Solving for k, we have:

$$k = ((2s+3) +/- \mathrm{sqrt}(s^2+9s+12))/3$$

We now show that $k_2$>s and 0<$k_1$<s for s>=1:

$3k_2 = (2s+3) + \mathrm{sqrt}(s^2+9s+12)$

$\Rightarrow 3(k_2-s) = (-s+3) + \mathrm{sqrt}(s^2+9s+12)$

$> (-s+3) + \mathrm{sqrt}(s^2+4s+4)$

$= (-s+3) + (s+2)$

$= 5$

$\Rightarrow k_2-s > 5/3 > 0$

$\Rightarrow k_2 > s$ for all $s>=0$

$3k_1 = (2s+3) - \mathrm{sqrt}(s^2+9s+12)$

$> (2s+3) - \mathrm{sqrt}(s^2+9s+20.25)$

$= (2s+3) - (s+4.5)$

$= s-1.5 > 0$ for all $s>=2$

Moreover, $3k_1 = 5-\mathrm{sqrt}(22) > 0$ for $s=1$; so $k_1 > 0$ for all $s>=1$.

$3k_1 = (2s+3) - \mathrm{sqrt}(s^2+9s+12)$

$\Rightarrow 3(s-k_1) = (s-3) + \mathrm{sqrt}(s^2+9s+12)$

$> (s-3) + \mathrm{sqrt}(s^2+4s+4)$

$= (s-3) + (s+2)$

$= 2s-1$

$\Rightarrow s-k_1 > (2s-1)/3 > 0$ for all $s>=1$

$\Rightarrow k_1 < s$ for all $s>=1$

It is important that $k_2>s$ and $0<k_1<s$ because this means that:

- the maximum value of $C_\Theta'(k)$ for $0 <= k < s$ is at $C_\Theta'(k_1)$
- $C_\Theta'(k)$ increases monotonically for $0 <= k < k_1$
- $C_\Theta'(k)$ decreases monotonically for $k_1 < k <= s$

It is reasonable to expect that $C_\Theta(k)$ will also have a maximum $k_1$ around s/3. To show that this is the case - and to find out the exact value of $k_1$ - we must work out values for both $C_\Theta(k)$ and $C_\Theta'(k)$ around $k = s/3$.

We can guarantee that we have found the true maximum $X_{opt} = C_\Theta(k_1)$ when:

a) $C_\Theta'(k_l) <= X_{opt}$, $C_\Theta'(k_u) <= X_{opt}$ and $k_l < k_1 < k_u$

b) $C_\Theta(k) <= X_{opt}$ for all $k_l < k < k_u$
(Note that $C_\Theta(k_l)<X_{opt}$ and $C_\Theta(k_u)<X_{opt}$ because $C_\Theta(k)<C_\Theta'(k)$ by definition)

[This is best understood by reference to Figure 401, where the crosses mark values of $C_\Theta(k)$; note that we have not proved that $C_\Theta(k)$ is monotonic before and after its maximum (although, in fact, it is).]

Unfortunately, there are nine cases to consider for $C_\Theta(k)$, and three cases for $C_\Theta'(k)$ - dependent on s mod 3 and (s-k) mod 3; the relevant formulae (developed by computer algebra (c0max)) follow:

a) $s=9t$: $k_l=3t-2$, $k_1=3t-1$, $k_u=3t+1$:

$$\begin{aligned}
C_\Theta'(3t-2) &= 36t^3 + 54t^2 + 17t - 5 \\
C_\Theta(3t-1) &= 36t^3 + 54t^2 + 22t && = X_{opt} \\
C_\Theta(3t) &= 36t^3 + 54t^2 + 22t - 2 \\
C_\Theta'(3t+1) &= 36t^3 + 54t^2 + 17t
\end{aligned}$$

b) $s=9t+1$: $k_l=3t-1$, $k_1=3t$, $k_u=3t+1$:

$$\begin{aligned}
C_\Theta'(3t-1) &= 36t^3 + 66t^2 + 35t + 2 \\
C_\Theta(3t) &= 36t^3 + 66t^2 + 36t + 4 && = X_{opt} \\
C_\Theta'(3t+1) &= 36t^3 + 66t^2 + 33t + 3 + 1/3
\end{aligned}$$

c) $s=9t+2$: $k_l=3t-1$, $k_1=3t$, $k_u=3t+1$:

$$\begin{aligned}
C_\Theta'(3t-1) &= 36t^3 + 78t^2 + 49t + 4 + 2/3 \\
C_\Theta(3t) &= 36t^3 + 78t^2 + 53t + 8 && = X_{opt}\ (t>=1) \\
C_\Theta'(3t+1) &= 36t^3 + 78t^2 + 51t + 8 + 2/3
\end{aligned}$$

d) $s=9t+3$: $k_l=3t$, $k_1=3t+1$, $k_u=3t+2$:

$$\begin{aligned}
C_\Theta'(3t) &= 36t^3 + 90t^2 + 71t + 15 \\
C_\Theta(3t+1) &= 36t^3 + 90t^2 + 71t + 15 && = X_{opt} \\
C_\Theta'(3t+2) &= 36t^3 + 90t^2 + 65t + 13
\end{aligned}$$

e) $s=9t+4$: $k_l=3t$, $k_1=3t+1$, $k_u=3t+2$:

$$\begin{aligned}
C_\Theta'(3t) &= 36t^3 + 102t^2 + 91t + 22 + 1/3 \\
C_\Theta(3t+1) &= 36t^3 + 102t^2 + 92t + 22 && = X_{opt}\ (t>=1) \\
C_\Theta'(3t+2) &= 36t^3 + 102t^2 + 89t + 23
\end{aligned}$$

f) $s=9t+5$: $k_l=3t$, $k_1=3t+1$, $k_u=3t+2$:

$$\begin{aligned}
C_\Theta'(3t) &= 36t^3 + 114t^2 + 113t + 31 \\
C_\Theta(3t+1) &= 36t^3 + 114t^2 + 116t + 36 && = X_{opt} \\
C_\Theta'(3t+2) &= 36t^3 + 114t^2 + 115t + 35 + 2/3
\end{aligned}$$

g) $s=9t+6$: $k_l=3t$, $k_1=3t+1$, $k_u=3t+3$:

$$\begin{aligned}
C_\Theta'(3t) &= 36t^3 + 126t^2 + 137t + 41 \\
C_\Theta(3t+1) &= 36t^3 + 126t^2 + 143t + 49 && = X_{opt}\ (t>=2) \\
C_\Theta(3t+2) &= 36t^3 + 126t^2 + 142t + 50 \\
C_\Theta'(3t+3) &= 36t^3 + 126t^2 + 137t + 46
\end{aligned}$$

h) $s=9t+7$: $k_l=3t+1$, $k_1=3t+2$, $k_u=3t+3$:

$$\begin{aligned}
C_\Theta'(3t+1) &= 36t^3 + 138t^2 + 171t + 65 + 1/3 \\
C_\Theta(3t+2) &= 36t^3 + 138t^2 + 173t + 68 && = X_{opt} \\
C_\Theta'(3t+3) &= 36t^3 + 138t^2 + 169t + 65 + 1/3
\end{aligned}$$

i) $s=9t+8$: $k_l=3t+1$, $k_1=3t+2$, $k_u=3t+3$:

$$\begin{aligned}
C_\Theta'(3t+1) &= 36t^3 + 150t^2 + 201t + 82 + 2/3 \\
C_\Theta(3t+2) &= 36t^3 + 150t^2 + 204t + 86 && = X_{opt}\ (t>=2) \\
C_\Theta'(3t+3) &= 36t^3 + 150t^2 + 203t + 88
\end{aligned}$$

Note: t>=0 wherever no constraint is explicitly stated.

This shows that $k_{opt}=k_1$ and $X_{opt}$ as given above define optimal OSG(s-k,1) as given in Table 103 for all s>=18.

Theorem 401

Let A = OSG($n_{opt}$,1) be the maximal optimal stride generator of order 1 for a given value of s>=18.

Then A is the stride generator underlying the cover C(A,3,s)=$X_{opt}$, where $n_{opt}$=s-$k_{opt}$ and $X_{opt}$ are given in Table 103.

[ Theorem 400 proved that A's potential cover is $X_{opt}$; this theorem shows that potential cover to be a real one. ]

Proof

There are nine cases to consider - see Table 103 for details.

In each case, we have only to consider values in stride 0 - $0<=x<a_3$ - since by Theorem 100 strides 1, 2, ... k can all be generated, because A has order 1.

Such a value x can be generated as:

$$x = c_2a_2 + c_1 \quad \text{where} \quad c_1 < a_2$$

Suppose $a_3 = C_2a_2 + C_1$ where $C_1 < a_2$

Then any value $0 <= x < a_3$ can be generated using at most $C_2+a_2-2$ stamps, because:

either a) $c_2 < C_2$ and $c_1 < a_2$

$$=> c_1+c_2 <= (C_2-1)+a_2-1 = C_2+a_2-2$$

or b) $c_1=C_2$ and $c_1 < C_1$

$$=> c_1+c_2 <= C_2+a_2-2$$

So we have only to determine under what conditions $C_2+a_2-2 <= s$:

e.g. From Table 105 for s = 9t we have:

$a_2 = 6t+3$, $a_3 = (2t+1)a_2 + (2t+1)$

So $C_2+a_2-2 = (2t+1) + (6t+3) - 2 = 8t+2$

=> $C_2+a_2-2 <= s$ iff $9t-(8t+2) >= 0$ => $t>=2$

The full results are shown in the following table:

| r | $a_2$ | $C_2$ | $s-(C_2+a_2-2)$ | >= 0 when: |
|---|---|---|---|---|
| 0) | 6t+3 | 2t+1 | (t-2 | t>=2 |
| 1) | | | (t-1 | t>=1 |
| 2) | 6t+5 | 2t+1 | (t-2 | t>=2 |
| 3) | | | (t-1 | t>=1 |
| 4 | 6t+7 | 2t+1 | t-2 | t>=2 |
| 5 | 6t+6 | 2t+2 | t-1 | t>=1 |
| 6) | 6t+8 | 2t+2 | (t-2 | t>=2 |
| 7) | | | (t-1 | t>=1 |
| 8 | 6t+10 | 2t+2 | t-2 | t>=2 |

Table 400

This shows that OSG(s-$k_{opt}$,1) generates a full cover for at least all s>=18.

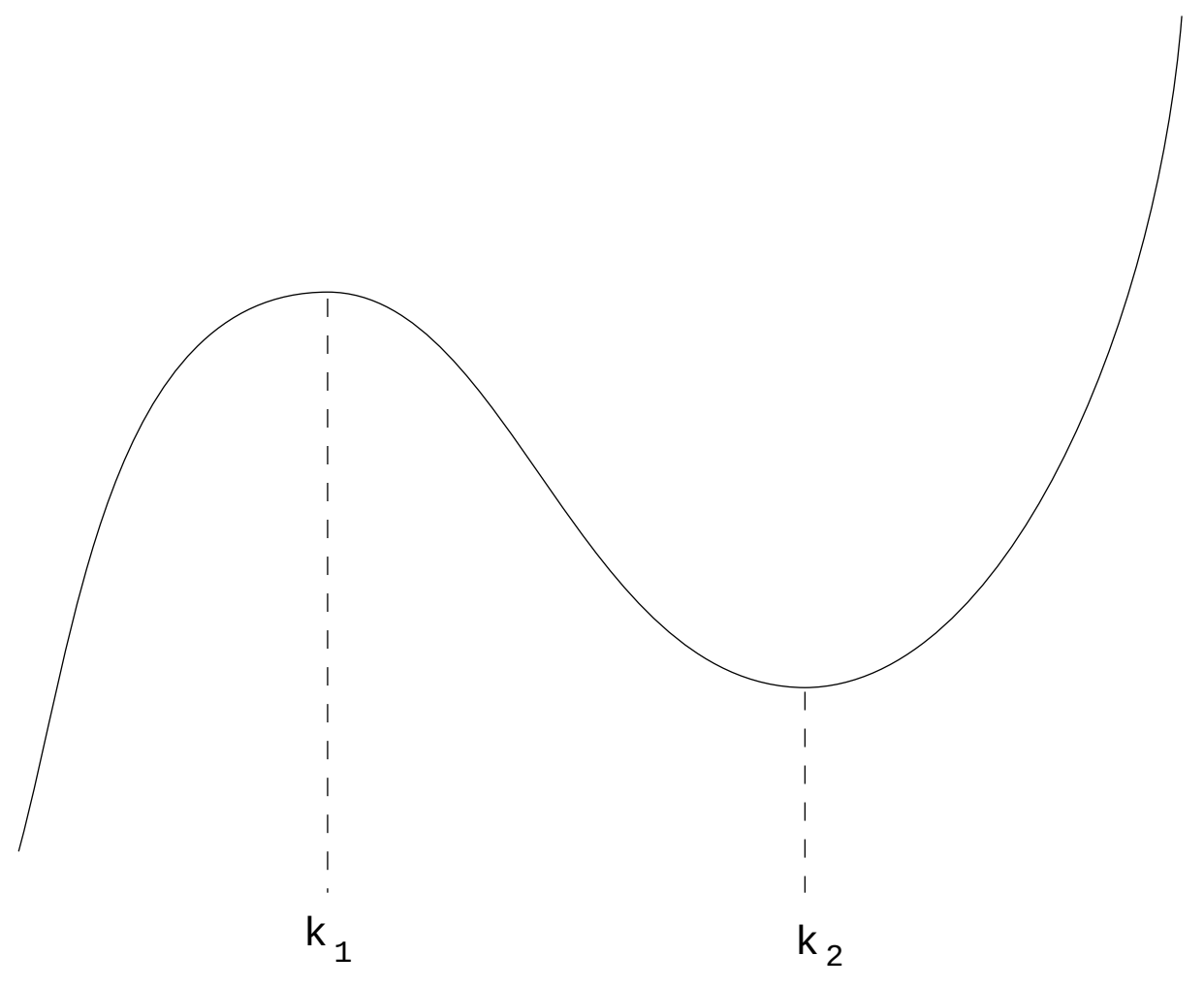

Figure 400

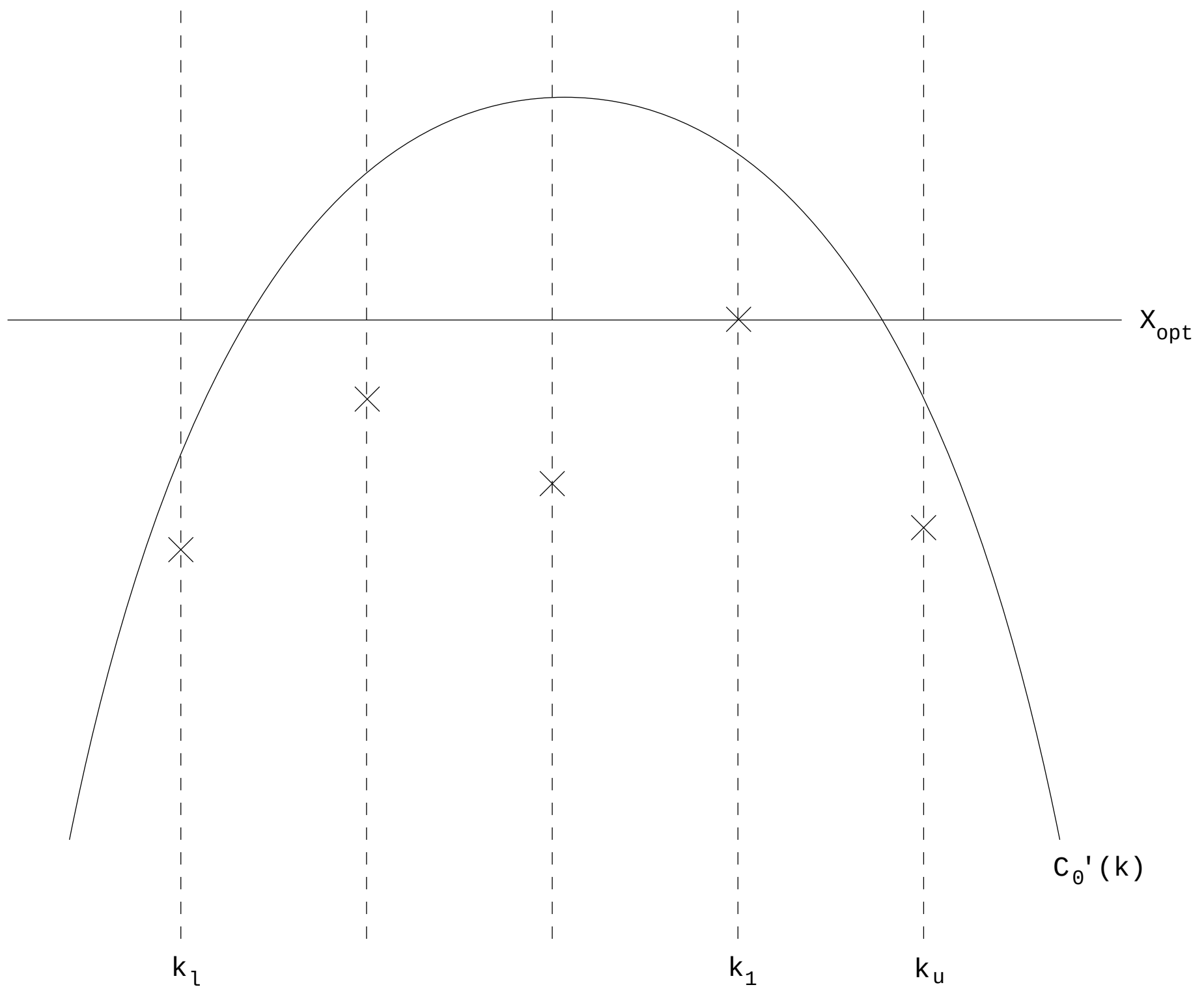

Figure 401

## 4.2 Sub-optimal SGi(s-k,p) for i>=2

### Definition 400

Let A = {1, $a_2$, $a_3$} be an n-stride generator SG(n,p), and {1, $A_2$, $A_3$} be the corresponding optimal n-stride generator OSG(n,1) of order 1.

Then A is called a "sub-optimal stride generator SGi(n,p)" if $a_3=A_3-i$.

If i=1, A is an "immediately sub-optimal stride generator SG1(n,p)".

### Theorem 402

Let A = {1, $a_2$, $a_3$} be a sub-optimal stride generator SGi(s-k,p) for some 0<=k<s and i>=2 which underlies a cover C(A,3,s).

Then:

$$C(A,3,s) < X_{opt} \quad \text{for all } s>=36$$

where $X_{opt}$ is the function of s given by Table 103.

### Proof

Note:
Some of the algebraic manipulations summarised in this proof have been confirmed using REDUCE; see section 6.3 for details.

Let OSG(s-k,1) = {1, $A_2$, $A_3$}; then by Definition 400:

$$a_3 = A_3-i$$

By Theorem 231 we know that:

$$C(A,3,s) = (k+1)(A_3-i) + y-1$$
$$< (k+2)(A_3-i)$$

We now define:

$$M_i(k) = (k+2)(A_3-i)$$

Now $M_i(k) - M_{i+1}(k) = (k+2)(A_3-i) - (k+2)(A_3-i-1)$

$$= k+2 > 0 \text{ for all } k>=0$$

=> $M_{i+1}(k) < M_i(k)$ for all i

=> $C(A,3,s) < M_i(k) <= M_2(k)$ for all i>=2 and for all k>=0

We now show that $M_2(k) < X_{opt}$ for all 0 <= k < s:

From Table 101:

$$M_2(k) = (k+2)(A_3-2)$$

where $A_3 = (n+3)(n+2)/3$ when n=s-k = 0 or 1 mod 3

or $A_3 = ((n+3)(n+2)+1)/3$ when n=s-k = 2 mod 3

We write:

$M_2'(k) = (k+2)(((n+3)(n+2)+1)/3 - 2)$; clearly, $M_2(k) <= M_2'(k)$

Substituting n=s-k and evaluating gives:

$M_2'(k) = (k^3 - (2s+3)k^2 + (s^2+s-9)k + 2s^2+10s+2)/3$

This function has the form of Figure 400, where the values $k_1$ and $k_2$ are the solutions of $dM_2'(k)/dk = 0$:

i.e. $3k^2 - 2(2s+3)k + s^2+s-9 = 0$

Solving for k, we have:

$k = (2(2s+3) +/- sqrt(4(2s+3)^2 - 12(s^2+s-9)))/6$

$=> k = ((2s+3) +/- sqrt((2s+3)^2 - 3(s^2+s-9)))/3$

$=> k = ((2s+3) +/- sqrt(s^2+9s+36))/3$

We now show that $k_2>s$ and $0<k_1<s$ for all s>=3, where $k_1<k_2$ are the two roots of this equation:

$3k_2 = (2s+3) + sqrt(s^2+9s+36)$

$> (2s+3) + sqrt(s^2+8s+16)$ for all s>=0

$= (2s+3) + (s+4)$

$= 3s+7$

$=> k_2>s+2$ for all s>=0

$3k_1 = (2s+3) - sqrt(s^2+9s+36)$

$> (2s+3) - sqrt(s^2+12s+36)$ for all s>=1

$= (2s+3) - (s+6)$

$= s-3$

$=> k_1>0$ for all s>=3

$3k_1 = (2s+3) - sqrt(s^2+9s+36)$

$< (2s+3) - sqrt(s^2+8s+16)$ for all s>=0

$= (2s+3) - (s+4)$

$= s-1$

$=> k_1<s$ for all s>=0

This is important, since it proves that the maximum value for $M_2'(k)$ in the range 0<=k<s lies at $M_2'(k_1)$, and that $M_2'(k)$ increases monotonically as k runs from 0 to $k_1$, and decreases monotonically as k runs from $k_1$ to s.

It is also clear that for large s, $k_1 \sim s/3$, since:

$$3k_1 = (2s+3) - \mathrm{sqrt}(s^2+9s+36)$$
$$\sim (2s+3) - \mathrm{sqrt}(s^2+9s+20.25)$$
$$= (2s+3) - (s+4.5)$$
$$= s-1.5$$
$$\Rightarrow k_1 \sim s/3 - 0.5$$

We now identify exactly where this maximum occurs and what it is by evaluating $M_2'(k)$ for values around k = s/3. There are three cases to consider, according to the value of s mod 3:

a) s=3u;    $k_1$ = u-1  for  u>=1:

$M_2'(u-2) = (4u^3 + 18u^2 + 15u)/3$
$M_2'(u-1) = (4u^3 + 18u^2 + 21u + 7)/3$
$M_2'(u) = (4u^3 + 18u^2 + 21u + 2)/3$

b) s=3u+1;  $k_1$ = u  for  u>=1:

$M_2'(u-1) = (4u^3 + 22u^2 + 33u + 15)/3$
$M_2'(u) = (4u^3 + 22u^2 + 35u + 14)/3$
$M_2'(u+1) = (4u^3 + 22u^2 + 31u + 3)/3$

c) s=3u+2;  $k_1$ = u  for  u>=0:

$M_2'(u-1) = (4u^3 + 26u^2 + 47u + 25)/3$
$M_2'(u) = (4u^3 + 26u^2 + 51u + 30)/3$
$M_2'(u+1) = (4u^3 + 26u^2 + 49u + 21)/3$

Next, for each of these three cases, we have to show that $M_2'(k_1) < X_{opt}$; since $X_{opt}$ depends on s mod 9, each case is further sub-divided into three:

aa)  s=9t,  u=3t

$M_2'(k_1) = 36t^3 + 54t^2 + 21t + 2 + 1/3$
$< 36t^3 + 54t^2 + 22t = X_{opt}$  (t>=3)

ab)  s=9t+3,  u=3t+1

$M_2'(k_1) = 36t^3 + 90t^2 + 69t + 16 + 2/3$
$< 36t^3 + 90t^2 + 71t + 15 = X_{opt}$  (t>=1)

ac)  s=9t+6,  u=3t+2

$M_2'(k_1) = 36t^3 + 126t^2 + 141t + 51$
$< 36t^3 + 126t^2 + 143t + 49 = X_{opt}$  (t>=2)

ba)  s=9t+1,  u=3t

$M_2'(k_1) = 36t^3 + 66t^2 + 35t + 4 + 2/3$
$< 36t^3 + 66t^2 + 36t + 4 = X_{opt}$  (t>=1)

bb)  s=9t+4,  u=3t+1

$M_2'(k_1) = 36t^3 + 102t^2 + 91t + 25$
$< 36t^3 + 102t^2 + 92t + 22 = X_{opt}$  (t>=4)

bc) $s=9t+7$, $u=3t+2$

$$M_2'(k_1) = 36t^3 + 138t^2 + 171t + 68$$
$$< 36t^3 + 138t^2 + 173t + 68 \quad = X_{opt} \quad (t>=1)$$

ca) $s=9t+2$, $u=3t$

$$M_2'(k_1) = 36t^3 + 78t^2 + 51t + 10$$
$$< 36t^3 + 78t^2 + 53t + 8 \quad = X_{opt} \quad (t>=2)$$

cb) $s=9t+5$, $u=3t+1$

$$M_2'(k_1) = 36t^3 + 114t^2 + 115t + 37$$
$$< 36t^3 + 114t^2 + 116t + 36 \quad = X_{opt} \quad (t>=2)$$

cc) $s=9t+8$, $u=3t+2$

$$M_2'(k_1) = 36t^3 + 150t^2 + 203t + 89 + 1/3$$
$$< 36t^3 + 150t^2 + 204t + 86 \quad = X_{opt} \quad (t>=4)$$

These results show that $M_2'(k)$ - and hence $M_2(k)$ and C(A,3,s) for i>=2 - is less than $X_{opt}$ for all t>=4 (ie s>=36).

## 4.3 Immediately sub-optimal SG1(s-k,1)

Theorem 403

Let A = {1, $a_2$, $a_3$} be an immediately sub-optimal stride generator SG1(s-k,1) for some 0<=k<s which underlies a cover C(A,3,s).

Then:

$$C(A,3,s) < X_{opt} \quad \text{for all } s>=9$$

where $X_{opt}$ is the function of s given by Table 103.

[ In other words, no stride generator SG(n,1) of order 1 whose length is just one less than that of OSG(n,1) can provide as good a cover. ]

Proof

Note:
Some of the algebraic manipulations summarised in this proof have been confirmed using REDUCE; see section 6.3 for details.

Let OSG(s-k,1) = {1, $A_2$, $A_3$}; then:

$$a_3 = A_3-1 \qquad \text{(see Definition 400)}$$

We define $C_0(k)$ to be the potential cover of OSG(s-k,1); by Theorems 231 and 301:

$$C_0(k) = (k+1)A_3 + y_0-1$$

where $A_3$ and $y_0$ are given by $a_3$ and y in Table 101.

We define $C_1(k)$ to be the potential cover of SG1(n,1); by Theorems 231 and 302:

$$C_1(k) = (k+1)(A_3-1) + y_1-1$$

where $y_1$ is given by y in Table 102.

We will show that:

$$C_1(k) < C_0(k) \quad \text{for all } 0 <= k < s$$

thus proving that $C_1(k) < X_{opt}$, since $X_{opt}$ is the maximum value of $C_0(k)$.

We write:

$$D = C_0(k) - C_1(k)$$

$$= (k+1)A_3 + (y_0-1) - (k+1)(A_3-1) - (y_1-1)$$

$$= y_0-y_1 + k+1$$

and show that D>0 for all k.

There are three cases to consider, according to the value of (s-k) mod 3:

a) (s-k) mod 3 = 0:

$$D = ((n^2+3n-9) - (n^2+3n))/3 + k+1$$
$$= k-2 > 0 \quad \text{for all } k>2$$

b) (s-k) mod 3 = 1:

$$D = ((n^2+3n-1) - (n^2+3n-16))/3 + k+1$$
$$= k+6 > 0 \quad \text{for all } k>=0$$

c) (s-k) mod 3 = 2:

There is no SG1(s-k,1) in this case, and so no case to answer!

There remains the issue of $C_1(k)$ for k<=2 in case (a): we must show that these values are less than $X_{opt}$ even though they may be greater than the corresponding $C_0(k)$. (Although this is obviously true for larger s - where $X_{opt} = C_0(k_{opt})$ for $k_{opt}$ ~ s/3 - it may not be true for low values of s.)

There are just three cases remaining to be considered:

$C_1(0)$ - with s mod 3 = 0
$C_1(1)$ - with (s-1) mod 3 = 0
$C_1(2)$ - with (s-2) mod 3 = 0

Substituting values from the r=0 row of Table 102, we have:

a) $C_1(0) = ((s^2+5s+3) + (s^2+3s) - 3)/3$
$\quad = (2s^2+8s)/3$

s=9t => $C_1(0) = 54t^2 + 24t$
$\quad < X_{opt} = 36t^3 + 54t^2 + 22t$ (t>=1)

s=9t+3 => $C_1(0) = 54t^2 + 60t + 14$
$\quad < X_{opt} = 36t^3 + 90t^2 + 71t + 15$ (t>=0)

s=9t+6 => $C_1(0) = 54t^2 + 96t + 40$
$\quad < X_{opt} = 36t^3 + 126t^2 + 143t + 149$ (t>=0)

b) $C_1(1) = (2((s-1)^2+5(s-1)+3) + ((s-1)^2+3(s-1)) - 3)/3$
$\quad = (3s^2+7s-7)/3$

s=9t+1 => $C_1(1) = 81t^2 + 39t + 1$
$\quad < X_{opt} = 36t^3 + 66t^2 + 36t + 4$ (t>=1)

s=9t+4 => $C_1(1) = 81t^2 + 93t + 23$
$\quad < X_{opt} = 36t^3 + 102t^2 + 92t + 22$ (t>=1)

s=9t+7 => $C_1(1) = 81t^2 + 147t + 63$
$\quad < X_{opt} = 36t^3 + 138t^2 + 173t + 68$ (t>=0)

c) $C_1(2) = (3((s-2)^2+5(s-2)+3) + ((s-2)^2+3(s-2)) - 3)/3$
$\quad = (4s^2+2s-14)/3$

s=9t+2 => $C_1(2) = 108t^2 + 54t + 2$
$\quad < X_{opt} = 36t^3 + 78t^2 + 53t + 8$ (t>=0)

s=9t+5 => $C_1(2) = 108t^2 + 126t + 32$
$\quad < X_{opt} = 36t^3 + 114t^2 + 116t + 36$ (t>=0)

s=9t+8 => $C_1(2) = 108t^2 + 198t + 86$
$\quad < X_{opt} = 36t^3 + 150t^2 + 204t + 86$ (t>=1)

Thus the case is proven for all t>=1 (ie s>=9).

# SECTION 5

## LIMITING THE ORDER

This Section collects together results which limit the order of stride generators which underlie good covers.

Section 5.1 shows that for s>=40 the order must be <=3; computer-aided algebra is of considerable assistance here!

Section 5.2 shows that for n>=52 there are always stride generators of order 1 that are longer than OSG(n,0), OSG(n,2) and OSG(n,3). To do this, further general properties of stride generators are developed, and the results strongly suggest that in the long run stride generators of order 1 are better than stride generators of any other order. This conjecture (although not part of the proof) is examined further in Section 5.3.

Finally, Section 5.4 shows that for s>=81 there is no need to consider any stride generators whose order is not equal to 1.

## 5.1 Order <= 3 for maximal cover for s>=40

### Theorem 500

Let SG(n,p) be the stride generator that underlies the maximal cover C(A,3,s).

Then p<=3 for all s>=40.

### Proof

We start with an overview of the proof.

a) We establish a lower bound X for the maximal cover M(3,s):

$$X = 4s^3/81 + 2s^2/3 + 22s/9 \quad \text{for } s>=35 \qquad - (1)$$

X is a function of s.

b) We establish an upper bound C on the value of $a_2$ for any stride generator A = {1, $a_2$, $a_3$} whose cover C(A,3,s) >= X:

$$C = s+3-a_3/s \quad \text{for } s>=40 \qquad - (2)$$

C is a function of s and $a_3$.

c) We establish an upper bound Q on the order of any stride generator A = SG(n,p) with C(A,3,s) >= X:

$$Q = ((n+1)C-a_3)/(a_3-C) \quad \text{for } s>=40 \qquad - (3)$$

Q is a function of n, s and $a_3$.

d) We establish an upper bound N for n for any stride generator A = SG(n,p) underlying a cover C(A,3,s) >= X:

$$N = (s+2) - X/a_3 \qquad - (4)$$

N is a function of s and $a_3$.

e) We substitute N for n in Q to obtain an upper bound P on the order of any stride generator A = SG(n,p) with C(A,3,s) >= X:

$$P = ((N+1)C-a_3)/(a_3-C) \quad \text{for } s>=40 \qquad - (5)$$

P is a function of s and $a_3$.

f) We show that the maximum value of P(s,$a_3$) as $a_3$ is varied is less than 4 for all s>=40.

The algebra involved in some of these steps is complex, and use is made of a computer algebraic manipulation package called REDUCE. Edited output from this package is included to support the proof, where input commands for REDUCE are preceded by "C: ", additional commentary is enclosed in square brackets [ ... ], and all other lines are the unedited output from the program itself.

Step (a):

We know from Theorem 401 that the formulae for $X_{opt}$ in Table 103 are less than or equal to M(3,s).

The following edited results from REDUCE show that X (as given in (1) above) is always less than or equal to those formulae for all s>=35.

[ Define the proposed lower bound X as a function of s: ]

```
C: for all s
C:   let X(s) = (4/81)*s*s*s + (2/3)*s*s + (22/9)*s;
```

[ Define the known lower bound $X_{opt}$ as a function of the polynomial coefficients a, b, c, d for convenience: ]

```
C: for all a,b,c,d
C:   let Xopt(a,b,c,d) = a*t*t*t + b*t*t + c*t + d;
```

[ Evaluate the difference between $X_{opt}$ and X for r = 0 ... 8 where r = s mod 9: ]

```
C: Xopt(36, 54, 22, 0) - X(9*t);

   0

C: Xopt(36, 66, 36, 4) - X(9*t+1);

   2/3*t + 68/81

C: Xopt(36, 78, 53, 8) - X(9*t+2);

   5/3*t + 4/81

C: Xopt(36, 90, 71,15) - X(9*t+3);

   t + 1/3

C: Xopt(36,102, 92,22) - X(9*t+4);

   2/3*t - 130/81

     [ This is >=0 provided:

             54t>=130

          =>   t>=3

          =>   s>=31  ]

C: Xopt(36,114,116,36) - X(9*t+5);

   2/3*t + 76/81

C: Xopt(36,126,143,49) - X(9*t+6);

   t - 1/3

     [ This is >=0 provided t>=1 => s>=15 ]

C: Xopt(36,138,173,68) - X(9*t+7);

   5/3*t + 104/81
```

C: Xopt(36,150,204,86) - X(9*t+8);

2/3*t - 122/81

[ This is >=0 provided:

54t>=122

=> t>=3

=> s>=35 ]

Step (b):

Let A = {1, $a_2$, $a_3$} be the stride generator underlying a cover C(A,3,s)>=X.

Then C({1,$a_2$},2,s) >= $a_3$-1

It is easy to establish a formula for C({1,$a_2$},2,s) (see my separate general document on the Postage Stamp Problem for details) and we have:

$a_2(s+3-a_2)-2 >= a_3-1$

=> $a_2^2-(s+3)a_2+(a_3+1) <= 0$

Values of $a_2$ which meet this condition must lie between the two roots of the quadratic equation, and so the larger root provides an upper bound for $a_2$; let us call this root $A_2$.

$A_2 = ((s+3) + sqrt((s+3)^2-4(a_3+1)))/2$

Looking at the value whose root is required, we will show that:

$(s+3-2a_3/s)^2 > (s+3)^2 - 4(a_3+1)$ for s>=40 - (6)

and so:

$C = ((s+3) + (s+3-2a_3/s))/2$

$= s+3-a_3/s$

$> A_2$ for s>=40

and hence C (as given in (2)) is an upper bound on $a_2$ for s>=40.

To show (6) we write $a_3$=ks and look at the difference between the lhs and the rhs:

$(s+3-2a_3/s)^2 - (s+3)^2 + 4(a_3+1)$

$= (s+3)^2 - 4(s+3)k + 4k^2 - (s+3)^2 + 4ks + 4$

$= 4k^2 - 4ks - 12k + 4ks + 4$

$= 4(k^2-3k+1)$

So (6) is true iff $k^2-3k+1 > 0$.

Solving the quadratic we see that this is true for:

k < (3-sqrt(5))/2 or k > (3+sqrt(5))/2

and so certainly for:

k < 0.38 or k > 2.62

Now for $\{1, a_2, a_3\}$ to be a stride generator with a cover >= X, we know that:

$sa_3 >= X$ - for otherwise there is no way to generate the value X

$=> sa_3 >= 4s^3/81 + 2s^2/3 + 22s/9$ by (1)

$=> a_3 >= 4s^2/81 + 2s/3 + 22/9$

$=> a_3 > 4s^2/81 + 2s/3$

$=> k > 4s/81 + 2/3$

So (6) is certainly true if:

$4s/81 + 2/3 > 2.62$

$=> 4s+54 > 212.22$

$=> s > 39.555$

Step (c):

Theorem 239 shows that for any stride generator:

$a_2 >= a_3(p+1)/(n+p+1)$

$=> a_2(n+p+1) >= a_3(p+1)$

$=> p(a_3-a_2) <= a_2(n+1)-a_3$

$=> p <= ((n+1)a_2-a_3)/(a_3-a_2)$

Now (2) gives an upper bound C on $a_2$ for any stride generator with C(A,3,s)>=X for s>=40, and so Q is an upper bound on p for s>=40:

$Q = ((n+1)C-a_3)/(a_3-C)$

Step (d):

Let $A = \{1, a_2, a_3\}$ be a stride generator SG(n,p) underlying a cover C(A,3,s)>=X.

Then, by Theorem 231, we know that:

$C(A,3,s) = (s-n+1)a_3 + y-1$

$< (s-n+2)a_3$

So:

$X < (s-n+2)a_3$

$=> na_3 < (s+2)a_3 - X$

$=> n < (s+2)-X/a_3$

So N is an upper bound for n for any stride generator SG(n,p) underlying a cover C(A,3,s)>=X where:

$N = (s+2)-X/a_3$

Step (e): needs no further explanation.

Step (f):

We start by using REDUCE to expand the function $P(a_3)$ given in (5):

[ We can read the statement below as:

"let $P(a_3) = ((N+1)C-a_3)/(a_3-C)$ substituting
$C = s+3-(a_3/s)$
and $N = (s+2)-X/a_3$ where $X = 4s^3/81 + ...$" ]

```
C: operator p;
C: for all a3 let p(a3) = sub(
C:   c = s+3-(a3/s),
C:   n = sub(  x = (4/81)*s*s*s + (2/3)*s*s + (22/9)*s,
C:                 (s+2)-x/a3
C:          ),
C:   ((n+1)*c-a3)/(a3-c)
C: );
```

[ We now evaluate $P(a_3)$: ]

```
C: p(a3);

                2            2           3            2                  5
    - (162*a3 *s + 243*a3  - 85*a3*s  - 540*a3*s  - 927*a3*s + 4*s

                  4        3        2                         2
        + 66*s  + 360*s  + 594*s )/(81*a3*(a3*s + a3 - s  - 3*s))
```

We can rewrite this as:

$P(a_3) = (D_2a_3^2+D_1a_3+D_0)/81a_3((s+1)a_3-s(s+3))$

where:

$D_2 = -(162s+243)$

$D_1 = 85s^3+540s^2+927s$

$D_0 = -(4s^5+66s^4+360s^3+594s^2)$

Looking at the denominator we see that $P(a_3)$ has singularities at:

$a_3 = 0$

and when:

$(s+1)a_3 = s(s+3)$

$\Rightarrow a_3 = s(s+3)/(s+1)$

Also, as $a_3 \to$ (plus or minus) infinity:

$P(a_3) \to -(162s+243)a_3^2/81a_3^2(s+1)$

$= -(2s+3)/(s+1)$

Figures 500 and 501 show typical curves for $P(a_3)$ - for s=5 and s=20.

We now show that for s>=23 the only portion of the graph that we need to consider is the right hand section - that is, values of $a_3 > A_3$ where:

$$A_3 = s(s+3)/(s+1)$$

Since $C(\{1, a_2, a_3\},3,s) >= X > 4s^3/81$ we know that:

$$sa_3 > 4s^3/81$$

So we need only consider the right hand section of the graph if:

$$sA_3 <= 4s^3/81$$

$$=> s(s+3) <= 4s^2(s+1)/81$$

$$=> 81(s+3) <= 4s(s+1)$$

$$=> 4s^2-77s-243 >= 0$$

Solving for s we have:

$$s = (77 +/- \mathrm{sqrt}(5929+3888))/8$$

$$=> s = 22.01 \quad \text{or} \quad s = -2.76$$

So our assumption is valid provided s>=23.

Now for values of $a_3 > A_3$ the denominator of $P(a_3)$ is positive, and so the shape of this part of the curve depends on the value of the numerator at $A_3$: if the numerator is positive we have the "s=5" case, and if it is negative we have the "s=20" case.

We now show that for s>=11 the numerator is negative at $a_3 = A_3$:

[ We first substitute $a_3=s(s+3)/(s+1)$ into the numerator of $P(a_3)$, and call the resulting expression xx: ]

```
C: xx := sub( a3=s*(s+3)/(s+1), num(p(a3)) );

                3     4       3        2
              s *(4*s  - 11*s  - 222*s  - 747*s - 864)
   xx :=  - -----------------------------------------
                                2
                               s  + 2*s + 1
```

We now show that xx<0 for all s>=11. Clearly the denominator is positive and $s^3$ is positive for all s>0, so we have only to show that y>0 for all s>=11 where:

$$y = 4s^4 - 11s^3 - 222s^2 - 747s - 864$$

We have:

$$s^4 >= 11s^3 >= 121s^2 >= 1331s >= 14641 \quad \text{for } s>=11$$

So:

$$11s^3 + 222s^2 + 747s + 864 <= s^4 + 2s^4 + 0.6s^4 + 0.1s^4 = 3.7s^4$$

$$=> y >= 4s^4 - 3.7s^4 = 0.3s^4 > 0 \text{ as required.}$$

This is a "sharp" bound: that is, xx changes sign at s=11:

[ We first multiply out the denominator and divide by $s^3$ for convenience: ]

```
C: xx := xx*(s+1)*(s+1)/(s*s*s);

                  4       3        2
   xx :=  - (4*s  - 11*s  - 222*s  - 747*s - 864)
```

[ Now we show that xx>0 for s=10, xx<0 for s=11: ]

```
C: sub(s=10, xx);

   1534

C: sub(s=11, xx);

   -7980
```

We now know that for s>=23:

- We need only look at values of $P(a_3)$ for $a_3>A_3$
- $P(a_3)$ -> -(2s+3)/(s+1) as $a_3$ -> infinity
- $P(a_3)$ -> -infinity as $a_3$ -> $A_3$ from above
- There are no singularities for $a_3>A_3$

How many maxima and minima can this part of the curve exhibit? Some possibilities are shown in Figure 502:

a) None
b) One maximum
c) One maximum and one minimum
d) Two maxima and one minimum

To determine which form the curve takes we need to look at the solutions to:

$$dP(a_3)/da_3 = 0$$

We will show that there are two solutions, and that one of these is greater than $A_3$. We also show that for s>=11 the middle section of the curve ($0<a_3<A_3$) must have at least one minimum, and so the other solution must lie in that range. This means that the form of the curve for $a_3>A_3$ is as shown in (b).

First, consider the middle section for s>=11:

As $a_3$ -> $A_3$ from below, $P(a_3)$ -> infinity.

As $a_3$ -> 0 from above, $P(a_3)$ -> infinity, because:

- the numerator -> $D_0$, which is negative for positive s
- the denominator is negative

There are no singularities in the range $0<a_3<A_3$, and at both ends of the range $P(a_3)$ -> infinity: so there must be at least one minimum value of $P(a_3)$ in this range, and hence one solution of $dP/da_3 = 0$.

Now we look at $dP/da_3$ itself:

[ We differentiate $P(a_3)$ with respect to $a_3$, and call the resulting expression dfp: ]

```
C: dfp := df(p(a3),a3);

                       2   3          2  2          2            2         5
   dfp :=  - (s*(85*a3 *s  + 463*a3 *s  + 738*a3 *s + 198*a3  - 8*a3*s

                     4            3             2
         - 140*a3*s  - 852*a3*s  - 1908*a3*s  - 1188*a3*s

                6       5        4         3         2          2
         + 4*s  + 78*s  + 558*s  + 1674*s  + 1782*s ))/(81*a3 *

              2  2       2       2         3         2            4
         (a3 *s  + 2*a3 *s + a3  - 2*a3*s  - 8*a3*s  - 6*a3*s + s

              3      2
         + 6*s  + 9*s ))
```

[ We see that the numerator is a quadratic in $a_3$ and so there are exactly two solutions to $dP/da_3 = 0$. These are given below as:

```
sol := { a3 = <first soln>, a3 = <second soln> }  ]
```

```
C: sol := solve(dfp,a3);

                                6       5         4          3          2
   sol := {a3= - (s*(sqrt(16*s  + 124*s  - 1062*s  - 13968*s  - 45846*s

                                          6       5        4
              - 57024*s)*s + 3*sqrt(16*s  + 124*s  - 1062*s

                       3          2                  4       3
              - 13968*s  - 45846*s  - 57024*s) - 4*s  - 70*s

                     2                    3        2
              - 426*s  - 954*s - 594))/(85*s  + 463*s  + 738*s

          + 198),

                            6        5         4          3          2
           a3=(s*(sqrt(16*s  + 124*s  - 1062*s  - 13968*s  - 45846*s

                                         6        5        4
            - 57024*s)*s + 3*sqrt(16*s  + 124*s  - 1062*s

                     3          2                  4       3
            - 13968*s  - 45846*s  - 57024*s) + 4*s  + 70*s

                   2                     3        2
            + 426*s  + 954*s + 594))/(85*s  + 463*s  + 738*s

          + 198)}
```

[ We now extract each solution and name them $a_{31}$, $a_{32}$ respectively.

(In detail, sol is a list, and so first(sol) is its first element "a3=-(...)"; the part function is then used to extract the second component of this expression which, in this case, is the right hand side of the equation.) ]

```
C: a31 := part(first(sol), 2);

                             6        5         4          3          2
   a31 :=  - (s*(sqrt(16*s  + 124*s  - 1062*s  - 13968*s  - 45846*s

                                          6        5         4
               - 57024*s)*s + 3*sqrt(16*s  + 124*s  - 1062*s

                      3          2                    4       3
               - 13968*s  - 45846*s  - 57024*s) - 4*s  - 70*s

                    2                     3        2
               - 426*s  - 954*s - 594))/(85*s  + 463*s  + 738*s + 198

      )

C: a32 := part(second(sol),2);

   a32 := (s*(sqrt(

                6        5         4          3          2
            16*s  + 124*s  - 1062*s  - 13968*s  - 45846*s  - 57024*s)

            *s + 3*sqrt(

                6        5         4          3          2
            16*s  + 124*s  - 1062*s  - 13968*s  - 45846*s  - 57024*s)

                  4       3        2                           3        2
             + 4*s  + 70*s  + 426*s  + 954*s + 594))/(85*s  + 463*s

       + 738*s + 198)
```

[ We now evaluate the differences $a_{31}$-$A_3$ and $a_{32}$-$A_3$: ]

```
C: aa3 := s*(s+3)/(s+1);

            s*(s + 3)
   aa3 := -----------
              s + 1

C: diff1 := a31 - aa3;

                               6        5         4          3          2
   diff1 :=  - (s*(sqrt(16*s  + 124*s  - 1062*s  - 13968*s  - 45846*s

                           2              6        5         4
                 - 57024*s)*s  + 4*sqrt(16*s  + 124*s  - 1062*s

                        3          2                                 6
                 - 13968*s  - 45846*s  - 57024*s)*s + 3*sqrt(16*s

                      5         4          3          2
                 + 124*s  - 1062*s  - 13968*s  - 45846*s  - 57024*s)

                      5       4        3        2                  4
                 - 4*s  + 11*s  + 222*s  + 747*s  + 864*s))/(85*s

                 3         2
          + 548*s  + 1201*s  + 936*s + 198)
```

```
     C: diff2 := a32 - aa3;

        diff2 := (s*(sqrt(

                     6        5         4          3          2
              16*s  + 124*s  - 1062*s  - 13968*s  - 45846*s  - 57024*s

                    2
              )*s  + 4*sqrt(

                     6        5         4          3          2
              16*s  + 124*s  - 1062*s  - 13968*s  - 45846*s  - 57024*s

              )*s + 3*sqrt(

                     6        5         4          3          2
              16*s  + 124*s  - 1062*s  - 13968*s  - 45846*s  - 57024*s

                        5       4        3        2                 4
              ) + 4*s  - 11*s  - 222*s  - 747*s  - 864*s))/(85*s

                    3         2
              + 548*s  + 1201*s  + 936*s + 198)
```

We see that:

diff1 = -s((s²+4s+3)x - y)/w

and diff2 =  s((s²+4s+3)x + y)/w

where:

x = sqrt(16s⁶ + 124s⁵ - 1062s⁴ - 13968s³ - 45846s² - 57024s)

y = 4s⁵ - 11s⁴ - 222s³ - 747s² - 864s

w = 85s⁴ + 548s³ + 1201s² + 936s + 198

We showed earlier that  4s⁴ - 11s³ - 222s² - 747s - 864 > 0 for all s>=11; so y>0 for all s>=11. We now show that x² (and hence x) is positive for all s>=11.

For s>=13 we have:

s⁶ >= 13s⁵ >= 169s⁴ >= 2197s³ >= 28561s² >= 371293s

So:

x² > 16s⁶ - (6.3s⁶ + 6.4s⁶ + 1.7s⁶ + 0.2s⁶) = 1.4s⁶ > 0

Numerical substitution confirms that xx>0 for s=11, s=12 but not for s=10:

```
  [ We first extract the formula for x from diff2: ]

     C: xx := part(diff2,1,2,1,1);

                         6        5         4          3          2
        xx := sqrt(16*s  + 124*s  - 1062*s  - 13968*s  - 45846*s  - 57024*s)

  [ and square it: ]

     C: xx := xx*xx;

                        5       4        3         2
        xx := 2*s*(8*s  + 62*s  - 531*s  - 6984*s  - 22923*s - 28512)
```

[ and show that it is negative for s=10, and positive for s=11 and s=12 as required: ]

```
C: sub(s=10, xx);

   -1342840

C: sub(s=11, xx);

   8000520

C: sub(s=12, xx);

   25186464
```

We now know that x, y and w are greater than zero for all s>=11, and so diff2>0 for all s>=11.

So for s>=11 $a_{32}$ is the solution to $dP/da_3=0$ that lies in the range $a_3>A_3$, and so $a_{31}$ must be the solution that lies in the range $0<a_3<A_3$.

We can now substitute $a_{32}$ back into $P(a_3)$ to obtain the formula which gives for s>=11 the maximum value of $P(a_3)$ in the range $a_3>A_3$ in terms of s:

[ Substitute $a_{32}$ back into $P(a_3)$ to define the function PP(s) - which provides an upper bound on p as a function of s: ]

```
C: operator pp;
C: for all s let pp(s) = sub(a3=a32,p(a3));

C: pp(s);

                     6        5         4          3          2
   (5929*sqrt(16*s  + 124*s  - 1062*s  - 13968*s  - 45846*s  - 57024*s)

      6
    *s  + 78418

               6        5         4          3          2           5
    *sqrt(16*s  + 124*s  - 1062*s  - 13968*s  - 45846*s  - 57024*s)*s  +

     401394

               6        5         4          3          2           4
    *sqrt(16*s  + 124*s  - 1062*s  - 13968*s  - 45846*s  - 57024*s)*s  +

     986922

               6        5         4          3          2           3
    *sqrt(16*s  + 124*s  - 1062*s  - 13968*s  - 45846*s  - 57024*s)*s  +

     1120203

               6        5         4          3          2           2
    *sqrt(16*s  + 124*s  - 1062*s  - 13968*s  - 45846*s  - 57024*s)*s  +

     299700

               6        5         4          3          2
    *sqrt(16*s  + 124*s  - 1062*s  - 13968*s  - 45846*s  - 57024*s)*s -

    315414
```

```
              6        5          4           3           2
*sqrt(16*s  + 124*s  - 1062*s  - 13968*s  - 45846*s  - 57024*s)

         9          8          7            6             5
 - 5184*s  - 79056*s  - 50544*s  + 6313140*s  + 54447552*s

              4              3              2
 + 215988120*s  + 467038224*s  + 533094372*s  + 249422976*s)/(81*(8

              6        5          4           3           2
*sqrt(16*s  + 124*s  - 1062*s  - 13968*s  - 45846*s  - 57024*s)

  6
*s  + 87

              6        5          4           3           2
*sqrt(16*s  + 124*s  - 1062*s  - 13968*s  - 45846*s  - 57024*s)

  5
*s  + 463

              6        5          4           3           2
*sqrt(16*s  + 124*s  - 1062*s  - 13968*s  - 45846*s  - 57024*s)

  4
*s  + 1455

              6        5          4           3           2
*sqrt(16*s  + 124*s  - 1062*s  - 13968*s  - 45846*s  - 57024*s)

  3
*s  + 2583

              6        5          4           3           2
*sqrt(16*s  + 124*s  - 1062*s  - 13968*s  - 45846*s  - 57024*s)

  2
*s  + 2646

              6        5          4           3           2
*sqrt(16*s  + 124*s  - 1062*s  - 13968*s  - 45846*s  - 57024*s)

*s + 1782

              6        5          4           3           2
*sqrt(16*s  + 124*s  - 1062*s  - 13968*s  - 45846*s  - 57024*s)

       9        8       7          6            5             4
 + 32*s  + 472*s  + 92*s  - 38796*s  - 316872*s  - 1194048*s

            3            2
 - 2425140*s  - 2535948*s  - 1026432*s))
```

Now we showed earlier that for $\{1, a_2, a_3\}$ to be a stride generator whose cover is greater than or equal to X, then $a_3 > A_3$ provided s>=23; but in any case the formula for $P(a_3)$ given by (5) is only valid for s>=40.

We now show that PP(s)<4 for all s>=40.

[ First, we use the structr function to help elucidate the structure of the function PP(s): ]

```
C: structr(pp(s));

                6                5                 4                 3
(5929*ans1*s  + 78418*ans1*s  + 401394*ans1*s  + 986922*ans1*s

                    2                                        9
  + 1120203*ans1*s  + 299700*ans1*s - 315414*ans1 - 5184*s

           8          7             6             5               4
  - 79056*s  - 50544*s  + 6313140*s  + 54447552*s  + 215988120*s

               3              2                                 6
  + 467038224*s  + 533094372*s  + 249422976*s)/(81*(8*ans1*s

                  5             4              3              2
       + 87*ans1*s  + 463*ans1*s  + 1455*ans1*s  + 2583*ans1*s

                                          9        8       7          6
       + 2646*ans1*s + 1782*ans1 + 32*s  + 472*s  + 92*s  - 38796*s

                 5              4              3              2
       - 316872*s  - 1194048*s  - 2425140*s  - 2535948*s  - 1026432*s

      ))
   where

      ans1 :=

       6        5         4          3          2            1/2
(16*s  + 124*s  - 1062*s  - 13968*s  - 45846*s  - 57024*s)
```

We see that PP(s) has the form:

$$PP(s) = \frac{G(s)p_1(s) - sp_2(s)}{81(G(s)p_3(s) + sp_4(s))}$$

where:

$$G(s) = \text{sqrt}(16s^6 + 124s^5 - 1062s^4 - 13968s^3 - 45846s^2 - 57024s)$$

$$p_1(s) = 5929s^6 + 78418s^5 + 401394s^4 + 986922s^3 + 1120203s^2 + 299700s - 315414$$

$$p_2(s) = 5184s^8 + 79056s^7 + 50544s^6 - 6313140s^5 - 54447552s^4 - 215988120s^3 - 467038224s^2 - 533094372s - 249422976$$

$$p_3(s) = 8s^6 + 87s^5 + 463s^4 + 1455s^3 + 2583s^2 + 2646s + 1782$$

$$p_4(s) = 32s^8 + 472s^7 + 92s^6 - 38796s^5 - 316872s^4 - 1194048s^3 - 2425140s^2 - 2535948s - 1026432$$

We now find values $k_1$, $k_2$ such that, for all s>=40:

$$G_1(s) = 4s^3 + k_1s^2 < G(s) < 4s^3 + k_2s^2 = G_2(s)$$

and values $c_1$, $c_2$, $c_3$, $c_4$ such that, for all s>=40:

$$p_1(s) \le 5929s^6 + c_1s^5$$

$$p_2(s) \ge 5184s^8 + c_2s^7$$

$$p_3(s) \ge 8s^6 + c_3s^5$$

$$p_4(s) \ge 32s^8 + c_4s^7$$

We then have:

$$PP(s) < \frac{(4s^3 + k_2s^2)(5929s^6 + c_1s^5) - s(5184s^8 + c_2s^7)}{81((4s^3 + k_1s^2)(8s^6 + c_3s^5) + s(32s^8 + c_4s^7))}$$

$$\Rightarrow PP(s) < \frac{(4s + k_2)(5929s + c_1) - s(5184s + c_2)}{81((4s + k_1)(8s + c_3) + s(32s + c_4))}$$

(Although I do not show this below, $k_1$, $k_2$, $c_1$, $c_2$, $c_3$ and $c_4$ are all "sharp" bounds.)

We find $k_1 = 10$:

$$G_1^2(s) = (4s^3 + 10s^2)^2 = 16s^6 + 80s^5 + 100s^4$$

$$\Rightarrow G^2(s) - G_1^2(s) = 44s^5 - (1162s^4 + 13968s^3 + 45846s^2 + 57024s)$$

Now for s>=40 we have:

$$s^5 >= 40s^4 >= 1600s^3 >= 64000s_2 >= 2560000s >= 102400000$$

So: $G^2(s) - G_1^2(s) >= 44s^5 - (29.1s^5 + 8.8s^5 + 0.8s^5 + 0.1s^5)$

$$= 44s^5 - 38.8s^5 = 5.2s^5 > 0 \text{ as required}$$

We find $k_2 = 16$:

$$G_2^2(s) = (4s^3 + 16s^2)^2 = 16s^6 + 128s^5 + 256s^4$$

$$\Rightarrow G_2^2(s) - G^2(s) = 4s^5 + 1318s^4 + 13968s^3 + 45846s^2 + 57024s > 0 \text{ as required}$$

We find $c_1 = 89088$:

$$(5929s^6 + 89088s^5) - p_1(s) =$$
$$10670s^5 - (401394s^4 + 986922s^3 + 1120203s^2 + 299700s) + 315414$$

$$>= 10670s^5 - (10034.9 + 616.9 + 17.6 + 0.2)s^5 \quad \text{for all } s>=40$$

$$= 0.4s^5 > 0 \text{ as required}$$

We find $c_2 = 75434$:

$$p_2(s) - (5184s^8 + 75434s^7) =$$
$$3622s^7 + 50544s^6 - (6313140s^5 + 54447552s^4 + 215988120s^3 + 467038224s^2 + 533094372s + 249422976)$$

$$> 3622s^7 + 50544s^6 - (157829 + 34030 + 3375 + 183 + 6 + 1)s^6 \text{ for } s>=40$$

$$= 3622s^7 - 144880s^6$$

$$>= 3622s^7 - 3622s^7 = 0 \quad \text{for all } s>=40 \text{ as required}$$

We find $c_3 = 87$:

$$p_3(s) - (8s^6 + 87s^5) =$$
$$463s^4 + 1455s^3 + 2583s^2 + 2646s + 1782 > 0 \text{ for all } s>=40 \text{ as required}$$

We find $c_4 = 444$:

$$p_4(s) - (32s^8 + 444s^7) = 28s^7 + 92s^6 - (38796s^5 + 316872s^4 + 1194048s^3 + 2425140s^2 + 2535948s + 1026432)$$

$$> 28s^7 + 92s^6 - (970 + 199 + 19 + 1 + 1 + 1)s^6 \quad \text{for all } s>=40$$

$$= 28s^7 - 1099s^6$$

$$>= 28s^7 - 27.5s^7 = 0.5s^7 > 0 \quad \text{for all } s>=40 \text{ as required}$$

Substituting back in our upper bound for PP(s) we find:

$$PP(s) < \frac{(4s + 16)(5929s + 89088) - s(5184s + 75434)}{81((4s + 10)(8s + 87) + s(32s + 444))}$$

$$= \frac{18532s^2 + 375782s + 1425408}{5184s^2 + 70632s + 70470}$$

$$< \frac{18532s^2 + 375782s + 1425408}{5184s^2 + 70632s}$$

$$< \frac{18532s^2 + 375782s + 35636s}{5184s^2 + 70632s} \qquad \text{for all } s>=40$$

$$= \frac{18532s + 411418}{5184s + 70632} = H(s), \text{ say}$$

We now show H(s) decreases as s increases:

We write $$H(s) = \frac{as + b}{cs + d}$$

Then $$H(s) - H(s+1) = \frac{as + b}{cs + d} - \frac{a(s+1) + b}{c(s+1) + d}$$

$$= \frac{(as + b)(c(s+1) + d) - (cs + d)(a(s+1) + b)}{(cs + d)(c(s+1) + d)}$$

Since c, d and s are positive, H(s) - H(s+1) is positive provided:

$$(as + b)(c(s+1) + d) - (cs + d)(a(s+1) + b) > 0$$

$$\iff asc(s+1) + asd + bc(s+1) + bd - asc(s+1) - bcs - ad(s+1) - db > 0$$

$$\iff asd + bcs + bc - bcs - ads - ad > 0$$

$$\iff bc - ad > 0$$

We find $411418*5184 - 18532*70632 = 823838688 > 0$ as required.

We see that H(59) = 3.9969..., so PP(s) < H(s) < 4 for all s>=59.

We calculate PP(s) explicitly for s=40 to s=58 to show PP(s)<4 for all s>=40:

[ We force numerical evaluation: ]

```
C: on bigfloat;
C: on numval;
```

[ And loop to calculate the required values: ]

```
C: for i := 40:58 do
C: <<
C:   write "pp(", i, ") = ", sub(s=i, pp(s));
C: >>;

   pp(40) = 3.895 87147
   pp(41) = 3.886 19646 7
   pp(42) = 3.877 08482 9
   pp(43) = 3.868 48889 6
   pp(44) = 3.860 36622 4
   pp(45) = 3.852 67889 6
   pp(46) = 3.845 39293
   pp(47) = 3.838 47778 4
   pp(48) = 3.831 90593 5
   pp(49) = 3.825 65250 9
   pp(50) = 3.819 69497 5
   pp(51) = 3.814 01287 6
   pp(52) = 3.808 58759 2
   pp(53) = 3.803 40214 3
   pp(54) = 3.798 44101 1
   pp(55) = 3.793 68998 5
   pp(56) = 3.789 13602 6
   pp(57) = 3.784 76715 3
   pp(58) = 3.780 57233 4
```

As additional confirmation of our conclusions, we calculate for s=40 the values of $a_{31}$, $a_{32}$ and the corresponding minimum and maximum values of $P(a_3)$:

[ We see that $A_3$ = 41.9512...: ]

```
C: sub(s=40, aa3);

   41.95 12195 1
```

[ and that the minimum value of $P(a_3)$ in the range $0<a_3<A_3$ occurs at $a_{31}$ = 23.8309...: ]

```
C: sub(s=40, a31);

   23.83 09969 5
```

[ and that the maximum value of $P(a_3)$ in the range $a_3>A_3$ occurs at $a_{32}$ = 175.0619...: ]

```
C: sub(s=40, a32);

   175.0 61965
```

[ and that the above-mentioned minimum value is 317.4529...: ]

```
C: sub(a3=sub(s=40,a31), sub(s=40,p(a3)));

   317.4 52965 7
```

[ and that the above-mentioned maximum value is 3.8958...: ]

```
C: sub(a3=sub(s=40,a32), sub(s=40,p(a3)));

   3.895 87147
```

Finally, we show that as s -> infinity the maximum value of P($a_3$) -> 3.5748...:

[ We first switch off numerical calculation: ]

```
C: off bigfloat;
C: off numval;
```

[ We now define lim as the limit of PP(s) as s -> infinity by simply copying out only those terms in the numerator and denominator of PP(s) with the highest power of s: ]

```
C: lim := (5929*4*s**9 - 5184*s**9)/(81*(8*4*s**9 + 32*s**9));

           4633
   lim := ------
           1296
```

[ Once again, we switch numerical calculation back on to show that the limit = 3.5748...:]

```
C: on bigfloat;
C: on numval;

C: lim;

   3.574 84567 9
```

This means that this formula does not yield any better bound on p for any larger s.

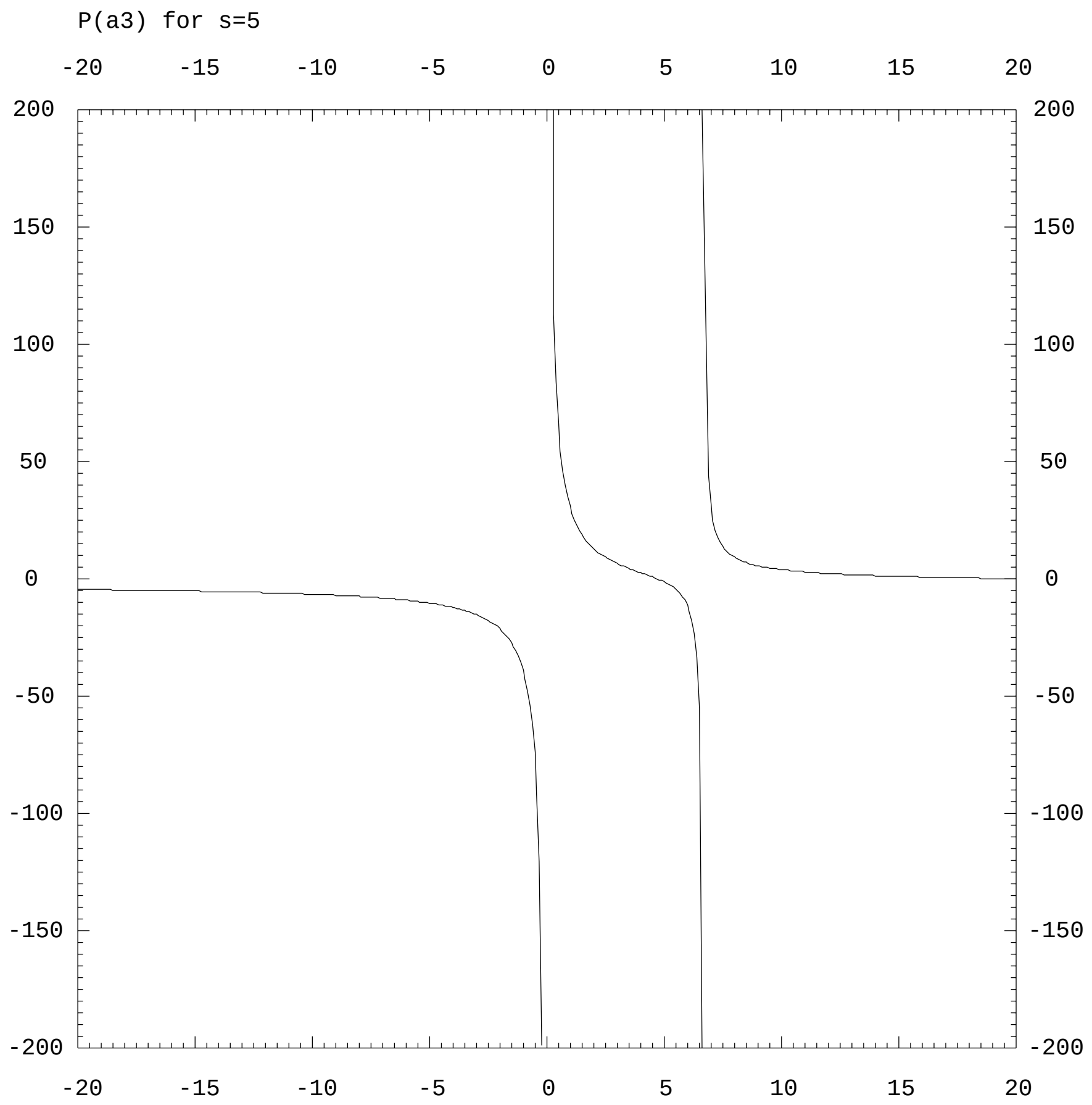

Figure 500

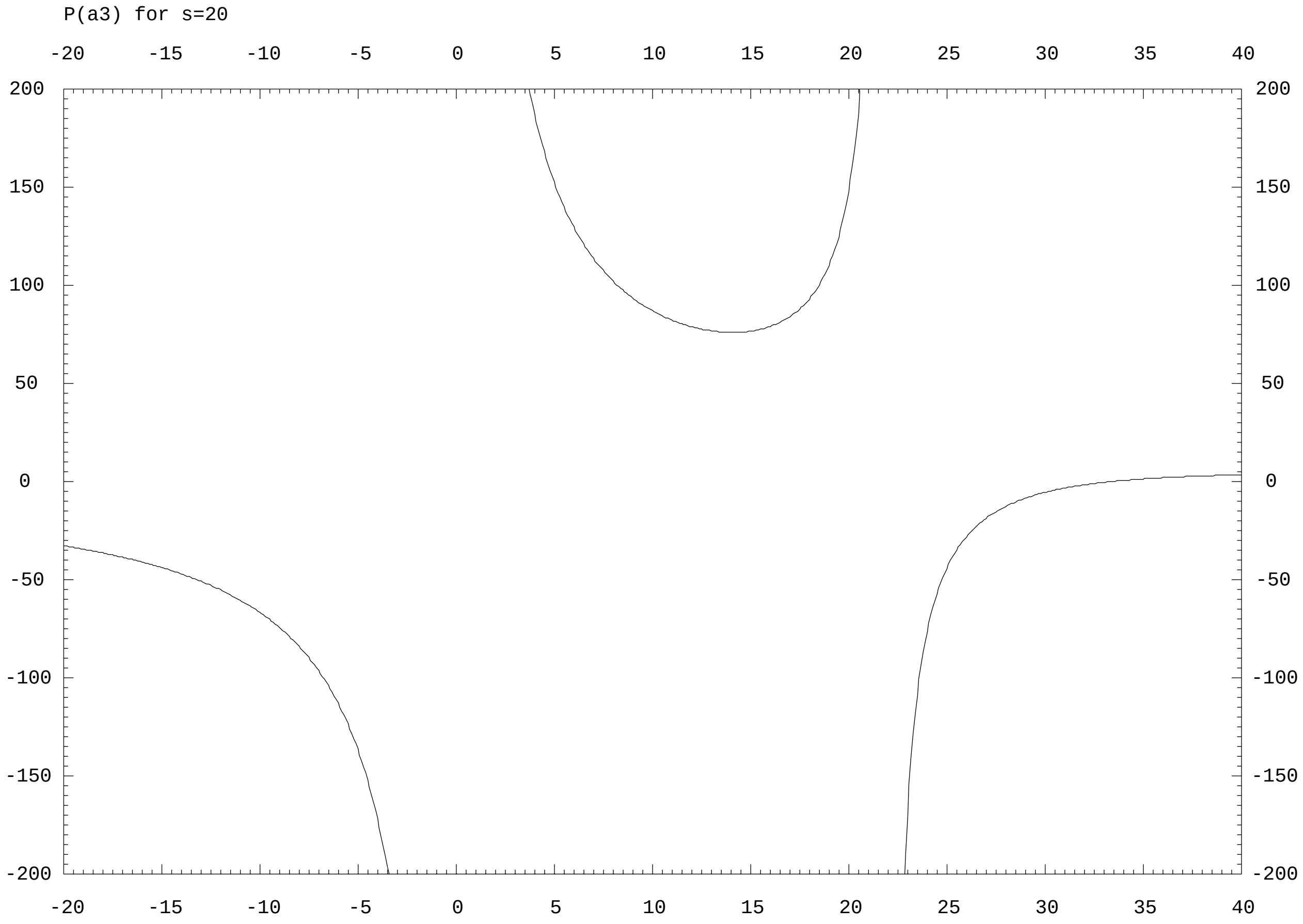
P(a3) for s=20
-20 -15 -10 -5 0 5 10 15 20 25 30 35 40
200 150 100 50 0 -50 -100 -150 -200

Figure 501

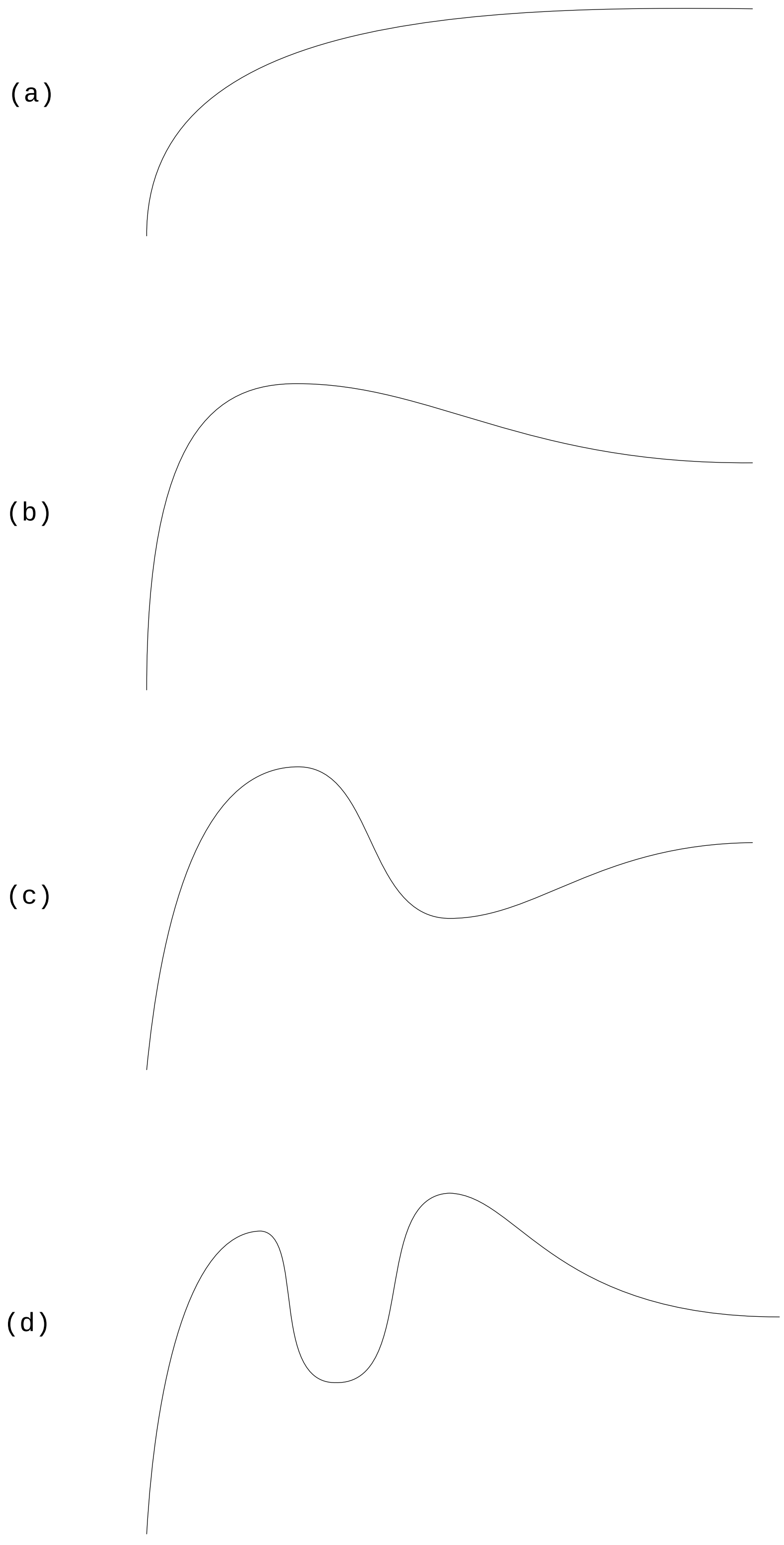

Figure 502

## 5.2 SG(n,1) is best for n>=52

### Definition 500

Any range of length $a_2$ within a stride generator must include the start of exactly one thread of order k for every 0<=k<=p (see Theorems 223 and 224).

If we choose the range 0 to $a_2$-1 inclusive, and write down the orders of each thread in the sequence in which they appear, we will have a permutation of the integers 0 to p in which the first number is 0: this is called the "signature" of the stride generator.

[ Examples:

{1,30,82} = SG(12,3) has signature (0,1,2,3) - see Figure 202

{1,30,38} = SG(4,10) has signature (0,7,3,10,6,2,9,5,1,8,4) - see Figure 219 ]

### Definition 501

The "relative position" of two threads $T_i(b_2)$ and $T_j(c_2)$ is defined by the:

"order gap" = j-i

and the "start gap" = $ST_j(c_2)-ST_i(b_2)$

We write the relative position as (order gap, start gap).

### Theorem 501

Let $T_i(b_2)$ and $T_j(c_2)$ be two threads with relative position (j-i, x), and let k be any positive or negative value such that 0 <= i+k, j+k <= p.

Then thread pairs $T_{i+k}(b_2+r)$, $T_{j+k}(c_2+r)$ will be in the same relative position for any value of r such that these threads both exist.

[ In other words, given a thread $T_m(e_2)$ and a pair of threads $T_i(b_2)$ and $T_j(c_2)$, we can derive a thread $T_{m+j-i}(e_2+c_2-b_2)$ provided only that the derived thread lies within the stride generator.

For example, suppose p=10 and we have:

$T_3$ starting at 6
$T_6$ starting at 12
$T_4$ starting at 28

Then we can deduce:

$T_9$ starting at 18, from $T_6$ via ($T_3$, $T_6$)
$T_7$ starting at 34, from $T_4$ via ($T_3$, $T_6$)
$T_2$ starting at 44, from $T_4$ via ($T_6$, $T_4$) ]

Proof

Clearly the "order gap" is the same, because (j+k)-(i+k) = j-i.

We now show that the "start gap" is the same:

$$ST_{j+k}(c_2+r) - ST_{i+k}(b_2+r) = ((c_2+r)-(b_2+r))a_2 - ((j+k)-(i+k))a_3$$
$$= (c_2-b_2)a_2 - (j-i)a_3$$
$$= ST_j(c_2) - ST_i(b_2)$$

The above theorem shows that "relative positions" must be replicated throughout a stride generator, and examination reveals that there are just two possible relative positions between successive threads: one arises when j>i, and the other when j<i.

[ Examples:

| | possible relative positions | |
|---|---|---|
| {1,30,82} = SG(12,3) | (1,8) | (-3,6) |
| {1,30,38} = SG(4,10) | (7,4) | (-4,2) ] |

If this is, indeed, so, then the corresponding order gaps are defined by the second and last elements of the corresponding signature (0,j, ... i): the "up" gap will be j, and the "down" gap will be -i. We first prove that i+j=p+1, and to do this we need to refer back to Theorem 225.

The proof of Theorem 225 shows that every stride generator includes related thread sequences $T_i$, $T_0$, $T_j$ and $T_{i-1}$, $T_p$, $T_{j-1}$ with a fundamental break at the end of $T_{i-1}$ and/or $T_p$; these triplets are illustrated in Figures 210 and 215. The following two theorems prove that i+j=p+1.

Theorem 502

Let (0,j, ... i) be the signature of a stride generator SG(n,p) for p>1.

Then i+j<=p+1.

Proof

From the proof of Theorem 225, we know there is a sequence of threads $T_{i-1}$, $T_p$, $T_{j-1}$ such that:

$$ST_{i-1} < ST_p < ST_{j-1}$$

and no other thread $T_k$ exists such that $ST_{i-1} < ST_k < ST_{j-1}$.

By Theorem 501, we can derive a thread $T_k$ from $T_{i-1}$ via the pair $(T_p, T_{j-1})$; $T_k$ will start beyond $T_{i-1}$, and will be of lower order - see Figure 503.

$$ST_k = ST_{i-1} + (ST_{j-1} - ST_p)$$

and so $ST_{i-1} < ST_k < ST_{j-1}$

So $T_k$ cannot be part of the stride generator, and so k must be less than 0:

=> k = (i-1) + (j-1-p) < 0

=> i+j-2 < p

=> i+j <= p+1

Theorem 503

Let (0,j, ... i) be the signature of a stride generator SG(n,p) for p>=1.

Then i+j=p+1.

Proof

If p=1, the signature must be (0,1), and 1+1 = 2 = p+1; so we may assume p>1.

From the proof of theorem 225, we know that there is a sequence of threads $T_{i-1}(b_2)$, $T_p(c_2)$, $T_{j-1}(d_2)$ that together cover the range from $(a_3-a_2)-1$ to $(a_3-a_2)+n+1$ inclusive.

Let us write:

$$x = ST_{i-1}(b_2)$$
$$y = ST_p(c_2)$$
$$z = ST_{j-1}(d_2)$$

Then we know from the proof of Theorem 225 that:

$$x < a_3-a_2$$
$$z > a_3-a_2$$

Now:

$$x + (i-1)a_3 = b_2a_2$$
$$z + (j-1)a_3 = d_2a_2$$

$$\Rightarrow (x+z) + (i+j-2)a_3 = (b_2+d_2)a_2$$

Subtracting $(a_3-a_2)$ from both sides we have:

$$(x+z)-(a_3-a_2) + (i+j-2)a_3 = (b_2+d_2)a_2 - (a_3-a_2)$$

$$\Rightarrow (x+z)-(a_3-a_2) + (i+j-1)a_3 = (b_2+d_2+1)a_2$$

Write $y' = (x+z)-(a_3-a_2)$:

$$y' + (i+j-1)a_3 = (b_2+d_2+1)a_2 \qquad - (1)$$

Now:

$$x < a_3-a_2 \Rightarrow (x+z)-(a_3-a_2) < z$$
$$\text{and } z > a_3-a_2 \Rightarrow (x+z)-(a_3-a_2) > x$$

So: $x < y' < z$

Theorem 502 above shows that i+j<=p+1. Furthermore, since p>1, i is not equal to j: for if it were, the signature would not be a permutation of the integers 0 ... p. This means that at most one of i and j is equal to 0, and so i+j-1>=0.

So (1) above shows $T_{i+j-1}(b_2+d_2+1)$ to be a thread of order <=p starting between $ST_{i-1}$ and $ST_{j-1}$ - and by hypothesis there is only one such thread: $T_p$.

So:

$$i+j-1 = p \Rightarrow i+j=p+1$$

and $c_2 = b_2+d_2+1$

We are now ready to prove that each stride generator contains only two "relative positions" between successive threads.

Theorem 504

Let i be an element in the signature of a stride generator SG(n,p) of order p>=1, and let j be the following element (or the first, if i is the last).

Choose $T_i$ and $T_j$ to be two successive threads in the stride generator.

Then all pairs $(T_i,T_j)$ for j>i have the same relative position, and all pairs $(T_i,T_j)$ for j<i have the same relative position.

Proof

We consider the case j>i only; the proof is analogous for the case j<i.

We complete the proof by showing that the relative position of $T_i$ and $T_j$ is the same as that of $T_0(C_2)$ and $T_J$ where J is the second element of the signature and $T_J$ is the thread immediately following $T_0(C_2)$.

To do this we consider the two threads on either side of $T_0$ and closest to it; let us call these $T_I$ and $T_J$ where $ST_I < ST_0 < ST_J$.

We now split the proof into two cases according to the value of i.

a) i <= p-J (refer to Figure 504)

In this case we derive a thread $T_{j'}$ from $T_i$ using theorem 501 and the pair $T_0$, $T_J$:

$j' = i + (J-0)$ - which is <= p by hypothesis

$ST_{j'}-ST_i = ST_J-ST_0 = x$, say

If j'=j, the theorem is proved.

If not, consider $ST_{j'}-ST_j$:

If $ST_{j'}=ST_j$, then one of $T_j$, $T_{j'}$ is covered by the other, which is not possible in a stride generator (see theorem 213).

If $ST_{j'}<ST_j$, then $T_{j'}$ would lie between $T_i$ and $T_j$ which is contrary to hypothesis.

If $ST_{j'}>ST_j$, we can derive a thread $T_{j''}$ from $T_0$ via the pair $T_i$, $T_j$:

$j'' = 0 + (j-i)$ - clearly $0<j''<=p$ since $0<=i<j<=p$ by hypothesis

$ST_{j''}-ST_0 = ST_j-ST_i < ST_{j'}-ST_i = x$

So $T_{j''}$ lies between $T_0$ and $T_J$ which is contrary to hypothesis.

b) i > p-J (refer to Figure 505)

By Theorem 503 we know that I+J = p+1, so:

$i > p-J \Rightarrow i > p - (p+1-I)$

$\Rightarrow i > I-1 \Rightarrow i >= I$

We derive a thread $T_{j'}$ from $T_i$ using the pair $T_I$, $T_0$:

$j' = i + (0-I) = i-I >= 0$

$ST_{j'} = ST_i + (ST_0-ST_I)$

Consider $ST_{j'}-ST_j$:

If $ST_{j'}=ST_j$, then one of $T_{j'}$, $T_j$ covers the other, which is not permitted in a stride generator by theorem 213.

If $ST_{j'}<ST_j$, then $T_{j'}$ would lie between $T_i$ and $T_j$ which is contrary to hypothesis.

If $ST_{j'}>ST_j$, we can derive a thread $T_{j''}$ from $T_J$ using the pair $T_j$, $T_i$:

$j'' = J + (i-j)$ - clearly $0<j''<J$ since $i-j<0$ and $i>p-J$ $\Rightarrow J+(i-j)>p-j >=0$

$$ST_{j''} = ST_J + (ST_i-ST_j)$$

$$= ST_J - (ST_j-ST_i) \qquad (\text{so } ST_{j''} < ST_J)$$

$$> ST_J - (ST_{j'}-ST_i)$$

$$= ST_J - (ST_0-ST_I)$$

$$> ST_I$$

So $T_{j''}$ ($j''>0$) lies between $T_I$ and $T_J$ which is contrary to hypothesis.

We have now shown that the order gap between successive threads in a stride generator is either j or p+1-j. Only one choice is valid for a given thread of order i because:

either a) $i+j>p$ in which case $i-(p+1-j)>=0$, and only the downward step is possible,

or b) $i+j<=p$ in which case $i-(p+1-j)<0$, and only the upward step is possible.

So all values in the signature are determined by the value j, and we can write the signature of a stride generator SG(n,p) as:

$(0, j, 2j \bmod (p+1), \ldots mj \bmod (p+1) \ldots pj \bmod (p+1))$

since if $i+j>p$, then $i-(p+1-j) = (i+j)-(p+1) = (i+j) \bmod (p+1)$.

Definition 502

The second element of the signature of a stride generator SG(n,p) is called its "key".

Theorem 505

Let SG(n,p) be a stride generator with key j.

Then j and p+1 are co-prime.

Proof

The signature $(0,j, \ldots pj \bmod (p+1))$ must be a permutation of the integers $0,1, \ldots p$, and so no integer can appear more than once in the sequence.

Suppose j and p+1 are not co-prime, but have a common factor f>1:

$$j = af$$
$$p+1 = bf \qquad 1 < b < p+1$$

Then:

$$bj \bmod (p+1) = baf \bmod (p+1) = 0$$

and so the value 0 would appear at position b as well as at position 0 in the signature.

So j and p+1 must be co-prime.

Theorem 505 proves that there can be at most p different signatures for stride generators of order p, and the possibilities for some values of p are listed below.

[ Note, however, that we have not proved that every possible signature has corresponding stride generators. ]

| | key | signature |
|---|---|---|
| p=2 | 1 | 0,1,2 |
| | 2 | 0,2,1 |
| p=3 | 1 | 0,1,2,3 |
| | 3 | 0,3,2,1 |
| p=4 | 1 | 0,1,2,3,4 |
| | 2 | 0,2,4,1,3 |
| | 3 | 0,3,1,4,2 |
| | 4 | 0,4,3,2,1 |
| p=5 | 1 | 0,1,2,3,4,5 |
| | 5 | 0,5,4,3,2,1 |
| | .... | |
| p=11 | 1 | 0,1,2,3,4,5,6,7,8,9,10,11 |
| | 5 | 0,5,10,3,8,1,6,11,4,9,2,7 |
| | 7 | 0,7,2,9,4,11,6,1,8,3,10,5 |
| | 11 | 0,11,10,9,8,7,6,5,4,3,2,1 |
| p=12 | 1 | 0,1,2,3,4,5,6,7,8,9,10,11,12 |
| | 2 | 0,2,4,6,8,10,12,1,3,5,7,9,11 |
| | 3 | 0,3,6,9,12,2,5,8,11,1,4,7,10 |
| | 4 | 0,4,8,12,3,7,11,2,6,10,1,5,9 |
| | 5 | 0,5,10,2,7,12,4,9,1,6,11,3,8 |
| | 6 | 0,6,12,5,11,4,10,3,9,2,8,1,7 |
| | 7 | 0,7,1,8,2,9,3,10,4,11,5,12,6 |
| | 8 | 0,8,3,11,6,1,9,4,12,7,2,10,5 |
| | 9 | 0,9,5,1,10,6,2,11,7,3,12,8,4 |
| | 10 | 0,10,7,4,1,11,8,5,2,12,9,6,3 |
| | 11 | 0,11,9,7,5,3,1,12,10,8,6,4,2 |
| | 12 | 0,12,11,10,9,8,7,6,5,4,3,2,1 |

Table 500

Note that the signature for key=(p+1-j) is the reverse of the signature for key=j.

For the purposes of this section, we need only consider stride generators of order 2 and 3, and in both cases the only possible keys are 1 and p.

The following theorem establishes an upper limit on the length of any stride generator SG(n,p) whose key is 1 or p*. We will then show that for sufficiently large n OSG(n,1) is longer, and finally determine by enumeration the lower limit for n for the cases p=2 and p=3.

[* These are just the "fundamental" stride generators whose thread diagrams appear as ascending or descending staircases, and which I defined in 1996 as part of my proof that $h_1,h_2<=h_0$. This later proof also shows that when A = SG(n,p) is a fundamental stride generator, there is no other stride generator SG(n',p') for A with n'>n - see Theorem 232. ]

Theorem 506

Let A = {1, $a_2$, $a_3$} be a stride generator SG(n,p) whose key is equal to 1 or p.

Then: $a_3 < (B+C)(B+3C)/4C$

where B = n(p+1)+1
and $C = p^2+p+1$

[ The proof also establishes all the conditions that must hold for A to be a stride generator with the specified characteristics. Not all of these conditions are needed to prove the limit result above, but they are useful when developing programs to list stride generators. ]

Proof

Let $a_3 = C_2a_2+C_1$

Figure 506 illustrates the threads from $T_0(C_2-1)$ to $T_0(C_2)$ for any stride generator SG(n,p) with key=1. By Theorem 225 we know that there must be a break at the end of $T_{p-1}$ or at the end of $T_p$.

Similarly Figure 507 illustrates the key=p case, and here a break must occur at the end of $T_0$ or at the end of $T_p$.

We now establish for each of these four possibilities conditions which are sufficient to ensure that A is a stride generator.

1) Key=1

We first prove that for the threads $T_i(k_i)$ in Figure 506:

$k_i = (i+1)C_2+i-1$

We observe that:

$ST_1(2C_2) = 2C_2a_2 - (C_2a_2+C_1)$

$= C_2a_2-C_1$

So $T_1(2C_2)$ is the only 1-thread that starts in the range $(C_2-1)a_2$ to $C_2a_2$ (since $C_1<a_2$), and so this must be the 1-thread in Figure 506.

So the relative position of $T_0,T_1$ is (1, $a_2-C_1$) and hence, by induction and Theorem 501, the relative position of $T_{i-1},T_i$ is (1, $a_2-C_1$) for 1<=i<=p.

This means that:

$$ST_i = ST_0(C_2-1) + i(a_2-C_1)$$

$$= (C_2-1)a_2 + ia_2 - iC_1$$

$$= (C_2-1)a_2 + ia_2 + iC_2a_2 - iC_2a_2 - iC_1$$
$$= (C_2-1+i+iC_2)a_2 - ia_3$$

$$= ((i+1)C_2+i-1)a_2 - ia_3$$

and so $ST_i$ is, in fact, $ST_i((i+1)C_2+i-1)$ for $0<=i<=p$ as shown in Figure 506.

Furthermore:

$$LT_i = n+i-(i+1)C_2-i+1+1$$

$$= n+2-(i+1)C_2 \qquad - (1)$$

and:

$$ST_i = (C_2+i-1)a_2-iC_1 \qquad - (2)$$

Next, we place limits on $C_2$.

For SG(n,p) to be a stride generator of order>0, there must be a gap between $T_0(C_2-1)$ and $T_0(C_2)$. This means that:

$$ET_0(C_2-1) + 1 < ST_0(C_2)$$

$$=> (C_2-1)a_2+n-(C_2-1)+1 < C_2a_2$$

$$=> C_2 > n-a_2+2 \qquad - (3)$$

On the other hand, an upper limit on $C_2$ can be established by considering the thread $T_p$; this must exist:

$$LT_p > 0$$

$$=> n+2-(p+1)C_2 > 0 \qquad \text{- by (1)}$$

$$=> C_2 < (n+2)/(p+1) \qquad - (4)$$

We have now ensured that there is a gap between $T_0(C_2-1)$ and $T_0(C_2)$ to be covered by threads $T_1$, $T_2$ ... $T_p$, and that all these threads exist. We now have to consider the cases where the first break y is at the end of $T_{p-1}$ and where it is at the end of $T_p$ separately.

1a) In this case y is at the end of $T_{p-1}$ and so:

$$ST_p = ET_{p-1} + 1$$

$$=> ST_p = ST_{p-1}+LT_{p-1}$$

$$=> LT_{p-1} = ST_p-ST_{p-1}$$

$$=> n+2-pC_2 = a_2-C_1$$

$$=> n = a_2+pC_2-C_1-2 \qquad - (5)$$

This ensures that we have a "staircase" of threads $T_0$, $T_1$ ... $T_p$ as shown in Figure 506 with $T_{p-1}$ and $T_p$ contiguous.

We now place limits on $C_1$.

First, $T_p$ must be long enough to reach $T_0$, and so:

$$ET_p >= ST_0 - 1$$

$$=> ST_p+LT_p >= ST_0$$

$$=> (C_2+p-1)a_2-pC_1+n+2-(p+1)C_2 >= C_2a_2$$

$$=> (p-1)a_2-pC_1+a_2+pC_2-C_1-2+2-(p+1)C_2 >= 0 \quad \text{- substituting (5)}$$

$$=> (p+1)C_1 <= pa_2-C_2$$

$$=> C_1 <= (pa_2-C_2)/(p+1) \qquad \text{- (6)}$$

On the other hand, $T_p$ must also start before $T_0$, and so:

$$ST_p < ST_0$$

$$=> (C_2-1)a_2+p(a_2-C_1) < C_2a_2$$

$$=> p(a_2-C_1) < a_2$$

$$=> C_1 > (p-1)a_2/p \qquad \text{- (7)}$$

We have now ensured that the threads do, indeed, appear as shown in Figure 506; however, these will only represent a stride generator of the given form if the value y = $ET_{p-1}$ is, indeed, a break.

This will be so iff:

$$y + ja_3 = c_2a_2 + c_1 \quad c_2+c_1<=n+j-1 \quad \text{is not soluble for any } j<=p+1$$

This is clearly not soluble for any j<=p, since y is covered only by the end of the (p-1)-thread.

It is also not soluble for j=p+1 because the (p+1)-thread starts after $ST_p$ and so cannot cover y*.

[* This is not a watertight argument - we must also show that the (p+1)-thread that starts before $ST_p$ does not also reach it - but this is of no consequence to the main proof.]

In summary, conditions (3), (4), (5), (6) and (7) are sufficient to guarantee that the threads in the configuration shown in Figure 506 define a stride generator SG(n,p) with key=1 and a break at y=$ET_{p-1}$.

1b) In this case y is at the end of $T_p$ and so:

$$ET_p = ST_0 - 1$$

$$=> (C_2+p-1)a_2-pC_1+n+2-(p+1)C_2 = C_2a_2$$

$$=> (p-1)a_2-pC_1+n+2-(p+1)C_2 = 0$$

$$=> n = (p+1)C_2+pC_1-(p-1)a_2-2 \qquad \text{- (8)}$$

This ensures that we have the staircase with $T_p$,$T_0$ contiguous.

We now place limits on $C_1$.

First, $T_p$ must cover the gap between $T_{p-1}$ and $T_0$:

$ST_p <= ET_{p-1} + 1$

$=> ST_p <= ST_{p-1} + LT_{p-1}$

$=> ST_p - ST_{p-1} <= LT_{p-1}$

$=> a_2-C_1 <= n+2-pC_2$

$=> a_2-C_1 <= (p+1)C_2+pC_1-(p-1)a_2-2+2-pC_2$ - substituting (8)

$=> a_2-C_1 <= C_2+pC_1-pa_2+a_2$

$=> (p+1)C_1 >= pa_2-C_2$

$=> C_1 >= (pa_2-C_2)/(p+1)$ - (9)

On the other hand, the break y at $ET_p$ must not be covered by the beginning of $T_{p+1}$ and so:

$ST_{p+1} >= ST_0$

$=> (C_2-1)a_2+(p+1)(a_2-C_1) >= C_2a_2$

$=> -a_2+(p+1)a_2 >= (p+1)C_1$

$=> C_1 <= pa_2/(p+1)$ - (10)

In summary, conditions (3), (4), (8), (9) and (10) are sufficient to guarantee that the threads in the configuration shown in Figure 506 define a stride generator SG(n,p) with key=1 and a break y at $ET_p$.

2) Key=p

We first prove that for the threads $T_i(k_i)$ in Figure 507:

$k_i = (i+1)C_2$

As before, we observe that:

$ST_1(2C_2) = C_2a_2-C_1$

and so $ST_1(C_2)$ must be the 1-thread in Figure 507.

So the relative position of $T_0,T_1$ is (1, $-C_1$) and hence, by induction and Theorem 501, the relative position of $T_{i-1},T_i$ is (1, $-C_1$) for 1<=i<=p.

This means that:

$ST_i = ST_0(C_2) - iC_1$

$= C_2a_2-iC_1$

$= C_2a_2+iC_2a_2-iC_2a_2-iC_1$

$= (i+1)C_2a_2-ia_3$

and so $ST_i$ is, in fact, $ST_i((i+1)C_2)$ for 0<=i<=p as shown in Figure 507.

Furthermore:

$LT_i = n+i-(i+1)C_2+1$ - (11)

and:

$ST_i = C_2a_2-iC_1$ - (12)

Next, we place limits on $C_2$.

The lower limit is the same as before:

$C_2 > n-a_2+2$ - (13)

The upper limit arises because $T_p$ must exist, and so:

$LT_p > 0$

$\Rightarrow n+p-(p+1)C_2+1 > 0$

$\Rightarrow C_2 < (n+p+1)/(p+1)$ - (14)

We now need to consider the two cases $y = ET_0$ and $y = ET_p$ separately.

2a) In this case y is at the end of $T_0$ and so:

$ST_p = ET_0 + 1$

$\Rightarrow C_2a_2-pC_1 = (C_2-1)a_2+n-(C_2-1)+1$

$\Rightarrow -pC_1 = -a_2+n-C_2+2$

$\Rightarrow n = a_2-pC_1+C_2-2$ - (15)

This ensures we have a descending staircase with $T_0,T_p$ contiguous.

We now place limits on $C_1$.

First, $T_p$ must be long enough to reach $T_{p-1}$:

$ET_p >= ST_{p-1} - 1$

$\Rightarrow ST_p + LT_p >= ST_{p-1}$

$\Rightarrow LT_p >= ST_{p-1}-ST_p$

$\Rightarrow n+p-(p+1)C_2+1 >= C_1$

$\Rightarrow C_1 <= a_2-pC_1+C_2-2+p-pC_2-C_2+1$ - substituting (15)

$\Rightarrow (p+1)C_1 <= a_2-pC_2+p-1$

$\Rightarrow C_1 <= (a_2-pC_2+p-1)/(p+1)$ - (16)

On the other hand, the break y at $ET_0$ must not be covered by any part of $T_{p+1}$ other than, possibly, its last element, so:

$ET_{p+1} <= ET_0$

$\Rightarrow C_2a_2-(p+1)C_1+n+p+1-(p+2)C_2 <= (C_2-1)a_2+n-(C_2-1)$

$\Rightarrow -(p+1)C_1+p+1-(p+2)C_2 <= -a_2-C_2+1$

$\Rightarrow -(p+1)C_1 <= -a_2-C_2-p+(p+2)C_2$

$\Rightarrow (p+1)C_1 >= a_2+p-(p+1)C_2$

$\Rightarrow C_1 >= (a_2+p)/(p+1) - C_2$ - (17)

In summary, conditions (13), (14), (15), (16) and (17) are sufficient to guarantee that the threads in Figure 507 define a stride generator SG(n,p) with key=p and a break at $y=ET_0$.

2b) In this case y is at the end of $T_p$ and so:

$ET_p = ST_{p-1} - 1$

$\Rightarrow LT_p = ST_{p-1} - ST_p$

$\Rightarrow n+p-(p+1)C_2+1 = C_1$

$\Rightarrow n = C_1+(p+1)C_2-p-1$ - (18)

This ensures we have a descending staircase with $T_p,T_{p-1}$ contiguous.

We now place limits on $C_1$.

First, $T_p$ must be long enough to reach $T_0$:

$ST_p <= ET_0 + 1$

$\Rightarrow C_2a_2-pC_1 <= (C_2-1)a_2+n-(C_2-1)+1$

$\Rightarrow -pC_1 <= -a_2+C_1+(p+1)C_2-p-1-C_2+1+1$ - substituting (18)

$\Rightarrow -(p+1)C_1 <= -a_2+pC_2-p+1$

$\Rightarrow (p+1)C_1 >= a_2-pC_2+p-1$

$\Rightarrow C_1 >= (a_2-pC_2+p-1)/(p+1)$ - (19)

Next, $T_{p-1}$ must start after $T_0$:

$ST_{p-1} > ET_0 + 1$

$\Rightarrow C_2a_2-(p-1)C_1 > (C_2-1)a_2+n-(C_2-1)+1$

$\Rightarrow -(p-1)C_1 > -a_2+n-C_2+2$

$\Rightarrow C_1 < (a_2+C_2-n-2)/(p-1)$ - (20)

Clearly the value $y = ET_p$ is not covered by any part of $T_{p+1}$ because $T_{p+1}$ does not even overlap $T_p$.

In summary, conditions (13), (14), (18), (19) and (20) are sufficient to guarantee that the threads in Figure 507 define a stride generator SG(n,p) with key=p and a break at $y=ET_p$.

We have now gathered together all the constraints associated with stride generators with key values of 1 or p, and are in a position to proceed to establish a limit on their length.

The next stage is to establish an upper bound on $C_2$ as a function of n, p and $a_2$ for each of the four cases.

1a) Key=1, break at $ET_{p-1}$

$$(5) \Rightarrow pC_2 = n-a_2+C_1+2$$

$$\Rightarrow p(p+1)C_2 \le (p+1)(n+2-a_2)+pa_2-C_2 \quad \text{- because of (6)}$$

$$\Rightarrow (p^2+p+1)C_2 \le (p+1)(n+2)-(p+1)a_2+pa_2$$

$$= (p+1)(n+2)-a_2$$

$$\Rightarrow C_2 \le ((p+1)n+2p+2-a_2)/(p^2+p+1) \qquad \text{- (21)}$$

1b) Key=1, break at $ET_p$

$$(8) \Rightarrow (p+1)C_2 = n-pC_1+(p-1)a_2+2$$

$$\Rightarrow (p+1)^2C_2 \le (p+1)n+(p^2-1)a_2+2(p+1)-p^2a_2+pC_2 \quad \text{- because of (9)}$$

$$\Rightarrow (p^2+p+1)C_2 \le (p+1)n+2p+2-a_2$$

$$\Rightarrow C_2 \le ((p+1)n+2p+2-a_2)/(p^2+p+1) \qquad \text{- (22)}$$

So both cases (1a) and (1b) result in the same condition on $C_2$ for key=1.

2a) Key=p, break at $ET_0$

$$(15) \Rightarrow C_2 = n-a_2+pC_1+2$$

$$\Rightarrow (p+1)C_2 \le (p+1)(n+2-a_2)+pa_2-p^2C_2+p^2-p \quad \text{- substituting (16)}$$

$$\Rightarrow (p^2+p+1)C_2 \le (p+1)(n+2)-a_2+p^2-p$$

$$= (p+1)n-a_2+2p+2+p^2-p$$

$$\Rightarrow C_2 \le ((p+1)n+p^2+p+2-a_2)/(p^2+p+1) \qquad \text{- (23)}$$

2b) Key=p, break at $ET_p$

$$(18) \Rightarrow (p+1)C_2 = n-C_1+p+1$$

$$\Rightarrow (p+1)^2C_2 \le (p+1)(n+p+1)-a_2+pC_2-p+1 \quad \text{- substituting (19)}$$

$$\Rightarrow (p^2+p+1)C_2 \le (p+1)n+(p+1)^2-a_2-p+1$$

$$= (p+1)n+p^2+p+2-a_2$$

$$\Rightarrow C_2 \le ((p+1)n+p^2+p+2-a_2)/(p^2+p+1) \qquad \text{- (24)}$$

Once again, cases (2a) and (2b) result in the same condition on $C_2$ for key=p.

Comparison of (21), (22) with (23), (24) shows that (23), (24) is the worst case because:

$$2p+2 \le p^2+p+2 \quad \text{for all } p \ge 1$$

We now optimise $C_2a_2$ with respect to $a_2$:

$$C_2a_2 \le (((p+1)n+p^2+p+2)a_2-a_2^2)/(p^2+p+1)$$

Simple differentiation shows this to be maximum when:

$$a_2 = ((p+1)n+p^2+p+2)/2$$

In which case:

$$C_2a_2 <= ((p+1)n+p^2+p+2)^2/4(p^2+p+1)$$

Now $C_1<a_2 \Rightarrow a_3<(C_2+1)a_2$, so:

$$a_3 < (((p+1)n+p^2+p+2)^2 + 2(p^2+p+1)((p+1)n+p^2+p+2))/4(p^2+p+1)$$

We now write:

$$B = n(p+1)+1$$
$$C = p^2+p+1$$

to show:

$$a_3 < ((B+C)^2 + 2C(B+C))/4C$$

$$\Rightarrow a_3 < (B+C)(B+3C)/4C$$

[ Programs have been written to investigate the behaviour of stride generators with key values of 1 or p, and to confirm the limit result above.

Firstly, programs were written to list all stride generators belonging to each of the cases (1a), (1b), (2a) and (2b) detailed in the proof above. These lists (for p=3 and n=2 to 10) were compared with a list of all stride generators of order 3 produced by an independent program, and were found to be the same - thus giving confidence that none have been missed in the proof!

It is interesting to note that the distribution of stride generators amongst the four cases is not even: (1a) and (2b) are more popular than (1b) and (2a):

[ My later proof - 1996 - shows that cases (1b) and (2a) are always canonical stride generators. In the notation of Lemma 14 of that document:

case (1a) is (A1)
case (1b) is (A2)
case (2a) is (D2)
case (2b) is (D1) ]

e.g. the following table gives the number of stride generators SG(n,3) of each kind for n=2 to 10:

| n | 1a | 1b | 2a | 2b |
|---|---|---|---|---|
| 2 | - | - | 1 | 2 |
| 3 | 1 | 1 | 1 | 3 |
| 4 | 2 | 2 | 1 | 4 |
| 5 | 3 | 2 | 2 | 6 |
| 6 | 4 | 2 | 3 | 8 |
| 7 | 6 | 3 | 3 | 10 |
| 8 | 8 | 4 | 3 | 12 |
| 9 | 10 | 5 | 4 | 15 |
| 10 | 12 | 5 | 5 | 18 |

Next, these four programs were combined and amended to list only the longest stride generator of each type found for a given n. The table below is an edited version of the output for p=3 and for n=2 to 60, in which the optimal SG(n,3) is indicated by an asterisk; I have also added for some lines the ratio of the limit given by the theorem above to the actual best SG(n,3).

Some interesting points to note:

a) In both cases the (A) and (B) columns are identical: in other words, the best SG(n,p) with key = 1 (or with key = p) is always one in which the p-thread fits exactly between the two threads on either side of it.

b) In the long term, key=p stride generators appear to be less successful than key=1 stride generators.

c) The ratio limit/actual decreases as n increases.

Optimal SG(n,3) with key = 1 and key = p

| | | Key = 1 | | | | Key = p | | | | Theoretical | | Theor/ |
|---|---|---|---|---|---|---|---|---|---|---|---|---|
| | | A | | B | | A | | B | | | | Actual |
| N = | 2 | 0 | 0 | 0 | 0 | 9 | 11* | 9 | 11 | 8.00 | 20.31 | 1.85 |
| N = | 3 | 7 | 12 | 7 | 12 | 13 | 16* | 13 | 16 | 10.00 | 26.00 | |
| N = | 4 | 11 | 19 | 11 | 19 | 17 | 21* | 17 | 21 | 12.00 | 32.31 | |
| N = | 5 | 15 | 26 | 15 | 26* | 21 | 26* | 21 | 26 | 14.00 | 39.23 | 1.51 |
| N = | 6 | 19 | 33 | 19 | 33* | 25 | 31 | 25 | 31 | 16.00 | 46.77 | |
| N = | 7 | 23 | 40 | 23 | 40* | 29 | 36 | 29 | 36 | 18.00 | 54.92 | |
| N = | 8 | 27 | 47 | 27 | 47* | 20 | 44 | 20 | 44 | 20.00 | 63.69 | |
| N = | 9 | 31 | 54 | 31 | 54* | 24 | 53 | 24 | 53 | 22.00 | 73.08 | |
| N = | 10 | 35 | 61 | 35 | 61 | 28 | 62* | 28 | 62 | 24.00 | 83.08 | 1.34 |
| N = | 11 | 26 | 71 | 26 | 71* | 32 | 71* | 32 | 71 | 26.00 | 93.69 | |
| N = | 12 | 30 | 82 | 30 | 82* | 36 | 80 | 36 | 80 | 28.00 | 104.92 | |
| N = | 13 | 34 | 93 | 34 | 93* | 40 | 89 | 40 | 89 | 30.00 | 116.77 | |
| N = | 14 | 38 | 104 | 38 | 104* | 31 | 99 | 31 | 99 | 32.00 | 129.23 | |
| N = | 15 | 42 | 115 | 42 | 115* | 35 | 112 | 35 | 112 | 34.00 | 142.31 | 1.24 |
| N = | 16 | 46 | 126 | 46 | 126* | 39 | 125 | 39 | 125 | 36.00 | 156.00 | |
| N = | 17 | 37 | 138 | 37 | 138* | 43 | 138* | 43 | 138 | 38.00 | 170.31 | |
| N = | 18 | 41 | 153 | 41 | 153* | 47 | 151 | 47 | 151 | 40.00 | 185.23 | |
| N = | 19 | 45 | 168 | 45 | 168* | 51 | 164 | 51 | 164 | 42.00 | 200.77 | |
| N = | 20 | 49 | 183 | 49 | 183* | 55 | 177 | 55 | 177 | 44.00 | 216.92 | 1.19 |
| N = | 21 | 53 | 198 | 53 | 198* | 46 | 193 | 46 | 193 | 46.00 | 233.69 | |
| N = | 22 | 57 | 213 | 57 | 213* | 50 | 210 | 50 | 210 | 48.00 | 251.08 | |
| N = | 23 | 61 | 228 | 61 | 228* | 54 | 227 | 54 | 227 | 50.00 | 269.08 | |
| N = | 24 | 52 | 246 | 52 | 246* | 58 | 244 | 58 | 244 | 52.00 | 287.69 | |
| N = | 25 | 56 | 265 | 56 | 265* | 62 | 261 | 62 | 261 | 54.00 | 306.92 | 1.16 |
| N = | 26 | 60 | 284 | 60 | 284* | 66 | 278 | 66 | 278 | 56.00 | 326.77 | |
| N = | 27 | 64 | 303 | 64 | 303* | 57 | 296 | 57 | 296 | 58.00 | 347.23 | |
| N = | 28 | 68 | 322 | 68 | 322* | 61 | 317 | 61 | 317 | 60.00 | 368.31 | |
| N = | 29 | 72 | 341 | 72 | 341* | 65 | 338 | 65 | 338 | 62.00 | 390.00 | |
| N = | 30 | 63 | 361 | 63 | 361* | 69 | 359 | 69 | 359 | 64.00 | 412.31 | 1.14 |
| N = | 31 | 67 | 384 | 67 | 384* | 73 | 380 | 73 | 380 | 66.00 | 435.23 | |
| N = | 32 | 71 | 407 | 71 | 407* | 77 | 401 | 77 | 401 | 68.00 | 458.77 | |
| N = | 33 | 75 | 430 | 75 | 430* | 81 | 422 | 81 | 422 | 70.00 | 482.92 | |
| N = | 34 | 79 | 453 | 79 | 453* | 72 | 446 | 72 | 446 | 72.00 | 507.69 | |
| N = | 35 | 83 | 476 | 83 | 476* | 76 | 471 | 76 | 471 | 74.00 | 533.08 | 1.12 |
| N = | 36 | 87 | 499 | 87 | 499* | 80 | 496 | 80 | 496 | 76.00 | 559.08 | |
| N = | 37 | 78 | 525 | 78 | 525* | 84 | 521 | 84 | 521 | 78.00 | 585.69 | |
| N = | 38 | 82 | 552 | 82 | 552* | 88 | 546 | 88 | 546 | 80.00 | 612.92 | |
| N = | 39 | 86 | 579 | 86 | 579* | 92 | 571 | 92 | 571 | 82.00 | 640.77 | |
| N = | 40 | 90 | 606 | 90 | 606* | 83 | 597 | 83 | 597 | 84.00 | 669.23 | 1.10 |
| N = | 41 | 94 | 633 | 94 | 633* | 87 | 626 | 87 | 626 | 86.00 | 698.31 | |
| N = | 42 | 98 | 660 | 98 | 660* | 91 | 655 | 91 | 655 | 88.00 | 728.00 | |
| N = | 43 | 89 | 688 | 89 | 688* | 95 | 684 | 95 | 684 | 90.00 | 758.31 | |
| N = | 44 | 93 | 719 | 93 | 719* | 99 | 713 | 99 | 713 | 92.00 | 789.23 | |
| N = | 45 | 97 | 750 | 97 | 750* | 103 | 742 | 103 | 742 | 94.00 | 820.77 | 1.09 |
| N = | 46 | 101 | 781 | 101 | 781* | 107 | 771 | 107 | 771 | 96.00 | 852.92 | |
| N = | 47 | 105 | 812 | 105 | 812* | 98 | 803 | 98 | 803 | 98.00 | 885.69 | |
| N = | 48 | 109 | 843 | 109 | 843* | 102 | 836 | 102 | 836 | 100.00 | 919.08 | |
| N = | 49 | 113 | 874 | 113 | 874* | 106 | 869 | 106 | 869 | 102.00 | 953.08 | |
| N = | 50 | 104 | 908 | 104 | 908* | 110 | 902 | 110 | 902 | 104.00 | 987.69 | 1.09 |
| N = | 51 | 108 | 943 | 108 | 943* | 114 | 935 | 114 | 935 | 106.00 | 1022.92 | |
| N = | 52 | 112 | 978 | 112 | 978* | 118 | 968 | 118 | 968 | 108.00 | 1058.77 | |
| N = | 53 | 116 | 1013 | 116 | 1013* | 109 | 1002 | 109 | 1002 | 110.00 | 1095.23 | |
| N = | 54 | 120 | 1048 | 120 | 1048* | 113 | 1039 | 113 | 1039 | 112.00 | 1132.31 | |
| N = | 55 | 124 | 1083 | 124 | 1083* | 117 | 1076 | 117 | 1076 | 114.00 | 1170.00 | 1.08 |
| N = | 56 | 115 | 1119 | 115 | 1119* | 121 | 1113 | 121 | 1113 | 116.00 | 1208.31 | |
| N = | 57 | 119 | 1158 | 119 | 1158* | 125 | 1150 | 125 | 1150 | 118.00 | 1247.23 | |
| N = | 58 | 123 | 1197 | 123 | 1197* | 129 | 1187 | 129 | 1187 | 120.00 | 1286.77 | |
| N = | 59 | 127 | 1236 | 127 | 1236* | 133 | 1224 | 133 | 1224 | 122.00 | 1326.92 | |
| N = | 60 | 131 | 1275 | 131 | 1275* | 124 | 1264 | 124 | 1264 | 124.00 | 1367.69 | 1.07 |

Table 501

All the points above are also confirmed for p=7, as shown in the following table:

Optimal SG(n,7) with key = 1 and key = p

| | Key = 1 | | | | Key = p | | | | Theoretical | | |
|---|---|---|---|---|---|---|---|---|---|---|---|
| | A | | B | | A | | B | | | | |
| N = 2 | 0 | 0 | 0 | 0 | 17 | 19* | 17 | 19 | 16.00 | 61.02 | 3.21 |
| N = 3 | 0 | 0 | 0 | 0 | 25 | 28* | 25 | 28 | 20.00 | 70.49 | |
| N = 4 | 0 | 0 | 0 | 0 | 33 | 37* | 33 | 37 | 24.00 | 80.53 | |
| N = 5 | 0 | 0 | 0 | 0 | 41 | 46* | 41 | 46 | 28.00 | 91.12 | |
| N = 6 | 0 | 0 | 0 | 0 | 49 | 55* | 49 | 55 | 32.00 | 102.28 | |
| N = 7 | 15 | 28 | 15 | 28 | 57 | 64* | 57 | 64 | 36.00 | 114.00 | 1.78 |
| N = 8 | 23 | 43 | 23 | 43 | 65 | 73* | 65 | 73 | 40.00 | 126.28 | |
| N = 9 | 31 | 58 | 31 | 58 | 73 | 82* | 73 | 82 | 44.00 | 139.12 | |
| N = 10 | 39 | 73 | 39 | 73 | 81 | 91* | 81 | 91 | 48.00 | 152.53 | |
| N = 11 | 47 | 88 | 47 | 88 | 89 | 100* | 89 | 100 | 52.00 | 166.49 | |
| N = 12 | 55 | 103 | 55 | 103 | 97 | 109* | 97 | 109 | 56.00 | 181.02 | |
| N = 13 | 63 | 118 | 63 | 118* | 105 | 118* | 105 | 118 | 60.00 | 196.11 | |
| N = 14 | 71 | 133 | 71 | 133* | 113 | 127 | 113 | 127 | 64.00 | 211.75 | |
| | | | | | | | | | | | |
| N = 21 | 127 | 238 | 127 | 238* | 112 | 237 | 112 | 237 | 92.00 | 337.02 | |
| N = 22 | 135 | 253 | 135 | 253 | 120 | 254* | 120 | 254 | 96.00 | 357.16 | |
| N = 23 | 143 | 268 | 143 | 268 | 128 | 271* | 128 | 271 | 100.00 | 377.86 | |
| N = 24 | 151 | 283 | 151 | 283 | 136 | 288* | 136 | 288 | 104.00 | 399.12 | |
| N = 25 | 159 | 298 | 159 | 298 | 144 | 305* | 144 | 305 | 108.00 | 420.95 | 1.38 |
| N = 26 | 110 | 316 | 110 | 316 | 152 | 322* | 152 | 322 | 112.00 | 443.33 | |
| N = 27 | 118 | 339 | 118 | 339* | 160 | 339* | 160 | 339 | 116.00 | 466.28 | |
| N = 28 | 126 | 362 | 126 | 362* | 168 | 356 | 168 | 356 | 120.00 | 489.79 | |
| | | | | | | | | | | | |
| N = 36 | 190 | 546 | 190 | 546* | 175 | 545 | 175 | 545 | 152.00 | 698.07 | |
| N = 37 | 198 | 569 | 198 | 569 | 183 | 570* | 183 | 570 | 156.00 | 726.63 | |
| N = 38 | 206 | 592 | 206 | 592 | 191 | 595* | 191 | 595 | 160.00 | 755.75 | |
| N = 39 | 214 | 615 | 214 | 615 | 199 | 620* | 199 | 620 | 164.00 | 785.44 | |
| N = 40 | 165 | 639 | 165 | 639 | 207 | 645* | 207 | 645 | 168.00 | 815.68 | 1.26 |
| N = 41 | 173 | 670 | 173 | 670* | 215 | 670* | 215 | 670 | 172.00 | 846.49 | |
| N = 42 | 181 | 701 | 181 | 701* | 223 | 695 | 223 | 695 | 176.00 | 877.86 | |
| | | | | | | | | | | | |
| N = 51 | 253 | 980 | 253 | 980* | 238 | 979 | 238 | 979 | 212.00 | 1185.44 | |
| N = 52 | 261 | 1011 | 261 | 1011 | 246 | 1012* | 246 | 1012 | 216.00 | 1222.42 | |
| N = 53 | 269 | 1042 | 269 | 1042 | 254 | 1045* | 254 | 1045 | 220.00 | 1259.96 | |
| N = 54 | 277 | 1073 | 277 | 1073 | 262 | 1078* | 262 | 1078 | 224.00 | 1298.07 | 1.20 |
| N = 55 | 228 | 1111 | 228 | 1111* | 270 | 1111* | 270 | 1111 | 228.00 | 1336.74 | |
| N = 56 | 236 | 1150 | 236 | 1150* | 278 | 1144 | 278 | 1144 | 232.00 | 1375.96 | |
| | | | | | | | | | | | |
| N = 66 | 316 | 1540 | 316 | 1540* | 301 | 1539 | 301 | 1539 | 272.00 | 1799.12 | |
| N = 67 | 324 | 1579 | 324 | 1579 | 309 | 1580* | 309 | 1580 | 276.00 | 1844.53 | |
| N = 68 | 332 | 1618 | 332 | 1618 | 317 | 1621* | 317 | 1621 | 280.00 | 1890.49 | 1.17 |
| N = 69 | 283 | 1662 | 283 | 1662* | 325 | 1662* | 325 | 1662 | 284.00 | 1937.02 | |
| N = 70 | 291 | 1709 | 291 | 1709* | 333 | 1703 | 333 | 1703 | 288.00 | 1984.11 | |
| | | | | | | | | | | | |
| N = 81 | 379 | 2226 | 379 | 2226* | 364 | 2225 | 364 | 2225 | 332.00 | 2539.12 | |
| N = 82 | 387 | 2273 | 387 | 2273 | 372 | 2274* | 372 | 2274 | 336.00 | 2592.95 | 1.14 |
| N = 83 | 338 | 2323 | 338 | 2323* | 380 | 2323* | 380 | 2323 | 340.00 | 2647.33 | |
| N = 84 | 346 | 2378 | 346 | 2378* | 388 | 2372 | 388 | 2372 | 344.00 | 2702.28 | |
| | | | | | | | | | | | |
| N = 96 | 442 | 3038 | 442 | 3038* | 427 | 3037 | 427 | 3037 | 392.00 | 3405.44 | |
| N = 97 | 393 | 3094 | 393 | 3094* | 435 | 3094* | 435 | 3094 | 396.00 | 3467.68 | 1.12 |
| N = 98 | 401 | 3157 | 401 | 3157* | 443 | 3151 | 443 | 3151 | 400.00 | 3530.49 | |

Table 502 ]

The theorem above gives an upper limit to the length of stride generators with keys equal to 1 or p. We now show that as n tends to infinity this limit tends to a value less than the length of OSG(n,1).

Theorem 507

The maximum length of a stride generator SG(n,p) whose key is equal to 1 or p is less than the length of OSG(n,1) for all sufficiently large n.

Proof

Theorem 506 gives an upper bound $L_p$ on the length of any stride generator SG(n,p) whose key is equal to 1 or p:

$$L_p = (B+C)(B+3C)/4C \quad \text{where} \quad B = n(p+1)+1 \quad \text{and} \quad C = p^2+p+1$$

As n -> infinity for fixed p:

$$L_p \to (n^2(p+1)^2)/4(p^2+p+1)$$

Theorem 301 shows that:

$$L_1 = (n^2+5n+7)/3$$

is an upper bound for the length of OSG(n,1).

As n -> infinity:

$$L_1 \to n^2/3$$

So:

$$L_p/L_1 \to 3(p+1)^2/4(p^2+p+1)$$

$$= 3(p^2+2p+1)/4(p^2+p+1)$$

$$< 1 \quad \text{for all } p>1$$

Since stride generators of order 2 or 3 must have a key equal to 1 or p, Theorem 506 gives upper limits for the length of OSG(n,2) and OSG(n,3). The next theorem determines specific values of n above which OSG(n,1) is always better than OSG(n,2) and OSG(n,3).

Theorem 508

Let $L_i$ be the length of OSG(n,i).

Then:

a) $L_1>L_2$ for all n>132

b) $L_1>L_3$ for all n>100

Proof

By Theorem 301:

$$L_1 >= (n^2+5n+6)/3$$

By Theorem 506:

$$L_2 < (B+C)(B+3C)/4C \quad \text{where } B = 3n+1 \text{ and } C = 7$$

$$=> L_2 < (3n+8)(3n+22)/28$$

$$= (9n^2+90n+176)/28$$

So $L_1$ is certainly greater than $L_2$ when:

$$28(n^2+5n+6) > 3(9n^2+90n+176)$$

$$=> 28n^2+140n+168 > 27n^2+270n+528$$

$$=> n^2-130n-360 > 0$$

Solving the quadratic we have:

$$n = (130+135.43)/2 \quad \text{or} \quad n = (130-135.43)/2$$

Clearly:

$$L_1>L_2 \quad \text{if} \quad n>132.71$$

Similarly, by Theorem 506:

$$L_3 < (B+C)(B+3C)/4C \quad \text{where } B = 4n+1 \text{ and } C = 13$$

$$=> L_3 < (4n+14)(4n+40)/52$$

$$= (16n^2+216n+560)/52$$

So $L_1$ is certainly greater than $L_3$ when:

$$52(n^2+5n+6) > 3(16n^2+216n+560)$$

$$=> 4n^2-388n-1368 > 0$$

Solving the quadratic we have:

$$n = (388+415.25)/8 \quad \text{or} \quad n = (388-415.25)/8$$

Clearly:

$$L_1>L_3 \quad \text{if} \quad n>100.41$$

We now show that OSG(n,1) is always better than OSG(n,0).

Theorem 509

Let $L_i$ be the length of OSG(n,i).

Then $L_1 > L_0$ for all n.

Proof

By Theorem 301:

$$L_1 >= (n^2+5n+6)/3$$

By Theorem 300:

$$L_0 <= (n^2+6n+5)/4$$

So:

$$12(L_1-L_0) >= 4(n^2+5n+6) - 3(n^2+6n+5)$$

$$= n^2+2n+9$$

$$> 0 \quad \text{for all } n>=1$$

Finally, we show by exhaustive enumeration that OSG(n,1) is actually better than OSG(n,2) and OSG(n,3) for all n>48, and that SG1(n,1) is better for all n>51.

Theorem 510

Let $L_i$ be the length of OSG(n,i).

Then:

a) $L_1 > L_i$ for i=0,2,3 for all n>48

b) $L_1 - 1 > L_i$ for i=0,2,3 for all n>51

Proof

Theorem 509 deals with i=0.

Theorem 508 proves that we need only investigate values of n<=132.

We use a program (exp15 - see section 6.1) to list all OSG(n,i) for i=1,2,3 and for n<=132, to show that in fact OSG(n,1) is longer than OSG(n,i) for i=2,3 for all n>48, and that SG1(n,1) is longer than OSG(n,i) for i=2,3 for all n>51.

The results are shown in Table 503 below, where the two columns on the right give limits for OSG(n,2) and OSG(n,3) from Theorem 506, and the stride generator of maximum length for each n is indicated by an asterisk.

Optimal stride generators of orders 1, 2 and 3 for n <= 133

| n | order 1 | | order 2 | | order 3 | | Limit for order 2 | Limit for order 3 |
|---|---|---|---|---|---|---|---|---|
| 1 | 3 | 4 | 4 | 5 | 5 | 6* | 9.82 | 15.23 |
| 2 | 5 | 7 | 7 | 9 | 9 | 11* | 14.00 | 20.31 |
| 3 | 7 | 10 | 10 | 13 | 13 | 16* | 18.82 | 26.00 |
| | | | 8 | 13 | | | | |
| 4 | 6 | 14 | 11 | 18 | 17 | 21* | 24.29 | 32.31 |
| 5 | 8 | 19 | 14 | 23 | 21 | 26* | 30.39 | 39.23 |
| | | | | | 15 | 26* | | |
| 6 | 10 | 24 | 17 | 28 | 19 | 33* | 37.14 | 46.77 |
| 7 | 9 | 30 | 15 | 34 | 23 | 40* | 44.54 | 54.92 |
| | | | 13 | 34 | | | | |
| 8 | 11 | 37 | 16 | 42 | 27 | 47* | 52.57 | 63.69 |
| 9 | 13 | 44 | 19 | 50 | 31 | 54* | 61.25 | 73.08 |
| 10 | 12 | 52 | 22 | 58 | 28 | 62* | 70.57 | 83.08 |
| 11 | 14 | 61 | 25 | 66 | 32 | 71* | 80.54 | 93.69 |
| | | | | | 26 | 71* | | |
| 12 | 16 | 70 | 21 | 76 | 30 | 82* | 91.14 | 104.92 |
| 13 | 15 | 80 | 24 | 87 | 34 | 93* | 102.39 | 116.77 |
| 14 | 17 | 91 | 27 | 98 | 38 | 104* | 114.29 | 129.23 |
| 15 | 19 | 102 | 30 | 109 | 42 | 115* | 126.82 | 142.31 |
| 16 | 18 | 114 | 33 | 120 | 46 | 126* | 140.00 | 156.00 |
| | | | 26 | 120 | | | | |
| 17 | 20 | 127 | 29 | 134 | 43 | 138* | 153.82 | 170.31 |
| | | | | | 37 | 138* | | |
| 18 | 22 | 140 | 32 | 148 | 41 | 153* | 168.29 | 185.23 |
| 19 | 21 | 154 | 35 | 162 | 45 | 168* | 183.39 | 200.77 |
| 20 | 23 | 169 | 38 | 176 | 49 | 183* | 199.14 | 216.92 |
| 21 | 25 | 184 | 34 | 191 | 53 | 198* | 215.54 | 233.69 |
| 22 | 24 | 200 | 37 | 208 | 57 | 213* | 232.57 | 251.08 |
| 23 | 26 | 217 | 40 | 225 | 61 | 228* | 250.25 | 269.08 |
| 24 | 28 | 234 | 43 | 242 | 52 | 246* | 268.57 | 287.69 |
| 25 | 27 | 252 | 46 | 259 | 56 | 265* | 287.54 | 306.92 |
| 26 | 29 | 271 | 42 | 278 | 60 | 284* | 307.14 | 326.77 |
| 27 | 31 | 290 | 45 | 298 | 64 | 303* | 327.39 | 347.23 |
| 28 | 30 | 310 | 48 | 318 | 68 | 322* | 348.29 | 368.31 |
| 29 | 32 | 331 | 51 | 338 | 72 | 341* | 369.82 | 390.00 |
| 30 | 34 | 352 | 54 | 358 | 63 | 361* | 392.00 | 412.31 |
| | | | 47 | 358 | | | | |
| 31 | 33 | 374 | 50 | 381 | 67 | 384* | 414.82 | 435.23 |
| 32 | 35 | 397 | 53 | 404 | 71 | 407* | 438.29 | 458.77 |
| 33 | 37 | 420 | 56 | 427 | 75 | 430* | 462.39 | 482.92 |
| 34 | 36 | 444 | 59 | 450 | 79 | 453* | 487.14 | 507.69 |
| 35 | 38 | 469 | 55 | 474 | 83 | 476* | 512.54 | 533.08 |
| 36 | 40 | 494 | 58 | 500* | 87 | 499 | 538.57 | 559.08 |
| 37 | 39 | 520 | 61 | 526* | 78 | 525 | 565.25 | 585.69 |
| 38 | 41 | 547 | 64 | 552* | 82 | 552* | 592.57 | 612.92 |
| 39 | 43 | 574 | 67 | 578 | 86 | 579* | 620.54 | 640.77 |
| 40 | 42 | 602 | 63 | 606* | 90 | 606* | 649.14 | 669.23 |
| 41 | 44 | 631 | 66 | 635* | 94 | 633 | 678.39 | 698.31 |
| 42 | 46 | 660 | 69 | 664* | 98 | 660 | 708.29 | 728.00 |
| 43 | 45 | 690 | 72 | 693* | 89 | 688 | 738.82 | 758.31 |
| 44 | 47 | 721 | 75 | 722* | 93 | 719 | 770.00 | 789.23 |
| | | | 68 | 722* | | | | |
| 45 | 49 | 752 | 71 | 754* | 97 | 750 | 801.82 | 820.77 |
| 46 | 48 | 784 | 74 | 786* | 101 | 781 | 834.29 | 852.92 |
| 47 | 50 | 817 | 77 | 818* | 105 | 812 | 867.39 | 885.69 |
| 48 | 52 | 850* | 80 | 850* | 109 | 843 | 901.14 | 919.08 |
| 49 | 51 | 884* | 76 | 883 | 113 | 874 | 935.54 | 953.08 |
| 50 | 53 | 919* | 79 | 918 | 104 | 908 | 970.57 | 987.69 |
| 51 | 55 | 954* | 82 | 953 | 108 | 943 | 1006.25 | 1022.92 |
| 52 | 54 | 990* | 85 | 988 | 112 | 978 | 1042.57 | 1058.77 |
| 53 | 56 | 1027* | 88 | 1023 | 116 | 1013 | 1079.54 | 1095.23 |
| 54 | 58 | 1064* | 84 | 1060 | 120 | 1048 | 1117.14 | 1132.31 |

| n | order 1 | order 2 | order 3 | Limit for order 2 | Limit for order 3 |
|---|---|---|---|---|---|
| 55 | 57 1102* | 87 1098 | 124 1083 | 1155.39 | 1170.00 |
| 56 | 59 1141* | 90 1136 | 115 1119 | 1194.29 | 1208.31 |
| 57 | 61 1180* | 93 1174 | 119 1158 | 1233.82 | 1247.23 |
| 58 | 60 1220* | 96 1212 | 123 1197 | 1274.00 | 1286.77 |
| | | 89 1212 | | | |
| 59 | 62 1261* | 92 1253 | 127 1236 | 1314.82 | 1326.92 |
| 60 | 64 1302* | 95 1294 | 131 1275 | 1356.29 | 1367.69 |
| 61 | 63 1344* | 98 1335 | 135 1314 | 1398.39 | 1409.08 |
| 62 | 65 1387* | 101 1376 | 139 1353 | 1441.14 | 1451.08 |
| 63 | 67 1430* | 97 1418 | 130 1395 | 1484.54 | 1493.69 |
| 64 | 66 1474* | 100 1462 | 134 1438 | 1528.57 | 1536.92 |
| 65 | 68 1519* | 103 1506 | 138 1481 | 1573.25 | 1580.77 |
| 66 | 70 1564* | 106 1550 | 142 1524 | 1618.57 | 1625.23 |
| 67 | 69 1610* | 109 1594 | 146 1567 | 1664.54 | 1670.31 |
| 68 | 71 1657* | 105 1640 | 150 1610 | 1711.14 | 1716.00 |
| 69 | 73 1704* | 108 1687 | 141 1654 | 1758.39 | 1762.31 |
| 70 | 72 1752* | 111 1734 | 145 1701 | 1806.29 | 1809.23 |
| 71 | 74 1801* | 114 1781 | 149 1748 | 1854.82 | 1856.77 |
| 72 | 76 1850* | 117 1828 | 153 1795 | 1904.00 | 1904.92 |
| | | 110 1828 | | | |
| 73 | 75 1900* | 113 1878 | 157 1842 | 1953.82 | 1953.69 |
| 74 | 77 1951* | 116 1928 | 161 1889 | 2004.29 | 2003.08 |
| 75 | 79 2002* | 119 1978 | 165 1936 | 2055.39 | 2053.08 |
| 76 | 78 2054* | 122 2028 | 156 1986 | 2107.14 | 2103.69 |
| 77 | 80 2107* | 118 2079 | 160 2037 | 2159.54 | 2154.92 |
| 78 | 82 2160* | 121 2132 | 164 2088 | 2212.57 | 2206.77 |
| 79 | 81 2214* | 124 2185 | 168 2139 | 2266.25 | 2259.23 |
| 80 | 83 2269* | 127 2238 | 172 2190 | 2320.57 | 2312.31 |
| 81 | 85 2324* | 130 2291 | 176 2241 | 2375.54 | 2366.00 |
| 82 | 84 2380* | 126 2346 | 167 2293 | 2431.14 | 2420.31 |
| 83 | 86 2437* | 129 2402 | 171 2348 | 2487.39 | 2475.23 |
| 84 | 88 2494* | 132 2458 | 175 2403 | 2544.29 | 2530.77 |
| 85 | 87 2552* | 135 2514 | 179 2458 | 2601.82 | 2586.92 |
| 86 | 89 2611* | 138 2570 | 183 2513 | 2660.00 | 2643.69 |
| | | 131 2570 | | | |
| 87 | 91 2670* | 134 2629 | 187 2568 | 2718.82 | 2701.08 |
| 88 | 90 2730* | 137 2688 | 191 2623 | 2778.29 | 2759.08 |
| 89 | 92 2791* | 140 2747 | 182 2681 | 2838.39 | 2817.69 |
| 90 | 94 2852* | 143 2806 | 186 2740 | 2899.14 | 2876.92 |
| 91 | 93 2914* | 139 2866 | 190 2799 | 2960.54 | 2936.77 |
| 92 | 95 2977* | 142 2928 | 194 2858 | 3022.57 | 2997.23 |
| 93 | 97 3040* | 145 2990 | 198 2917 | 3085.25 | 3058.31 |
| 94 | 96 3104* | 148 3052 | 202 2976 | 3148.57 | 3120.00 |
| 95 | 98 3169* | 151 3114 | 193 3036 | 3212.54 | 3182.31 |
| 96 | 100 3234* | 147 3178 | 197 3099 | 3277.14 | 3245.23 |
| 97 | 99 3300* | 150 3243 | 201 3162 | 3342.39 | 3308.77 |
| 98 | 101 3367* | 153 3308 | 205 3225 | 3408.29 | 3372.92 |
| 99 | 103 3434* | 156 3373 | 209 3288 | 3474.82 | 3437.69 |
| 100 | 102 3502* | 159 3438 | 213 3351 | 3542.00 | 3503.08 |
| | | 152 3438 | | | |
| 101 | 104 3571* | 155 3506 | 217 3414 | 3609.82 | 3569.08 |
| 102 | 106 3640* | 158 3574 | 208 3480 | 3678.29 | 3635.69 |
| 103 | 105 3710* | 161 3642 | 212 3547 | 3747.39 | 3702.92 |
| 104 | 107 3781* | 164 3710 | 216 3614 | 3817.14 | 3770.77 |
| 105 | 109 3852* | 160 3779 | 220 3681 | 3887.54 | 3839.23 |
| 106 | 108 3924* | 163 3850 | 224 3748 | 3958.57 | 3908.31 |
| 107 | 110 3997* | 166 3921 | 228 3815 | 4030.25 | 3978.00 |
| 108 | 112 4070* | 169 3992 | 219 3883 | 4102.57 | 4048.31 |
| 109 | 111 4144* | 172 4063 | 223 3954 | 4175.54 | 4119.23 |
| 110 | 113 4219* | 168 4136 | 227 4025 | 4249.14 | 4190.77 |
| 111 | 115 4294* | 171 4210 | 231 4096 | 4323.39 | 4262.92 |
| 112 | 114 4370* | 174 4284 | 235 4167 | 4398.29 | 4335.69 |
| 113 | 116 4447* | 177 4358 | 239 4238 | 4473.82 | 4409.08 |

| n | order 1 | order 2 | order 3 | Limit for order 2 | Limit for order 3 |
|---|---|---|---|---|---|
| 114 | 118 4524* | 180 4432 | 243 4309 | 4550.00 | 4483.08 |
| | | 173 4432 | | | |
| 115 | 117 4602* | 176 4509 | 234 4383 | 4626.82 | 4557.69 |
| 116 | 119 4681* | 179 4586 | 238 4458 | 4704.29 | 4632.92 |
| 117 | 121 4760* | 182 4663 | 242 4533 | 4782.39 | 4708.77 |
| 118 | 120 4840* | 185 4740 | 246 4608 | 4861.14 | 4785.23 |
| 119 | 122 4921* | 181 4818 | 250 4683 | 4940.54 | 4862.31 |
| 120 | 124 5002* | 184 4898 | 254 4758 | 5020.57 | 4940.00 |
| 121 | 123 5084* | 187 4978 | 245 4834 | 5101.25 | 5018.31 |
| 122 | 125 5167* | 190 5058 | 249 4913 | 5182.57 | 5097.23 |
| 123 | 127 5250* | 193 5138 | 253 4992 | 5264.54 | 5176.77 |
| 124 | 126 5334* | 189 5220 | 257 5071 | 5347.14 | 5256.92 |
| 125 | 128 5419* | 192 5303 | 261 5150 | 5430.39 | 5337.69 |
| 126 | 130 5504* | 195 5386 | 265 5229 | 5514.29 | 5419.08 |
| 127 | 129 5590* | 198 5469 | 269 5308 | 5598.82 | 5501.08 |
| 128 | 131 5677* | 201 5552 | 260 5390 | 5684.00 | 5583.69 |
| | | 194 5552 | | | |
| 129 | 133 5764* | 197 5638 | 264 5473 | 5769.82 | 5666.92 |
| 130 | 132 5852* | 200 5724 | 268 5556 | 5856.29 | 5750.77 |
| 131 | 134 5941* | 203 5810 | 272 5639 | 5943.39 | 5835.23 |
| 132 | 136 6030* | 206 5896 | 276 5722 | 6031.14 | 5920.31 |
| 133 | 135 6120* | 202 5983 | 280 5805 | 6119.54 | 6006.00 |

Table 503

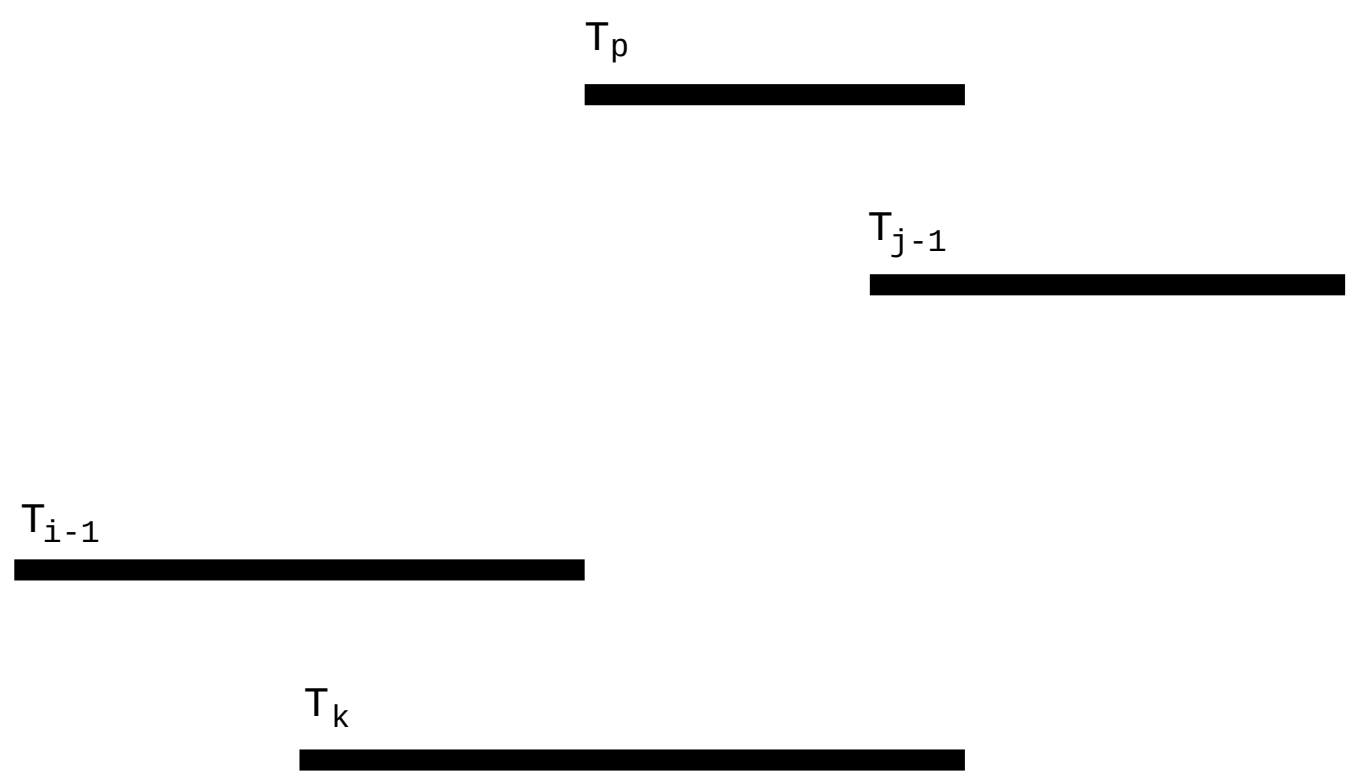

Figure 503

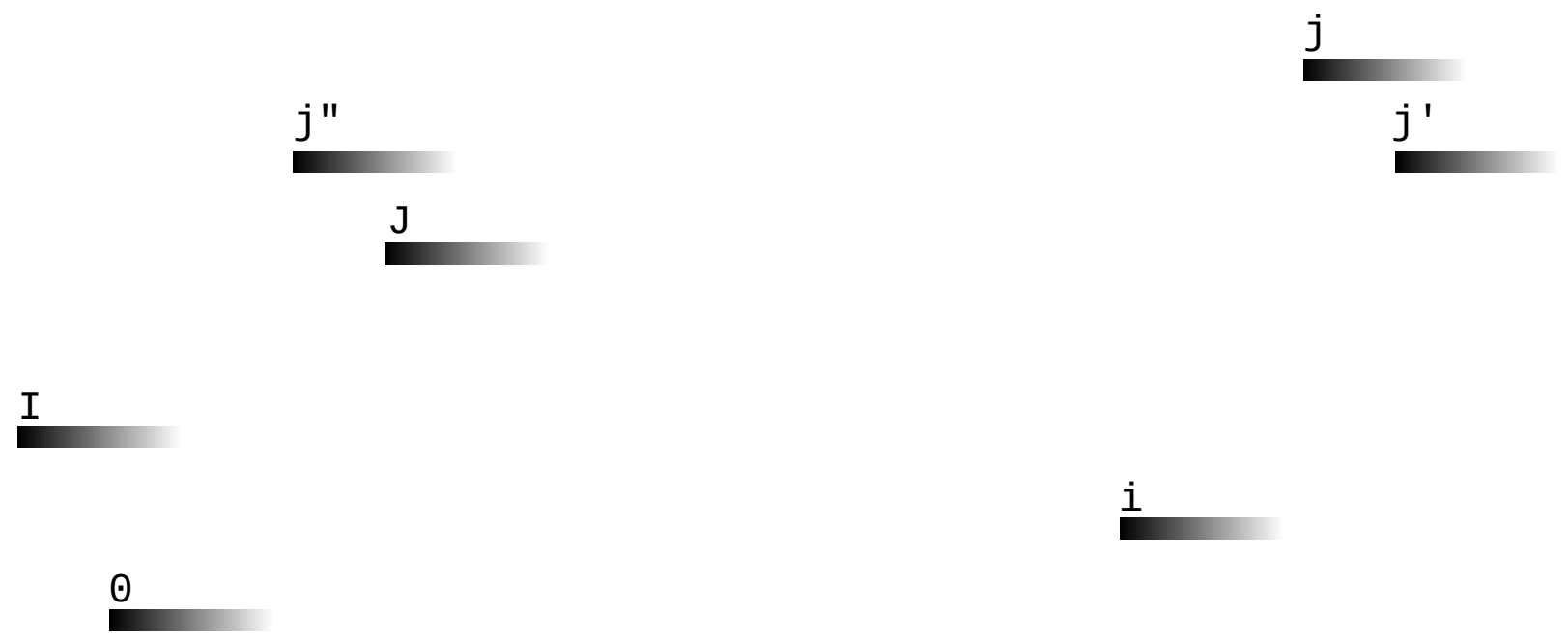

Figure 504

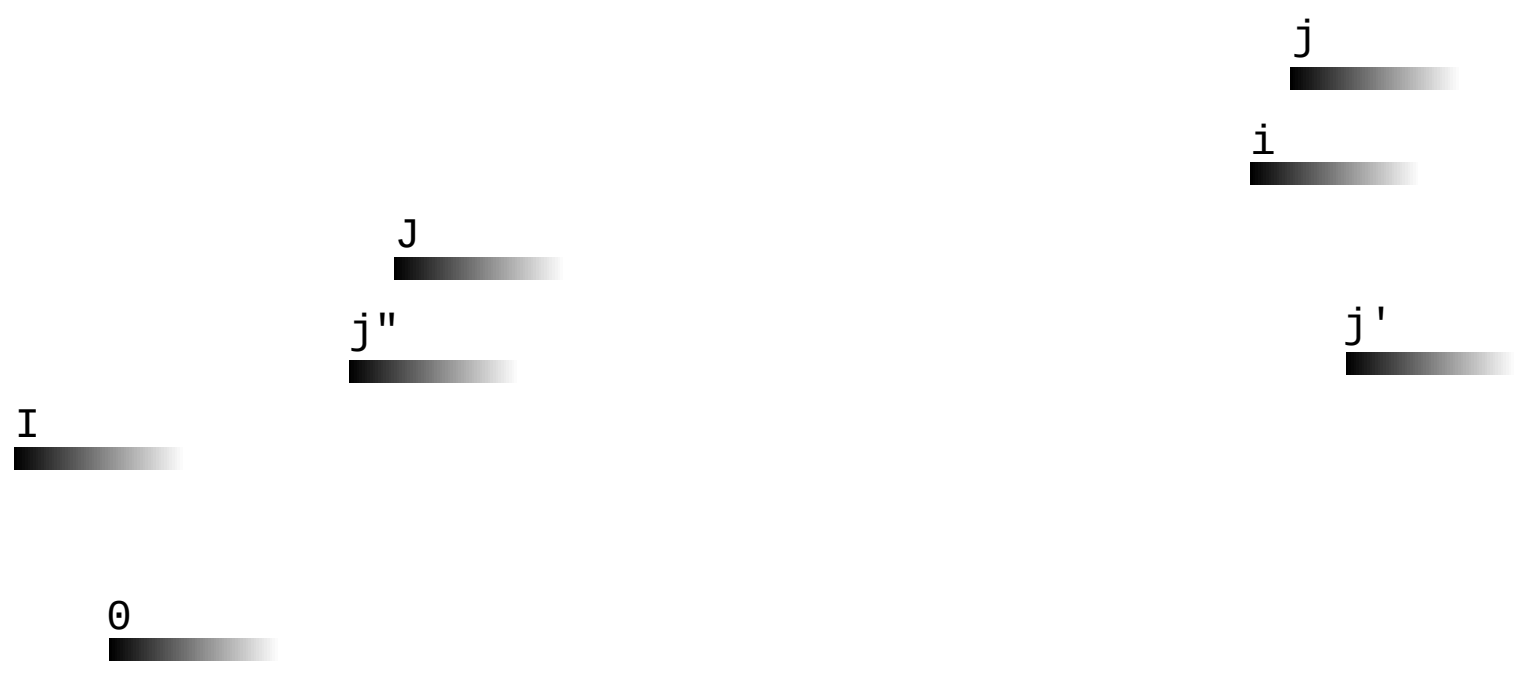

Figure 505

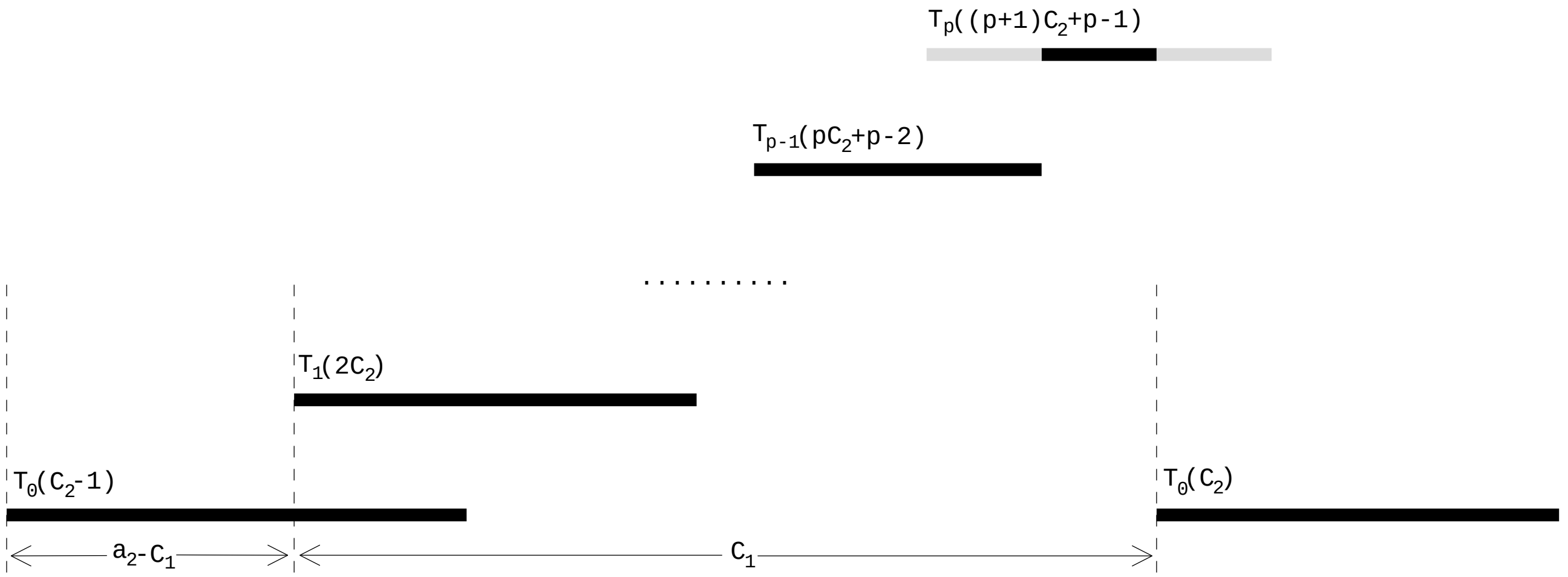

Figure 506

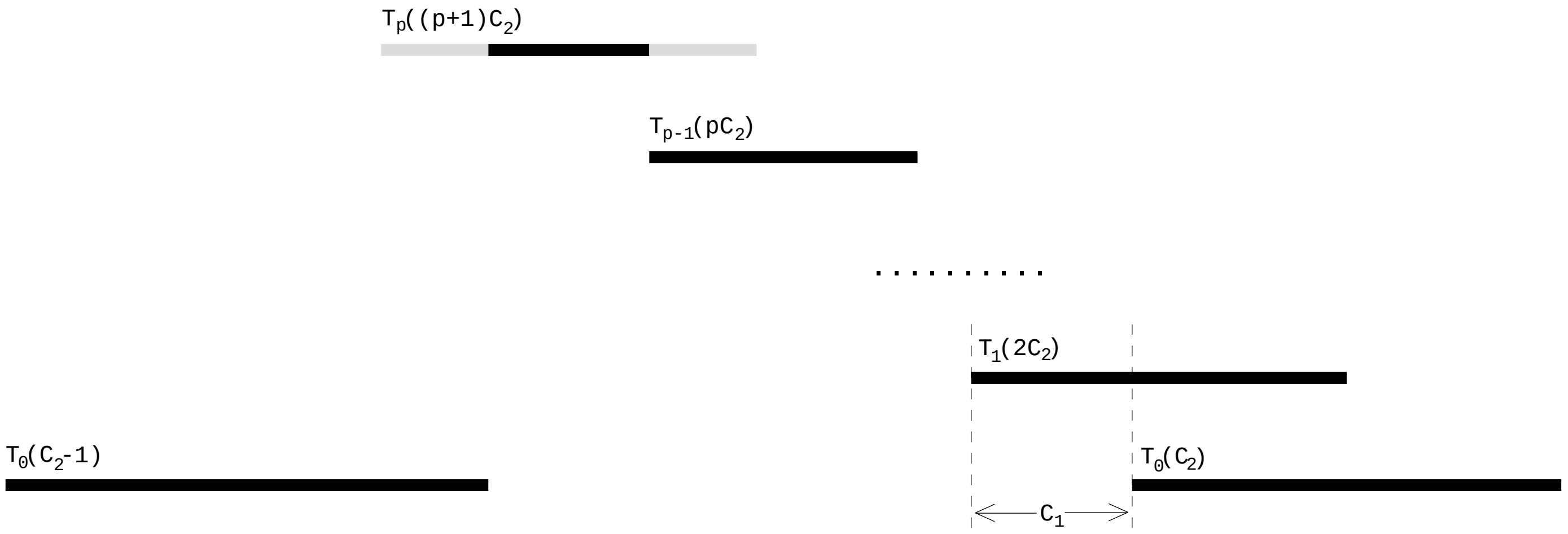

Figure 507

## 5.3 Is OSG(n,1) always best in the long run?

Theorem 507 proves that for sufficiently large n, OSG(n,1) is superior to any SG(n,p) with key = 1 or p; to be precise:

Given p, we can find N such that:

$L_1 > L_p$ for all $n >= N$

where $L_1$ = length of OSG(n,1)
and $L_p$ = maximum length of SG(n,p) with key = 1 or p

Is Theorem 507 true for all SG(n,p) - regardless of key value?

The available evidence strongly suggests that this is the case, but I have not yet succeeded in finding a proof.

This section presents some of the evidence, and some of the lines of investigation that I have followed. None of the results are necessary for the overall proof.

We start by considering the case p=11, key=5. We follow part (a) of the proof of Theorem 506 (the break is just before $T_{11}$), and show that in this case the conjecture holds. At the same time, we try to present the steps of the proof in a general way to show how the Theorem might be extended.

Theorem 511

Let $L_{11}$ be the maximum length of any stride generator SG(n,11) with key=5 whose first break occurs at $ST_{11}$-1.

Then $L_1 > L_{11}$ for sufficiently large n.

Proof

By Theorem 301 we know that:

$L_1 \rightarrow n^2/3$ as $n \rightarrow$ infinity

We prove the result by showing that:

$L_{11} \rightarrow 36n^2/109 < n^2/3$ as $n \rightarrow$ infinity

Step 1: Identify the threads surrounding $T_p$

We need to know not only the orders of the threads, but also their start positions and lengths.

In the case of Theorem 506 this was easy, since the threads formed a staircase such that $ST_i = ST_0 + iG$, where G was the "up gap"; and in its turn the "up gap" was easily shown to be $a_2$-$C_1$ because the thread immediately following $T_0$ was $T_1$.

In this case the pattern is not so obvious (see Figure 508) and we have to derive the information by enumeration.

We know that the sequence of threads from $T_0(C_2\text{-}1)$ to $T_0(C_2)$ is:

$T_0(C_2\text{-}1), T_5, T_{10}, T_3, T_8, T_1, T_6, T_{11}, T_4, T_9, T_2, T_7, T_0(C_2)$

and that:

$ST_i((i+1)C_2+m) - C_1 = ST_{i+1}((i+2)C_2+m)$

So $ST_0(C_2) - C_1 = ST_1(2C_2)$; so $T_1$ must be $T_1(2C_2)$.

Then:

$ST_1(2C_2) - C_1 = ST_2(3C_2)$

But there is no 2-thread between $ST_0(C_2\text{-}1)$ and $ST_1(2C_2)$ and so $ST_2(3C_2) < ST_0(C_2\text{-}1)$. This means that the 2-thread that follows $T_1$ in the stride generator must be the one following $T_2(3C_2)$, which is $T_2(3C_2+1)$.

Then:

$ST_2(3C_2+1) - C_1 = ST_3(4C_2+1)$

There is a 3-thread between $ST_0(C_2\text{-}1)$ and $ST_2(3C_2+1)$ - so it must be $T_3(4C_2+1)$.

Then:

$ST_3(4C_2+1) - C_1 = ST_4(5C_2+1)$

There is no 4-thread between $ST_0(C_2\text{-}1)$ and $ST_3(4C_2+1)$ - so the 4-thread we are looking for must be $T_4(5C_2+2)$.

We repeat this process until we reach $T_{11}$:

$T_0(C_2\text{-}1)$
$T_1(2C_2)$
$T_2(3C_2+1)$
$T_3(4C_2+1)$
$T_4(5C_2+2)$
$T_5(6C_2+2)$
$T_6(7C_2+3)$
$T_7(8C_2+4)$
$T_8(9C_2+4)$
$T_9(10C_2+5)$
$T_{10}(11C_2+5)$
$T_{11}(12C_2+6)$

Because key=5, we know that $T_{11}$ is surrounded by $T_6$ and $T_4$, and we have:

$ST_6(7C_2+3) < ST_{11}(12C_2+6) < ST_4(5C_2+2)$

Step 2: Obtain formula for n in terms of $a_2$, $C_2$, $C_1$

We use the fact that $T_p$ and its predecessor are contiguous (because the break lies just before $T_p$):

$ST_{11} = ET_6+1$

Now:

$ST_{11} = (12C_2+6)a_2 - 11(C_2a_2+C_1)$

$= C_2a_2+6a_2-11C_1$

$$ET_6 = ST_6(7C_2+3) + n+6-(7C_2+3)$$

$$= (7C_2+3)a_2 - 6(C_2a_2+C_1) + n-7C_2+3$$

$$= C_2a_2+3a_2-6C_1+n-7C_2+3$$

So:

$$ST_{11} = ET_6+1$$

$$\Rightarrow 3a_2-5C_1 = n-7C_2+4$$

$$\Rightarrow n = 3a_2+7C_2-5C_1-4 \qquad - (1)$$

Step 3: Establish an upper bound on $C_1$ in terms of $a_2$, $C_2$

We use the fact that $T_p$ must be long enough to cover the gap between its predecessor and successor:

$$ET_{11} >= ST_4-1$$

Now:

$$ET_{11} = ST_{11} + n+11-(12C_2+6)$$

$$= C_2a_2+6a_2-12C_2-11C_1+n+5$$

$$ST_4 = (5C_2+2)a_2 - 4(C_2a_2+C_1)$$

$$= C_2a_2+2a_2-4C_1$$

So:

$$ET_{11} >= ST_4-1$$

$$\Rightarrow 4a_2-12C_2-7C_1+n+6 >= 0$$

$$\Rightarrow 4a_2-12C_2-7C_1+3a_2+7C_2-5C_1-4+6 >= 0 \quad - \text{ using (1)}$$

$$\Rightarrow 7a_2-5C_2-12C_1+2 >= 0$$

$$\Rightarrow C_1 <= (7a_2-5C_2+2)/12 \qquad - (2)$$

Step 4: Establish an upper bound on $C_2$ in terms of $a_2$, n

We use (1) to determine $C_2$ in terms of n, $a_2$, $C_1$ and then substitute worst case $C_1$ from (2):

$$(1) \Rightarrow 7C_2 = n-3a_2+5C_1+4$$

$$\Rightarrow 84C_2 = 12n-36a_2+60C_1+48$$

$$<= 12n-36a_2+5(7a_2-5C_2+2)+48 \quad - \text{ using (2)}$$

$$\Rightarrow 109C_2 <= 12n-a_2+58$$

$$\Rightarrow C_2 <= (12n-a_2+58)/109$$

Step 5: Maximise $C_2a_2$ and let n -> infinity

$$C_2a_2 <= ((12n+58)a_2-a_2^2)/109$$

This is a maximum when $a_2 = (12n+58)/2$

$$\Rightarrow C_2a_2 \le (12n+58)^2/436$$

$$= (6n+29)^2/109$$

Now $a_3 < (C_2+1)a_2$, so:

$$L_{11} = (6n+29)^2/109 + (6n+29)$$

And as n -> infinity:

$$L_{11} \to 36n^2/109 < n^2/3 \quad \text{as required}$$

It is interesting to see how close to $n^2/3$ this limit is: it is much closer than the limits for key=1 and key=11:

$$(B+C)(B+3C)/4C \quad \text{where } B = n(p+1)+1 \text{ and } C = p^2+p+1$$

$$\to (n^2(p+1)^2)/4(p^2+p+1) \quad \text{as } n \to \text{infinity}$$

$$= 144n^2/4.133$$

$$= 36n^2/133$$

We observe that for key=5 the p-thread is close to the middle of the range $ST_0(C_2-1)$ to $ST_0(C_2)$, whereas for key=1 or key=p the p-thread is at either end. We suggest that the limit value approaches $n^2/3$ as the p-thread approaches the middle of the range; this is supported by the fact that for p=1 the 1-thread is, by definition, in the middle.

Furthermore, it is only with p=1 that the p-thread can be exactly in the middle, because:

a) It is only possible with odd values of p

b) It requires key=(p+1)/2 - and this is not possible for p>1 because ((p+1)/2,p+1) are not then co-prime as required by Theorem 505.

The difficulty with generalising Theorem 511 arises in Step 1: the key defines the "order gap" between successive threads, and from this we have to find their start positions. What we need is a way of determining the "start gap" from the "order gap". If we had such a formula, we might be able to use it directly to derive the start positions of the two threads surrounding $T_p$.

The only "hint" in this direction is the observation that if the p-thread is the (m+1)th thread in the sequence, then it is $T_p((p+1)C_2+m)$.

eg Consider SG(n,11) with key=5.

$T_{11}$ is the 7th thread in the sequence (counting $T_0(C_2-1)$ as the 0th) and it is $T_{11}(12C_2+6)$.

This means that p = (m+1)j mod (p+1) where j is the key, and so might give a clue to determining m from p and j. However, this does not immediately suggest a way of deriving information about the threads on either side of $T_p$.

There is, however, an alternative approach to the categorisation of stride generators of order p which might help.

Let SG(n,p) be an order p stride generator, and suppose that the p-thread that lies between $ST_0(C_2-1)$ and $ST_0(C_2)$ is $T_p((p+1)C_2+m)$.

Then:

$$ST_p = ((p+1)C_2+m)a_2 - p(C_2a_2+C_1)$$

$$= C_2a_2+ma_2-pC_1$$

So $ST_p > (C_2-1)a_2 \;\Rightarrow\; ma_2-pC_1+a_2 > 0$

$$\Rightarrow\; m > (pC_1/a_2) - 1$$

and $ST_p < C_2a_2 \;\Rightarrow\; ma_2-pC_1 < 0$

$$\Rightarrow m < (pC_1/a_2)$$

So:

$$(pC_1/a_2) - 1 < m < pC_1/a_2$$

$$\Rightarrow\; m = \text{intpt}(pC_1/a_2)$$

There are therefore p possible values for m: 0,1, ... p-1, according as:

$$ma_2 < pC_1 < (m+1)a_2$$

So the number of different values of m is exactly the same as the number of different possible key values - which suggests there may be a correspondence between the two. This correspondence is not obvious, as the following tables show (but there is a suggestion that maybe the only valid forms of m are those such that (m+1,p+1) are co-prime):

| | key | m |
|---|---|---|
| p=11 | 1 | 10 |
| | 5 | 6 |
| | 7 | 4 |
| | 11 | 0 |
| p=12 | 1 | 11 |
| | 2 | 5 |
| | 3 | 3 |
| | 4 | 2 |
| | 5 | 4 |
| | 6 | 1 |
| | 7 | 10 |
| | 8 | 7 |
| | 9 | 9 |
| | 10 | 8 |
| | 11 | 6 |
| | 12 | 0 |

The idea is that it may be easier to prove limit theorems about the length of stride generators SG(n,p) where:

$$ma_2 < pC_1 < (m+1)a_2 \quad \text{for} \quad 0<=m<=p-1$$

instead of limit theorems about SG(n,p) where:

$$\text{key}=j \quad \text{for} \quad 1<=j<=p$$

One final note about the size of the "start gaps".

Consider the "up" gap between $y=ST_j(c_2)$ and $x=ST_i(b_2)$:

$$y-x = (c_2-b_2)a_2 - (j-i)a_3$$

Suppose without loss of generality that $i=0$, $x=0$ and $b_2=0$:

$$y = c_2a_2-ja_3$$

$$\Rightarrow y = (c_2-jC_2)a_2 - jC_1$$

So to minimise y we must minimise:

$$(-jC_1) \bmod a_2 \quad \text{over } 0<j<=p$$

Similarly the "down" gap is found by minimising:

$$(jC_1) \bmod a_2 \quad \text{over } 0<j<=p$$

[Note that $x \bmod a_2 = r \Rightarrow -x \bmod a_2 = a_2-r$]

SG(20,11) = {1,80,206}

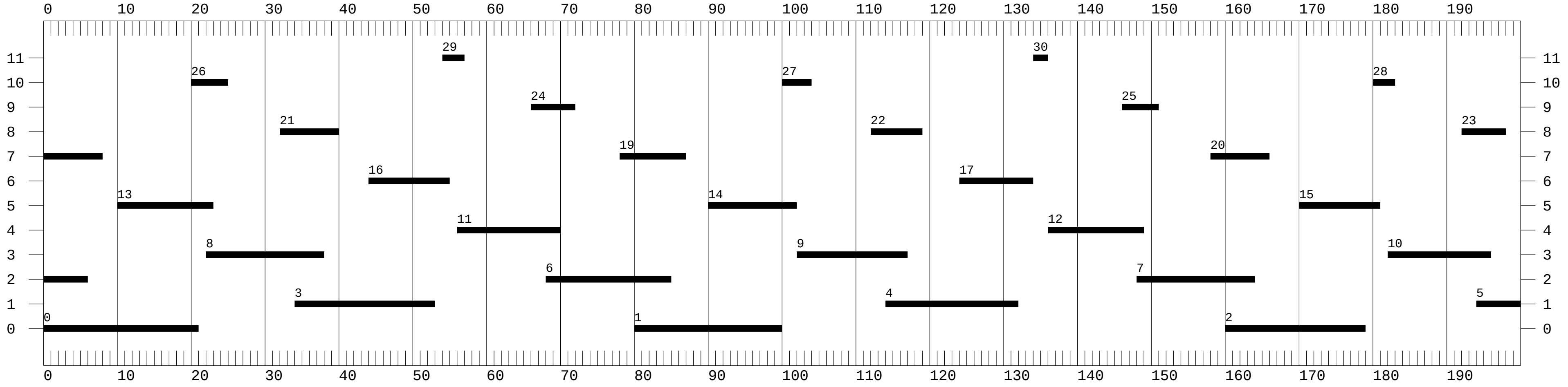

Figure 508

## 5.4 Order = 1 for maximal cover for s>=81

Theorem 512

Let A = {1, $a_2$, $a_3$} be a stride generator SG(n,p) where:

p = 0, 2 or 3
and n < 52

whose length is greater than or equal to one less than the length of OSG(n,1); that is:

$a_3$ >= $A_3$-1 where OSG(n,1) = {1, $A_2$, $A_3$}

Then if A is the stride generator underlying a cover C(A,3,s) for s>=81:

C(A,3,s) < $X_{opt}$

where $X_{opt}$ is given in Table 103.

[ This theorem finally shows that if s>=81, only order 1 stride generators matter:

Theorem 500 shows that any stride generator A = SG(n,p) that underlies a cover C(A,3,s) >= $X_{opt}$ must have p<=3 if s>=40.

Theorem 510 shows that there are no stride generators A = SG(n,p) with $a_3$>=$A_3$-1 and p=0,2,3 for n>=52.

So this theorem shows that even those stride generators whose length is at least as great as SG1(n,1) cannot improve on the cover derived from OSG(s-$k_{opt}$,1) for s>=81.

Figure 509 may also help to understand what this theorem is about.

The solid line indicates:

X = (k+1)$a_3$ + y-1 for OSG(s-k,1) = {1, $a_2$, $a_3$} as k varies.

The maximum, $X_{opt}$, occurs at approximately k=s/3.

We know that there are no stride generators of order <= 3 that can improve upon X in the range 0<=k<s-52, but such stride generators might exist in the range s-52<=k<s - because we know there are SG(n,p) for n<52 whose length >= $A_3$-1 where $A_3$ is the length of OSG(n,1).

This theorem shows that although such stride generators may improve on OSG(n,1) for values of n in the range 0<=n<52 - as indicated by the dotted line in Figure 509 - they can never improve on $X_{opt}$ provided s>=81. ]

Proof

By Theorem 231 we know that:

$$C(A,3,s) = (s-n+1)a_3 + y-1$$
$$< (s-n+2)a_3$$

We now show that for s>=81, (s-n+2)$a_3$ < $X_{opt}$, and so C(A,3,s)<$X_{opt}$ as required.

Table 504 below tabulates for s=81 values of:

$$X' = (s-n+2)a_3 \quad \text{for } 1<=n<52$$

where $a_3$ is the length of the longer of OSG(n,2) and OSG(n,3) taken from Table 502; note that Theorem 509 shows that we do not need to consider any stride generators of order 0.

| s-n | n | $a_3$ | X' |
|---|---|---|---|
| 30 | 51 | 953 | 30496 |
| 31 | 50 | 918 | 30294 |
| 32 | 49 | 883 | 30022 |
| 33 | 48 | 850 | 29750 |
| 34 | 47 | 818 | 29448 |
| 35 | 46 | 786 | 29082 |
| 36 | 45 | 754 | 28652 |
| 37 | 44 | 722 | 28158 |
| 38 | 43 | 693\ | |
| ...... | | ) | < 43*693 = 29799 |
| | | ) | |
| 41 | 40 | 606/ | |
| 42 | 39 | 579\ | |
| ...... | | ) | < 52*579 = 30108 |
| | | ) | |
| 50 | 31 | 384/ | |
| 51 | 30 | 361\ | |
| ...... | | ) | < 82*361 = 29602 |
| | | ) | |
| 80 | 1 | 6/ | |

Table 504

Now $X_{opt}$ = 30816 for s=81, and so the table shows that $X'<X_{opt}$ for s=81.

We now have only to show that what holds for s=81 holds also for s>81:

a) Consider how X' increases for a particular value of n as s increases by 1:

$$\text{delta}(X') = ((s-n+1)+2)a_3 - (s-n+2)a_3 = a_3$$

So for any 1 <= n < 52, delta(X') <= 953 - from the table above.

b) From Table 103, we see that $X_{opt}$ increases by at least $12t^2+14t+4$ as s increases by 1 (where s = 9t+r, 0 <= r < 9):

$$\text{delta}(X_{opt}) >= 1102$$

This shows that the rate of increase of X' for any 1 <= n < 52 is less than the rate of increase of $X_{opt}$, and so the result holds for all s>=81.

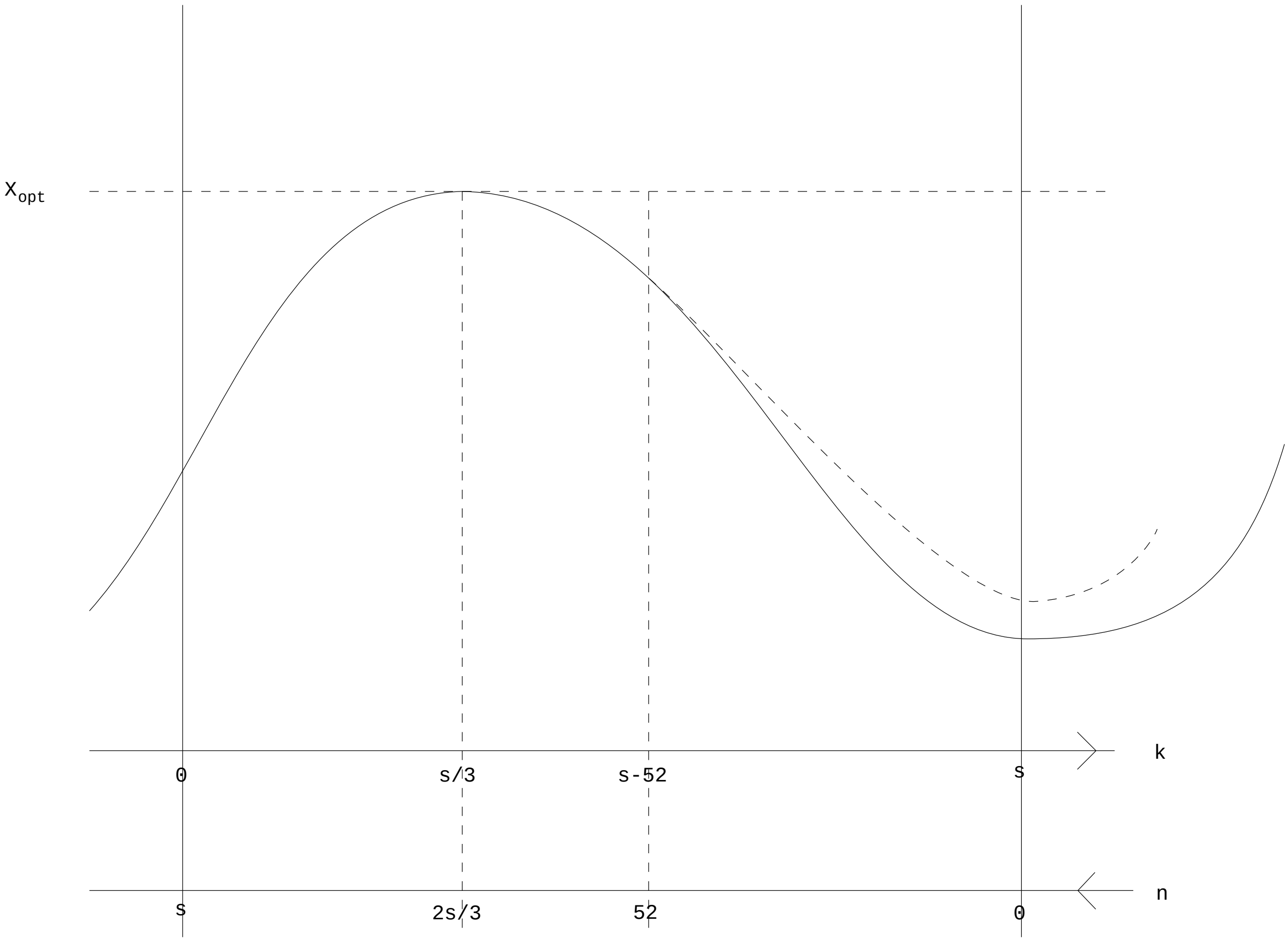

Figure 509

# SECTION 6

## SUPPORTING PROGRAMS

This Section gives further details of the various programs - both numerical and algebraic - that have been referenced within the proof.

## 6.1 Stride generator programs

Section 5.2 refers to two programs used to list stride generators. Both of these are written in BASIC for the Acorn Archimedes range of computers.

The first of these, exp14, was used to list the longest stride generators of a given order with key = 1 or p; the results confirm Theorem 506. A listing of the program follows:

```
 10REM  Lists actual best SG(n,p) for given p and a range of n with key
 20REM   value of 1 or p.
 30
 40INPUT "p = " p%
 50INPUT "n1 = " n1%
 60INPUT "n2 = " n2%
 70
 80PRINT "                 Key = 1               Key = p          Theoretical"
 90PRINT "              A          B          A          B"'
100
110FOR n% = n1% TO n2%
120
130PRINT "N = ";
140@%=3: PRINT n%;
150
160
170REM  Type (1a)
180
190a3max%=0
200a2opt%=0: a2opt2%=0
210
220FOR a2% = 2 TO n%*(p%+1)+1
230
240  C2max% = (n%+2)/(p%+1)
250  C2min% = n%-a2%+3: IF C2min%<1 C2min%=1
260
270  C2%=C2min%
280  WHILE C2%<=C2max%
290
300    C1% = a2%+p%*C2%-2-n%
310
320    IF (C1%>((p%-1)*a2%)/p%) AND (C1%<=(p%*a2%-C2%)/(p%+1)) THEN
                                              PROCprocess(a2%, C2%*a2%+C1%)
330
340    C2%+=1
350  ENDWHILE
360
370NEXT
380
390PRINT " ";: @%=5: PRINT a2opt%;: PRINT a3max%;
400IF a2opt2%<>0 THEN
410  PRINT '"        ";
420  @%=5: PRINT a2opt2%;: PRINT a3max%;
430ENDIF
440
450
```

```
 460REM  Type (1b)
 470
 480a3max%=0
 490a2opt%=0: a2opt2%=0
 500
 510FOR a2% = 2 TO n%*(p%+1)+1
 520
 530  C2max% = (n%+2)/(p%+1)
 540  C2min% = n%-a2%+3: IF C2min%<1 THEN C2min%=1
 550
 560  C2%=C2min%
 570  WHILE C2%<=C2max%
 580
 590    pC1% = (p%-1)*a2%+n%+2-(p%+1)*C2%
 600
 610    IF pC1% MOD p% = 0 THEN
 620
 630      C1% = pC1%/p%
 640      IF (C1%>=(p%*a2%-C2%)/(p%+1)) AND (C1%<=p%*a2%/(p%+1)) THEN
                                                 PROCprocess(a2%, C2%*a2%+C1%)
 650
 660    ENDIF
 670
 680    C2%+=1
 690  ENDWHILE
 700
 710NEXT
 720
 730PRINT " ";: @%=5: PRINT a2opt%;: PRINT a3max%;
 740IF a2opt2%<>0 THEN
 750  PRINT '"                    ";
 760  @%=5: PRINT a2opt2%;: PRINT a3max%;
 770ENDIF
 780
 790
 800REM  Type (2a)
 810
 820a3max%=0
 830a2opt%=0: a2opt2%=0
 840
 850FOR a2% = 2 TO n%*(p%+1)+1
 860
 870  C2max% = (n%+p%+1)/(p%+1)
 880  C2min% = n%-a2%+3: IF C2min%<1 THEN C2min%=1
 890
 900  C2%=C2min%
 910  WHILE C2%<=C2max%
 920
 930    pC1% = (a2%+C2%-n%-2)
 940
 950    IF pC1% MOD p% = 0 THEN
 960
 970      C1% = pC1%/p%
 980      IF (C1%>=((a2%+p%)/(p%+1)-C2%)) AND (C1%<=(a2%-p%*C2%+p%-1)/(p%+1)) THEN
                                                 PROCprocess(a2%, C2%*a2%+C1%)
 990
1000    ENDIF
1010
1020    C2%+=1
1030  ENDWHILE
1040
1050NEXT
1060
```

```
1070PRINT " ";: @%=5: PRINT a2opt%;: PRINT a3max%;
1080IF a2opt2%<>0 THEN
1090  PRINT '"                              ";
1100  @%=5: PRINT a2opt2%;: PRINT a3max%;
1110ENDIF
1120
1130
1140REM  Type (2b)
1150
1160a3max%=0
1170a2opt%=0: a2opt2%=0
1180
1190FOR a2% = 2 TO n%*(p%+1)+1
1200
1210  C2max% = (n%+p%+1)/(p%+1)
1220  C2min% = n%-a2%+3: IF C2min%<1 C2min%=1
1230
1240  C2%=C2min%
1250  WHILE C2%<=C2max%
1260
1270    C1% = n%+p%+1-(p%+1)*C2%
1280
1290    IF (C1%>=((a2%-p%*C2%+p%-1)/(p%+1))) AND (C1%<((a2%+C2%-n%-2)/(p%-1)))
                                          THEN PROCprocess(a2%, C2%*a2%+C1%)
1300
1310    C2%+=1
1320  ENDWHILE
1330
1340NEXT
1350
1360PRINT " ";: @%=5: PRINT a2opt%;: PRINT a3max%;
1370IF a2opt2%<>0 THEN
1380  PRINT '"                                         ";
1390  @%=5: PRINT a2opt2%;: PRINT a3max%;
1400ENDIF
1410
1420
1430a2theor = ((p%+1)*n%+2*p%+2)/2
1440alpha = n%*(p%+1)+1
1450beta = p%*p%+p%+1
1460a3theor = (alpha+beta)*(alpha+3*beta)/(4*beta)
1470PRINT " ";: @%=&020208: PRINT a2theor;: PRINT a3theor
1480
1490NEXT
1500
1510END
1520
1530
1540
1550
1560DEF PROCprocess(a2%, a3%)
1570REM  @%=3: PRINT n%, ":", a2%, a3%
1580  IF a3%=a3max% THEN
1590    IF a2opt2%<>0 THEN PRINT "*** too many best ones ***" ELSE a2opt2%=a2%
1600  ELSE
1610    IF a3%>a3max% THEN a2opt%=a2%: a3max%=a3%: a2opt2%=0
1620  ENDIF
1630ENDPROC
1640
1650
```

An independent program, exp15, was used to generate optimal stride generators in the proof of Theorem 510. This program works from "first principles": the main functions are FNtry and FNcan which together determine whether $\{1, a_2, a_3\}$ is a stride generator SG(n,i) for some i by first determining the maximum order necessary to generate all values $0<x<a_3$, and then checking that there is also at least one break value in the same range.

The range of values of $a_3$ - and hence $a_2$ - that need to be considered are determined by applying the limit theorems of section 2.12. FNmaxa3 applies Theorem 241, FNa2min applies Theorem 239, and FNa2max uses Theorem 240.

The results confirm those of exp14 for p=2 and 3.

Note that it is easy to modify exp15 so that all stride generators SG(n,p) for given n and p are listed.

```
   10REM Lists all optimal stride generators of order between ord1% and ord2% for
   20REM  values of n between n1% and n2%
   30
   40INPUT "n1 = " n1%
   50INPUT "n2 = " n2%
   60INPUT "ord1 = " ord1%
   70INPUT "ord2 = " ord2%
   80@%=5
   90PRINT "Optimal stride generators"'
  100DIM gota%(ord2%), got%(ord2%)
  110
  120FOR n% = n1% TO n2%
  130  a3% = FNmaxa3(n%)
  140  PRINT "N =", n%;
  150  FOR i% = ord1% TO ord2%: got%(i%)=FALSE: NEXT
  160  REPEAT
  170    FOR i% = ord1% TO ord2%: gota%(i%)=FALSE: NEXT
  180    a2min% = FNa2min(n%,a3%): a2max% = FNa2max(n%,a3%)
  190    FOR a2% = a2max% TO a2min% STEP -1
  200      maxorder% = FNtry
  210      IF maxorder%>0 THEN
  220        IF NOT got%(maxorder%) THEN PROCout(maxorder%): gota%(maxorder%)=TRUE
  230      ENDIF
  240    NEXT
  250    alldone%=TRUE
  260    FOR i% = ord1% TO ord2%
  270      IF gota%(i%) THEN got%(i%)=TRUE
  280      IF NOT got%(i%) alldone%=FALSE
  290    NEXT
  300    a3%-=1
  310  UNTIL alldone%
  320  PRINT
  330NEXT
  340
  350END
  360
  370
  380
  390
  400DEF PROCout(i%)
  410  PRINT a2%, a3%, "(";i%;")  ";
  420ENDPROC
  430
  440
  450
  460
```

```
 470REM On entry, a2%, a3% and n% are set.
 480REM  If {1,a2%,a3%} is an SG(n%,i%) for i%<=ord2%, then the result is i%; if
 490REM  not, the result is -1
 500
 510
 520DEF FNtry
 530LOCAL order%, x%, p%, q%, pp%, sgorder%
 540  order%=0
 550  p%=a3%/a2%: q%=a3% MOD a2%
 560  pp%=p%: x%=q%-1
 570  WHILE x%>=0
 580    IF NOT FNcan THEN =-1
 590  ENDWHILE
 600  pp%-=1: x%=a2%-1
 610  WHILE x%>=q%
 620    IF NOT FNcan THEN =-1
 630  ENDWHILE
 640  sgorder%=order%: n%-=1: ord2%+=1:  REM check there is a break!
 650  pp%=p%: x%=q%-1
 660  WHILE x%>=0
 670    IF NOT FNcan THEN n%+=1: ord2%-=1: =sgorder%
 680  ENDWHILE
 690  pp%-=1: x%=a2%-1
 700  WHILE x%>=q%
 710    IF NOT FNcan THEN n%+=1: ord2%-=1: =sgorder%
 720  ENDWHILE
 730  n%+=1: ord2%-=1
 740=-1
 750
 760
 770
 780
 790REM On entry:
 800REM    pp% = p%-1 or p%
 810REM    x% < a2%
 820REM    order% >= 0
 830REM
 840REM This function determines whether pp%*a2%+x% can be generated with a
 850REM  generation of order <= ord2%; the result is TRUE or FALSE accordingly.
 860REM  If TRUE, order% is set to the greater of order% and the order of
 870REM  the generation, and x% is set to indicate the next gap below.
 880
 890
 900DEF FNcan
 910  IF (pp%+x%)<=n% THEN x%=-1: =TRUE
 920  LOCAL na2%,nu%,i%
 930  i%=1: na2%=pp%: nu%=x%
 940  REPEAT
 950    na2%+=p%: nu%+=q%
 960    IF nu%>=a2% na2%+=1: nu%-=a2%
 970    IF (na2%+nu%)<=(n%+i%) THEN
 980      x%-=(nu%+1)
 990      IF order%<i% order%=i%
1000      UNTIL TRUE: =TRUE
1010    ENDIF
1020    i%+=1
1030  UNTIL i%>ord2%
1040=FALSE
1050
1060
1070
1080
```

```
1090DEF FNmaxa3(n%)
1100=n%*(n%+ord2%+1)+(n%+ord2%+1)/(ord2%+1)
1110
1120
1130
1140
1150DEF FNa2min(n%,a3%)
1160=((ord1%+1)*a3%)/(n%+ord1%+1)
1170
1180
1190
1200
1210DEF FNa2max(n%,a3%)
1220LOCAL t%
1230  t%=n%*(ord2%+1)+1
1240  IF t%>=a3% t%=a3%-1
1250=t%
1260
1270
1280
1290
```

Finally, mention should be made of the interactive program exp3 which was developed to investigate the behaviour of threads (see Section 2.4); routines were later added to generate thread diagrams in the form of "drawfiles" which can be displayed on the screen or printed out, and many of these appear as figures in Section 2. No listing is included here, since the program is quite long and mostly concerned with screen and graphics management which are of no relevance to this proof.

## 6.2 Algebra associated with Section 5.1

We list here the two REDUCE programs used in the proof of Theorem 500 in Section 5.1; for details of their operation, see the proof itself.

The program (xopt) used in step (a) of the proof of Theorem 500 is as follows:

```
off acorn;
off allfac;
off mcd;

operator X, Xopt;

for all s
  let X(s) = (4/81)*s*s*s + (2/3)*s*s + (22/9)*s;

for all a,b,c,d
  let Xopt(a,b,c,d) = a*t*t*t + b*t*t + c*t + d;

Xopt(36, 54, 22, 0) - X(9*t);
Xopt(36, 66, 36, 4) - X(9*t+1);
Xopt(36, 78, 53, 8) - X(9*t+2);
Xopt(36, 90, 71,15) - X(9*t+3);
Xopt(36,102, 92,22) - X(9*t+4);
Xopt(36,114,116,36) - X(9*t+5);
Xopt(36,126,143,49) - X(9*t+6);
Xopt(36,138,173,68) - X(9*t+7);
Xopt(36,150,204,86) - X(9*t+8);

;END;
```

The program (minord) used in step (f) of the proof of Theorem 500 is as follows:

```
off acorn;

operator p;
for all a3 let p(a3) = sub(
  c = s+3-(a3/s),
  n = sub(  x = (4/81)*s*s*s + (2/3)*s*s + (22/9)*s,
                (s+2)-x/a3
         ),
  ((n+1)*c-a3)/(a3-c)
);

p(a3);

xx := sub( a3=s*(s+3)/(s+1), num(p(a3)) );
xx := xx*(s+1)*(s+1)/(s*s*s);
sub(s=10, xx);
sub(s=11, xx);

dfp := df(p(a3),a3);

sol := solve(dfp,a3);
a31 := part(first(sol), 2);
a32 := part(second(sol),2);
```

```
aa3 := s*(s+3)/(s+1);
diff1 := a31 - aa3;
diff2 := a32 - aa3;

xx := part(diff2,1,2,1,1);
xx := xx*xx;
sub(s=10, xx);
sub(s=11, xx);
sub(s=12, xx);

operator pp;
for all s let pp(s) = sub(a3=a32,p(a3));

pp(s);
structr(pp(s));

on bigfloat;
on numval;

for i := 40:58 do
<<
  write "pp(", i, ") = ", sub(s=i, pp(s));
>>;

sub(s=40, aa3);
sub(s=40, a31);
sub(s=40, a32);
sub(a3=sub(s=40,a31), sub(s=40,p(a3)));
sub(a3=sub(s=40,a32), sub(s=40,p(a3)));

off bigfloat;
off numval;

structr(pp(s));
lim := (5929*4*s**9 - 5184*s**9)/(81*(8*4*s**9 + 32*s**9));

on bigfloat;
on numval;

lim;

;end;
```

## 6.3 Algebra associated with Section 4

This section gives details of the REDUCE programs used to verify the algebra of Section 4.

Theorem 400 in Section 4.1 gives formulae for $C_\Theta'(k)$ and for the two roots of $dC_\Theta'(k)/dk$. These are verified by the following program (c0d):

```
off acorn;

operator n, c0d;

for all k let n(k) = s-k;
for all k let c0d(k) = (k+1)*(((n(k)+3)*(n(k)+2)+1)/3) + (n(k)*(n(k)+3)-4)/3;

c0d(k);

dfc0d := df(c0d(k),k);

k12 := solve(dfc0d, k);

;end;
```

with annotated output as follows:

[ Define $C_\Theta'(k)$ as c0d(k): ]

```
C: for all k let n(k) = s-k;
C: for all k let c0d(k) = (k+1)*(((n(k)+3)*(n(k)+2)+1)/3)+(n(k)*(n(k)+3)-4)/3;
```

[ Evaluate $C_\Theta'(k)$: ]

```
C: c0d(k);

  3      2        2      2                 2
 k  - 2*k *s - 3*k  + k*s  + k*s - k + 2*s  + 8*s + 3
------------------------------------------------------
                         3
```

[ Evaluate its derivative with respect to k: ]

```
C: dfc0d := df(c0d(k),k);

                 2                   2
              3*k  - 4*k*s - 6*k + s  + s - 1
dfc0d := ---------------------------------
                         3
```

[ Solve for k: ]

```
C: k12 := solve(dfc0d, k);

                          2
                    sqrt(s  + 9*s + 12) - 2*s - 3
k12 := {k= - ------------------------------,
                              3

                       2
                 sqrt(s  + 9*s + 12) + 2*s + 3
         k=------------------------------}
                           3
```

A similar program (m2d) confirms the expansion of $M_2$'(k) and the roots of its derivative given in the proof of Theorem 402 in Section 4.2; the output is as follows:

```
C: for all k let n(k) = s-k;
C: for all k let m2d(k) = (k+2)*(((n(k)+3)*(n(k)+2)+1)/3 - 2);

C: m2d(k);

  3      2        2      2                  2
 k  - 2*k *s - 3*k  + k*s  + k*s - 9*k + 2*s  + 10*s + 2
---------------------------------------------------------
                            3

C: dfm2d := df(m2d(k),k);

              2                  2
           3*k  - 4*k*k*s - 6*k + s  + s - 9
dfm2d := ---------------------------------
                         3

C: k12 := solve(dfm2d, k);

                       2
                 sqrt(s  + 9*s + 36) - 2*s - 3
k12 := {k= - ------------------------------,
                            3

                    2
              sqrt(s  + 9*s + 36) + 2*s + 3
        k=------------------------------}
                         3
```

The proof of Theorem 400 also has nine cases to consider in order to be able to identify $k_{opt}$ correctly. The following annotated program (cover) was used to confirm the algebra involved:

```
    off acorn;

  [ This procedure expects exp = k*var + n, and returns the value of
    n mod const:

    e.g. if exp = 6t-1, var = t, const = 3, then the result is 2]

    procedure mod(exp, var, const);
    begin
      while not numberp exp do exp := exp - var;
      while exp<0 do exp := exp + const;
      while exp>=const do exp := exp - const;
      return exp;
    end;
```

[ This procedure determines the integer and fractional parts of the numerical expression exp:

e.g. if exp = 29/3, the result is {9, 2/3} ]

```
procedure dispn(exp);
begin
local p, q;
  on mcd;
  off rational;
  p := 0;
  q := den(exp);
  while den(exp) > 1 do
  <<  p := p + 1;
      exp := exp - 1/q;
  >>;
  off mcd;
  on rational;
  return {exp, p/q};
end;

operator n, c0d, c00, c01, c02;
```

[ We now define $C_\Theta$'(k) as c0d(k), and c00(k), c01(k), c02(k) as the three possibilities for $C_\Theta$(k): ]

```
for all k let n(k) = s-k;
for all k let c0d(k) = (k+1)*(((n(k)+3)*(n(k)+2)+1)/3)+(n(k)*(n(k)+3)-4)/3;
for all k let c00(k) = (k+1)*((n(k)+3)*(n(k)+2)/3) + (n(k)*(n(k)+3)-12)/3;
for all k let c01(k) = (k+1)*((n(k)+3)*(n(k)+2)/3) + (n(k)*(n(k)+3)-4)/3;
for all k let c02(k) = (k+1)*(((n(k)+3)*(n(k)+2)+1)/3)+(n(k)*(n(k)+3)-7)/3;
```

[ This procedure defines $C_\Theta$(k) as c0(k) - by choosing amongst c00(k), c01(k), c02(k) according to the value of (s-k) mod 3: ]

```
procedure c0(k);
begin
local c0n;
  c0n := mod(n(k), t, 3);
  if c0n=0 then return c00(k);
  if c0n=1 then return c01(k);
  if c0n=2 then return c02(k);
end;
```

[ This procedure prints out values of $C_\Theta$'(k) and $C_\Theta$(k) around k = s/3 in terms of t where s = 9t+r.

On entry: ss - is the value of s (between 9t and 9t+8)
min, max - determine the range around 3t to be evaluated: $C_\Theta$'(k) and $C_\Theta$(k) will be printed for 3t+min <= k <= 3t+max

As an aid to checking the results, numerical values of $C_\Theta$'(k) and $C_\Theta$(k) are also printed for t=6: ]

```
procedure dos(ss, min, max);
begin
local kk, x;
  s := ss;   [ to enable the function n(k) and the procedure c0(k)
               to evaluate correctly ]
  write " ";
  write "### s = ", s;
  kk := 3*t+min;
  for i := 1:(max-min+1) do
  <<  x := dispn(sub(t=6,c0d(kk)));
      if second(x) = 0 then
        write "c0d(", kk, ") = ", c0d(kk),
              "   c0d(", sub(t=6,kk), ") = ", first(x)
      else
        write "c0d(", kk, ") = ", c0d(kk),
              "   c0d(", sub(t=6,kk), ") = ", first(x), " + ", second(x);
      kk := kk + 1;
  >>;
  write " ";
  kk := 3*t+min;
  for i := 1:(max-min+1) do
  <<  x := dispn(sub(t=6,c0(kk)));
      if second(x) = 0 then
          write "c0(", kk, ") = ", c0(kk),
                "   c0(", sub(t=6,kk), ") = ", first(x)
      else
          write "c0(", kk, ") = ", c0(kk),
                "   c0(", sub(t=6,kk), ") = ", first(x), " + ", second(x);
      kk := kk + 1;
  >>;
end;

off allfac;
off mcd;
on rational;

[ Note that the ranges have been chosen by previous experiment! ]

dos(9*t, -2, 1);
dos(9*t+1, -1, 1);
dos(9*t+2, -1, 1);
dos(9*t+3, 0, 2);
dos(9*t+4, 0, 2);
dos(9*t+5, 0, 2);
dos(9*t+6, 0, 3);
dos(9*t+7, 1, 3);
dos(9*t+8, 1, 3);

;end;
```

The results of the run are as follows, and the conclusions of Theorem 300 (a) - (i) can be checked against them:

### s = 9*t

$c0d(3*t - 2) = 36*t^3 + 54*t^2 + 17*t - 5$ c0d(16) = 9817

$c0d(3*t - 1) = 36*t^3 + 54*t^2 + 23*t$ c0d(17) = 9858

$c0d(3*t) = 36*t^3 + 54*t^2 + 23*t + 1$ c0d(18) = 9859

$c0d(3*t + 1) = 36*t^3 + 54*t^2 + 17*t$ c0d(19) = 9822

$c0(3*t - 2) = 36*t^3 + 54*t^2 + 17*t - 6$ c0(16) = 9816

$c0(3*t - 1) = 36*t^3 + 54*t^2 + 22*t$ c0(17) = 9852

$c0(3*t) = 36*t^3 + 54*t^2 + 22*t - 2$ c0(18) = 9850

$c0(3*t + 1) = 36*t^3 + 54*t^2 + 17*t - 1$ c0(19) = 9821

### s = 9*t + 1

$c0d(3*t - 1) = 36*t^3 + 66*t^2 + 35*t + 2$ c0d(17) = 10364

$c0d(3*t) = 36*t^3 + 66*t^2 + 37*t + 13/3$ c0d(18) = 10378 + 1/3

$c0d(3*t + 1) = 36*t^3 + 66*t^2 + 33*t + 10/3$ c0d(19) = 10353 + 1/3

$c0(3*t - 1) = 36*t^3 + 66*t^2 + 35*t + 1$ c0(17) = 10363

$c0(3*t) = 36*t^3 + 66*t^2 + 36*t + 4$ c0(18) = 10372

$c0(3*t + 1) = 36*t^3 + 66*t^2 + 32*t$ c0(19) = 10344

### s = 9*t + 2

$c0d(3*t - 1) = 36*t^3 + 78*t^2 + 49*t + 14/3$ c0d(17) = 10882 + 2/3

$c0d(3*t) = 36*t^3 + 78*t^2 + 53*t + 9$ c0d(18) = 10911

$c0d(3*t + 1) = 36*t^3 + 78*t^2 + 51*t + 26/3$ c0d(19) = 10898 + 2/3

$c0(3*t - 1) = 36*t^3 + 78*t^2 + 48*t + 2$ c0(17) = 10874

$c0(3*t) = 36*t^3 + 78*t^2 + 53*t + 8$ c0(18) = 10910

$c0(3*t + 1) = 36*t^3 + 78*t^2 + 50*t + 8$ c0(19) = 10892

### s = 9*t + 3

$c0d(3*t) = 36*t^3 + 90*t^2 + 71*t + 15 \quad c0d(18) = 11457$

$c0d(3*t + 1) = 36*t^3 + 90*t^2 + 71*t + 16 \quad c0d(19) = 11458$

$c0d(3*t + 2) = 36*t^3 + 90*t^2 + 65*t + 13 \quad c0d(20) = 11419$

$c0(3*t) = 36*t^3 + 90*t^2 + 70*t + 12 \quad c0(18) = 11448$

$c0(3*t + 1) = 36*t^3 + 90*t^2 + 71*t + 15 \quad c0(19) = 11457$

$c0(3*t + 2) = 36*t^3 + 90*t^2 + 64*t + 12 \quad c0(20) = 11412$

### s = 9*t + 4

$c0d(3*t) = 36*t^3 + 102*t^2 + 91*t + 67/3 \quad c0d(18) = 12016 + 1/3$

$c0d(3*t + 1) = 36*t^3 + 102*t^2 + 93*t + 76/3 \quad c0d(19) = 12031 + 1/3$

$c0d(3*t + 2) = 36*t^3 + 102*t^2 + 89*t + 23 \quad c0d(20) = 12005$

$c0(3*t) = 36*t^3 + 102*t^2 + 90*t + 22 \quad c0(18) = 12010$

$c0(3*t + 1) = 36*t^3 + 102*t^2 + 92*t + 22 \quad c0(19) = 12022$

$c0(3*t + 2) = 36*t^3 + 102*t^2 + 89*t + 22 \quad c0(20) = 12004$

### s = 9*t + 5

$c0d(3*t) = 36*t^3 + 114*t^2 + 113*t + 31 \quad c0d(18) = 12589$

$c0d(3*t + 1) = 36*t^3 + 114*t^2 + 117*t + 110/3 \quad c0d(19) = 12618 + 2/3$

$c0d(3*t + 2) = 36*t^3 + 114*t^2 + 115*t + 107/3 \quad c0d(20) = 12605 + 2/3$

$c0(3*t) = 36*t^3 + 114*t^2 + 113*t + 30 \quad c0(18) = 12588$

$c0(3*t + 1) = 36*t^3 + 114*t^2 + 116*t + 36 \quad c0(19) = 12612$

$c0(3*t + 2) = 36*t^3 + 114*t^2 + 114*t + 32 \quad c0(20) = 12596$

```
 ### s = 9*t + 6
                     3        2
  c0d(3*t) = 36*t  + 126*t  + 137*t + 41   c0d(18) = 13175

                         3        2
  c0d(3*t + 1) = 36*t  + 126*t  + 143*t + 50   c0d(19) = 13220

                         3        2
  c0d(3*t + 2) = 36*t  + 126*t  + 143*t + 51   c0d(20) = 13221

                         3        2
  c0d(3*t + 3) = 36*t  + 126*t  + 137*t + 46   c0d(21) = 13180

                    3        2
  c0(3*t) = 36*t  + 126*t  + 136*t + 38   c0(18) = 13166

                        3        2
  c0(3*t + 1) = 36*t  + 126*t  + 143*t + 49   c0(19) = 13219

                        3        2
  c0(3*t + 2) = 36*t  + 126*t  + 142*t + 50   c0(20) = 13214

                        3        2
  c0(3*t + 3) = 36*t  + 126*t  + 136*t + 42   c0(21) = 13170

### s = 9*t + 7
                         3        2
  c0d(3*t + 1) = 36*t  + 138*t  + 171*t + 196/3   c0d(19) = 13835 + 1/3

                         3        2
  c0d(3*t + 2) = 36*t  + 138*t  + 173*t + 69   c0d(20) = 13851

                         3        2
  c0d(3*t + 3) = 36*t  + 138*t  + 169*t + 196/3   c0d(21) = 13823 + 1/3

                        3        2
  c0(3*t + 1) = 36*t  + 138*t  + 170*t + 62   c0(19) = 13826

                        3        2
  c0(3*t + 2) = 36*t  + 138*t  + 173*t + 68   c0(20) = 13850

                        3        2
  c0(3*t + 3) = 36*t  + 138*t  + 168*t + 64   c0(21) = 13816

### s = 9*t + 8
                         3        2
  c0d(3*t + 1) = 36*t  + 150*t  + 201*t + 248/3   c0d(19) = 14464 + 2/3

                         3        2
  c0d(3*t + 2) = 36*t  + 150*t  + 205*t + 269/3   c0d(20) = 14495 + 2/3

                         3        2
  c0d(3*t + 3) = 36*t  + 150*t  + 203*t + 88   c0d(21) = 14482

                        3        2
  c0(3*t + 1) = 36*t  + 150*t  + 200*t + 82   c0(19) = 14458

                        3        2
  c0(3*t + 2) = 36*t  + 150*t  + 204*t + 86   c0(20) = 14486

                        3        2
  c0(3*t + 3) = 36*t  + 150*t  + 203*t + 87   c0(21) = 14481
```

Algebra connected with the optimisation of $M_2$'(k) in the proof of Theorem 402 in section 4.2 was confirmed with the help of the following program (m2d2):

```
off acorn;

operator n, m2d;

for all k let n(k) = s-k;
for all k let m2d(k) = (k+2)*(((n(k)+3)*(n(k)+2)+1)/3 - 2);

off allfac;
off mcd;

for i := 0:2 do
<< s := 3*u+i;
   write "s = ", s;
   for j := -2:2 do
     write "  3*m2d(", u+j, ") = ", 3*m2d(u+j);
>>;

procedure doit(ss, u, k1);
begin
  s := ss;
  write "s = ", s, "; u = ", u, "; m2d(", k1, ") = ", m2d(k1);
end;

doit(9*t, 3*t, 3*t-1);
doit(9*t+3, 3*t+1, 3*t);
doit(9*t+6, 3*t+2, 3*t+1);
doit(9*t+1, 3*t, 3*t);
doit(9*t+4, 3*t+1, 3*t+1);
doit(9*t+7, 3*t+2, 3*t+2);
doit(9*t+2, 3*t, 3*t);
doit(9*t+5, 3*t+1, 3*t+1);
doit(9*t+8, 3*t+2, 3*t+2);

;end;
```

The results are as follows, and can be checked against (a), (b), (c) and (aa) to (cc) as marked:

```
   (a): s = 3*u
                            3       2
          3*m2d(u - 2) = 4*u  + 18*u  + 15*u

                            3       2
          3*m2d(u - 1) = 4*u  + 18*u  + 21*u + 7

                        3       2
          3*m2d(u) = 4*u  + 18*u  + 21*u + 2

                            3       2
          3*m2d(u + 1) = 4*u  + 18*u  + 15*u - 9

                            3       2
          3*m2d(u + 2) = 4*u  + 18*u  + 3*u - 20
```

(b): $s = 3*u + 1$

$$3*m2d(u - 2) = 4*u^3 + 22*u^2 + 25*u$$

$$3*m2d(u - 1) = 4*u^3 + 22*u^2 + 33*u + 15$$

$$3*m2d(u) = 4*u^3 + 22*u^2 + 35*u + 14$$

$$3*m2d(u + 1) = 4*u^3 + 22*u^2 + 31*u + 3$$

$$3*m2d(u + 2) = 4*u^3 + 22*u^2 + 21*u - 12$$

(c): $s = 3*u + 2$

$$3*m2d(u - 2) = 4*u^3 + 26*u^2 + 37*u$$

$$3*m2d(u - 1) = 4*u^3 + 26*u^2 + 47*u + 25$$

$$3*m2d(u) = 4*u^3 + 26*u^2 + 51*u + 30$$

$$3*m2d(u + 1) = 4*u^3 + 26*u^2 + 49*u + 21$$

$$3*m2d(u + 2) = 4*u^3 + 26*u^2 + 41*u + 4$$

(aa): $s = 9*t;\ u = 3*t;\ m2d(3*t - 1) = 36*t^3 + 54*t^2 + 21*t + 7/3$

(ab): $s = 9*t + 3;\ u = 3*t + 1;\ m2d(3*t) = 36*t^3 + 90*t^2 + 69*t + 50/3$

(ac): $s = 9*t + 6;\ u = 3*t + 2;\ m2d(3*t + 1) = 36*t^3 + 126*t^2 + 141*t + 51$

(ba): $s = 9*t + 1;\ u = 3*t;\ m2d(3*t) = 36*t^3 + 66*t^2 + 35*t + 14/3$

(bb): $s = 9*t + 4;\ u = 3*t + 1;\ m2d(3*t + 1) = 36*t^3 + 102*t^2 + 91*t + 25$

(bc): $s = 9*t + 7;\ u = 3*t + 2;\ m2d(3*t + 2) = 36*t^3 + 138*t^2 + 171*t + 68$

(ca): $s = 9*t + 2;\ u = 3*t;\ m2d(3*t) = 36*t^3 + 78*t^2 + 51*t + 10$

(cb): $s = 9*t + 5;\ u = 3*t + 1;\ m2d(3*t + 1) = 36*t^3 + 114*t^2 + 115*t + 37$

(cc): $s = 9*t + 8;\ u = 3*t + 2;\ m2d(3*t + 2) = 36*t^3 + 150*t^2 + 203*t + 268/3$

In Section 4.3, the proof of Theorem 403 involves examination of $C_1(k)$ for the case (s-k) mod 3 = 0 with k<=2. The following program confirms the detailed algebraic manipulations (c1small):

```
off acorn;
off allfac;

operator n, c1;
for all k let n(k) = s-k;
for all k let c1(k) = (k+1)*((n(k)*n(k)+5*n(k)+3)/3) +
                      ((n(k)*n(k)+3*n(k))/3-1);

c1(0);
  sub(s=9*t, c1(0));
  sub(s=9*t+3, c1(0));
  sub(s=9*t+6, c1(0));
c1(1);
  sub(s=9*t+1, c1(1));
  sub(s=9*t+4, c1(1));
  sub(s=9*t+7, c1(1));
c1(2);
  sub(s=9*t+2, c1(2));
  sub(s=9*t+5, c1(2));
  sub(s=9*t+8, c1(2));
```

The results are as follows:

```
C: c1(0);

      2
  2*s  + 8*s
 ------------
      3

  C: sub(s=9*t, c1(0));

      2
 54*t  + 24*t

  C: sub(s=9*t+3, c1(0));

      2
 54*t  + 60*t + 14

  C: sub(s=9*t+6, c1(0));

      2
 54*t  + 96*t + 40
```

```
C: c1(1);

       2
 3*s  + 7*s - 7
----------------
       3

  C: sub(s=9*t+1, c1(1));

       2
81*t  + 39*t + 1

  C: sub(s=9*t+4, c1(1));

       2
81*t  + 93*t + 23

  C: sub(s=9*t+7, c1(1));

       2
81*t  + 147*t + 63

C: c1(2);

       2
 4*s  + 2*s - 14
-----------------
        3

  C: sub(s=9*t+2, c1(2));

        2
108*t  + 54*t + 2

  C: sub(s=9*t+5, c1(2));

        2
108*t  + 126*t + 32

  C: sub(s=9*t+8, c1(2));

        2
108*t  + 198*t + 86
```

# SECTION 7

## SPECIAL CASES

M(3,s) and the corresponding maximal sets for s<81 are listed, and their characteristics examined.

## 7.1 Results for s<81

Before work on this proof was started, values of M(3,s) and the corresponding maximal sets were computed using a variety of techniques for all values of s<=284; details of these techniques are given in a separate document.

These results confirm and extend Theorem 101, showing that the the limit of s>=81 can be extended to s>=23.

Results for s<23 are tabulated below in Table 700, together with details of the underlying stride generators.

| s | $a_1$ | $a_2$ | $a_3$ | M(3,s) | SG |
|---|---|---|---|---|---|
| 1 | 1 | 2 | 3 | 3 | OSG(1,0) |
| 2 | 1 | 3 | 4 | 8 | OSG(1,1) |
| 3 | 1 | 4 | 5 | 15 | OSG(1,2) |
| 4 | 1 | 5 | 8 | 26 | SG(2,2) |
| 5 | 1 | 6 | 7 | 35 | OSG(1,4) |
| 6 | 1 | 7 | 12 | 52 | SG(3,3) |
| 7 | 1 | 8 | 13 | 69 | OSG(3,2) |
| 8 | 1 | 9 | 14 | 89 | SG(4,2) |
| 9 | 1 | 9 | 20 | 112 | SG(5,2) |
| 10 | 1 | 10 | 26 | 146 | SG(6,2) |
| 11 | 1 | 9 | 30 | 172 | OSG(7,1) |
| | 1 | 10 | 26 | 172 | SG(6,2) |
| 12 | 1 | 11 | 37 | 212 | OSG(8,1) |
| 13 | 1 | 13 | 34 | 259 | OSG(7,2) |
| 14 | 1 | 12 | 52 | 302 | OSG(10,1) |
| 15 | 1 | 12 | 52 | 354 | OSG(10,1) |
| 16 | 1 | 15 | 54 | 418 | SG(10,2) |
| 17 | 1 | 14 | 61 | 476 | OSG(11,1) |
| 18 | 1 | 15 | 80 | 548 | OSG(13,1) |
| 19 | 1 | 18 | 65 | 633 | SG(11,2) |
| 20 | 1 | 17 | 91 | 714 | OSG(14,1) |
| 21 | 1 | 17 | 91 | 805 | OSG(14,1) |
| 22 | 1 | 19 | 102 | 902 | OSG(15,1) |
| | 1 | 20 | 92 | 902 | SG(14,2) |

Table 700

# SECTION 8

# INDICES

Section 8.1 provides a short glossary of the terms used in the proof.

Sections 8.2 to 8.5 contain summary listings of all the definitions, theorems, figures and tables that form part of the proof.

Section 8.6 lists some relevant references.

## 8.1 Glossary

**break** [Def. 200, 201]

A "break" in a stride generator SG(n,p) is a value $0<y<a_3$ such that:

$$y + ja_3 = c_2a_2 + c_1 \quad c_2+c_1<=n+j-1 \quad c_2,c_1>=0$$

is not soluble for any j<=p+1. The "order of the break" is the smallest value of j>p+1, if any, for which this equation is soluble.

If the stride generator underlies a cover C(A,3,s), then the smallest break value corresponds to the first value that cannot be generated by the cover.

**canonical** [Def. 201, 202; Sect. 1.1]

a) A "canonical break" is one whose equation (see "break") is not soluble for any j>p+1.

b) A "canonical stride generator" is one in which all breaks are canonical. We show also that canonical breaks can only exist in canonical stride generators.

c) The "canonical generation" is the one with smallest order.

**cover** [Def. 100, 206]

a) A set's "cover" X with respect to a value s is one less than the smallest value that cannot be generated; we write X = C(A,3,s).

b) One thread "covers" another if all the values in the second are included in the first.

**final stride** [Sect. 1.1]

The "final stride" is the last complete stride that can be generated by a set A with respect to a value s. It is the kth stride where:

$$C(A,3,s) = (k+1)a_3 + Y \quad \text{where} \quad 0<=Y<a_3-1$$

**fundamental break** [Def. 207]

A "fundamental break" y in a stride generator SG(n,p) is one which lies in the range $a_3-a_2<=y<a_3-a_2+n$. We show that every stride generator has at least one and at most two fundamental breaks.

**generation** [Def. 100; Sect. 1.1, 2.1]

a) A set A = {1, $a_2$, $a_3$} "generates" a value x with respect to s if:

$$x = c_3a_3 + c_2a_2 + c_1 \quad \text{where} \quad c_3+c_2+c_1<=s \quad \text{and} \quad c_3,c_2,c_1>=0$$

If $x = b_3a_3 + b_2$ where $0<=b_2<a_3$, this generation is said to have "order" $b_3-c_3$.

b) A stride generator SG(n,p) = {1, $a_2$, $a_3$} "generates" a value x if:

$$x + ia_3 = c_2a_2 + c_1 \quad \text{where} \quad c_2+c_1<=n+i, \quad i<=p \quad \text{and} \quad c_2,c_1>=0$$

This generation is said to have "order" i.

key [Def. 502]

We show that the signature of any stride generator SG(n,p) can be generated as the sequence {kj mod (p+1): 0<=k<=p} for some value 1<=j<=p. The value j is called the "key".

length [Def. 200]

The "length" of a stride generator SG(n,p) = {1, $a_2$, $a_3$} is $a_3$.

maximal [Def. 100]

a) The "maximal cover" M(3,s) with respect to a value s is the highest cover value C(A,3,s) over all sets A.

b) A "maximal set" is one which has maximal cover; there may be more than one maximal set for a given value of s.

non-trivial cover [Def. 203]

A set A has a "non-trivial cover" if C(A,3,s) >= $a_3$.

optimal stride generator [Sect. 1.1]

An "optimal stride generator" OSG(n,p) of order p is one that is at least as long as any other stride generator of order p. We show that for large enough s, maximal sets are always optimal stride generators of order 1.

order

a) "Stride generator order" - see "stride generator".

b) "Generation order" - see "generation".

c) "Thread order" - see "thread".

d) "Break order" - see "break".

potential cover [Def. 208]

The "potential cover" of a stride generator A = SG(n,p) with respect to a value s is defined as:

$$(s-n+1)a_3 + y-1$$

where y is the first break in the stride generator. We prove that the cover of a set A is equal to the potential cover of its underlying stride generator.

relative position [Def. 501]

The "relative position" of two threads in a thread diagram is represented by (j-i, y-x) where the thread orders are j,i and their first elements are y,x respectively.

signature [Def. 500]

The threads of a stride generator SG(n,p) appear in a particular order sequence from left to right in the thread diagram. Each order appears exactly once in the range $0<=x<a_2$, and the sequence is then repeated. This sequence of p+1 integers which starts with 0 is called the "signature" of the stride generator.

stride [Sect. 1.1]

The set of values $ja_3<=x<(j+1)a_3$ is called the jth "stride".

stride generator [Def. 200]

An "n-stride generator of order p" is a set A = {1, $a_2$, $a_3$}, and is denoted SG(n,p). The concept is introduced informally in Section 1.1 and defined formally in Section 2.1.

Roughly speaking, a stride generator is just able to generate all values $0<x<a_3$, but fails to generate at least one value $a_3<x<2a_3$.

We show that there is a correspondence between stride generator generations and cover generations, and that every cover has an underlying stride generator.

sub-optimal stride generator [Def. 400]

A "sub-optimal stride generator" SGi(n,p) is one whose length is i less than that of the corresponding first order optimal stride generator OSG(n,1).

thread [Def. 205]

A "thread of order i" (or "i-thread") $T_i(c_2)$ represents the contiguous set of order i generations x to y where:

$x = c_2a_2$
$y = c_2a_2 + c_1$ and $c_2+c_1=n+i$

thread diagram [Sect. 2.4]

Every generation in a stride generator SG(n,p) belongs to one or more threads of order <= p. The "thread diagram" for a stride generator represents each i-thread graphically as a horizontal line at height i running from the first to the last generation of the thread. These diagrams illustrate many of the properties of stride generators in a pictorial form.

underlying stride generator [Def. 204]

We prove that every set A with a non-trivial cover C(A,3,s) is also a stride generator SG(s-k,p) for some 0<=p<=k<s; this stride generator is said to "underlie" the cover.

## 8.2 Definitions

| No. | Section | Description |
|---|---|---|
| 100 | 1.1 | The postage stamp problem - cover, generation, maximal set |
| 200 | 2.1 | Definition of a stride generator - length, order, breaks |
| 201 | 2.2 | Order of a break in a stride generator |
| 202 | 2.2 | Canonical stride generator |
| 203 | 2.3 | Non-trivial cover |
| 204 | 2.3 | The stride generator underlying a cover |
| 205 | 2.4 | Thread notation |
| 206 | 2.7 | One thread covers another iff ... |
| 207 | 2.9 | Fundamental breaks |
| 208 | 2.10 | Potential cover |
| 400 | 4.2 | Sub-optimal stride generators SGi(n,p) |
| 500 | 5.2 | Signature of a stride generator |
| 501 | 5.2 | Relative position of two threads |
| 502 | 5.2 | The key of a stride generator's signature |

## 8.3 Theorems

| No. | Section | Description |
|---|---|---|
| 100 | 1.1 | If stride k can be generated, then so can strides p to k |
| 101 | 1.2 | OSG(s-$k_{opt}$,1) is the maximal set with cover $X_{opt}$ for s>=81 |
| 200 | 2.1 | Alternative definition of a stride generator (SG1, SG2, SG3) |
| 201 | 2.3 | $C(A,3,s)=(k+1)a_3+Y$ where $0<=Y<a_3-1$ |
| 202 | 2.3 | Every cover defines an underlying stride generator |
| 203 | 2.3 | If the underlying SG is canonical, then $C(A,3,s)=(k+1)a_3+y-1$ |
| 204 | 2.5 | $y+ia_3=c_2a_2+c_1$  $c_2+c_1=n+i$: breaks lie at the end of threads |
| 205 | 2.5 | $(y+1)+ka_3=c_2a_2$: the value after a break starts a thread |
| 206 | 2.5 | Equivalent definition of break order in terms of threads |
| 207 | 2.5 | Equivalent definition of stride generators in terms of threads |
| 208 | 2.6 | $T_i(c_2)$ at x <=> $T_i(c_2-1)$ at $x-a_2$ |
| 209 | 2.6 | $T_i(c_2)$ at x <=> $T_{i+1}(c_2+C_2)$ at $x-C_1$ |
| 210 | 2.6 | $T_i(c_2)$ at x <=> $T_{i+1}(c_2+1)$ at $x-(a_3-a_2)$ |
| 211 | 2.7 | i>j, $ST_i>ST_j$ => $L_i<L_j$ and $b_2<=c_2+1$ |
| 212 | 2.7 | i-thread covered by j-thread => all k-threads for k>=i are covered |
| 213 | 2.7 | No thread can be covered by another thread in a stride generator |
| 214 | 2.8 | SG(n,0) => $C_1=0$ or $a_2-C_2<=C_1<a_2$ and then $n=a_2+C_2-2$ |
| 215 | 2.8 | $a_3>=(a_2-1)a_2$ means the stride generator is order 0 |
| 216 | 2.8 | SG(n,0) with two breaks have $a_3=(C_2+1)a_2-1$ and $y_2=a_3-1$ |
| 217 | 2.8 | All SG(n,0) are canonical |
| 218 | 2.8 | SG(n,0) cannot also be SG(n',p') for p'>0 |
| 219 | 2.9 | $y>=a_3-a_2$ |
| 220 | 2.9 | $(a_3-a_2)<=y<(a_3-a_2)+n$ |
| 221 | 2.9 | $x>=a_3-a_2$ requires an order p generation |
| 222 | 2.9 | $(a_3-a_2)<=x<=(a_3-a_2)+n$ requires an order p generation |
| 223 | 2.9 | There is exactly one k-thread that starts in the range $a_3-a_2$ to $a_3$ |
| 224 | 2.9 | All possible positions for threads in a SG are filled |
| 225 | 2.9 | One of the threads surrounding a break must be of order p |
| 226 | 2.9 | No SG of order>0 can have a fundamental break at $y=a_3-1$ |
| 227 | 2.9 | There are at most two fundamental breaks |
| 228 | 2.9 | If there are two fundamental breaks, both are canonical |
| 229 | 2.9 | In a break series, all breaks are canonical or the order decreases |
| 230 | 2.9 | SG(n,p) p>0 has all canonical breaks or the break order decreases |
| 231 | 2.10 | Potential cover is real for an underlying stride generator |
| 232 | 2.11 | Each set defines a SG series with n decreasing, p increasing |
| 233 | 2.11 | Every set A defines a stride generator |
| 234 | 2.11 | J exists st $y+ja_3=c_2a_2+c_1$, $c_2+c_1<=n+j-1$ is not soluble for j>=J |
| 235 | 2.11 | n'<n => p'>p+1 for SG(n,p), SG(n',p') |
| 236 | 2.11 | Every set A has a canonical stride generator |
| 237 | 2.11 | The canonical stride generator for a set A is unique |
| 238 | 2.11 | SG(n',p') is canonical => SG(n,p) non-canonical, n>n' |
| 239 | 2.12 | $a_2>=a_3(p+1)/(n+p+1)$ |
| 240 | 2.12 | $a_2<=n(p+1)+1$ |
| 241 | 2.12 | $a_3<=n(n+p+1)+(n+p+1)/(p+1)$ |
| 242 | 2.13 | SG(n,p) can be as long as you like |
| 300 | 3.1 | Formulae for OSG(n,0) |
| 301 | 3.2 | Formulae for OSG(n,1) |
| 302 | 3.3 | Formulae for SG1(n,1) |
| 400 | 4.1 | OSG(s-k,1)'s maximum potential cover is $X_{opt}$ at $k_{opt}$ |
| 401 | 4.1 | OSG(s-$k_{opt}$,1) defines a real cover = $X_{opt}$ |
| 402 | 4.2 | Sub-optimal SGi(n,p) cannot better OSG(n,1) for s>=36 |
| 403 | 4.3 | SG1(n,1) cannot better OSG(n,1) for s>=9 |
| 500 | 5.1 | The order of the SG underlying the maximal cover <=3 for s>=40 |
| 501 | 5.2 | Relative positions replicate |
| 502 | 5.2 | i+j<=p+1 |
| 503 | 5.2 | i+j=p+1 |
| 504 | 5.2 | There are only two relative positions in a stride generator |
| 505 | 5.2 | (j,p+1) are coprime for SG(n,p) with key j |
| 506 | 5.2 | Limit on SG(n,p) with key = 1 or key = p |
| 507 | 5.2 | OSG(n,1) betters SG(n,p) with key=1,p for large n |
| 508 | 5.2 | OSG(n,1) betters OSG(n,2) for n>132, OSG(n,3) for n>100 |
| 509 | 5.2 | OSG(n,1) betters OSG(n,0) for all n |
| 510 | 5.2 | OSG(n,1) betters OSG(n,2), OSG(n,3) for n>48 and SG1(n,1) for n>51 |
| 511 | 5.3 | OSG(n,1) better than SG(n,11) with key=5, y=$ST_{12}$-1, large n |
| 512 | 5.4 | Only OSG(n,1) matters for s>=81 |

## 8.4 Figures

## 8.5 Tables